\documentclass[trsc,nonblindrev]{informs3} 

\OneAndAHalfSpacedXI 

\usepackage{natbib}

 \bibpunct[, ]{(}{)}{,}{a}{}{,}%
 \def\bibfont{\small}%
\TheoremsNumberedThrough     

\EquationsNumberedThrough    

   \usepackage[hidelinks]{hyperref} 
\usepackage{subcaption}
\usepackage{mathtools}
\usepackage{dsfont}
\usepackage{color}
\usepackage{natbib}
    \usepackage{multirow}
\usepackage{ulem}

\newcommand{\calS}{\mathcal{S}}
\newcommand{\calX}{\mathcal{X}}
 \newcommand{\calF}{\mathcal{F}}

\newcommand{\WL}{\underline{W}}
\newcommand{\WU}{\overline{W}}
\newcommand{\Wb}{\pmb{W}}
\newcommand{\E}{\mathbb{E}}
\newcommand{\Prob}{\mathbb{P}}
\newcommand{\R}{\mathbb{R}}
\newcommand{\mub}{\pmb{\mu}}
\newcommand{\xbb}{\pmb{x}}
\newcommand{\yb}{\pmb{y}}
\newcommand{\etab}{\pmb{\eta}}
\newcommand{\lambdab}{\pmb{\lambda}}
\newcommand{\vb}{\pmb{v}}
\newcommand{\pib}{\pmb{\pi}}
\newcommand{\calP}{\mathcal{P}}

\newcommand{\F}{\mathbb{F}}	         
\newcommand{\CVAR}{\text{CVaR}}	              
     
\usepackage[ruled, commentsnumbered]{algorithm2e}

\begin{document}


 \RUNAUTHOR{Shehadeh, K.S.}

 \RUNTITLE{DRO Approaches for MFRSP}

\TITLE{Distributionally Robust Optimization Approaches for a Stochastic Mobile Facility \textcolor{black}{Fleet Sizing,} Routing and Scheduling Problem}

\ARTICLEAUTHORS{%
\AUTHOR{Karmel S. Shehadeh}
\AFF{Department of Industrial and Systems Engineering, Lehigh University, Bethlehem, PA, \EMAIL{kas720@lehigh.edu}}
} 

\ABSTRACT{

\noindent \textcolor{black}{We propose two distributionally robust optimization (DRO) models for a mobile facility (MF) fleet sizing, routing, and scheduling problem (MFRSP) with time-dependent and random demand,  as well as methodologies for solving these models. Specifically, given a set of MFs, a planning horizon, and a service region, our models aim to find the number of MFs to use  (i.e., fleet size) within the planning horizon and a route and time schedule for each MF in the fleet. The objective is to minimize the fixed cost of establishing the MF fleet  plus a risk measure (expectation or mean conditional value-at-risk) of the operational cost over all demand distributions defined by an ambiguity set.}  In the first model, we use an ambiguity set based on the demand's mean, support, and mean absolute deviation. In the second model, we use an ambiguity set that incorporates all distributions within a 1-Wasserstein distance from a reference distribution. \textcolor{black}{To solve the proposed DRO models, we propose a decomposition-based algorithm. In addition, we derive valid lower bound inequalities that efficiently strengthen the master problem in the decomposition algorithm, thus improving convergence. We also derive two families of symmetry breaking constraints that improve the solvability of the proposed models. Finally, we present extensive computational experiments comparing the operational and computational performance of the proposed models and a stochastic programming model, demonstrating where significant performance improvements could be gained and derive insights into the MFRSP.}}

\color{black}

\KEYWORDS{Facility location, mobile facility,  demand uncertainty, scheduling  and routing, mixed-integer programming, distributionally robust optimization.}

\maketitle

\section{Introduction}

\noindent A \textit{mobile facility} (MF) is a facility capable of moving from one place to another, providing real-time service to customers in the vicinity of its location when it is stationary \citep{halper2011mobile}. In this paper, we study a mobile facility \textcolor{black}{fleet sizing, }routing, and scheduling problem (MFRSP) with stochastic demand. Specifically, in this problem, we aim to find the number of MFs (i.e., fleet size) to use in a given service region over a specified planning horizon and the route and schedule for each MF in the fleet. The demand level of each customer in each time period is random.  The probability distribution of the demand is unknown, and only partial information about the demand (e.g., mean and range) may be available.  The objective is to find the MF fleet size, routing, and scheduling decisions that minimize the sum of the fixed cost of establishing the MF fleet, the cost of assigning demand to the MFs (e.g., transportation cost), and the cost of unsatisfied demand (i.e., shortage cost).

The concept of MF routing and scheduling is very different than conventional static facility location (FL) and conventional vehicle routing (VR) problems. In static FL problems, we usually consider opening facilities at fixed locations. Conventional VR problems aims at handling the movement of items between facilities (e.g., depots) and customers. \textcolor{black}{A mobile facility is a \textit{facility-like vehicle} that functions as a traditional facility when it is stationary, except that it can move from one place to another if necessary \citep{lei2014multicut}}. Thus, the most evident advantage of MFs over fixed facilities is their flexibility in moving to accommodate the change in the demand over time and location \citep{halper2011mobile, lei2014multicut, lei2016two}.

 MFs are used in many applications ranging from cellular services, healthcare services, to humanitarian relief logistics. For example, light trucks with portable cellular stations can provide cellular service in areas where existing cellular network of base stations temporarily fails \citep{halper2011mobile}.  Mobile clinics (i.e., customized MFs fitted with medical equipment and staffed by health professionals) can travel to rural and urban areas to provide various (prevention, testing, diagnostic) health services. Mobile clinics also offer alternative healthcare (service) delivery options when a disaster, conflict, or other events cause stationary healthcare facilities to close or stop operations \citep{blackwell2007use, brown2014mobile, du2007mobile, gibson2011households, oriol2009calculating, song2013mobile}. For example, mobile clinics played a significant role in providing drive-through COVID-19 testing sites or triage locations during the COVID-19 pandemic. In 2019, the mobile health clinic market was valued at nearly 2  billion USD and is expected to increase to $\sim$12 billion USD by 2028 \citep{MFMarket}. In humanitarian relief logistics, MFs give relief organizations the ability to provide aid to populations dispersed in remote and dense  areas. These examples motivate the need for computationally efficient optimization tools to support decision-making in all areas of the MF industry.

MF operators often seek a  \textcolor{black}{strategic} and tactical plan, including the size of the MF fleet (\textcolor{black}{strategic}) and a routing plan for each MF in the fleet (tactical and operational) that minimize their fixed operating costs and maximize demand satisfaction. Determining the fleet size, in particular, is very critical as it is a major fixed investment for starting any MF-based business. The fleet sizing problem depends on the MF operational performance, which depends on the routing and scheduling decisions. The allocation of the demand to the MFs is also very important for the entire system performance \citep{lei2016two}. For example, during the COVID-19 pandemic, Latino Connection, a community health leader, has established Pennsylvania's first COVID-19 Mobile Response Unit, CATE (i.e., Community-Accessible Testing \& Education). The goal of CATE is to provide affordable and accessible COVID-19 education, testing, and vaccinations to low-income, vulnerable communities across Pennsylvania to ensure the ability to stay safe, informed, and healthy \citep{CATEMF}. \textcolor{black}{During  COVID-19, CATE published an online schedule consisting of the mobile unit stops and the schedule at each stop.} A model that optimizes CATE's fleet size and schedules considering demand uncertainty could help improve CATE's operational performance and achieve better access to health services.


Unfortunately, the MFRSP is a challenging optimization problem for two primary reasons. First, customers' demand is random and hard to predict in advance, especially with limited data during the planning process. Second, even in a perfect world in which we know with certainty the amount of demand in each period, the deterministic MFRSP is challenging because  it is similar to the classical  FL problem \citep{halper2011mobile, lei2014multicut}.  \textcolor{black}{Thus, the incorporation of demand variability increases the overall complexity of the MFRSP.  However, ignoring demand uncertainty may lead to sub-optimal decisions and, consequently, the inability to meet customer demand (i.e., shortage).} Failure to meet customer demand may lead to adverse outcomes, especially in healthcare, as it impacts population health. It also impacts customers' satisfaction and thus the reputation of the service providers and may increase their operational cost (due to, e.g., outsourcing the excess demand to other providers).

To model uncertainty, \cite{lei2014multicut} assumed that the probability distribution of the demand is known and accordingly proposed the first a \textit{priori} two-stage stochastic optimization model (SP) for a closely related MFRSP. Although attractive, the applicability of the SP approach is limited to the case in which we know the distribution of the demand or we have sufficient data to model it. \textcolor{black}{In practice, however, one might not have access to a sufficient amount of high-quality data to estimate the demand distribution accurately.} This is especially true in application domains where the use of mobile facilities to deliver services is relatively new (e.g., mobile COVID-19 testing clinics).  Moreover, it is challenging for MF companies to obtain data from other companies (competitors) due to privacy issues. Finally, various studies show that different distributions can typically explain raw data of uncertain parameters, indicating distributional ambiguity \citep{esfahani2018data, vilkkumaa2021causes}.

Suppose we model uncertainty using a data sample from a potentially biased distribution or an assumed distribution (as in SP). In this case, the resulting nominal decision problem evaluates the cost only at this training sample, and thus the resulting decisions may be overfitted (optimistically biased). Accordingly, SP solutions may demonstrate disappointing out-of-sample performance (`\textit{black swans}') under the true distribution (or unseen data). In other words, solutions of SP decision problems often display an optimistic in-sample risk, which cannot be realized in out-of-sample settings. This phenomenon is known as the \textit{Optimizers' Curse} (i.e., an attempt to optimize based on imperfect estimates of distributions leads to biased decisions with disappointing performance) and is reminiscent of the overfitting effect in statistics \citep{smith2006optimizer}.

\textcolor{black}{Alternatively, one can construct an ambiguity set of all distributions that possess certain partial information about the demand. Then, using this ambiguity set, one can formulate a distributionally robust optimization (DRO) problem to minimize a risk measure  (e.g., expectation or conditional value-at-risk (CVaR)) of the operational cost over all distributions residing within the ambiguity set. In particular, in the DRO approach, the optimization is based on the worst-case distribution within the ambiguity set, which effectively means that the distribution of the demand is a decision variable.}

 DRO \textcolor{black}{has received substantial attention recently} in various application domains due to the following striking benefits. First, as pointed out by \cite{esfahani2018data}, DRO models are more ``\textit{honest}'' than their SP counterparts as they acknowledge the presence of distributional uncertainty. \textcolor{black}{Therefore, DRO solutions often faithfully anticipate the possibility of black swan (i.e., out-of-sample disappointment). Moreover, depending on the ambiguity set used, DRO often guarantees an out-of-sample cost that falls below the worst-case optimal cost}.  Second, DRO alleviates the unrealistic assumption of the decision-maker's complete knowledge of distributions. Third, several studies have proposed DRO models for real-world problems that are more computationally tractable than their SP counterparts, see, e.g., \cite{basciftci2019distributionally, luo2018distributionally, saif2020data,  ShehadehSanci, shehadeh2020distributionallyTucker, tsang2021distributionally, wang2020distributionally, wang2021two,wu2015approximation}. In this paper, we propose tractable DRO approaches for the MFRSP.

The ambiguity set is a key ingredient of DRO models that must (1) capture the true distribution with a high degree of certainty, and  (2) be computationally manageable (i.e., allow for a tractable DRO model or solution method). There are several methods to construct the ambiguity set.  Most applied DRO literature employs moment-based ambiguity \citep{delage2010distributionally, zhang2018ambiguous}, consisting of all distributions sharing particular moments (e.g., mean-support ambiguity). The main advantage of the mean-support ambiguity set, for example, is that it incorporates intuitive statistics that a decision-maker may easily approximate and change.  Moreover, various techniques have been developed to derive tractable moment-based DRO models. However,  asymptotic properties of the moment-based DRO model cannot often be guaranteed because the moment information represents descriptive statistics.

Recent DRO approaches define the ambiguity set by choosing a distance metric (e.g., $\phi$--divergence \citep{jiang2016data}, Wasserstein distance \citep{esfahani2018data, gao2016distributionally}) to describe the deviation from a reference (often empirical) distribution. The main advantage of Wasserstein ambiguity, for example, is that it enable decision-makers to incorporate possibly small-size data in the ambiguity set and optimization, enjoys asymptotic properties, and often offers a strong out-of-sample performance guarantee \citep{esfahani2018data,mevissen2013data}. Recent results indicate that Wasserstein's ambiguity centered around a given empirical distribution contains the unknown true distribution with a high probability and is richer than other divergence-based ambiguity sets (in particular, they contain discrete and continuous distributions as compared to, e.g., $\phi$-divergence ball centered at the empirical distribution which does not contain any continuous distribution, and Kullback-Leibler divergence ball, which must be absolutely continuous with respect to the nominal distribution).

Despite the potential advantages, there are no moment-based, Wasserstein-based, or any other DRO approaches for the specific MFRSP that we study in this paper (see Section~\ref{Sec:LitRev}).  This inspires this paper's central question: \textit{what are the computational and operational performance values of employing DRO to address demand uncertainty and ambiguity compared to the classical SP approach for the MFRSP}. To answer this question, we design and analyze two DRO models based on the demand's mean, support, and mean absolute deviation ambiguity and Wasserstein ambiguity and compare the performance of these models with the classical SP approach.

\subsection{\textcolor{black}{Contributions}}

\noindent In this paper, we present two distributionally robust MF fleet sizing, routing, and scheduling (DMFRS) models for the MFRSP, as well as methodologies for solving these models. We summarize our main contributions as follows.

\begin{enumerate}
\item \textbf{Uncertainty Modeling and Optimization Models.}  We propose the first two-stage DRO models for the MFRSP. These models aim to find the optimal  (1)  number of MFs to use within a planning horizon, (2) a routing plan and a schedule for the selected MFs, i.e., the node that each MF is located at in each time period, (3) assignment of MFs to customers. Decisions (1)-(2) are planning (first-stage) decisions, which cannot be changed in the short run. Conversely, the assignments of the demand are decided based on the demand realization, and thus are second-stage decisions. The objective is to minimize the fixed cost (i.e., cost of establishing the MF fleet and traveling inconvenience cost) plus the maximum of a risk measure (expectation or mean CVaR) of the operational cost (i.e., transportation and unsatisfied demand costs) over all possible distributions of the demand defined by an ambiguity set. In the first model (MAD-DRO), we use an ambiguity set based on the demand's mean, support, and mean absolute deviation (MAD). In the second model (W-DRO), we use an ambiguity set that incorporates all distributions within a 1-Wasserstein distance from a reference distribution.  To the best of our knowledge, and according to our literature review in Section~\ref{Sec:LitRev}, our paper is the first to address the distributional ambiguity of the demand in the MFRSP using DRO.

\item \textbf{Solution Methods.} We derive equivalent solvable reformulations of the proposed mini-max nonlinear DRO models. We propose a computationally efficient decomposition-based algorithm to solve the reformulations. In addition, we derive valid  lower bound inequalities that efficiently strengthen the master problem in the decomposition algorithm, thus improving convergence.

\item \textbf{Symmetry-Breaking Constraints.}  We derive two families of new symmetry breaking constraints, which break symmetries in the solution space of the first-stage routing and scheduling decisions and thus improve the solvability of the proposed models.  These constraints are independent of the method of modeling uncertainty. Hence, they are valid for any (deterministic and stochastic) formulation that employ the first-stage decisions of the MFRSP.  Our paper is the first to attempt to break the symmetry in the solution space of these planning decisions of the MFRSP.

\item \textbf{Computational Insights.} We conduct extensive computational experiments comparing the proposed DRO models and a classical SP model empirically and theoretically, demonstrating where significant performance improvements can be gained. Specifically, our results show  (1) how the DRO approaches have superior operational performance in terms of satisfying customers demand as compared to the SP approach; (2) the MAD-DRO model is more computationally efficient than the W-DRO model; (3) the MAD-DRO model yield more conservative decisions than the W-DRO model, which often have a higher fixed cost but significantly lower operational cost; (4) how mobile facilities can move from one location to another to accommodate the change in demand over time and location; (5) efficiency of the proposed symmetry breaking constraints and lower bound inequalities; (6) the trade-off between cost, number of MFs, MF capacity, and operational performance; and (7) the trade-off between the risk-neutral and risk-averse approaches. Most importantly, our results show the value of modeling uncertainty and distributional ambiguity.



\end{enumerate}

\subsection{Structure of the paper}

\color{black}
The remainder of the paper is structured as follows. In Section~\ref{Sec:LitRev}, we review the relevant literature.  Section~\ref{sec3:DMFRSformulation} details our problem setting. In Section~\ref{sec:SP}, we present our SP. In Section~\ref{sec:DRO_Models}, we present and analyze our proposed DRO models. In Section~\ref{sec:solutionmethods},  we present our decomposition algorithm and strategies to improve convergence. In Section~\ref{sec:computational},  we present our numerical experiments and corresponding insights. Finally, we draw conclusions and discuss future directions in Section~\ref{Sec:Conclusion}. 

\color{black}

\section{Relevant Literature}\label{Sec:LitRev}

\noindent \textcolor{black}{In this section, we review recent literature that is most relevant to our work, mainly studies that propose stochastic optimization approaches for closely related problems to the MFRSP.} There is limited literature on MF as compared to stationary facilities. However, as pointed out by \cite{lei2014multicut}, the MFRSP share some features with several well-studied problems, including Dynamic Facility Location Problem (DFLP), Vehicle Routing Problem (VRP), and the Covering Tour Problem (CTP). First, let us briefly discuss the similarities and differences between the MFRSP and these problems. Given that we consider making decisions over a planning period, then the MFRSP is somewhat similar to DFLP, which seeks to locate/re-locate facilities over a planning horizon. To mitigate the impact of demand fluctuation along the planning period,  decision-makers may  open new facilities and close or relocate existing facilities at a relocation cost (\cite{albareda2009multi, antunes2009location, contreras2011dynamic, drezner1991facility,jena2015dynamic, jena2017lagrangian,  van1982dual}). Most DFLPs assume that the relocation time is relatively short as compared to the planning horizon. In contrast, the MFRSP takes into account the relocation time of MFs.  In addition, each MF needs to follow a specific route during the entire planning horizon, which is not a requirement in DFLP.

In CTP, one seeks to select a subset of nodes to visit that can cover other nodes within  a particular coverage (\cite{current1985maximum, flores2017multi, gendreau1997covering, hachicha2000heuristics, tricoire2012bi}). In contrast to the MFRSP, CTP does not consider the variations of demand over time and assumes that the amount of demand to be met by vehicles is not related to the length of time the MF is spending at the stop. The VRP is one of the most extensively studied problem in operations research. The VRP also has numerous applications and variants \citep{subramanyam2020robust}. Both the MFRSP and the VRP consider the routing decisions of vehicles. However, the MFRSP is different than the VRP in the following ways \citep{lei2014multicut}. First, in the MFRSP, we can meet customer demand by a nearby  MF (e.g., cellular stations). In the VRP, vehicles visit customers to meet their demand. Second, the amount of demand that an MF can serve at each location depends on the duration of the MF stay, which is a decision variable. In contrast, VRPs often assume a fixed service time. Finally, most VRPs require that each customer has to be visited exactly once in each route. In contrast, in the MFRSP, some customers may not be visited, and some may be visited multiple times.

Next, we review studies that proposed stochastic optimization approaches to problems similar to the MFRSP.  \cite{halper2011mobile} introduced the concept of MF and proposed a continuous-time formulation to model the maximum covering mobile facility routing problem under deterministic settings. To solve their model, \cite{halper2011mobile} proposed several computationally effective heuristics. \cite{lei2014multicut} and \cite{lei2016two} are two closely related (and only) papers that proposed stochastic optimization approaches for MF routing and scheduling. 

\cite{lei2014multicut} assumed that the distribution of the demand is known and accordingly proposed the first a \textit{priori} two-stage SP for MFRSP. \cite{lei2014multicut}'s SP  seeks optimal first-stage routing and scheduling decisions to minimize the total expected system-wide cost, where the expectation is taken with respect to the known distribution of the demand. A priori optimization has a managerial advantage since it guarantees the regularity of service, which is beneficial for both customer and service provider. That is, a prior plan allow the customers to know when and where to receive service and enable MF service providers to be familiar with routes and better manage their time schedule during the day. The applicability of the SP approach is limited to the case in which the distribution of the demand is fully known, or we have sufficient data to model it.

Robust optimization (RO) and distributionally robust optimization (DRO) are alternative techniques to model, analyze and optimize decisions under uncertainty and ambiguity (where the underlying distributions are unknown).  RO assumes that the uncertain parameters can take any value from a pre-specified uncertainty set of possible outcomes with some structure \citep{bertsimas2004price, ben2015deriving, soyster1973convex}.  In RO, optimization is based on the worst-case scenario within the uncertainty set.

\textcolor{black}{Notably,  \cite{lei2016two} are the first to motivate the importance of handling demand uncertainty using RO. They argue that RO is useful because it only requires moderate information about the uncertain demand rather than a detailed description of the probability distribution or a large data set. Specifically, \cite{lei2016two}  proposed the first two-stage RO approach for MF feet sizing and routing problem with demand uncertainty. \cite{lei2016two}'s model aims to find the fleet size and routing decisions that minimize the fixed cost of establishing the MF fleet (first-stage) and a penalty cost for the unmet demand (second-stage).}  Optimization in \cite{lei2016two}'s RO model is based on the worst-case scenario of the demand occurring within a polyhedral uncertainty set. By focusing the optimization on the worst-case scenario, RO may lead to overly conservative and suboptimal  decisions for other more-likely scenarios \citep{chen2019robust, delage2018value}.

DRO models the uncertain parameters as random variables whose underlying probability distribution can be any distribution within a pre-defined ambiguity set.   The ambiguity set is a family of all possible distributions characterized by some known properties of random parameters  \citep{esfahani2018data}. In DRO, optimization is based on the worst-case distribution within this set. DRO is an attractive approach to model uncertainty with ambiguous distributions because: (1) it alleviates the unrealistic assumption of the decision-makers' complete knowledge of the distribution governing the uncertain parameters, (2) it is usually more computationally tractable than its SP and RO counterparts \citep{delage2018value, rahimian2019distributionally}, and (3) one can use minimal distributional information or a small sample to construct the ambiguity set and then build DRO models. \cite{rahimian2019distributionally} provide a comprehensive survey of the DRO literature.

\color{black}

The computational tractability of DRO models depends on the ambiguity sets. These sets are often based on moment information \citep{delage2010distributionally, mehrotra2014models,zhang2018ambiguous} or statistical measures such as the Wasserstein distance \citep{esfahani2018data}. To derive tractable DRO models for the MFRSP, we construct two ambiguity sets of demand, one based on 1-Wasserstein distance and one using the demand's support, mean, and mean absolute deviation (MAD). As mentioned in the introduction, we use the Wasserstein ambiguity because it is richer than other divergence-based ambiguity sets. We use the MAD as a dispersion measure instead of the variance because it allows tractable reformulation and better captures outliers and small deviations. In addition, the MAD exists for some distributions while the second moment does not \citep{ben1985approximation, postek2018robust}. We refer to \cite{postek2018robust}  and references therein for rigorous  discussions on properties of MAD.

\color{black}

\textcolor{black}{Next, we discuss some relevant results on the mean-support-MAD ambiguity set (henceforth denoted as  MAD ambiguity).  \cite{ben1972more} derived tight upper and lower bounds on the expectation of a general convex function of a random variable under MAD ambiguity. In particular, when the random variable is one-dimensional, \cite{ben1972more} show that the worst-case distribution under MAD ambiguity is a three-point distribution on the mean, support, and MAD. Recently, \cite{postek2018robust} used the results of \cite{ben1972more} to treat ambiguous expected feasibility constraints to obtain exact reformulations for both functions that are convex and concave in the components of the random variable under MAD ambiguity. These reformulations require independence of the random variables and  involve an exponential number of terms. However, for the special case of linearly aggregated random variables, \cite{postek2018robust} derived  polynomial-sized upper bounds on the worst-case expectations of convex functions.  Finally, under the assumption of independent random variables, they derived tractable approximations of ambiguous chance constraints under mean, support, and MAD information.}

\color{black}

A reviewer of this paper brought our attention to the results in a working paper by \cite{longsupermodularity}  on supermodularity in two-stage DRO problems. Specifically,   \cite{longsupermodularity} identified  a tractable class of two-stage DRO problems based on the scenario-based ambiguity set proposed by \cite{chen2019robust}. They showed that any two-stage DRO problem with mean, support, and upper  bounds of MAD  has a computationally tractable reformulation whenever the second-stage cost function is supermodular in the random parameter. Furthermore, they proposed an algorithm to compute the worst-case distribution for this reformulation.  They argued that using the computed worst-case distribution in the reformulation can make the two-stage DRO problem tractable. In addition, they provided a necessary and sufficient condition to check whether any given two-stage optimization problem has the property of supermodularity. \textcolor{black}{In Appendix~\ref{Appx:Points_MAD}, we show that even if our recourse is supermodular in demand realization, the number of support points in the worst-case distribution of the demand is large, which renders our two-stage MAD-DRO model computationally challenging to solve using \cite{longsupermodularity}'s approach. In contrast, we can efficiently solve an equivalent reformulation of our MAD-DRO model using our proposed decomposition algorithm.}

\color{black}

\begin{figure}[t!]
 \centering
        \includegraphics[width=\textwidth, height=60mm]{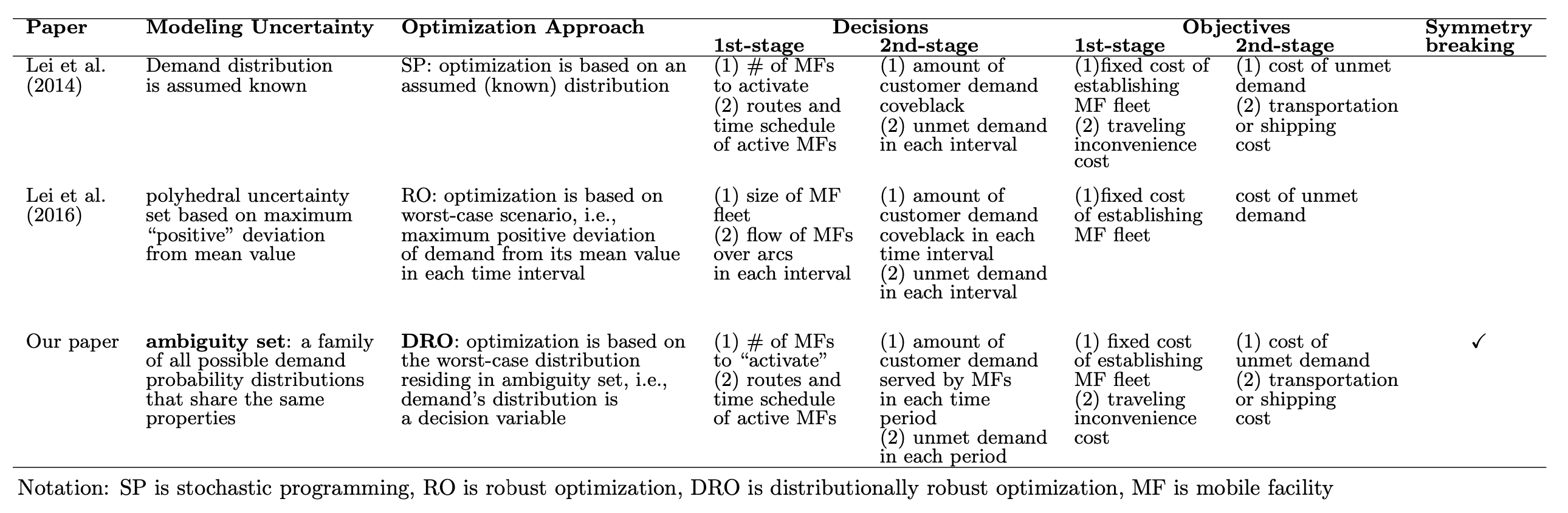}
    \caption{Comparison between \cite{lei2014multicut}, \cite{lei2016two}, and our paper. }\label{Fig:Compare}
\end{figure}

Despite the potential advantages, there are no DRO approaches for the specific  MFRSP that we study in this paper. Therefore, our paper is the first to  propose and analyze DRO approaches for the MFRSP. In Figure~\ref{Fig:Compare}, we provide a comparison between  \cite{lei2014multicut}, \cite{lei2016two}, and our approach based on the assumption made on uncertainty distribution, proposed stochastic optimization approach, decision variables, objectives, and addressing symmetry. We note that our paper and these papers share the common goal of deriving generic optimization models that can be used in any application of MF where one needs to determine the same sets of decisions under the same criteria/objective considered in each paper. 

We make the following observations from Figure~\ref{Fig:Compare}.  In contrast to \cite{lei2016two},  we additionally incorporate the MF traveling inconvenience cost in the first-stage objective and the random transportation cost in the second-stage objective.  In contrast to \cite{lei2016two} and \cite{lei2014multicut}, we model both uncertainty and distributional ambiguity and optimize the system performance over all demand distributions residing within the ambiguity sets. Our master and sub-problems and lower bound inequalities have a different structure than those of \cite{lei2016two}  due to the differences in the decision variables and objectives. We also propose two families of symmetry-breaking constraints, which break symmetries in the solution space of the routing and scheduling decisions. These constraints can improve the solvability of any formulation that uses the same routing and scheduling decisions of the MFRSP. \cite{lei2014multicut} and \cite{lei2016two} did not address the issue of symmetry in the MFRSP.   Finally, to model decision makers' risk-averse attitudes, we propose both risk-neutral (expectation) and risk-averse (mean-CVaR-based) DRO models for the MFRSP.  \cite{lei2014multicut} and \cite{lei2016two} models are risk-neutral. 


Finally, it is worth mentioning that our work uses similar reformulation techniques in recent DRO static FL literature (see, e.g., \cite{basciftci2019distributionally, luo2018distributionally, saif2020data,  ShehadehSanci,  shehadeh2020distributionallyTucker, tsang2021distributionally, wang2020distributionally,wang2021two, wu2015approximation} and references therein).

\section{\textcolor{black}{Problem Setting}}\label{sec3:DMFRSformulation}

\noindent As in \cite{lei2014multicut}, we consider a fleet of $M$ mobile facilities and define the MFRSP on a directed network $G(V, E)$ with node set $V:=\lbrace v_1, \ldots, v_n\rbrace$ and edge set $E:=\lbrace e_1, \ldots, e_m \rbrace$. The sets $I \subseteq V$ and $J \subseteq V$ are the set of all customer points and the subset of nodes where MFs can be located, respectively. The distance matrix $D=(d_{i,j})$ is defined on $E$ and satisfies the triangle inequality, where $d_{i,j}$ is a deterministic and time-invafriant distance between any pair of nodes $i$ and $j$. For simplicity in modeling, we consider a planning horizon of $T$ identical time periods, and we assume that the length of each period $t\in T$ is sufficiently short such that, without loss of generality, all input parameters are the same from one time period to another (this is the same assumption made in \cite{lei2014multicut} and \cite{lei2016two}). The demand, $W_{i,t}$, of each customer $i$ in each time period $t$ is random. The probability distribution of the demand is unknown, and only a possibly small data on the demand may be available. We assume that we know the mean $\mub:=[\mu_{1,1}, \ldots, \mu_{|I|, |T|}]^\top$ and range [$\pmb{\WL}$, $\pmb{\WU}$] of $\Wb$.  Mathematically,  we make the following assumption on the support  of $\Wb$.

\noindent \textbf{Assumption 1.} The support set $\calS$ of $\Wb$ in \eqref{support} is nonempty, convex, and compact.
\begin{align} \label{support}
&  \calS:=\left\{ \Wb  \geq 0: \begin{array}{l} \underline{W}_{i,t} \leq W_{i,t} \leq \overline{W}_{i,t},  \ \forall i \in I, \ \forall t \in T \end{array} \right\}.
\end{align}

\noindent We consider the following basic features as in \cite{lei2014multicut}: (1) each MF has all the necessary service equipment and can move from one place to another; (2) all MFs are homogeneous, providing the same service, and traveling at the same speed; (3) we explicitly account for the travel time of the MF in the model, and service time are only incurred when the MF is not in motion; (4) the travel time $t_{j,j^\prime}$ from location $j$ to $j^\prime$ is an integer multiplier of a single time period (\cite{lei2014multicut, lei2016two}; and (5) the amount of demand to be served is proportional to the duration of the service time at the location serving the demand.

We consider a cost $f$ for using an MF, which represents the expenses associated with purchasing or renting an MF, staffing cost, equipment, etc. Each MF has a capacity limit $C$, which represents the amount of demand that an MF can serve in a single time unit.  Due to the random fluctuations of the demand and the limited capacity of each MF, there is a possibility that the MF fleet will fail to satisfy customers' demand fully.  To minimize shortage, we consider a penalty cost $\gamma$ for each unit of unmet demand.  This penalty cost can represent  the opportunity cost for the loss of demand or expense for outsourcing the excess demand to other companies \citep{basciftci2019distributionally, lei2016two}.  Thus, maximizing demand satisfaction is an important objective that we incorporate in our model \citep{lei2014multicut}.

Given that an MF cannot provide service when in motion, it is not desirable to keep it moving for a long time to avoid losing potential benefits. On the other hand, it is not desirable to keep the MF stationary all the time because this may  lead to losing the potential benefits of making a strategic move to locations with higher demand. Thus, the trade-off of the problem includes the decision to move or keep the MF stationery. Accordingly, we consider a traveling inconvenience cost $\alpha$ to discourage unnecessary moving in cases where moving would neither improve nor degrade the total performance. As in \cite{lei2014multicut}, we assume that $\alpha$ is much lower than other costs such that its impact over the major trade-off is negligible.

We assume that the quality of service a customer receives from a mobile facility  is inversely proportional to the distance between the two to account for the ``access cost'' (this assumption is common in practice and in the literature, see, e.g., \cite{ahmadi2017survey, reilly1931law, drezner2014review, lei2016two, berman2003locating,lei2014multicut}). Accordingly, we consider a demand assignment cost that is linearly proportional to the distance between the customer point and the location of an MF, i.e., $\beta d_{i,j}$, where $\beta\geq 0$ represents the assignment cost factor per demand unit and per distance unit.  Table~\ref{table:notation} summarizes these notation.

Given a set of MFs, $M$, $T$, $I$, and $J$, our models aim to find: (1) the number of MFs to use within $T$; (2) a routing plan and a schedule for the selected MFs, i.e., the node that each MF is located at in each time period; and (3) assignment of MFs to customers. Decisions (1)--(2) are first-stage decisions that we make before realizing the demand. The assignment decisions (3) represent the recourse (second-stage) actions in response to the first-stage decisions \textcolor{black}{and demand realizations} (i.e., you cannot assign demand to the MFs before realizing the demand).  \textcolor{black}{The objective is to minimize  the fixed cost (i.e., cost of establishing the MF fleet and traveling inconvenience cost)  plus a risk measure (expectation or mean CVaR) of the operational cost (i.e., transportation and unsatisfied demand costs).}   We refer to \cite{lei2014multicut} for an excellent visual representation of MF operations and some of the basic features mentioned above.

\noindent \textbf{\textit{Additional Notation}}: For $a,b \in \mathbb{Z}$, we define $[a]:=\lbrace1,2,\ldots, a \rbrace$ as the set of positive integer running indices to $a$. Similarly, we define $[a,b]_\mathbb{Z}:=\lbrace c \in \mathbb{Z}: a \leq c \leq b \rbrace$ as the set of positive running indices from $a$ to $b$. The abbreviations ``w.l.o.g.'' and ``w.l.o.o.'' respectively represent ``without loss of generality'' and ``without loss of optimality.''

\section{\textbf{Stochastic Programming Model}}\label{sec:SP}

\noindent \textcolor{black}{In this section, we present a two-stage SP formulation of the MFRSP that assumes that the probability distribution of the demand is known. A complete listing of the parameters and decision variables of the model can be found in Table~\ref{table:notation}.}

\color{black}
First, let us introduce the variables and constraints defining the first-stage of this SP model. For each $m \in M$, we define a binary decision variable $y_m$ that equals 1 if MF $m$ is used, and is 0 otherwise.  For all $j \in J$, $m \in M$, and $t \in T$, we define a binary decision variable $x_{j,m}^t$ that equals 1 if MF $m$ stays at location $j$ at period $t$, and is 0 otherwise. The feasible region $\mathcal{X}$ of variables $\xbb$ and $\yb$ is defined in \eqref{eq:RegionX}. 
\color{black}
 \begin{align}
\mathcal{X}&=\left\{ (\xbb, \yb) :  \begin{array}{l} x_{j,m}^t+x_{j^\prime,m}^{t^\prime} \leq y_m, \ \ \forall t, m, j, \ j\neq j^\prime, \ t^\prime \in \lbrace t, \ldots, \min \lbrace t+t_{j,j^\prime}, T\rbrace \rbrace \\
  x_{j,m}^t \in \lbrace 0, 1\rbrace, \ y_m \in \lbrace 0, 1\rbrace,  \ \forall j, m , t \end{array} \right\}.  \label{eq:RegionX}
\end{align} 
As defined in \cite{lei2014multicut}, region $\calX$ represents: (1) the requirement that an MF can only be in service when it is stationary; (2) MF $m$ at location $j$ in period $t$ can only be available at location $j^\prime\neq j$ after a certain period of time depending on the time it takes to travel from location $j$ to $j^\prime$, $t_{j,j^\prime}$ (i.e., $x_{j,m}^{t^\prime}=0$ for all $j^\prime \neq j$ and $t^\prime \in \lbrace t, \ldots, \min \lbrace t+t_{j,j^\prime}, T\rbrace \rbrace$); and (3) MF $m$ has to be in an active condition before providing service. We refer the reader to Appendix~\ref{Appx:FirstStageDec} for a detailed derivation of region $\mathcal{X}$.

\color{black}
Now, let us introduce the variables defining our second-stage problem. For all $(i,\ j, \ m, \ t)$, we define a nonnegative continuous variable  $z_{i,j,m}^t$ to represent the amount of node $i$ demand served by MF $m$ located at $j$ in period $t$. For each $t \in T$, we define a nonnegative continuous variable $u_{i,t}$ to represent the amount of unmet demand of node \textit{i} in period $t$.   Finally, we define a random vector $\Wb:=[W_{1,1}. \ldots, W_{|I|,|T|}]^\top$. Our SP model can now be stated as follows:
\color{black}
\begin{table}[t]  
\small
\center
   \renewcommand{\arraystretch}{0.6}
  \caption{Notation.} 
\begin{tabular}{ll}
\hline
\multicolumn{2}{l}{\textbf{Indices}} \\
$m$& index of MF, $m=1,\ldots, M$\\
$i$ & index of customer location, $i=1,\ldots,I$\\
$j$ & index of MF location, $j=1,...,J$\\
\multicolumn{2}{l}{\textbf{Parameters and sets}} \\
$T$& planning horizon\\
$M$ & number, or set, of MFs \\
$J$& number, or, set of locations \\
$f$ & fixed operating cost  \\
$d_{i,j}$ & distance between any pair of nodes $i$ and $j$ \\
$t_{j,j'}$ & travel time from $j$ to $j'$\\
$C$ & the amount of demand that can be served by an MF in a single time unit, i.e., MF capacity \\
$W_{i,t}$& demand at customer site  $i$ for each period $t$\\
$\underline{W}_{i,t}/\overline{W}_{i,t}$ & lower/upper bound of demand at customer location $i$ for each period $t$ \\
$\gamma$ &  penalty cost for each unit of unmet demand\\
 \multicolumn{2}{l}{\textbf{First-stage decision variables } } \\
$y_{m}$ &   $\left\{\begin{array}{ll}
1, & \mbox{if MF \textit{m} is permitted  to use,} \\
0, & \mbox{otherwise.}
\end{array}\right.$ \\
$x_{jm}^t$ &   $\left\{\begin{array}{ll}
1, & \mbox{if MF \textit{m} stays at location \textit{j} at period \textit{t},} \\
0, & \mbox{otherwise.}
\end{array}\right.$ \\
 \multicolumn{2}{l}{\textbf{Second-stage decision variables } } \\
$z_{i,j,m}^t$ & amount of demand of node $i$ being served by MF $m$ located at $j$ in period $t$  \\
$u_{i,t}$ & total amount of unmet demand of node \textit{i} in period $t$ \\
\hline
\end{tabular}\label{table:notation}
\end{table} 
\begin{subequations}\label{DMFRS_SP}
\begin{align}
(\text{SP}) \quad Z^*=& \min \limits_{(\xbb, \yb)\in \mathcal{X}}   \left\lbrace \sum_{m \in M} f y_{m} -\sum_{t \in T} \sum_{j \in J} \sum_{m \in M } \alpha x_{j,m}^t+ \textcolor{black}{\varrho
\big(Q(\xbb, \Wb)\big)} \ \right\rbrace, \label{ObjSP}
\end{align}
\end{subequations} 
where for a feasible $(\xbb, \yb) \in \calX$ and a realization of $\Wb$
\allowdisplaybreaks
\begin{subequations}\label{2ndstage}
\begin{align}
  Q(\xbb, {\Wb} ):=  & \min_{\pmb{z, u}} \Big (\sum_{j \in J} \sum_{i \in I} \sum_{m \in M} \sum_{t \in T} \beta d_{i,j} z_{i,j,m}^t+ \gamma \sum_{t \in T} \sum_{i \in I} u_{i,t}\Big)\label{2ndstageObj}\\
&  \ \  \text{s.t.} \ \ \  \  \sum_{j \in J} \sum_{m \in M} z_{i,j,m}^t+u_{i,t}=  W_{i,t}, \qquad  \forall i \in I, \ t \in T, \label{2ndstage_const1}\\
&\quad  \   \quad  \quad  \sum_{i \in I} z_{i,j,m}^t \leq C x_{j,m}^t, \qquad  \forall j \in J, \ m \in M, \ t \in T, \label{2ndstage_const2}\\
&\quad  \   \quad  \quad u_t \geq 0, \ z_{i,j,m}^t \geq 0, \qquad   \forall i \in I, \ j \in J, m \in M, \ t \in T. \label{2ndstage_const3}
\end{align} 
\end{subequations}
\noindent \textcolor{black}{Formulation \eqref{DMFRS_SP} aims to find  first-stage decisions $(\xbb,\yb) \in \calX$ that minimize the sum of (1) the fixed cost of establishing the MF fleet (first term); (2)  the traveling inconvenience cost\footnote{Minimizing the traveling inconvenience cost is equivalent to maximizing the profit of keeping the MF stationary whenever possible. Parameter $\alpha$ is the profit weight factor as detailed in \cite{lei2014multicut})} (second term); and (3) a risk measure $\varrho (\cdot) $ of the random second-stage function $Q(\xbb, {\Wb} )$ (third term).  A risk-neutral decision-maker may opt to set $\varrho(\cdot)=\E(\cdot)$, whereas a risk-averse decision-maker might set $\varrho(\cdot)$ as the CVaR or mean-CVaR.  Classically, the MFRSP literature 
employs $\varrho(\cdot)=\E(\cdot)$,  which might be more intuitive for MF providers. For brevity, we relegate further discussion of the mean-CVaR-based SP model to Appendix~\ref{Appex:SP_CVAR}.}

\section{\textcolor{black}{Distributionally Robust Optimization (DRO) Models}}\label{sec:DRO_Models}

\noindent In this section, we present our proposed DRO models for the MFRSP that do not assume that the probability distribution of the demand is known. In Sections~\ref{sec:MAD-DRO} and \ref{sec:WDMFRS_model}, we respectively present and analyze the risk-neutral MAD-DRO and W-DRO models. For brevity, we relegate the formulations and discussions of the risk-averse mean-CVaR-based DRO models to Appendix~\ref{Sec:meanCVAR}.


\subsection{\textbf{The DRO Model with MAD Ambiguity (MAD-DRO)}}\label{sec:MAD-DRO}

\noindent  \textcolor{black}{In this section, we present our proposed MAD-DRO model, which is based on an ambiguity set that incorporates the demand's mean ($\mub$), MAD ($\etab$), and support ($\calS$)}. As mentioned earlier, we use the MAD as a dispersion or variability measure because it allows us to derive a computationally attractive reformulation   \citep{postek2018robust, wang2019distributionally, wang2020distributionally}.

\textcolor{black}{First, let us introduce some additional sets and notation defining our MAD ambiguity set and MAD-DRO model}. We define $\E_\Prob$ as the expectation under distribution $\Prob$. We let $\mu_{i,t}=\E_\Prob[W_{i,t}]$ and $\eta_{i,t}=\E_{\Prob}(|W_{i,t}-\mu_{i,t}|)$ respectively represent the mean value and MAD of  $W_{i,t}$, for all $i \in I$ and $t \in T$. \textcolor{black}{Using this notation and the support $\calS$ defined in \eqref{support}, we construct the following MAD ambiguity set}:
\begin{align}\label{eq:ambiguityMAD}
\calF(\calS, \mub, \etab) := \left\{ \Prob\in \calP(\calS) \middle|\:  \begin{array}{l} \Prob (\Wb \in \calS)=1,\\ \mathbb{E_P}[W_{i,t}] = \mu_{i,t}, \ \forall i \in I, \ t \in T, \\
\E_{\Prob}(|W_{i,t}-\mu_{i,t}|) \leq \eta_{i,t}, \ \forall i \in I, \ t \in T.
 \end{array} \right\},
\end{align}

\noindent where $\calP(S)$ in $\calF(\calS, \mub, \etab)$ represents the set of distributions supported on $\calS$ with mean $\pmb \mu$ and dispersion measure $ \leq \pmb \eta$.   Using $\calF(\calS, \mub, \etab)$ defined in \eqref{eq:ambiguityMAD}, we formulate our MAD-DRO model as
\begin{align}
(\text{MAD-DRO}) & \min \limits_{(\xbb, \yb)\in \mathcal{X}}   \left\lbrace \sum_{m \in M} f y_{m} -\sum_{t \in T} \sum_{j \in J} \sum_{m \in M } \alpha x_{j,m}^t+\ \Bigg[\sup \limits_{\Prob \in \calF(\calS, \mub, \etab)} \E_\Prob[Q(\xbb,{\Wb})]  \Bigg] \ \right\rbrace. \label{MAD-DMFRS}
\end{align}
\noindent The MAD-DRO formulation in \eqref{MAD-DMFRS} seeks first-stage decisions ($\xbb,\yb)$ that minimize the first-stage cost and the worst-case expectation of the second-stage (recourse) cost, \textcolor{black}{where the expectation is taken over all distributions
residing in $\calF(\calS, \mub, \etab)$}. Note that we do not incorporate higher moments (e.g., co-variance) in $\calF(\calS, \mub, \etab)$ for the following primary reasons. First, the mean and range are intuitive statistics that a decision-maker may approximate and change in the model (e.g.,  the mean may be estimated from limited data or approximated by subject matter experts, and the range may represent the error margin in the estimates).  Second, it is not straightforward for decision-makers to approximate or accurately estimate the correlation between uncertain parameters.   Third, mathematically speaking,  various studies have demonstrated that incorporating higher moments in the ambiguity set often undermines the computational tractability of DRO models and, therefore, their applicability in practice.  Indeed, as we will show next, using  $\calF(\calS, \mub, \etab)$ allows us to derive a tractable equivalent reformulation of the MAD-DRO model and an efficient solution method to solve the reformulation (see Sections~\ref{sec3:reform}, \ref{sec:Alg}, and \ref{sec5:CPU}).

\textcolor{black}{Finally, note that parameters $W_{i,t}$, $\mu_{i,t}$, $\WL_{i,t}$, $\WU_{i,t}$, $\eta_{i,t}$ are all indexed by time period $t$ and location $i$. Thus, if in any application, there is a relationship (e.g., correlation) between the demand of a subset of locations in a subset of periods, one can easily adjust $\calS$, $\mub$, and $\etab$ in $\calF(\calS, \mub, \etab)$ to reflect this relationship. For example, if urban cities have higher demand, then we can adjust the mean and range of the demand of  these cities to reflect such a relationship. Similarly, if there is a correlation between the time period $t$ and the demand, then we can define the mean and the support based on this correlation. For example, the morning service hours may  have lower demand on Monday.  In this case, we can adjust  $\mub$, $\calS$, and $\etab$ to reflect this relationship. Similarly, if the demand in a subset of periods and locations is correlated, we can adjust $\calS$, $\mub$, and $\etab$ in $\calF(\calS, \mub, \etab)$ to reflect this relationship}. Nevertheless, we acknowledge that not incorporating higher moments and complex relations may be a limitation of our work and thus is worth future investigation.

\subsubsection{\textbf{Reformulation of the MAD-DRO model}}\label{sec3:reform}
\textcolor{white}{ }

\textcolor{black}{Recall that $Q(\xbb,{\Wb})$ is defined by a minimization problem; hence, in \eqref{MAD-DMFRS}, we have an inner max-min problem. As such, it is not straightforward to solve formulation \eqref{MAD-DMFRS} in its presented form}. In this section, we derive an equivalent formulation of \eqref{MAD-DMFRS} that is solvable. \textcolor{black}{First,  in Proposition~\ref{Prop1:DualMinMax}, we present an equivalent reformulation of the inner problem in \eqref{MAD-DMFRS} (see Appendix \ref{Appx:ProofofProp1} for a proof).}

\begin{proposition}\label{Prop1:DualMinMax}
For any fixed $(\xbb, \yb) \in \mathcal{X}$,  problem $\sup \limits_{\Prob \in \calF(\calS, \mub, \etab)} \E_\Prob[Q(\xbb,{\Wb})]$ in \eqref{MAD-DMFRS} is equivalent to 
\begin{align}  \label{eq:FinalDualInnerMax-1}
& \min_{\pmb{\rho}, \pmb{\psi} \geq 0}  \ \left \lbrace  \sum \limits_{t \in T} \sum \limits_{i \in I} ( \mu_{i,t} \rho_{i,t}+\eta_{i,t} \psi_{i,t})+\max \limits_{\pmb{W}\in \calS} \Big\lbrace Q(\xbb,{\Wb}) +   \sum \limits_{t \in T} \sum \limits_{i \in I} -(W_{i,t} \rho_{i,t}+ |W_{i,t}-\mu_{i,t}| \psi_{i,t})\Big \rbrace  \right  \rbrace.
\end{align}
\end{proposition}

\noindent \textcolor{black}{Again, the problem in \eqref{eq:FinalDualInnerMax-1} involves an inner max-min problem that is not straightforward to solve in its presented form. However, we next derive an equivalent reformulation of the inner problem in \eqref{eq:FinalDualInnerMax-1}  that is solvable.} 
First, we observe that $Q(\xbb,{\Wb})$ is a feasible linear program  (LP) for a given first-stage solution $(\xbb, \yb) \in \calX$ and a realization of ${\Wb}$. The dual of $Q(\xbb, {\Wb})$ is as follows.
\begin{subequations}\label{DualOfQ}
\begin{align}
Q(\xbb,{\Wb})=& \max \limits_{\pmb{\lambda, v}} \ \sum_{t \in T} \sum_{i \in I} \lambda_{i,t} W_{i,t}+ \sum_{t \in T} \sum_{j \in J} \sum_{m \in M} Cx_{j,m}^t v_{j,m}^t \label{Obj:Qdual}\\
& \ \text{s.t. } \ \lambda_{i,t}+ v_{j,m}^t \leq  \beta d_{i,j}, \ \qquad  \forall i \in I,   j \in J,  m \in M,   t \in T, \label{Const1:QudalX}\\
& \qquad \ \ \lambda_{i,t} \leq \gamma,  \qquad    \qquad    \qquad  \ \  \forall i \in I, \  t \in T, \label{Const2:QdualU}\\
& \qquad \ \ v_{j,m}^t \leq 0,   \qquad    \qquad    \qquad  \ \  \forall j \in J, \ t \in T, \label{Const3:Qdualv}
\end{align}
\end{subequations}
\noindent where $\pmb{\lambda}$ and $\pmb{v}$ are the dual variables associated with constraints \eqref{2ndstage_const1} and \eqref{2ndstage_const2}, respectively. It is easy to see that w.l.o.o, $\lambda_{i,t} \geq 0$ for all $i \in I$ due to constraints \eqref{Const2:QdualU} and the objective of maximizing $W_{i,t}\geq 0$ times $\lambda_{i,t}$. Additionally, $ v_{j,m}^t \leq \min\{  \min_i\{\beta d_{i,j}-\lambda_{i,t}\},0  \}$ by constraints \eqref{Const1:QudalX} and \eqref{Const3:Qdualv}. Given the objective of maximizing a nonnegative term $Cx_{j,m}^t $ multiplied by $v_{j,m}^t$,  $ v_{j,m}^t = \min\{  \min_i\{\beta d_{i,j}-\lambda_{i,t}\},0  \}$ in the optimal solution. Given that $\pmb{\beta}$, $\pmb d$, and $\pmb \lambda$ are finite,  $\pmb{v} $ is finite. It follows that problem \eqref{DualOfQ} is a feasible and bounded LP. Note that $W_{i,t} \in [\underline{W}_{i,t}, \overline{W}_{i,t}]$ and ($\underline{W}_{i,t}, \overline{W}_{i,t})\geq 0$ by definition, in view of dual formulation \eqref{DualOfQ},  we can rewrite the inner maximization problem $\max \{ \cdot \} $ in \eqref{eq:FinalDualInnerMax-1} as
\begin{subequations}\label{DualQAndOuter}
\begin{align}
 \max \limits_{\pmb{\lambda, v, W, k}} & \Bigg\{ \sum_{t \in T} \sum_{i \in I} \lambda_{i,t} W_{i,t}+ \sum_{t \in T} \sum_{j \in J} \sum_{m \in M} Cx_{j,m}^t v_{j,m}^t +\sum \limits_{t \in T} \sum \limits_{i \in I} -(W_{i,t} \rho_{i,t}+k_{i,t}\psi_{i,t}) \Bigg \} \\
  \ \text{s.t. }  &\eqref{Const1:QudalX}-\eqref{Const3:Qdualv}, \ \ W_{i,t} \in [\underline{W}_{i,t}, \overline{W}_{i,t}], \  \forall i \in I, \ \forall t \in T, \\
& k_{i,t} \geq W_{i,t}-\mu_{i,t}, \ k_{i,t} \geq \mu_{i,t}-W_{i,t},  \ \forall i \in I, \ \forall t \in T . \label{Const1:DualQAndOuter}
\end{align}
\end{subequations}
Note that the objective function in \eqref{DualQAndOuter} contains the interaction term $\lambda_{i,t} W_{i,t}$. To linearize formulation \eqref{DualQAndOuter}, we define $\pi_{i,t}=\lambda_{i,t} W_{i,t}$ for all $i \in I$ and $t \in T$. Also, we introduce the following McCormick inequalities for variables $\pi_{i,t}$:
\begin{subequations}
\begin{align}
&\pi_{i,t} \geq \lambda_{i,t} \underline{W}_{i,t},  \qquad \qquad \qquad   \ \ \   \pi_{i,t} \leq \lambda_{i,t} \overline{W}_{i,t}, \ \ \forall i \in I, \ \forall t \in T, \label{McCormick1}\\
& \pi_{i,t} \geq \gamma W_{i,t}+  \overline{W}_{i,t} (\lambda_{i,t}- \gamma ), \qquad  \pi_{i,t} \leq \gamma W_{i,t}+\underline{W}_{i,t}( \lambda_{i,t}-\gamma), \ \forall i \in I, \ \forall t \in T. \label{McCormick2}
\end{align}
\end{subequations}
 Accordingly, for a fixed $(\xbb \in \mathcal{X}, \pmb{\rho}, \pmb{\psi})$, problem \eqref{DualQAndOuter} is equivalent to the following mixed-integer linear program (MILP):
\begin{subequations}\label{eq:FinalDualInnerMax-2}
\begin{align}
\max \limits_{\pmb{\lambda, v, W, \pi, k} }&   \ \Bigg\{ \sum_{t \in T} \sum_{i \in I} \pi_{i,t}+ \sum_{t \in T} \sum_{j \in J} \sum_{m \in M} Cx_{j,m}^t v_{j,m}^t +\sum \limits_{t \in T} \sum \limits_{i \in I} -(W_{i,t} \rho_{i,t}+k_{i,t}\psi_{i,t})  \Bigg\}\\
\text{s.t. } &   \ (\pmb{\lambda, v}) \in \{\eqref{Const1:QudalX}-\eqref{Const3:Qdualv}\}, \ \pmb{\pi} \in \{ \eqref{McCormick1}-\eqref{McCormick2} \}, \label{Const1_finalinner} \\
& \ W_{i,t} \in [\underline{W}_{i,t}, \overline{W}_{i,t}], \ k_{i,t} \geq W_{i,t}-\mu_{i,t}, \ k_{i,t} \geq \mu_{i,t}-W_{i,t},  \ \forall i \in I, \ \forall t \in T. \label{Const2_finalinner}
\end{align}
\end{subequations}
Combining the inner problem in the form of  \eqref{eq:FinalDualInnerMax-2} with the outer minimization problems in \eqref{eq:FinalDualInnerMax-1} and  \eqref{MAD-DMFRS}, we derive the following equivalent  reformulation of  the MAD-DRO model in \eqref{MAD-DMFRS}:
\begin{subequations}\label{FinalDR}
\begin{align}
&\min_{\pmb{x,y,\rho, \psi, \delta}}  \Bigg \{ \sum_{m \in M} f y_{m} -\sum_{t \in T} \sum_{j \in J} \sum_{m \in M } \alpha x_{j,m}^t+  \sum \limits_{t \in T} \sum \limits_{i \in I} \Big[\mu_{i,t} \rho_{i,t}+\eta_{i,t} \psi_{i,t}\Big] + \delta  \Bigg \}  \\
& \ \ \text{s. t.} \ \ (\xbb, \yb) \in \mathcal{X}, \ \psi \geq 0, \\
& \qquad  \ \ \delta \geq \textcolor{black}{h(\xbb, \pmb \rho, \pmb \psi)}, \label{const1:FinalDR}
 \end{align}
\end{subequations}
where $\textcolor{black}{h(\xbb, \pmb \rho, \pmb \psi)}=\max \limits_{\pmb{\lambda, v, W, \pi, k} } \ \Big\{ \sum \limits_{t \in T} \sum\limits_{i \in I} \pi_{i,t}+ \sum\limits_{t \in T} \sum\limits_{j \in J} \sum \limits_{m \in M} Cx_{j,m}^t v_{j,m}^t +\sum \limits_{t \in T} \sum \limits_{i \in I} -(W_{i,t} \rho_{i,t}+k_{i,t}\psi_{i,t}):  \eqref{Const1:QudalX}-\eqref{Const3:Qdualv}, \eqref{McCormick1}-\eqref{McCormick2}, \eqref{Const2_finalinner} \Big\} $.

\begin{proposition}\label{Prop2:Convexh} \textcolor{black}{For any fixed values of variables $(\xbb, \yb) \in \calX$, $\pmb \rho,$ and $\pmb \psi$, $h(\xbb, \pmb \rho, \pmb \psi)< +\infty$. Furthermore, function $(\xbb, \pmb \rho, \pmb \psi) \mapsto\textcolor{black}{h(\xbb, \pmb \rho, \pmb \psi)}$ is a convex piecewise linear function in $\xbb$, $\pmb\rho$, and $\pmb \psi$ with a finite number of pieces (see Appendix \ref{Appx:Prop2} for a detailed proof).}
\end{proposition}
\vspace{1mm}

\subsection{\textit{\textbf{The DRO Model with 1-Wasserstein Ambiguity (W-DRO)}}}\label{sec:WDMFRS_model}

\noindent In this section, we consider the case that $\Prob$ may be observed via a possibly small finite set $\{\hat{\Wb}^1, \ldots, \hat{\Wb}^N\}$ of $N$ i.i.d. samples, which may come from the limited historical realizations of the demand or a reference empirical distribution. Accordingly, we construct an ambiguity set based on 1-Wasserstein distance, which often admits tractable reformulation in most real-life applications (see, e.g., \cite{Daniel2020, hanasusanto2018conic, jiang2019data, tsang2021distributionally, saif2020data}). 

\textcolor{black}{First, let us define the 1-Wasserstein distance.} Suppose that random vectors $\pmb{\Wb}_1$ and $\pmb{\Wb}_2$ follow $\F_1$ and $\F_2$, respectively, where probability distributions $\F_1$ and $\F_2$ are defined over the common support $\calS$. The 1-Wasserstein distance \textcolor{black}{$\text{dist}(\F_1, \F_2)$ between $\F_1$ and $\F_2$ is the minimum  transportation  cost of moving from  $\F_1$ to $\F_2$, where the cost of moving masses $\Wb_1$ to $\Wb_2$ is the norm $||\Wb_1-\Wb_2||$}. \textcolor{black}{ Mathematically,}
\begin{equation}\label{W_distance}
 \textcolor{black}{\text{dist}}(\mathbb{F}_1, \mathbb{F}_2):=\inf_{\Pi \in \calP(\mathbb{F}_1, \mathbb{F}_2)}  \Bigg \{\int_{\calS} ||\Wb_1-\Wb_2|| \ \Pi(\text{d}\Wb_1, \text{d}\Wb_2) \Bigg | \begin{array}{ll} & \Pi \text{ is a joint distribution of }\Wb_1 \text{ and } \Wb_2\\
& \text{with marginals } \mathbb{F}_1 \text{ and }  \mathbb{F}_2, \text{ respectively} 
 \end{array} \Bigg \},
\end{equation}

\noindent where $\calP(\mathbb{F}_1, \mathbb{F}_2) $  is the set of all joint distributions of $(\Wb_1, \Wb_2)$ supported on $\calS$ with  marginals ($\mathbb{F}_1$, $\mathbb{F}_2$).  Accordingly, we construct the following $1$-Wasserstein ambiguity set:
\begin{align}\label{W-ambiguity}
\calF (\hat{\Prob}^N, \epsilon)= \left\{ \Prob \in \calP(\calS) \middle|\begin{array}{l} \textcolor{black}{\text{dist}}(\Prob, \hat{\Prob}^N) \leq \epsilon\end{array} \right\},
\end{align}
\noindent where $ \calP(\calS)$ is the set of all distributions supported on $\calS$, $\hat{\Prob}^N=\frac{1}{N} \sum_{n=1}^N \delta_{\hat \Wb^n}$ is the empirical distribution of $\Wb$ based on $N$ i.i.d samples, and $\epsilon >0$ is the radius of the ambiguity set. The set $\calF (\hat{\Prob}^N, \epsilon)$ represents a Wasserstein ball of radius $\epsilon$ centered at the empirical distribution $\hat{\Prob}^N$. \textcolor{black}{Using  the ambiguity set $\calF(\hat{\Prob}^N, \epsilon)$ defined  in \eqref{W-ambiguity}, we formulate our W-DRO model as}
\begin{align}
  (\text{W-DRO})   \ \hat{Z}(N, \epsilon)&= \min \limits_{(\xbb, \yb)\in \mathcal{X}}   \Bigg\lbrace \sum_{m \in M} f y_{m} -\sum_{t \in T} \sum_{j \in J} \sum_{m \in M } \alpha x_{j,m}^t+\Bigg [\sup_{\Prob \in \calF (\hat{\Prob}^N, \epsilon) }\E_{\Prob} [Q(\xbb, \Wb) ] \Bigg] \Bigg \}. \label{W-DMFRS}
\end{align}
\noindent Formulation \eqref{W-DMFRS} finds first-stage decisions ($\xbb,\yb) \in \mathcal{X}$ that minimize the first-stage cost and the maximum expectation of the second-stage cost over all distributions residing in $\calF (\hat{\Prob}^N, \epsilon)$. 

The W-DRO model in \eqref{W-DMFRS} can be used to model uncertainty in general and distributional ambiguity when there is a possibly small finite data sample on uncertainty. As detailed in \cite{esfahani2018data} and discussed earlier, if we have a small sample and we optimize using this sample, then the optimizer's curse cannot be avoided. To mitigate the optimizer's curse (estimation error), we robustify the nominal decision problem (the MFRSP optimization problem) against all distributions $\Prob$ under which the estimated  distribution $\hat{\Prob}^N$ based on the $N$ data points has a small estimation error (i.e., with $ \textcolor{black}{\text{dist}}(\Prob, \hat{\Prob}^N) \leq \epsilon$).  Therefore, in some sense one can think of Wasserstein ball as the set of all distributions under which our estimation error is below $\epsilon$, where $\epsilon$ is the maximum estimation error against which we seek protection. When $\epsilon=0$, the ambiguity set contains the empirical distribution and the W-DRO problem in \eqref{W-DMFRS} reduces to the SP problem.  A larger radius $\epsilon$ indicates that we seek more robust solutions  (see Appendix \ref{WDRO_Radius}).

In the next section, we show that using $\ell_1$-norm instead of the \textcolor{black}{$\ell_p$}-norm  ($1<p<\infty$) in our  Wasserstein ambiguity set allows us to derive a linear and tractable reformulation of the W-DRO model in \eqref{W-DMFRS}. Note that $\ell_1$-norm (i.e., the sum of the magnitudes of the vectors in space) is the most intuitive and natural way to measure the distance between vectors.  In contrast, the \textcolor{black}{$\infty$-norm}-based Wasserstein ball is an extreme case. That is, the $\infty$-norm gives the largest magnitude among each element of a vector. Thus, from the perspective of the Wasserstein DRO framework, the \textcolor{black}{$\infty$-norm}-based distance metric only picks the most influential value to determine the closeness between data points \citep{chen2018robust}, which, in our case, may not be reasonable since every demand point plays a role. Deriving and comparing DRO models with different Wasserstein sets is out of the scope of this paper but is worth future investigation in more comprehensive MF optimization problems.

\subsubsection{\textbf{\textit{Reformulation of the W-DRO model}}}\label{sec:Wreformualtion}
\textcolor{white}{ }

\noindent In this section, we derive an equivalent solvable reformulation of the W-DRO model in \eqref{W-DMFRS}.  First,  in Proposition~\ref{Prop1}  we present an equivalent dual formulation of  the inner maximization problem $\sup\  [\cdot ]$ in \eqref{W-DMFRS} (see  Appendix~\ref{Proof_Prop1} for a detailed proof).
 \begin{proposition}\label{Prop1} 
For for a fixed $(\xbb, \yb) \in \calX$, problem $\sup \limits_{\Prob\in \calF (\hat{\Prob}^N, \epsilon) }\E_{\Prob} [Q(\xbb,\Wb) ]$ in \eqref{W-DMFRS} is equivalent to
\begin{align}
& \inf_{\rho \geq 0} \Bigg\{ \epsilon \rho + \Big[ \frac{1}{N} \sum_{n=1}^N \sup_{\Wb \in \calS}  \big \{ Q(\xbb, \Wb)-\rho || \Wb -\Wb^n || \big \} \Big] \Bigg\}. \label{DualOfInner}
\end{align}
\end{proposition}
\noindent Formulation \eqref{DualOfInner} is potentially challenging to solve because it require solving $N$ non-convex optimization problems. Fortunately, given that the support of $\Wb$ is rectangular and finite (see Assumption 1) and $Q(\xbb,\Wb)$ is feasible and bounded for every $\xbb$ and $\Wb	$, we next recast these inner problems as LPs for each $\rho \geq 0$ and $\xbb \in \calX$. First, using the dual formulation of $Q(\xbb, \Wb)$ in \eqref{DualOfQ}, we rewrite the inner problem $\sup \{ \cdot  \} $ in \eqref{DualOfInner} for each $n$ as
\begin{subequations}\label{InnerInnerW}
\begin{align}
 & \max \limits_{\lambdab, \vb, \Wb} \Big \{ \sum_{t \in T} \sum_{i \in I} \lambda_{i,t} W_{i,t}-\rho |W_{i,t}-W_{i,t}^n|+ \sum_{t \in T} \sum_{j \in J} \sum_{m \in M} Cx_{j,m}^t v_{j,m}^t \Big \}\\
& \text{s.t. } \ \ (\lambdab, \vb) \in \{  \eqref{Const1:QudalX}-\eqref{Const3:Qdualv} \}, \Wb \in [\pmb{\WL}, \pmb{\WU} ].
\end{align}
\end{subequations}
Second, using the same techniques in Section~\ref{sec3:reform}, we define an epigraphical random variable $\eta_{i,t}^n$ for the term $ |W_{i,t}-W_{i,t}^n|$. Then, using variables $\pmb{\eta}$, $\pi_{i,t}=\lambda_{i,t} W_{i,t}$, and inequalities \eqref{McCormick1}-\eqref{McCormick2} for variables $\pi_{i,t}$, we derive the following equivalent reformulation of \eqref{InnerInnerW} (for each $n \in [N]$):
\begin{subequations}\label{InnerInnerW2}
\begin{align}
 \max \limits_{\lambdab, \vb, \Wb, \pib, \pmb{\eta}}   &\big \{ \sum_{t \in T} \sum_{i \in I}  \pi_{i,t}-\rho  \eta_{i,t}^n+ \sum_{t \in T}\sum_{j \in J} \sum_{m \in M} Cx_{j,m}^t v_{j,m}^t \big \} \label{Obj:InnerInnerW2}\\
\text{s.t. } &  \ (\pmb{\lambda, v}) \in \{  \eqref{Const1:QudalX}-\eqref{Const3:Qdualv} \}, \label{C1:InnerInnerW2} \\
& \ \ \pi_{i,t} \in \{ \eqref{McCormick1}- \eqref{McCormick1}\},  \ \  W_{i,t} \in [\WL_{i,t}, \WU_{i,t}] , && \forall i \in I , \forall t \in T, \label{C2:InnerInnerW2} \\
&   \ \   \eta_{i,t}^n \geq W_{i,t}-W_{i,t}^n,  \  \eta_{i,t}^n \geq W_{i,t}^n- W_{i,t}, && \forall i \in I , \forall t \in T. \label{C_InnerInner}
\end{align}
\end{subequations}
Third, combining the inner problem in the form of \eqref{InnerInnerW2} with the outer minimization problems in \eqref{DualOfInner} and \eqref{W-DMFRS}, we derive the following equivalent reformulation of the W-DRO model in \eqref{W-DMFRS}
\begin{align}
\hat{Z}(N, \epsilon) &=  \min_{ (\xbb, \yb) \in \mathcal{X}, \ \rho  \geq 0 }  \Bigg\lbrace \sum_{m \in M} f y_{m} -\sum_{t \in T} \sum_{j \in J} \sum_{m \in M } \alpha x_{j,m}^t+  \epsilon \rho  \nonumber\\
  & \quad  \ \ \  \ +\frac{1}{N} \sum_{n=1}^N  \max \limits_{\lambdab, \vb, \Wb, \pib, \pmb{\eta}} \Big \{ \sum_{t \in T} \sum_{i \in I}  \pi_{i,t}- \rho  \eta_{i,t}^n+ \sum_{t \in T}\sum_{j \in J} \sum_{m \in M} Cx_{j,m}^t v_{j,m}^t: \eqref{C1:InnerInnerW2}-\eqref{C_InnerInner} \Big \} \Bigg\}.\label{Final_W_DMFRS}
\end{align}

\noindent Using the same techniques in the proof of Proposition~\ref{Prop2:Convexh}, one can easily  verify that function $[\max \limits_{\pmb{\lambda, v, \pi, W, \eta}} \big \{ \sum\limits_{t \in T} \sum\limits_{i \in I} \pi_{i,t}-\rho  \eta_{i,t}^n+ \sum \limits_{t \in T}\sum\limits_{j \in J} \sum\limits_{m \in M} Cx_{j,m}^t v_{j,m}^t \big\}]<\infty $ and is a convex piecewise linear function in $\xbb \in \calX$ and $\rho$.

\section{Solution Method}\label{sec:solutionmethods}
\noindent In this section, we present a decomposition-based algorithm to solve the MAD-DRO formulation in \eqref{FinalDR}, and strategies to improve the solvability of the formulation. The algorithmic steps for solving the W-DRO in \eqref{Final_W_DMFRS} are similar. In Section~\ref{sec:Alg}, we present our decomposition algorithm. In Section~\ref{sec:Improv}, we derive valid lower bound inequalities to strengthen the master problem in the decomposition algorithm. In Section~\ref{sec:symm}, we derive two families of symmetry breaking constraints that improve the solvability of the proposed models.
\subsection{\textbf{Decomposition Algorithm}}\label{sec:Alg}
\noindent Proposition~\ref{Prop2:Convexh} suggests that constraint \eqref{const1:FinalDR} describes the epigraph of a convex and piecewise linear function of decision variables in formulation \eqref{FinalDR}. Therefore, given the two-stage characteristics of MAD-DRO in \eqref{FinalDR}, it is natural to attempt to solve problem \eqref{FinalDR} via a separation-based decomposition algorithm. Algorithm~\ref{Alg1:CAG} presents our proposed decomposition algorithm, and the algorithm for the W-DRO model  in  \eqref{Final_W_DMFRS} has the same steps. Algorithm~\ref{Alg1:CAG}  is finite because we identify a new piece of the function $\max \limits_{\pmb{\lambda, v, W, \pi, k}} \ \big\{ \sum \limits_{t \in T} \sum\limits_{i \in I} \pi_{i,t}+ \sum\limits_{t \in T} \sum\limits_{j \in J} \sum \limits_{m \in M} Cx_{j,m}^t v_{j,m}^t +\sum \limits_{t \in T} \sum \limits_{i \in I} -(W_{i,t} \rho_{i,t}+k_{i,t}\psi_{i,t}) \big\}$ each time when the set $ \lbrace L (\xbb, \delta)\geq 0  \rbrace$ is augmented in step 4, and the function has a finite number of pieces according to Proposition~\ref{Prop2:Convexh}. Note that this algorithm is based on the same theory and art of cutting plane-based decomposition algorithms employed in various other papers using decomposition to solve problems with similar structure. Nevertheless, we customized Algorithm~\ref{Alg1:CAG} to solve our proposed DRO models. In addition, in the following subsections, we derive problem-specific valid inequalities to strengthen the master problem, thus improving convergence.

\begin{algorithm}[t!]
\small
   \renewcommand{\arraystretch}{0.3}
\caption{Decomposition algorithm.}
\label{Alg1:CAG}
\noindent \textbf{1. Input.} Feasible region $\mathcal{X}$; support $\calS$; set of cuts $ \lbrace L (\xbb, \delta)\geq 0 \rbrace=\emptyset $; $LB=-\infty$ and $UB=\infty.$

\vspace{2mm}

\noindent \textbf{2. Master Problem.} Solve the following master problem
\begin{subequations}\label{Master}
\begin{align}
Z=&\min_{\xbb, \yb, \pmb \rho, \pmb \psi} \Bigg \{ \sum_{m \in M} f y_{m} -\sum_{t \in T} \sum_{j \in J} \sum_{m \in M } \alpha x_{j,m}^t+  \sum \limits_{t \in T} \sum \limits_{i \in I} \Big[\mu_{i,t} \rho_{i,t}+\eta_{i,t} \psi_{i,t}\Big] + \delta  \Bigg \}   \label{MasterObj} \\
& \ \ \text{s. t.} \qquad  (\xbb,\yb ) \in \mathcal{X},  \ \ \pmb{\psi} \geq 0, \ \ \{ L (\xbb, \delta)\geq 0 \},
 \end{align}
\end{subequations}

and record the optimal solution $(\xbb^	*, \pmb{\rho}^*, \pmb{\psi}^*,  \delta^*)$ and optimal value $Z^*$. Set $LB=Z^*$.

\noindent \textbf{3. Sub-problem.} 
\begin{itemize}
\item[3.1.] With $(\xbb, \pmb{\rho}, \pmb{\psi})$ fixed to $(\xbb^*, \pmb{\rho}^*, \pmb{\psi}^*)$, solve the following problem 
\begin{subequations}\label{MILPSep}
\begin{align}
h(\xbb, \pmb \rho, \pmb \psi)=  &\max \limits_{\pmb{\lambda, v, W, \pi, k} } \ \Big\{ \sum \limits_{t \in T} \sum\limits_{i \in I} \pi_{i,t}+ \sum\limits_{t \in T} \sum\limits_{j \in J} \sum \limits_{m \in M} Cx_{j,m}^t v_{j,m}^t +\sum \limits_{t \in T} \sum \limits_{i \in I} -(W_{i,t} \rho_{i,t}+k_{i,t}\psi_{i,t})  \Big\}\\
& \ \ \ \text{s.t. } \  \eqref{Const1:QudalX}-\eqref{Const3:Qdualv}, \eqref{McCormick1}-\eqref{McCormick2}, \eqref{Const2_finalinner},
\end{align}
\end{subequations}

$\qquad \ $ and record optimal solution $(\pmb{\pi^*, \lambda^*, W^*, v^*, k^*})$ and $h(\xbb, \pmb \rho, \pmb \psi)^*$. 
\item[3.2.] Set $UB=\min \{ UB, \ h(\xbb, \pmb \rho, \pmb \psi)^*+ (LB-\delta^*) \}$.
\end{itemize}

\noindent \textbf{4. if} $\delta^* \geq  \sum \limits_{t \in T} \sum\limits_{i \in I} \pi_{i,t}^*+ \sum\limits_{t \in T} \sum\limits_{j \in J} \sum \limits_{m \in M} Cx_{j,m}^{t*} v_{j,m}^{t*} +\sum \limits_{t \in T} \sum \limits_{i \in I} -(W_{i,t}^* \rho_{i,t}^*+k_{i,t}^*\psi_{i,t}^*) $ \textbf{then}

$\qquad \ \ $ stop and return $\xbb^*$ and $\yb^*$ as the optimal solution to problem \eqref{Master} (i.e., MAD-DRO in \eqref{FinalDR}).

\noindent $\ \ $ \textbf{else} add the cut $\delta \geq   \sum \limits_{t \in T} \sum \limits_{i \in I} \pi_{i,t}^*+ \sum\limits_{t \in T} \sum \limits_{j \in J} \sum\limits_{m \in M} Cx_{j,m} v_{j,m}^{t*}+\sum \limits_{t \in T} \sum \limits_{i \in I} -(W_{i,t}^* \rho_{i,t}+k_{i,t}^*\psi_{i,t} )$  to the set of

$\qquad \quad$  cuts   $ \{L (\xbb, \delta) \geq 0 \}$ and go to step 2.

\noindent $\ \ $ \textbf{end if}

\end{algorithm}

\subsection{\textbf{Multiple Optimality Cuts and Lower Bound Inequalities}}\label{sec:Improv}

\noindent In this section, we aim to incorporate more second-stage information into the first-stage without adding optimality cuts into the master problem by exploiting the structural properties of the recourse problem. We first observe that once the first-stage solutions and the demand are known, the second-stage problem can be decomposed into independent sub-problems with respect to time periods. Accordingly, we can construct cuts for each sub-problem in step 4. Let $\delta_t$ represent the optimality cut for each period $t$, we replace $\delta$ in  \eqref{MasterObj} with $\sum_{t} \delta_t$ and add constraints
\begin{align}
\delta_t\geq   \sum \limits_{i \in I} \pi_{i,t}^*+  \sum \limits_{j \in J} \sum\limits_{m \in M} Cx_{j,m} v_{j,m}^{t*}+ \sum \limits_{i \in I} -(W_{i,t}^* \rho_{i,t}+k_{i,t}^*\psi_{i,t} ), \qquad \forall t \in T.\label{Cut1}
\end{align}

\noindent The original single cut is the summation of multiple cuts of the form, i.e., $\delta=\sum_{t \in T} \delta_t$. Hence, in each iteration, we incorporate more or at least an equal amount of information into the master problem using \eqref{Cut1} as compared with the original single cut approach.  In this manner, the optimality cuts become more specific, which may result in better lower bounds and, therefore, a faster convergence. In Proposition \ref{Prop3:LowerB1}, we further identify valid lower bound inequalities for each time period to tighten the master problem (see Appendix~\ref{Appx:Prop3} for a proof).
\begin{proposition}\label{Prop3:LowerB1}
Inequalities \eqref{Vineq1} are valid lower bound inequalities on the recourse of the MFRSP.
\begin{align}
\sum \limits_{i \in I} \min \lbrace \gamma,  \min \limits_{j \in J} \lbrace \beta d_{i,j} \rbrace \rbrace \underline{W}_{i,t}, \qquad \forall t \in T. \label{Vineq1}
\end{align}

\end{proposition}
It follows that $\sum \limits_{i \in I} \mu_{i,t} \rho_{i,t}+\eta_{i,t} \psi_{i,t}+\delta_t\geq \sum \limits_{i \in I} \min \lbrace \gamma,  \min \limits_{j \in J} \lbrace \beta d_{i,j} \rbrace \rbrace \underline{W}_{i,t} $, for all $t \in T$, are valid.


\subsection{\textbf{Symmetry-Breaking Constraints}}\label{sec:symm}
\noindent Suppose there are three homogeneous MFs. As such, solutions $y=[1,1,0]^\top$, $y=[0,1,1]^\top$, and $y=[1, 0,1 ]^\top$ are equivalent (i.e., yield the same objective) in the sense that they all permit 2 out of 3 MFs to be used in the planning period. To avoid \textcolor{black}{wasting time} exploring such equivalent solutions, we assume that MFs are numbered sequentially and add constraints \eqref{Sym1} to the first-stage.
\begin{align}\label{Sym1}
y_{m+1} \leq y_m, \qquad \forall m <M.
\end{align}

\noindent \textcolor{black}{Constraints \eqref{Sym1} enforce} arbitrary ordering or scheduling of MFs. Second, recall that in the first period, $t=1$, we decide the initial locations of the MFs.  Therefore, it doesn't matter which MF is assigned to location $j$. For example, suppose that we have three candidate locations, and MFs 1 and 2 are active. Then, feasible solutions $(x_{1,1}^1=1$, $x_{3,2}^1=1)$ and $(x_{1,2}^1=1$, $x_{3,1}^1=1)$ yield the same objective. To avoid exploring such equivalent solutions, we define a dummy location $J+1$ and \textcolor{black}{add constraints \eqref{Ineq:Symm1}--\eqref{Ineq:Symm2} to the first-stage.}
\begin{subequations}\label{Sym2}
\begin{align}
& x_{j,m}^1 - \sum \limits_{j^\prime=j}^{J+1}x_{j^\prime,m+1}^1\leq 0, && \forall m <M, \forall j \in J, \label{Ineq:Symm1}\\
& x_{J+1, m}^1= 1- \sum \limits_{j=1}^Jx_{j,m}^1, && \forall m \in M.\label{Ineq:Symm2}
\end{align}
\end{subequations}
\noindent Constraints  \eqref{Ineq:Symm1}--\eqref{Ineq:Symm2} are valid for any formulation that uses the same sets of first-stage routing and scheduling decisions and constraints.  We derived constraints \eqref{Sym1}--\eqref{Sym2}  based on similar symmetry breaking principles  in  \cite{ostrowski2011orbital} and \cite{shehadeh2019analysis}. Although breaking symmetry is very important and standard in integer programming problems, our paper is the first to attempt to break the symmetry in the solution space of the first-stage planning decisions of the MFRSP. In Section~\ref{sec5:symmetry}, we demonstrate the computational advantages that could be gained by incorporating these inequalities.


\section{Computational Experiments}\label{sec:computational}

\noindent  In this section, we conduct extensive computational experiments comparing the proposed DRO models and a sample average approximation (SAA) of the SP model computationally and operationally, demonstrating where significant performance improvements could be gained.  The sample average model solves model \eqref{DMFRS_SP} with $\Prob$ replaced by an empirical distribution based on $N$ samples of $\Wb$  (see Appendix~\ref{Appx:SAA} for the formulation).  In Section~\ref{sec5.1:instancegen}, we describe the set of problem  instances and discuss other experimental setups.   In Section~\ref{sec5:CPU}, we compare solution time of the proposed models. In Section~\ref{sec5:symmetry}, we demonstrate efficiency of the proposed lower bound inequalities and symmetry breaking  constraints. \textcolor{black}{We compare optimal solutions of the proposed models and their out-of-sample performance in Sections~\ref{sec:Optimal Solutions} and \ref{sec5:OutSample}, respectively}. We analyze the sensitivity of the DRO expectation models to different parameter settings in Section~\ref{sec5:sensitivity}. \textcolor{black}{We close by comparing the risk-neutral and risk-averse models under critical parameter settings in Section~\ref{sec:CVaRExp}.}

\subsection{\textbf{Experimental Setup}}\label{sec5.1:instancegen}

\begin{figure}
\begin{center}

\includegraphics[scale=0.8]{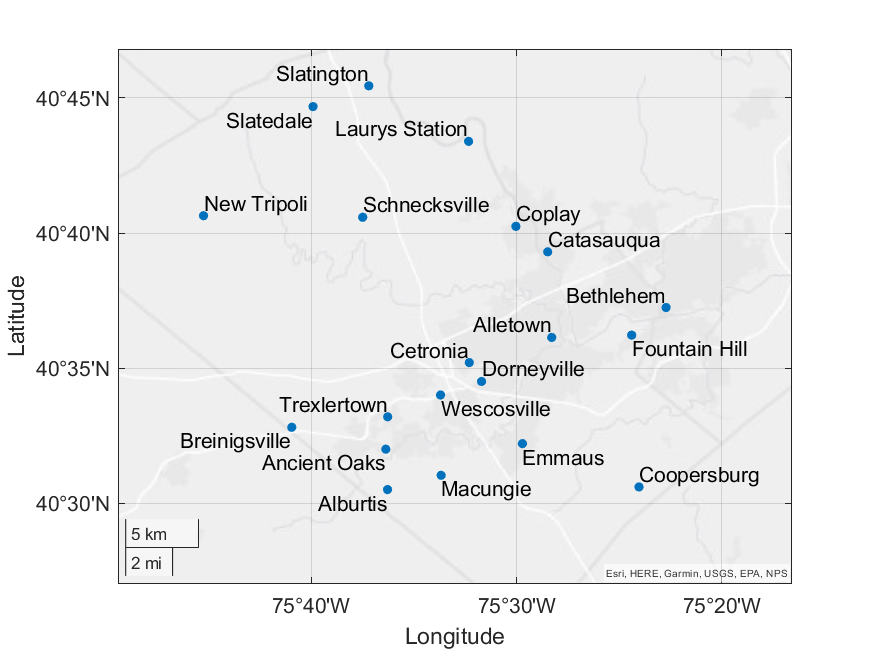}
\caption{Map of 20 cities in Lehigh County. We created this map using the \texttt{geoscatter} function (MATLAB).}\label{LehighMap}
\end{center}

\end{figure}

\noindent  We constructed 10 MFRSP instances, in part based on the same parameters settings and assumptions made in \cite{lei2014multicut} and \cite{lei2016two}. We summarize these instances in Table~\ref{table:DMFRSInstances}.  Each of the 10  instances is characterized by the number of customers locations $I$, number of candidate locations $J$, and the number of periods $T$.  Instances 1--4 are from \cite{lei2014multicut} and Instances 5--10 are from \cite{lei2016two}.  These benchmark instances represent a wide range of potential service regions in terms of problem size as a function of the number of demand nodes/locations and time periods. For example, if we account for the scale of the problem in the sense of a static facility location problem, instance 10 consists of  $30 \times 20= 600$ customers, which is relatively large for many practical applications. In addition, we constructed a service region based on 20 selected nodes (see \textcolor{black}{Figure~\ref{LehighMap}}) in Lehigh County of Pennsylvania (USA). Then,  as detailed below, we constructed two instances (denoted as Lehigh 1 and Lehigh 2) based on this region.
\begin{table}[t]
 \center 
 \footnotesize
   \renewcommand{\arraystretch}{0.3}
  \caption{MFRSP instances. Notation: $I$ is \# of customers, $J$ is \# of locations, and $T$ is  \# of periods.}
\begin{tabular}{cccccccccc}
\hline
\textbf{Inst} & \textbf{$\pmb{I}$}  &  \textbf{$\pmb{J}$}   & \textbf{$\pmb{T}$} \\
\hline
1&  10 & 10 & 10& \\
2 & 10 & 10 & 20& \\
3 & 15 & 15 & 10& \\
4 &15 & 15 & 20& \\
5 & 20 & 20 & 10\\
6 &20 & 20 & 20\\
7& 25 & 25& 10\\
8& 25 & 25& 20\\
9& 30 & 30& 10\\
10& 30 & 30& 20\\
\hline																									
\end{tabular} 
\label{table:DMFRSInstances}
\end{table}


For each instance in Table~\ref{table:DMFRSInstances}, we generated a total of $I$ vertices as uniformly distributed random numbers on a 100 by 100 plane and computed the distance between each pair of nodes  in Euclidean sense as in \cite{lei2014multicut}.  For Lehigh county instances, we first extracted the latitude and longitude of each node and used Bing Maps Developer API to compute the travel time in minutes between each pair of nodes.

We followed the same procedures in the DRO scheduling and facility location literature to generate random parameters as follows. \textcolor{black}{For instances 1--10,  we  generated $\mu_{i,t}$ from a uniform distribution $U[\underline{W}, \overline{W}]=[20, 60]$, and set the standard deviation $\sigma_{i,t}=0.5\mu_{i,t}$, for all $i \in I$ and $t \in T$.} For Lehigh county instances, we used the population estimate for each node based on the most updated information posted on the 2010 US Census Bureau (see Appendix~\ref{AppexLehigh}) to construct the following two demand structures.  In Lehigh 1, we generated the demand's mean as follows: if the population $\geq$10,000, we set $\mu_{i,t}=40$ (i.e., the mean of  $U[20,60]$); if the population $\in$ [5,000, 10,000), we set the mean to $\mu=30$; if the population $\in$ [1,000, 5,000),  we set $\mu_{i,t}=20$; and if the population $<1,000$, we set $\mu_{i,t}=15$. In Lehigh 2, we use the population percentage (weight) at each node to generate the demand's mean as $\mu_{i,t}=\min(60, \texttt{population\%} \times 1000)$. To a certain extent, these structures reflect what may be observed in real life, i.e., locations with more population may potentially create greater demand. We refer to Appendix~\ref{AppexLehigh} for the details related to Lehigh 1 and Lehigh 2. 


For each instance, we sample $N$ realizations $W_{i,t}^n, \ldots, W_{I,T}^n$, for $n=1, \ldots, N$, by following lognormal (LogN) distributions with the generated  $\mu_{i,t}$ and $\sigma_{i,t}$.   We round each  parameter to the nearest integer. We solve the SAA and W-DRO models using the $N$ sample and the MAD-DRO model with the corresponding mean, MAD, and range. The Wasserstein ball's radius $\epsilon$ in the W-DRO model is an input parameter, where different values of $\epsilon$ may result in a different robust solution $\xbb (\epsilon, N)$ with a very different out-of-sample performance $\hat{\E} [Q(\xbb (\epsilon, N),\Wb)]$. To estimate $\epsilon^{\mbox{\tiny opt}}$ that minimizes $\hat{\E} [Q(\xbb (\epsilon, N),\Wb)]$, we employed a widely used cross-validation method (see Appendix~\ref{WDRO_Radius}).

We assume that all cost parameters are calculated in terms of present monetary value.  Specifically,  as in \cite{lei2014multicut}, for each instance, we randomly generate (1) the fixed cost from a uniform distribution $U[a,b]$ with $a = 1000$ and $b = 1500$; (2) the assignment cost factor per unit distance per unit demand $\beta$ from $U[0.0001a,0.0001b]$; and (3) the penalty cost per unit demand $\gamma$ form $U[0.01a, 0.01b]$. Finally, we set the traveling inconvenience cost factor to $0.0001a$, and unless stated otherwise, we use a capacity parameter $C=100$. We implemented all models and the decomposition algorithm using AMPL2016 Programming language calling CPLEX V12.6.2 as a solver with default settings. We run all experiments on a MacBook Pro with Apple M1 Max Chip, 32GB of memory, and 10-core CPUs.   Finally, we imposed a solver time limit of 1 hour.

\subsection{\textbf{CPU Time}}\label{sec5:CPU}

\noindent   In this section, we analyze solution times of the proposed DRO models. \textcolor{black}{We consider two ranges of the demand:  $\Wb \in [20, 60]$ (base-case) and $\Wb \in [50, 100]$ (higher demand variability and volume). We also consider two MF capacities: $C=60$ (relatively small capacity) and $C=100$ (relatively large capacity).} We focus on solving problem instances with a small sample size, which is often seen in most real-world applications (especially in healthcare) and is the primary motivation for using DRO. Specifically, we use $N=10$ and $N=50$ as a sample size for the W-DRO model. For each of the 10 instances in Table~\ref{table:DMFRSInstances}, demand range, $C$, and $N$, we generated and solved 5 instances using each model as described in Section~\ref{sec5.1:instancegen}.

\color{black}
Let us first analyze solution times of the risk-neutral DRO models. Tables~\ref{table:MAD_CPU_2060}--\ref{table:MAD_CPU_50100} present the computational details (i.e., CPU time in seconds and the number of iterations in the decomposition algorithm before it converges to the optimum) of solving the  MAD-DRO model. Tables~\ref{table:WASS_CPU_2060_N10}--\ref{table:WASS_CPU_50100_N10} and Tables~\ref{table:WASS_CPU_2060_N50}--\ref{table:WASS_CPU_50100_N50} present the computational details of solving the W-DRO model with $N=10$ and $N=50$, respectively.  We observe the following from these tables. First, the computational effort (i.e., solution time per iteration)  increases with the size of instance ($I \times J \times T$). Second, solution times are shorter under tight capacity ($C=60$) than under large capacity ($C=100$). In addition, the decomposition algorithm takes fewer iterations to converge to the optimum under tight capacity. Intuitively, when $C=100$, each facility can satisfy more demand than $C=60$. Thus, there are more feasible choices for the MF fleet size and schedule when $C=100$. That said, the search space when $C=100$ is larger, potentially leading to a longer computational time.

Third, using the MAD-DRO model, we were able to solve all instances within the time limit. The average solution time for Instances 1--7 using the MAD-DRO model ranges from 2 to 914 seconds. The average solution time of larger instances (Instances 8--10) ranges from 241 to 1,901 seconds. In contrast, using the W-DRO model, we were able to solve Instances 1--8 with $N=10$ and Instances 1--7 with $N=50$. The average solution time of the W-DRO model with $N=10$ and $N=50$ ranges from 5 to 715 and from 14 to 1,725 seconds, respectively.  Note that solution times of the W-DRO model increase with the sample size. This makes sense because the number of variables and constraints of the W-DRO model increases as $N$ increases. In addition, solution times of the W-DRO model with $N=50$ are generally longer than the MAD-DRO model.  This also makes sense because the MAD-DRO model is a smaller deterministic model (i.e., it has fewer variables and constraints).

\begin{table}[t!]
\color{black}
\center 
      \small 
\caption{Computational details of solving the MAD-DRO model ($ \Wb \in [20, 60]$).}
   \renewcommand{\arraystretch}{0.6}
\begin{tabular}{lllllllllllllllllllllllllllllllllll}
 \hline
&  &  & \multicolumn{6}{c}{$C=60$ }  &&   \multicolumn{6}{c}{$C=100$ } \\  \cline{5-10} \cline{12-17}
Inst & $I$  &  $T$  & &  \multicolumn{3}{c}{CPU time} &   \multicolumn{3}{c}{iteration}  &&  \multicolumn{3}{c}{CPU time}   & \multicolumn{3}{c}{iteration} \\ \cline{5-10} \cline{12-17}
& &  & &  Min & Avg & Max & Min & Avg & Max & & Min & Avg & Max & Min & Avg & Max\\
\hline
1	&	10	&	10	&	&	1	&	2	&	4	&	3	&	5	&	12	&	&	5	&	6	&	8	&	18	&	23	&	30	\\
2	&	10	&	20	&	&	3	&	7	&	13	&	21	&	30	&	51	&	&	2	&	3	&	4	&	5	&	9	&	15	\\
3	&	15	&	10	&	&	2	&	5	&	10	&	3	&	7	&	16	&	&	15	&	24	&	37	&	27	&	41	&	62	\\
4	&	15	&	20	&	&	4	&	7	&	12	&	6	&	11	&	22	&	&	18	&	25	&	43	&	34	&	45	&	72	\\
5	&	20	&	10	&	&	6	&	10	&	19	&	6	&	11	&	22	&	&	19	&	26	&	45	&	34	&	45	&	72	\\
6	&	20	&	20	&	&	9	&	41	&	65	&	6	&	28	&	45	&	&	91	&	308	&	640	&	53	&	68	&	80	\\
7	&	25	&	10	&	&	23	&	169	&	578	&	8	&	25	&	44	&	&	503	&	703	&	1004	&	110	&	138	&	190	\\
8	&	25	&	20	&	&	57	&	241	&	352	&	5	&	33	&	53	&	&	573	&	872	&	1272	&	147	&	196	&	245	\\
9	&	30	&	10	&	&	99	&	416	&	696	&	16	&	61	&	80	&	&	185	&	1860	&	1871	&	113	&	147	&	161	\\
10	&	30	&	20	&	&	463	&	905	&	1785	&	7	&	38	&	91	&	&	1098	&	1592	&	1887	&	6	&	125	&	191	\\
\hline 
\end{tabular}\label{table:MAD_CPU_2060}
\end{table}

\begin{table}[t!]
\color{black}
\center 
      \small 
\caption{Computational details of solving the MAD-DRO model ($ \Wb \in [50, 100]$).}
   \renewcommand{\arraystretch}{0.6}
\begin{tabular}{lllllllllllllllllllllllllllllllllll}
 \hline
&  &  & \multicolumn{6}{c}{$C=60$ }  &&   \multicolumn{6}{c}{$C=100$ } \\  \cline{5-10} \cline{12-17}
Inst & $I$  &  $T$  & &  \multicolumn{3}{c}{CPU time} &   \multicolumn{3}{c}{iteration}  &&  \multicolumn{3}{c}{CPU time}   & \multicolumn{3}{c}{iteration} \\ \cline{5-10} \cline{12-17}
& &  & &  Min & Avg & Max & Min & Avg & Max & & Min & Avg & Max & Min & Avg & Max\\
\hline
1	&	10	&	10	&	&	1	&	2	&	2	&	1	&	2	&	4	&	&	10	&	11	&	13	&	20	&	22	&	24	\\
2	&	10	&	20	&	&	1	&	2	&	3	&	2	&	4	&	9	&	&	5	&	8	&	13	&	15	&	26	&	41	\\
3	&	15	&	10	&	&	4	&	9	&	17	&	2	&	2	&	3	&	&	29	&	47	&	65	&	28	&	34	&	45	\\
4	&	15	&	20	&	&	3	&	9	&	23	&	2	&	5	&	13	&	&	43	&	54	&	70	&	44	&	54	&	63	\\
5	&	20	&	10	&	&	27	&	68	&	177	&	3	&	5	&	14	&	&	215	&	292	&	392	&	70	&	88	&	124	\\
6	&	20	&	20	&	&	6	&	14	&	38	&	3	&	5	&	14	&	&	292	&	322	&	347	&	74	&	94	&	107	\\
7	&	25	&	10	&	&	20	&	385	&	796	&	2	&	9	&	29	&	&	553	&	686	&	914	&	175	&	205	&	258	\\
8	&	25	&	20	&	&	63	&	334	&	835	&	3	&	15	&	36	&	&	122	&	622	&	1635	&	3	&	14	&	36	\\
9&	30	&	10	&	&	25	&	579	&	1692	&	3	&	39	&	143	&	&	1851	&	1860	&	1871	&	113	&	147	&	161	\\
10&	30	&	20	&	&	81	&	300	&	1240	&	3	&	10	&	41	&	&	1821	&	1901	&	1937	&	123	&	147	&	156	\\
\hline 
\end{tabular}\label{table:MAD_CPU_50100}
\end{table}

\begin{table}[t!]
\color{black}
\center 
      \small 
\caption{Computational details of solving the W-DRO model ($ \Wb \in [20, 60]$, $N=10$).}
   \renewcommand{\arraystretch}{0.6}
\begin{tabular}{lllllllllllllllllllllllllllllllllll}
 \hline
&  &  & \multicolumn{6}{c}{$C=60$ }  &&   \multicolumn{6}{c}{$C=100$ } \\  \cline{5-10} \cline{12-17}
Inst & $I$  &  $T$  & &  \multicolumn{3}{c}{CPU time} &   \multicolumn{3}{c}{iteration}  &&  \multicolumn{3}{c}{CPU time}   & \multicolumn{3}{c}{iteration} \\ \cline{5-10} \cline{12-17}
& &  & &  Min & Avg & Max & Min & Avg & Max & & Min & Avg & Max & Min & Avg & Max\\
\hline
1	&	10	&	10	&	&	7	&	9	&	18	&	7	&	9	&	11	&	&	15	&	20	&	22	&	14	&	18	&	20	\\
2	&	10	&	20	&	&	8	&	10	&	12	&	7	&	8	&	10	&	&	18	&	24	&	28	&	15	&	18	&	20	\\
3	&	15	&	10	&	&	35	&	39	&	45	&	16	&	18	&	20	&	&	12	&	15	&	18	&	6	&	8	&	9	\\
4	&	15	&	20	&	&	30	&	43	&	62	&	5	&	6	&	78	&	&	115	&	134	&	147	&	23	&	25	&	29	\\
5	&	20	&	10	&	&	41	&	48	&	64	&	6	&	8	&	9	&	&	99	&	448	&	682	&	21	&	28	&	33	\\
6	&	20	&	20	&	&	49	&	139	&	177	&	6	&	8	&	10	&	&	142	&	207	&	273	&	24	&	35	&	44	\\
7	&	25	&	10	&	&	81	&	86	&	97	&	8	&	9	&	10	&	&	296	&	410	&	506	&	32	&	41	&	50	\\
8	&	25	&	20	&	&	100	&	678	&	1820	&	7	&	9	&	11	&	&	534	&	715	&	1068	&	47	&	51	&	56	\\
\hline 
\end{tabular}\label{table:WASS_CPU_2060_N10}
\end{table}

\begin{table}[t!]
\color{black}
\center 
      \small 
\caption{Computational details of solving the W-DRO model ($ \Wb \in [50, 100]$, $N=10$).}
   \renewcommand{\arraystretch}{0.6}
\begin{tabular}{lllllllllllllllllllllllllllllllllll}
 \hline
&  &  & \multicolumn{6}{c}{$C=60$ }  &&   \multicolumn{6}{c}{$C=100$ } \\  \cline{5-10} \cline{12-17}
Inst & $I$  &  $T$  & &  \multicolumn{3}{c}{CPU time} &   \multicolumn{3}{c}{iteration}  &&  \multicolumn{3}{c}{CPU time}   & \multicolumn{3}{c}{iteration} \\ \cline{5-10} \cline{12-17}
& &  & &  Min & Avg & Max & Min & Avg & Max & & Min & Avg & Max & Min & Avg & Max\\
\hline
1	&	10	&	10	&	&	4	&	5	&	7	&	3	&	3	&	4	&	&	7	&	9	&	10	&	6	&	7	&	8	\\
2	&	10	&	20	&	&	9	&	12	&	17	&	3	&	3	&	4	&	&	15	&	16	&	18	&	6	&	6	&	7	\\
3	&	15	&	10	&	&	16	&	32	&	55	&	3	&	3	&	3	&	&	35	&	48	&	61	&	6	&	7	&	9	\\
4	&	15	&	20	&	&	30	&	34	&	43	&	3	&	3	&	3	&	&	71	&	81	&	85	&	8	&	9	&	10	\\
5	&	20	&	10	&	&	29	&	42	&	76	&	3	&	3	&	3	&	&	75	&	79	&	81	&	9	&	10	&	11	\\
6	&	20	&	20	&	&	82	&	92	&	108	&	3	&	3	&	3	&	&	110	&	182	&	229	&	7	&	11	&	13	\\
7	&	25	&	10	&	&	108	&	123	&	135	&	9	&	10	&	11	&	&	151	&	192	&	243	&	9	&	10	&	12	\\
\hline 
\end{tabular}\label{table:WASS_CPU_50100_N10}
\end{table}

\begin{table}[t!]
\color{black}
\center 
      \small 
\caption{Computational details of solving the W-DRO model ($ \Wb \in [20, 60]$, $N=50$).}
   \renewcommand{\arraystretch}{0.6}
\begin{tabular}{lllllllllllllllllllllllllllllllllll}
 \hline
&  &  & \multicolumn{6}{c}{$C=60$ }  &&   \multicolumn{6}{c}{$C=100$ } \\  \cline{5-10} \cline{12-17}
Inst & $I$  &  $T$  & &  \multicolumn{3}{c}{CPU time} &   \multicolumn{3}{c}{iteration}  &&  \multicolumn{3}{c}{CPU time}   & \multicolumn{3}{c}{iteration} \\ \cline{5-10} \cline{12-17}
& &  & &  Min & Avg & Max & Min & Avg & Max & & Min & Avg & Max & Min & Avg & Max\\
\hline
1	&	10	&	10	&	&	18	&	28	&	35	&	6	&	10	&	12	&	&	53	&	60	&	71	&	17	&	19	&	21	\\
2	&	10	&	20	&	&	30	&	35	&	42	&	8	&	9	&	11	&	&	89	&	104	&	114	&	18	&	20	&	21	\\
3	&	15	&	10	&	&	28	&	32	&	81	&	6	&	6	&	7	&	&	104	&	124	&	139	&	21	&	23	&	25	\\
4	&	15	&	20	&	&	40	&	43	&	48	&	6	&	6	&	7	&	&	184	&	214	&	231	&	21	&	23	&	25	\\
5	&	20	&	10	&	&	33	&	57	&	95	&	4	&	6	&	7	&	&	271	&	341	&	463	&	29	&	32	&	36	\\
6	&	20	&	20	&	&	89	&	306	&	723	&	6	&	8	&	10	&	&	455	&	1461	&	2141	&	22	&	29	&	35	\\
7	&	25	&	10	&	&	109	&	127	&	169	&	6	&	7	&	8	&	&	1006	&	1725	&	2095	&	37	&	40	&	42	\\
\hline 
\end{tabular}\label{table:WASS_CPU_2060_N50}
\end{table}

\begin{table}[t!]
\color{black}
\center 
      \small 
\caption{Computational details of solving the W-DRO model ($ \Wb \in [50, 100]$, $N=50$).}
   \renewcommand{\arraystretch}{0.6}
\begin{tabular}{lllllllllllllllllllllllllllllllllll}
 \hline
&  &  & \multicolumn{6}{c}{$C=60$ }  &&   \multicolumn{6}{c}{$C=100$ } \\  \cline{5-10} \cline{12-17}
Inst & $I$  &  $T$  & &  \multicolumn{3}{c}{CPU time} &   \multicolumn{3}{c}{iteration}  &&  \multicolumn{3}{c}{CPU time}   & \multicolumn{3}{c}{iteration} \\ \cline{5-10} \cline{12-17}
& &  & &  Min & Avg & Max & Min & Avg & Max & & Min & Avg & Max & Min & Avg & Max\\
\hline
1	&	10	&	10	&	&	12	&	14	&	17	&	3	&	3	&	4	&	&	30	&	32	&	33	&	7	&	7	&	8	\\
2	&	10	&	20	&	&	17	&	22	&	25	&	3	&	3	&	4	&	&	44	&	57	&	63	&	7	&	9	&	9	\\
3	&	15	&	10	&	&	65	&	81	&	92	&	6	&	6	&	7	&	&	71	&	105	&	145	&	6	&	6	&	7	\\
4	&	15	&	20	&	&	42	&	53	&	87	&	3	&	3	&	7	&	&	48	&	65	&	84	&	4	&	5	&	7	\\
5	&	20	&	10	&	&	50	&	60	&	66	&	3	&	3	&	3	&	&	111	&	144	&	223	&	7	&	7	&	7	\\
6	&	20	&	20	&	&	580	&	1174	&	2024	&	3	&	3	&	3	&	&	319	&	1009	&	3000	&	8	&		10&11		\\
7	&	25	&	10	&	&	106	&	113	&	129	&	3	&	3	&	3	&	&	370	&	402	&	464	&	9	&	9	&	10	\\
\hline 
\end{tabular}\label{table:WASS_CPU_50100_N50}
\end{table}

\color{black}

\textcolor{black}{Let us now analyze solution times of the risk-averse models. We use MAD-CVaR (W-CVaR) to denote the mean-CVaR-based DRO model with MAD (1-Wasserstein) ambiguity. Since we observe similar computational performance with different values of $\Theta \in [0, 1)$, we present results with $\Theta=0$. Tables  \ref{table:MADCVAR_2060}--\ref{table:MADCVAR__50100}  and Tables \ref{table:Wass_CVAR_2060_N10}--\ref{table:Wass_CVAR_50100_N50} in Appendix~\ref{Appx:CPU2} present the computational details of solving the MAD-CVaR model and W-CVaR model, respectively. Using the MAD-CVaR model, we were able to solve Instances 1--8 with $\Wb \in [20,60]$ and Instances 1--7 with $\Wb \in [50,100]$. In addition, solution times of MAD-CVaR model are longer than the risk-neutral MAD-DRO model. In contrast, using the W-CVaR model, we were able to solve Instances 1--5. Solution times of these instances are generally longer than the risk-neutral model, especially when $N=50$.  It is not surprising that the CVaR models are more computationally challenging to solve than the risk-neutral models because the former models are larger (have more variables and constraints). In particular, the master problem of the CVaR models in the decomposition algorithm is larger than the expectation models (see Algorithm~\ref{Alg2:Decomp2} in Appendix~\ref{Sec:meanCVAR}).}

Finally, it is worth mentioning that using an \textit{enhanced multicut L-shaped} (E-LS) method to solve their SP model, \cite{lei2014multicut} were able to solve Instance 1--4. The average solution time of Instance 4 using E-LS is 3000 seconds obtained at 5\% optimality gap.  The CVaR-based SP model is more challenging to solve than the risk-neutral SP model.

\subsection{\textbf{Efficiency of  Inequalities  \eqref{Vineq1}--\eqref{Sym2} }}\label{sec5:symmetry}
\noindent In this section, we study the efficiency of symmetry breaking constraints~\eqref{Sym1}--\eqref{Sym2} and lower bounding inequalities \eqref{Vineq1}. Given the challenges of solving large instances without  \eqref{Vineq1}--\eqref{Sym2}, we use Instance 1 with $\Wb \in [20, 60]$ and $C=100$ in this experiment.

\begin{figure}[t!]
     \begin{subfigure}[b]{0.5\textwidth}
 \centering
        \includegraphics[width=\textwidth]{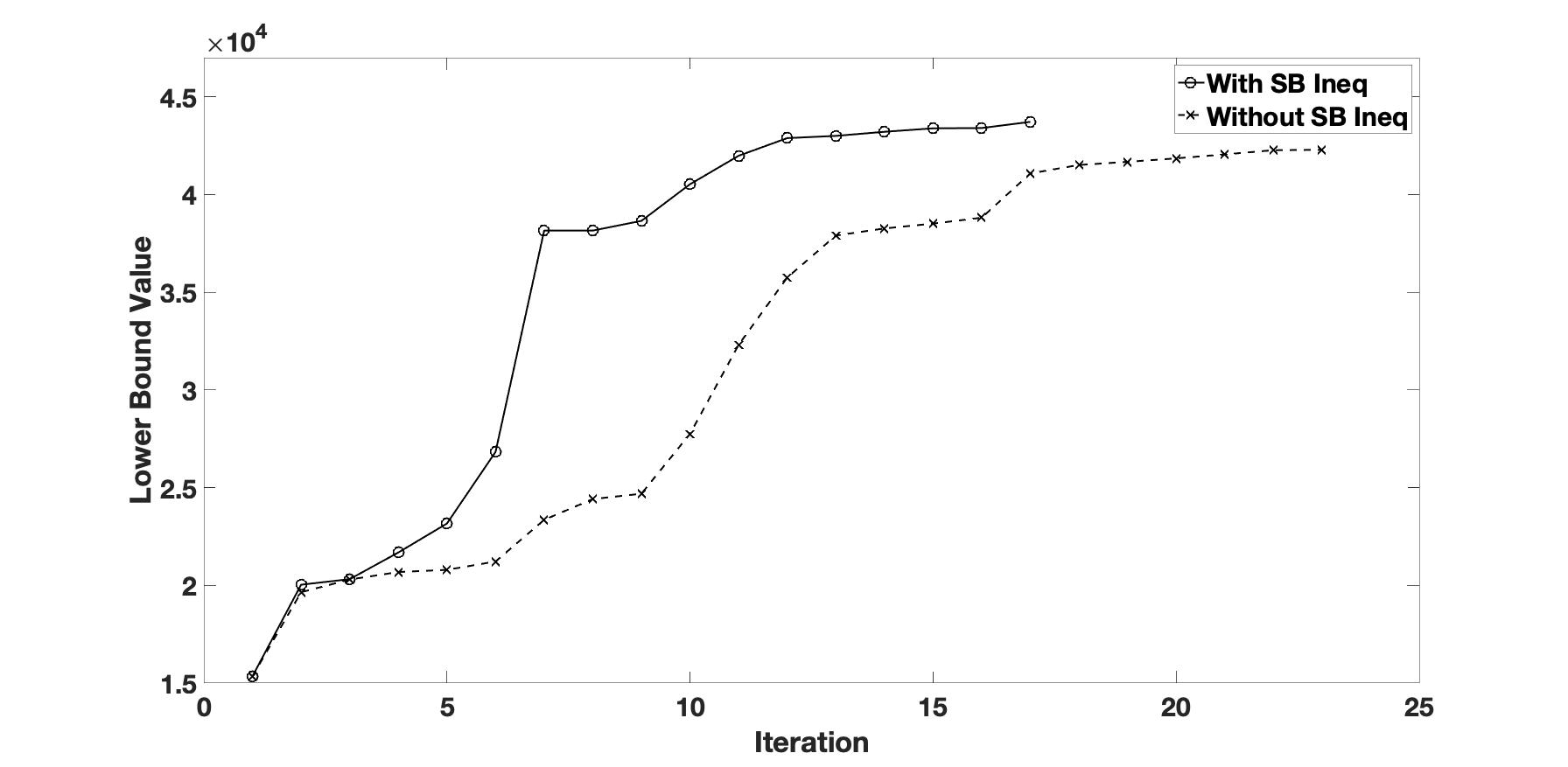}
        \caption{LB values, W-DRO}
        \label{FigSBa}
    \end{subfigure}%
    \begin{subfigure}[b]{0.5\textwidth}
            \includegraphics[width=\textwidth]{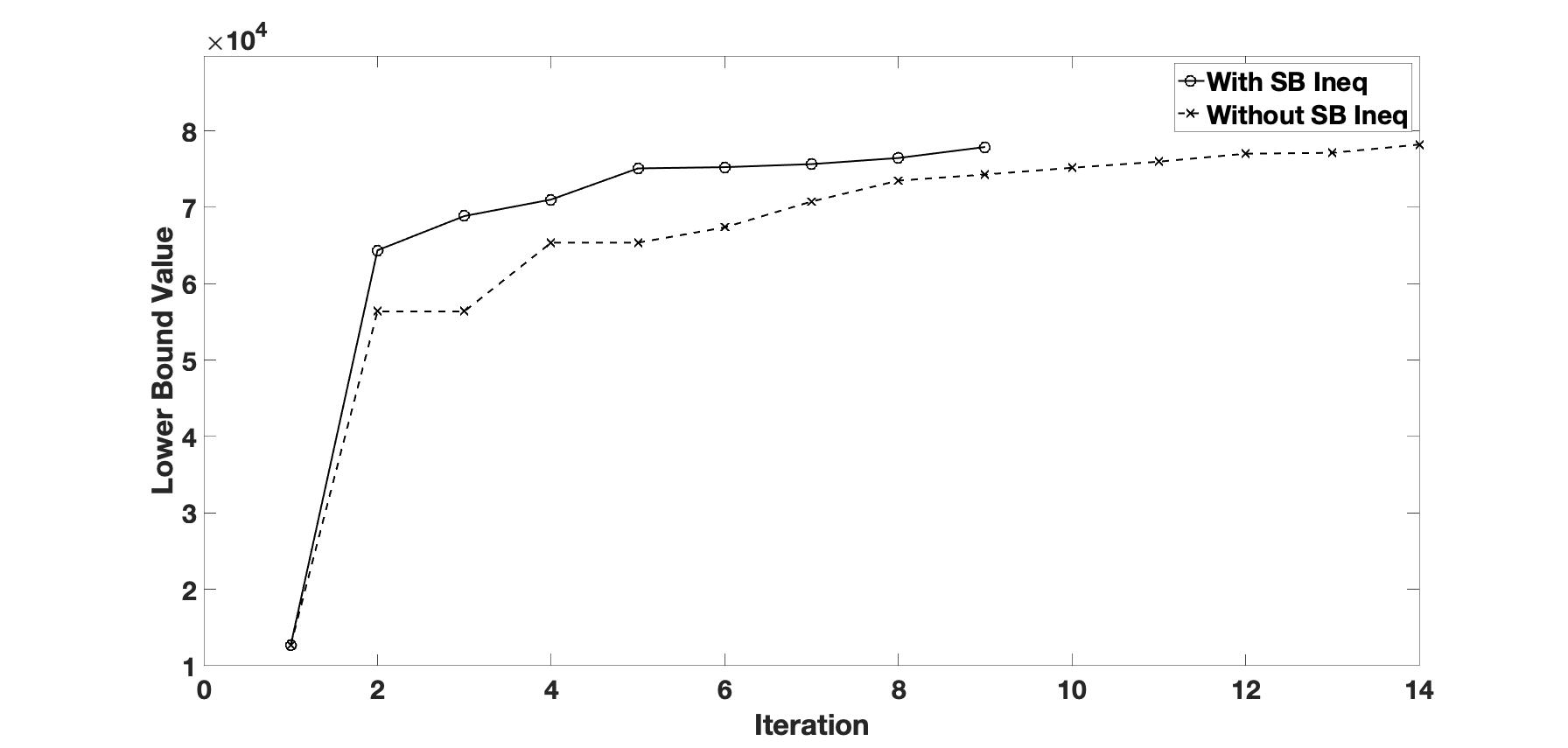}
      \caption{LB value, MAD-DRO}
      \label{FigSBb}
    \end{subfigure}%
    
         \begin{subfigure}[b]{0.5\textwidth}
 \centering
        \includegraphics[width=\textwidth]{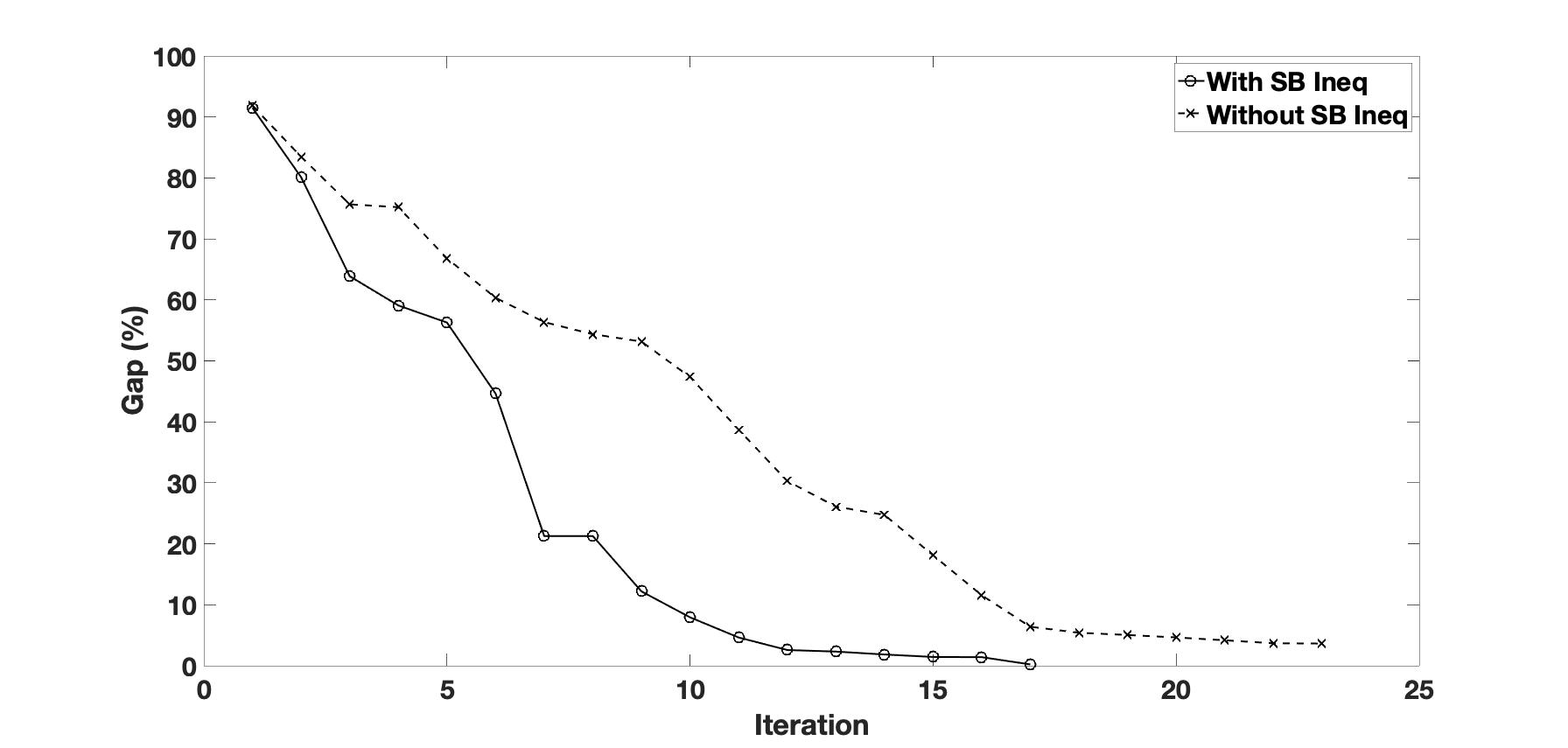}
        \caption{Gap values, W-DRO}
        \label{FigSBGAPa}
    \end{subfigure}%
    \begin{subfigure}[b]{0.5\textwidth}
            \includegraphics[width=\textwidth]{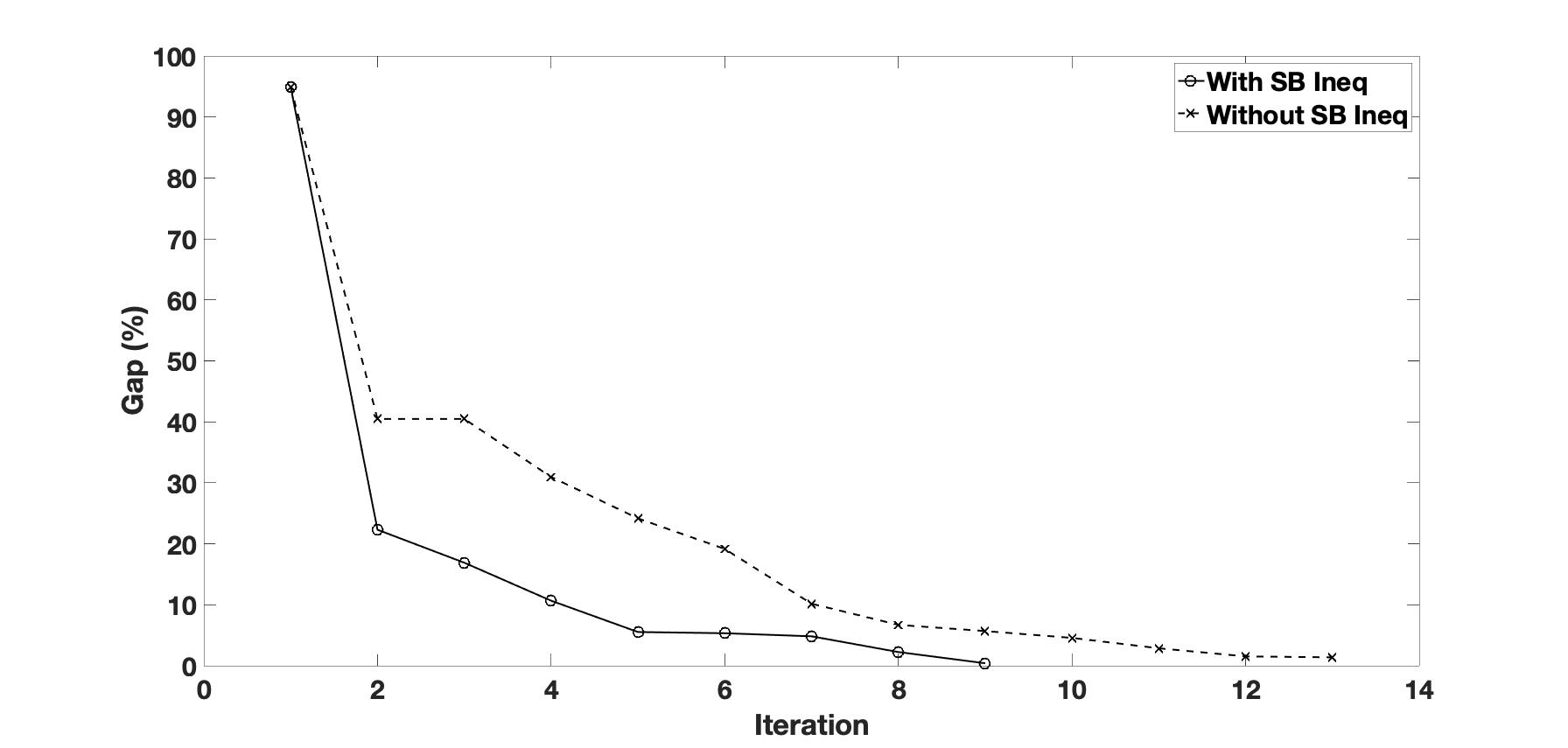}
      \caption{Gap values, W-DRO}
      \label{FigSBGAPb}
    \end{subfigure}%
            \caption{Comparisons of lower bound and gap values with and without SB constraints \eqref{Sym1}-\eqref{Sym2}.}\label{FigSB}
\end{figure}

First, we separately solve the proposed models with and without  symmetry-breaking (SB) constraints \eqref{Sym1}--\eqref{Sym2}. First, we observe that without these SB constraints, solution times of Instance 1 using (W-DRO, MAD-DRO, SP) significantly increase from (20, 6, 70) to (1,765, 1,003, 3,600) seconds. Instances 3--10 terminated with a large gap after one hour without these SB constraints. Second, as shown in Figure~\ref{FigSB}, both the lower bound and gap (i.e., the relative difference between the upper and lower bounds on the objective value)  converge faster when we include constraints \eqref{Sym1}--\eqref{Sym2} in the master problem. Moreover, constraints \eqref{Sym1}--\eqref{Sym2} lead to a stronger bound in each iteration.   These results demonstrate the importance of breaking the symmetry in the first-stage decisions and the effectiveness of our SB constraints.

\begin{figure}[t!]
     \begin{subfigure}[b]{0.5\textwidth}
 \centering
        \includegraphics[width=\textwidth]{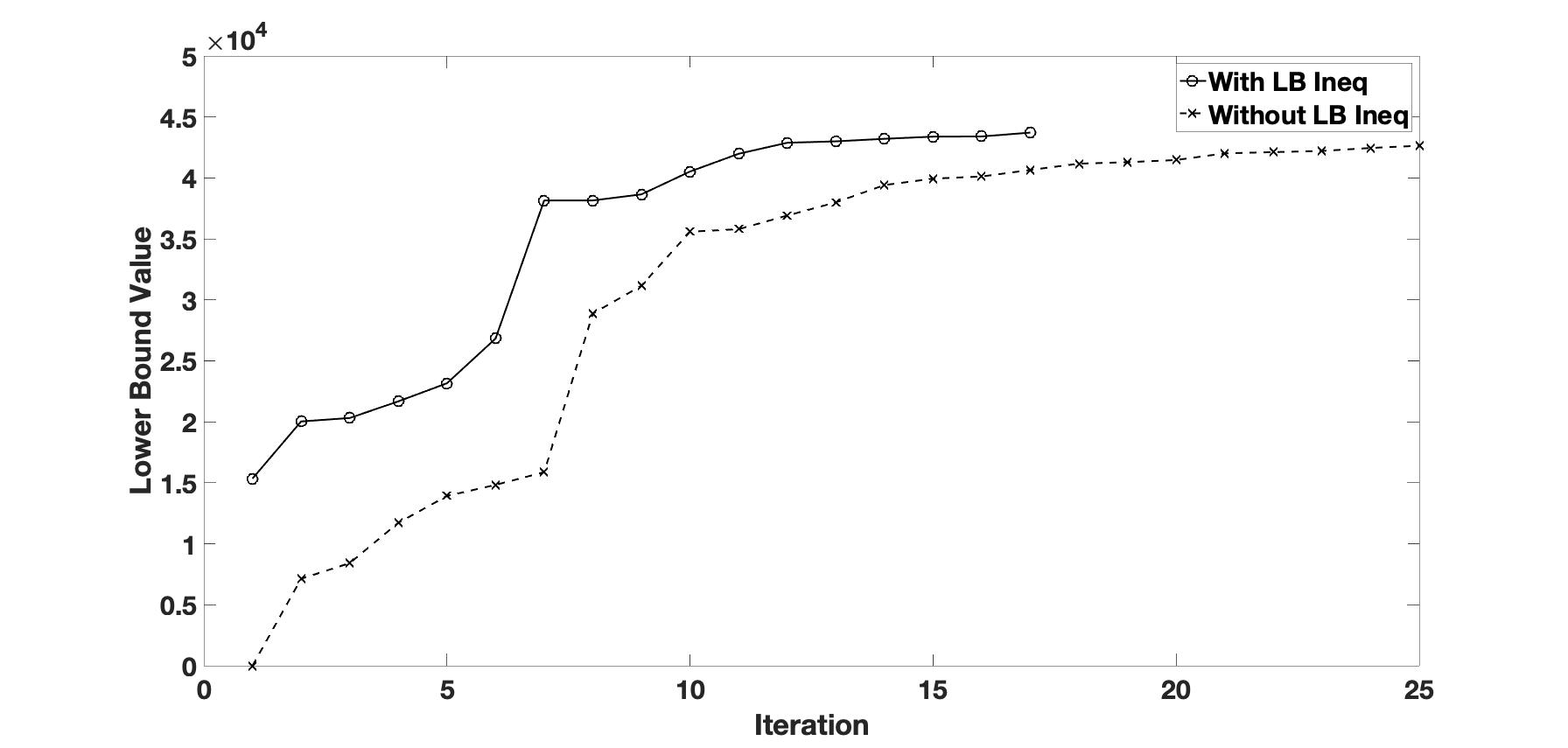}
        \caption{LB values, W-DRO}
        \label{FigVIa}
    \end{subfigure}%
    \begin{subfigure}[b]{0.5\textwidth}
            \includegraphics[width=\textwidth]{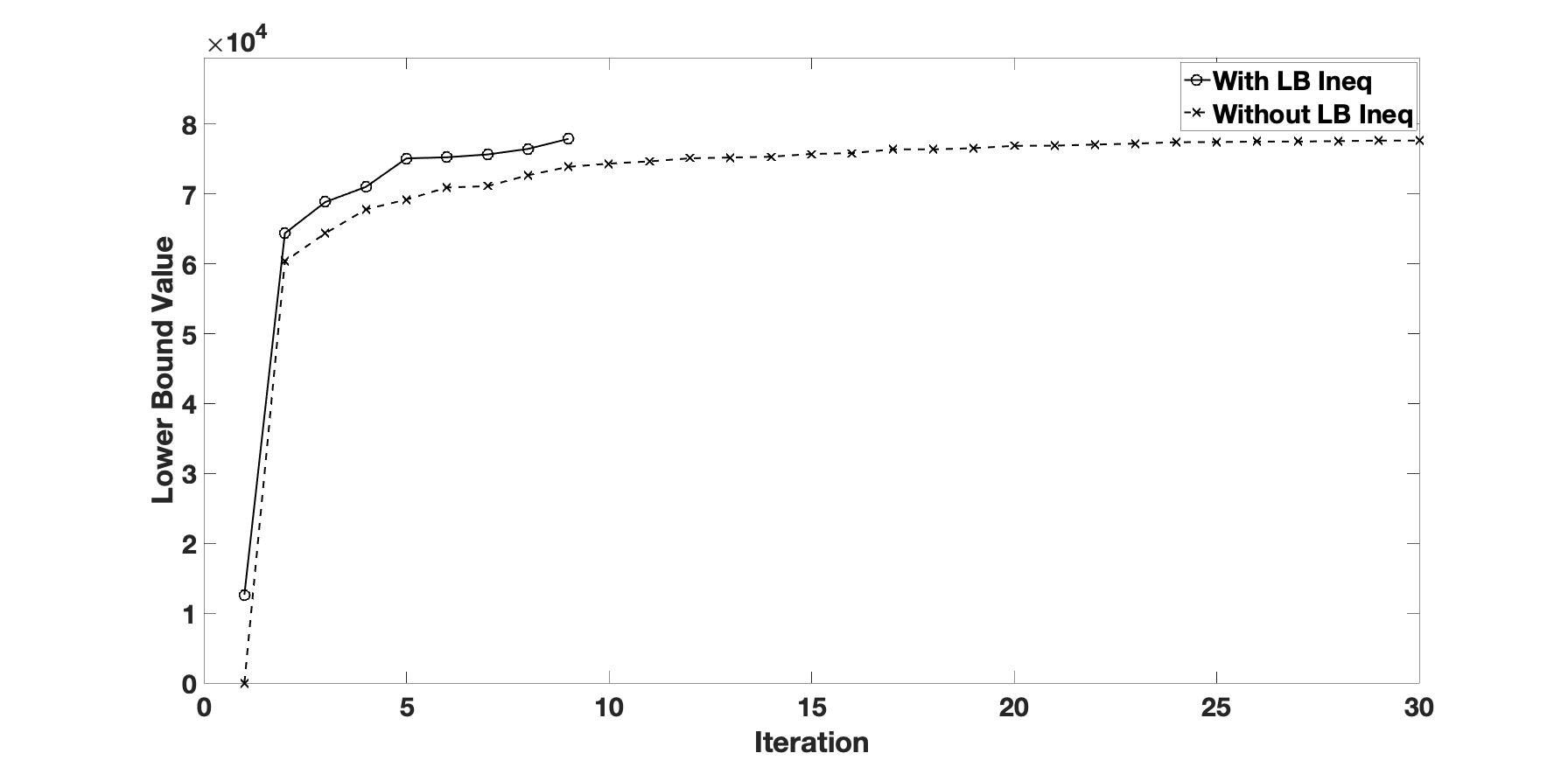}
      \caption{LB values, MAD-DRO}
      \label{FigVIb}
    \end{subfigure}%
    
     \begin{subfigure}[b]{0.5\textwidth}
 \centering
        \includegraphics[width=\textwidth]{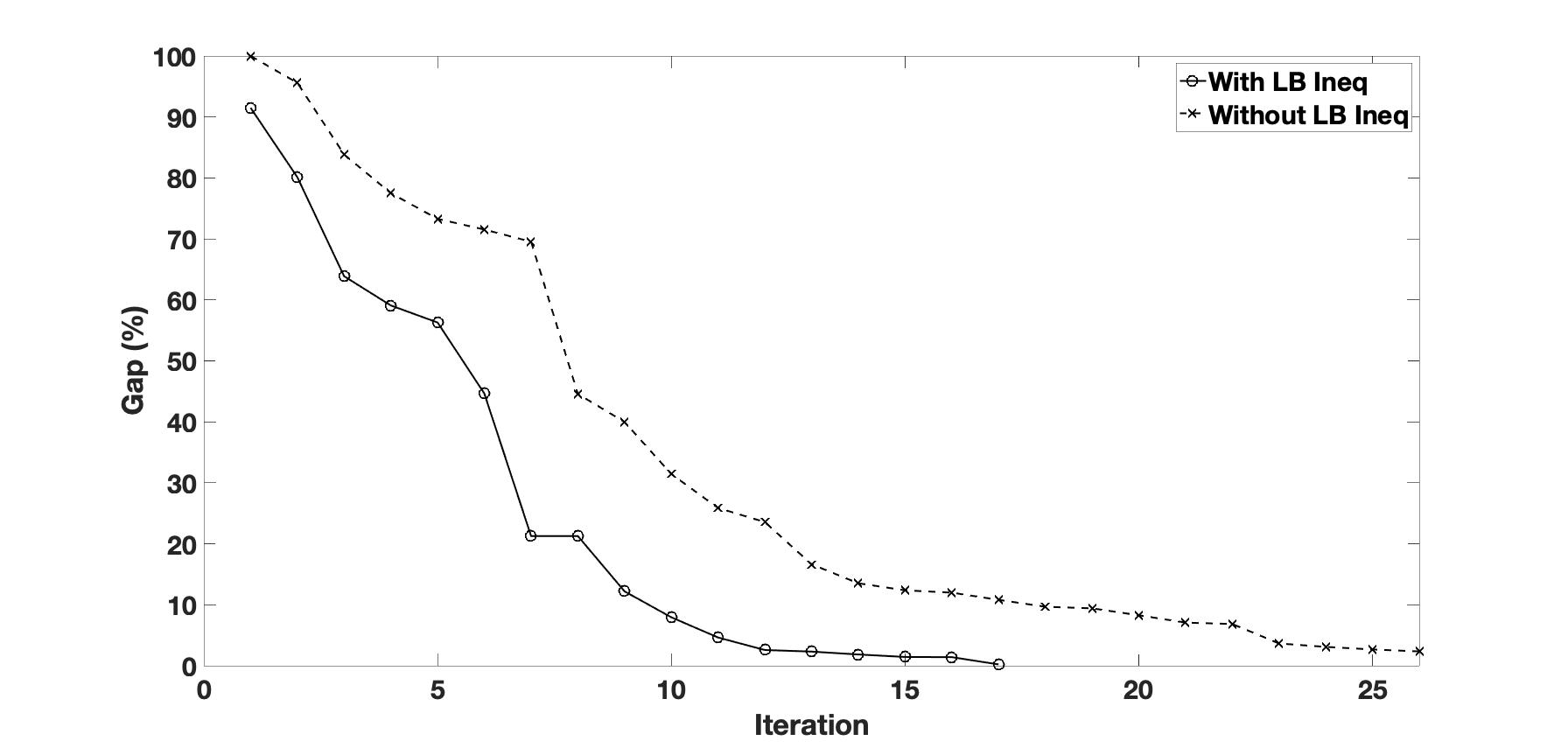}
        \caption{Gap values, W-DRO}
        \label{FigVIaGap}
    \end{subfigure}%
    \begin{subfigure}[b]{0.5\textwidth}
            \includegraphics[width=\textwidth]{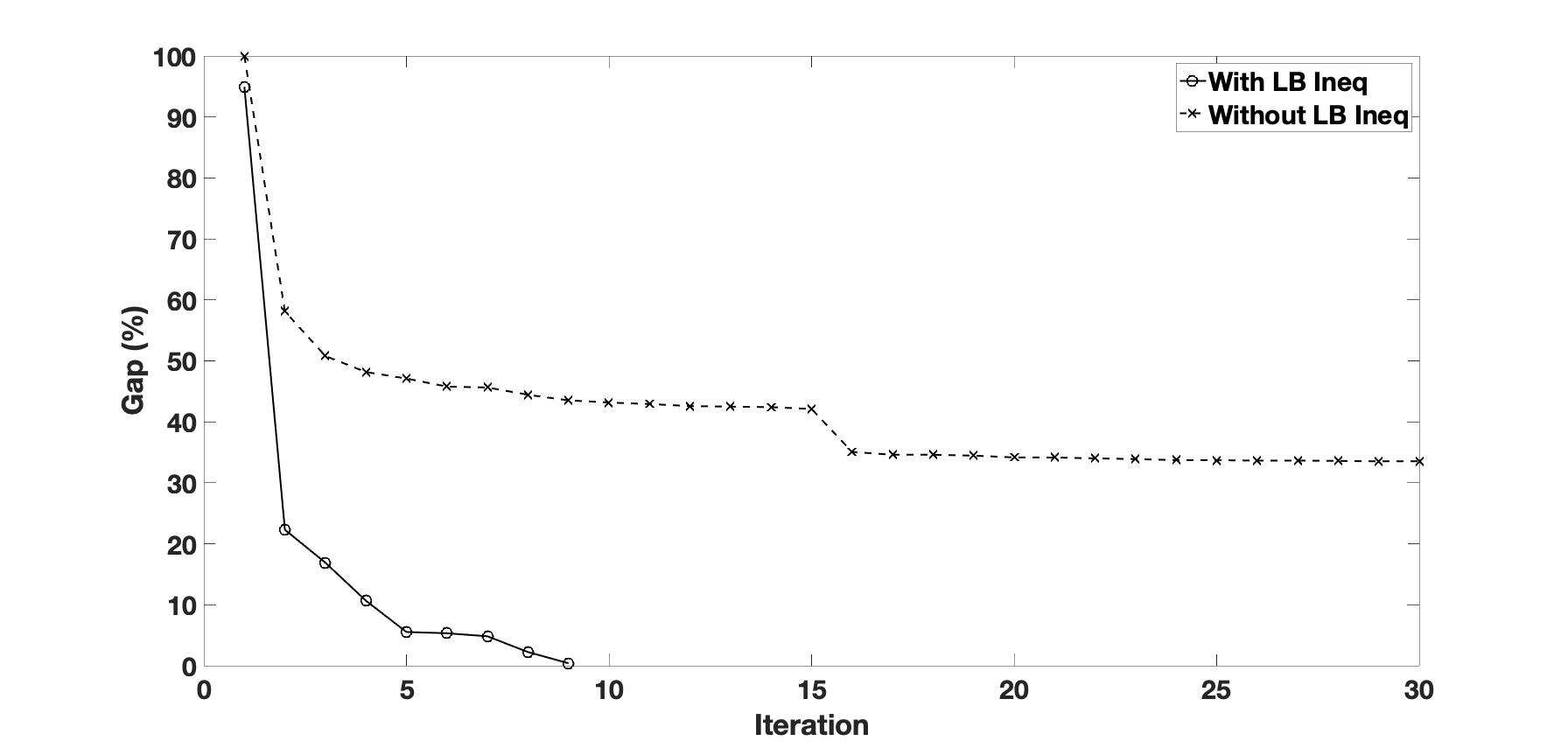}
      \caption{Gap values, MAD-DRO}
      \label{FigVIbGap}
    \end{subfigure}%

        \caption{Comparisons of lower bound and gap values with and without valid LB inequalities \eqref{Vineq1}.}\label{FigVI}
\end{figure}

Next, we analyze the impact of including the valid lower bounding (LB) inequalities \eqref{Vineq1} in the master problem of decomposition algorithm. We first observe that the algorithm takes a very large number of iterations and a longer time until convergence without these LB inequalities. Therefore, in Figure~\ref{FigVI}, we present the LB and gap values with and without inequalities \eqref{Vineq1} from the first 25 iterations. It is obvious that both the lower bound and gap values converge faster when we introduce inequalities  \eqref{Vineq1} into the master problem. Moreover, because of the better bonding effect, the algorithm converges to the optimal solution in fewer iterations and shorter solution times.  
For example, the algorithm takes 10 seconds and 9 iterations to solve the MAD-DRO instance with these inequalities and terminates with a 33\% gap after an hour without these inequalities. The results in this section demonstrate the importance and efficiency of inequalities \eqref{Vineq1}--\eqref{Sym2}.

\subsection{\textbf{Analysis of Optimal Solutions}}\label{sec:Optimal Solutions}

\textcolor{black}{In this section, we compare the optimal solutions of the SP, MAD-DRO, and W-DRO models.}  Given that the SP model can only solve small instances to optimality, for a fair comparison and brevity, we use Instance 3  ($I=15$, $J=15$, and $T=10$) as an example of an average-sized instance. In addition, we present results for Lehigh 1 and Lehigh 2. Table~\ref{table:Optimal} presents the number of MFs (i.e.,  fleet size) for each instance.

\begin{table}[t!]
 \center 
 \footnotesize
   \renewcommand{\arraystretch}{0.3}
  \caption{Optimal number of MFs.}
\begin{tabular}{lccccccccc}
\hline
 \multicolumn{4}{c}{Instance 3 (\textcolor{black}{$\Wb\in [20, 60]$)}} \\
 \hline 
\textbf{Model} & \textbf{$N=10$}  &  \textbf{$N=50$}   & \textbf{$N=100$} \\
\hline
SP&  7 & 8 & 8 \\
W-DRO & 9 & 8 & 8\\
MAD-DRO & 10 & 10 & 10 \\
\\
\hline 
 \multicolumn{4}{c}{\textcolor{black}{Instance 3 ($\Wb\in [50, 100]$)}} \\
 \hline
 SP&  12 & 12 & 12 \\
W-DRO &13 & 13 & 13\\
MAD-DRO & 15 & 15& 15 \\
\hline
 \multicolumn{4}{c}{Lehigh 1} \\
 \hline 
\textbf{Model} & \textbf{$N=10$}  &  \textbf{$N=50$}   & \textbf{$N=100$} \\
\hline
SP&  6 & 7 & 7 \\
W-DRO & 9 & 8 & 8\\
MAD-DRO & 10 & 10 & 10 \\
\\
\hline
 \multicolumn{4}{c}{Lehigh 2} \\
\hline
\textbf{Model} & \textbf{$N=10$}  &  \textbf{$N=50$}   & \textbf{$N=100$} \\
\hline
SP&  5 & 5 & 5 \\
W-DRO & 6 & 6 & 6\\
MAD-DRO & 7 & 7 & 7 \\
\hline																									
\end{tabular} 
\label{table:Optimal}
\end{table}

We observe the following from Table~\ref{table:Optimal}.  First, the MAD-DRO model always activates (schedules) a higher number of MFs than the SP model, and a larger number of MFs than the W-DRO model. By scheduling more MFs, the MAD-DRO  model tends to conservatively mitigate the ambiguity of the demand (reflected by lower shortage and transportation costs reported later in Section~\ref{sec5:OutSample}).  Second, the W-DRO model schedules a larger number of MFs than the SP model, and the difference is significant when the sample size is small ($N=10$). This makes sense as a small sample does not have sufficient distributional information. Thus, in this case, the W-DRO model makes conservative decisions to hedge against ambiguity.  As $N$ increases (i.e., more information becomes available), the W-DRO model often makes less conservative decisions. Consider Instance 3 with $\Wb \in [20, 60]$, for example. The W-DRO model schedules 9 and 8 MFs when $N=10$ and $N=50$, respectively.  \textcolor{black}{Third, we observe that all models scheduled more MFs when we increased the demand's range from $\Wb\in [20, 60]$ to $\Wb\in [50, 100]$ to hedge against the increase in the demand's volume and variability.}

\color{black}

\begin{figure}
\begin{subfigure}[b]{\textwidth}
\includegraphics[width=\textwidth]{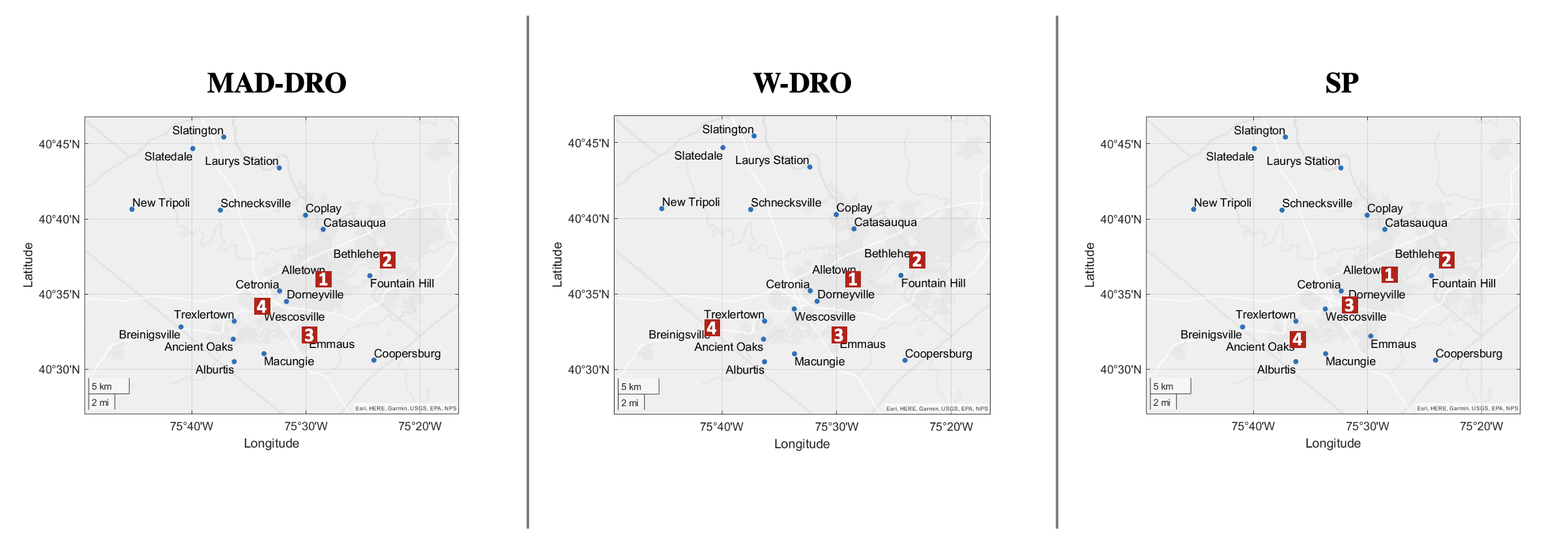}
\caption{Initial locations of the MFs in period 1.}\label{period1}
\end{subfigure}

\begin{subfigure}[b]{\textwidth}
\center 
\includegraphics[width=\textwidth]{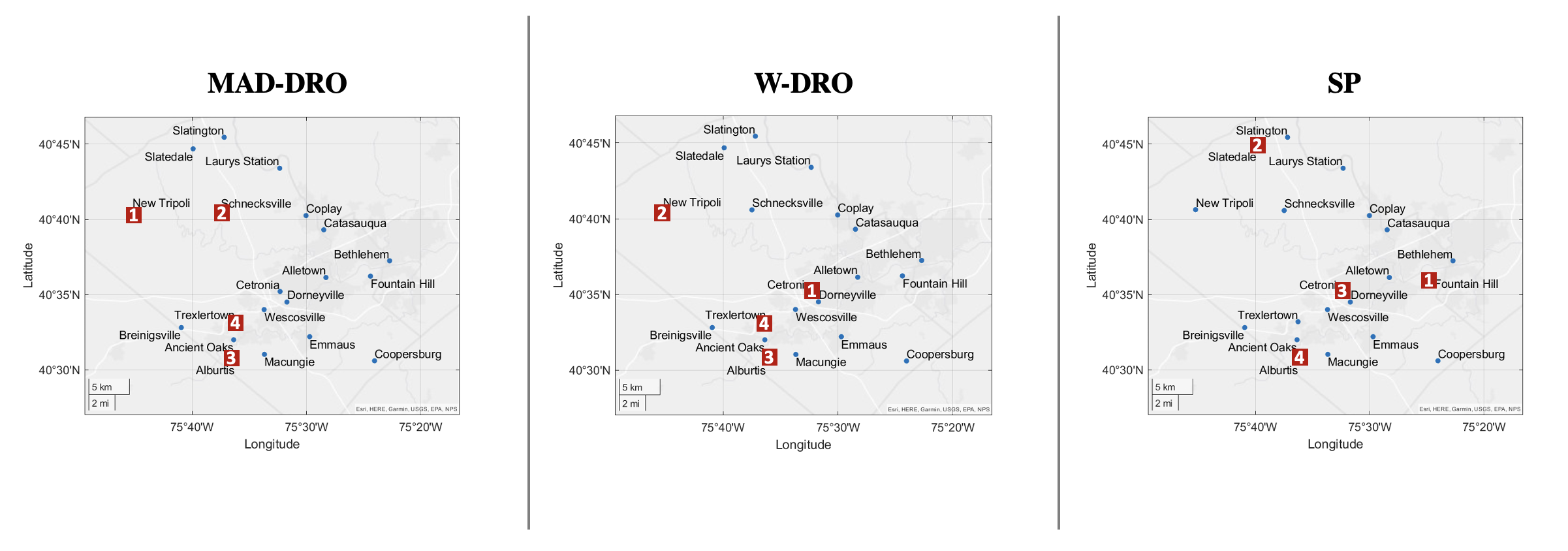}
\caption{Locations of the MFs in period 2.}\label{period2}
\end{subfigure}

\caption{Optimal MF locations in period 1 and period 2 (Lehigh 2 instance). Color code: the red square is MF, and the black circle is a demand node/city.}\label{base_swap}
\end{figure}

Let us now analyze the optimal locations of the MFs. For illustrative purposes and brevity, we use  Lehigh 2 in this analysis. Recall that the SP and DRO models yield different fleet sizes and thus different routing decisions. Therefore, to facilitate the analysis, we first fixed the fleet size to 4 in the three models. Second, to demonstrate how MFs can move to accommodate the change in the demand over time and location, we consider two periods (two days) with the following demand structure. In period 1, we kept the demand structure as described in Section~\ref{sec5.1:instancegen}.  In the second period, we swapped the average demand of the following nodes: Allentown and Alburtis, Bethlehem and Cetronia, Emmaus and Trexlertown, Ancient Oaks and Laurys Station, Catasauqua and  New Tripoli, and Wescosvill and Slatedale. That is, in the second period, we decreased the demand of the 6 nodes with the highest demand (Allentown--Wescosvill) to that of the nodes that generate the lowest demand (Alburtis--Slatedale) and increased the demand of (Alburtis--Slatedale) to that of (Allentown--Wescosvill). We refer to Table~\ref{table:period} in Appendix~\ref{AppexLehigh} for a summery of the average demand of each node in period 1 and period 2. Figure~\ref{base_swap} illustrates the MFs' locations in period 1 (Figure~\ref{period1}) and period 2 (Figure~\ref{period2}). We provide a summary of these results in  Table~\ref{table:OptimalLocations} in Appendix~\ref{AppexLehigh}.

  We observe the following about the initial locations in period 1 (Figure~\ref{period1}). First, all models scheduled  MF  \#1 and \#2 at Allentown and Bethlehem, respectively. This makes sense because, by construction, these nodes generate greater demand than the remaining nodes in period 1. Second, we do not see any MF at or near any of the nodes in the top left of the map (New Tripoli, Slatedale, Slatington, Laury Station, Schnecksville). This makes sense because these nodes generate significantly lower demand than the remaining nodes in period 1.  MF  \#3 and MF\#4  are scheduled at nodes that generate higher demand or near nodes that generate higher demand than the nodes in the top left of the map.  For example, the DRO models scheduled MF \#3 at Emmaus (which has the greatest demand after Allentown and Bethlehem in period 1).  The MAD-DRO model scheduled MF \#4 at Wescoville, while the W-DRO model scheduled this MF at Breinigsville. The SP model scheduled MF \#3 and MF \#4 at Dorneyville and Ancient Oaks. Note that Ancient Oaks generates the highest demand after Allentown, Bethlehem, and Emmaus in period 1. Moreover, Dorneyville, Breinigsville, and Wescoville are closer to demand nodes that generate higher demand in period 1 than the remaining nodes on the top left of the map  (see Table~\ref{table:period} in Appendix~\ref{AppexLehigh}).
  

We make the following observations from the results in period 2 presented in Figure~\ref{period2}. First, it is clear that all MFs moved from their initial locations to other locations in period 2 to accommodate the change in the demand. Second, all models scheduled one MF at Alburtis, where the average demand increased from 9 to 60 (average demand of Allentown in period 1). This makes sense because, in period 2, Alburtis generates the highest demand. Third, the DRO models scheduled one MF at  Trexlertown and one at New Tripoli, where the average demand increased from 8 and 3 to 43 (average demand of Emmaus in period 1) and 23 (average demand of Catasauqua in period 1), respectively. The SP and W-DRO models scheduled one MF at Cetronia, where the demand increased from 8 to 60. 

\color{black}

\subsection{\textbf{Analysis of Solutions Quality}}\label{sec5:OutSample}

\noindent  In this section, we compare the operational performance of the optimal solutions to Instance 3, Lehigh 1, and Lehigh 2 via out-of-sample simulation.  First, we fix the optimal first-stage decisions yielded by each model in the second-stage of the SP.  Then, we solve the second-stage problem in \eqref{2ndstage} with the fixed first-stage decisions and the following sets of $N'=10,000$ out-of-sample data of $W_{i,t}^n$, for all $i \in I, t \in T,$ and $n \in [N']$, to compute the corresponding out-of-sample second-stage cost.

\begin{enumerate}\itemsep0em
\item[Set 1.] \textit{Perfect distributional information}. We use the same settings and distribution (LogN) that we use for generating the $N$ data  in the optimization to generate $N^\prime$ data. This is to simulate the performance when the true distribution is the same as the one used in the optimization.
\item[Set 2.] \textit{Misspecified distributional information}. We follow the same out-of-sample simulation procedure described in \cite{wang2020distributionally} and employed in \cite{shehadeh2020distributionallyTucker} to generate the $N^\prime=10,000$ data.  Specifically, we perturb the distribution of the demand by a parameter $\Delta$ and use a parameterized uniform distribution $U$[$(1-\Delta)\WL,  (1+\Delta)\WU$ ] for which a higher value of $\Delta$ corresponds to a higher variation level.  We apply $\Delta \in \{ 0, 0.25, 0.5\}$ with $\Delta=0$ indicating that we only vary the demand distribution from LogN to Uniform. This is to simulate the performance when the true distribution is different from the one we used in the optimization.  In addition, we generate $N'$ correlated data points with 0.2 and 0.6 correlation coefficients.
\end{enumerate}

\color{black}

\begin{figure}[t!]
 \centering
 \begin{subfigure}[b]{0.5\textwidth}
          \centering
        \includegraphics[width=\textwidth]{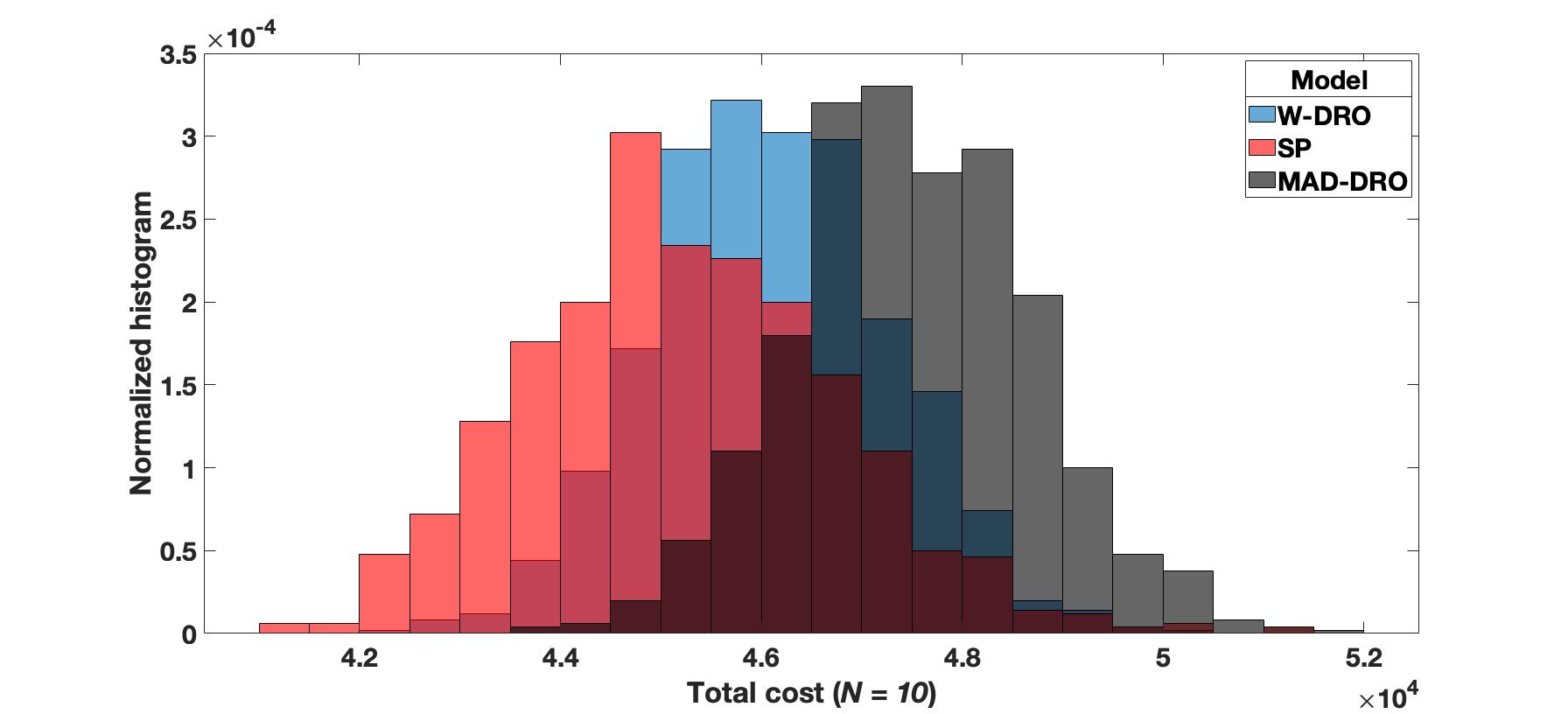}
        \caption{TC (Set 1, LogN)}\label{Inst3_LogNTC}
    \end{subfigure}%
  \begin{subfigure}[b]{0.5\textwidth}
        \centering
        \includegraphics[width=\textwidth]{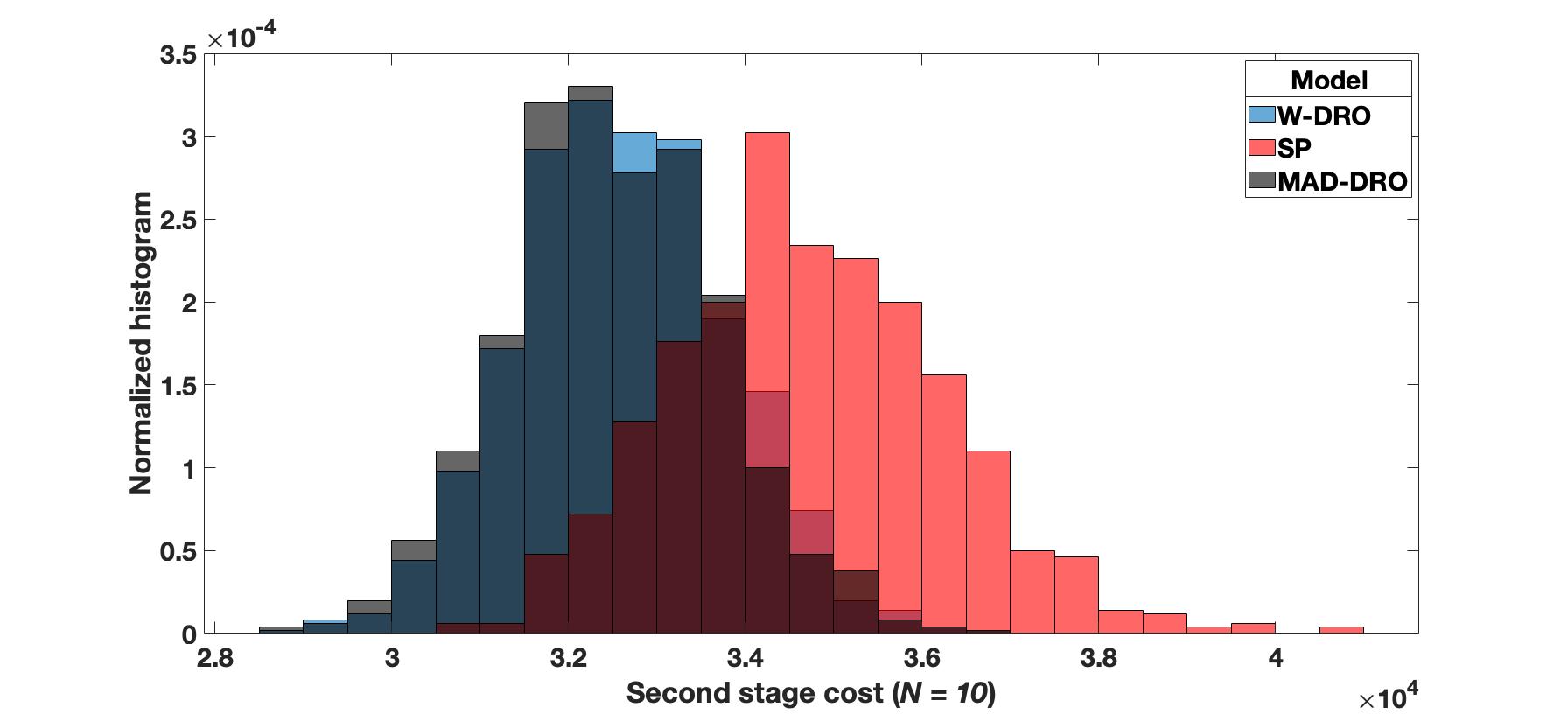}
        \caption{2nd (Set 1, LogN)}\label{Inst3_LogN2nd}
    \end{subfigure}%

  \begin{subfigure}[b]{0.5\textwidth}
          \centering
        \includegraphics[width=\textwidth]{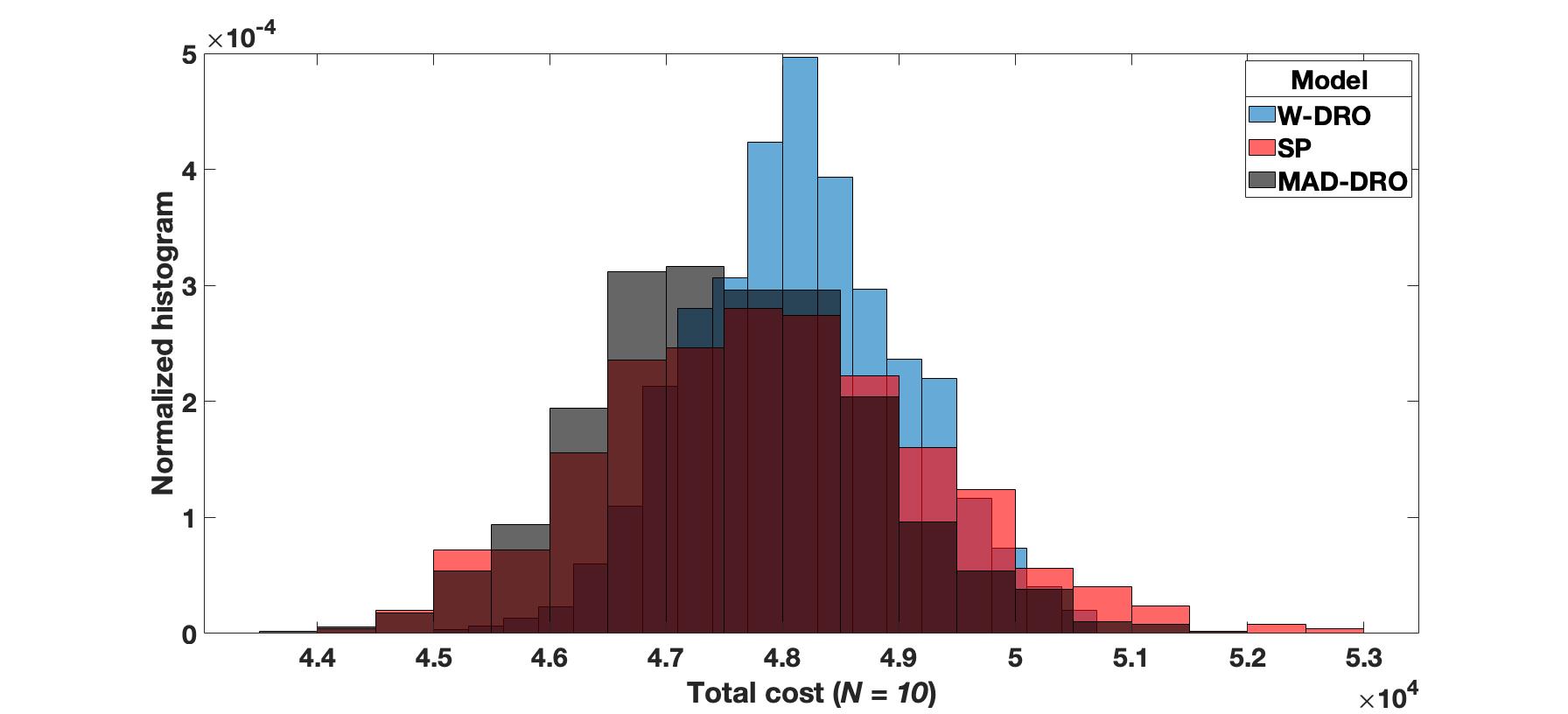}
        \caption{TC (Set 2, $\Delta=0$)}\label{Inst3_Uni0_TC}
    \end{subfigure}%
      \begin{subfigure}[b]{0.5\textwidth}
          \centering
        \includegraphics[width=\textwidth]{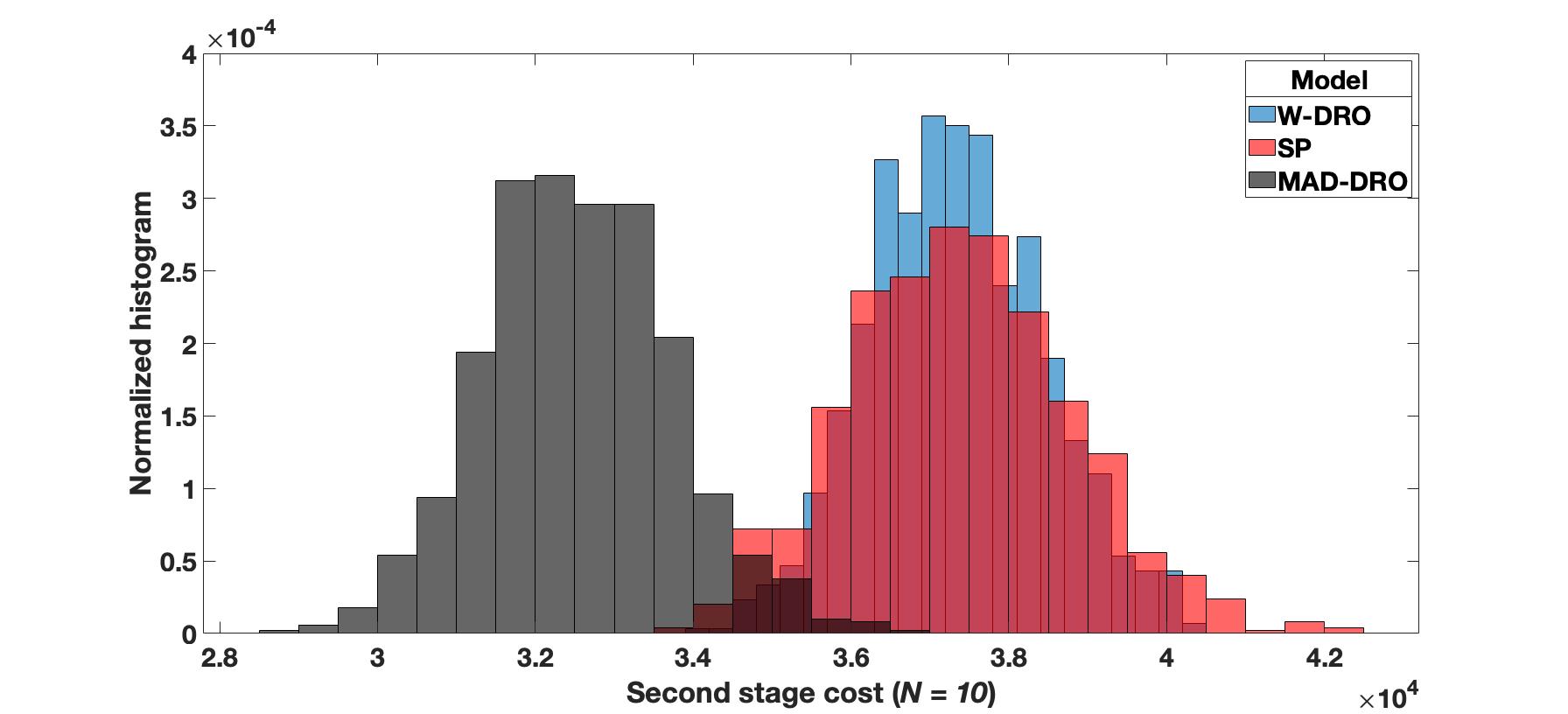}
        \caption{2nd (Set 2, $\Delta=0$)}\label{Inst3_Uni0_2nd}
    \end{subfigure}
 

  \begin{subfigure}[b]{0.5\textwidth}
        \centering
        \includegraphics[width=\textwidth]{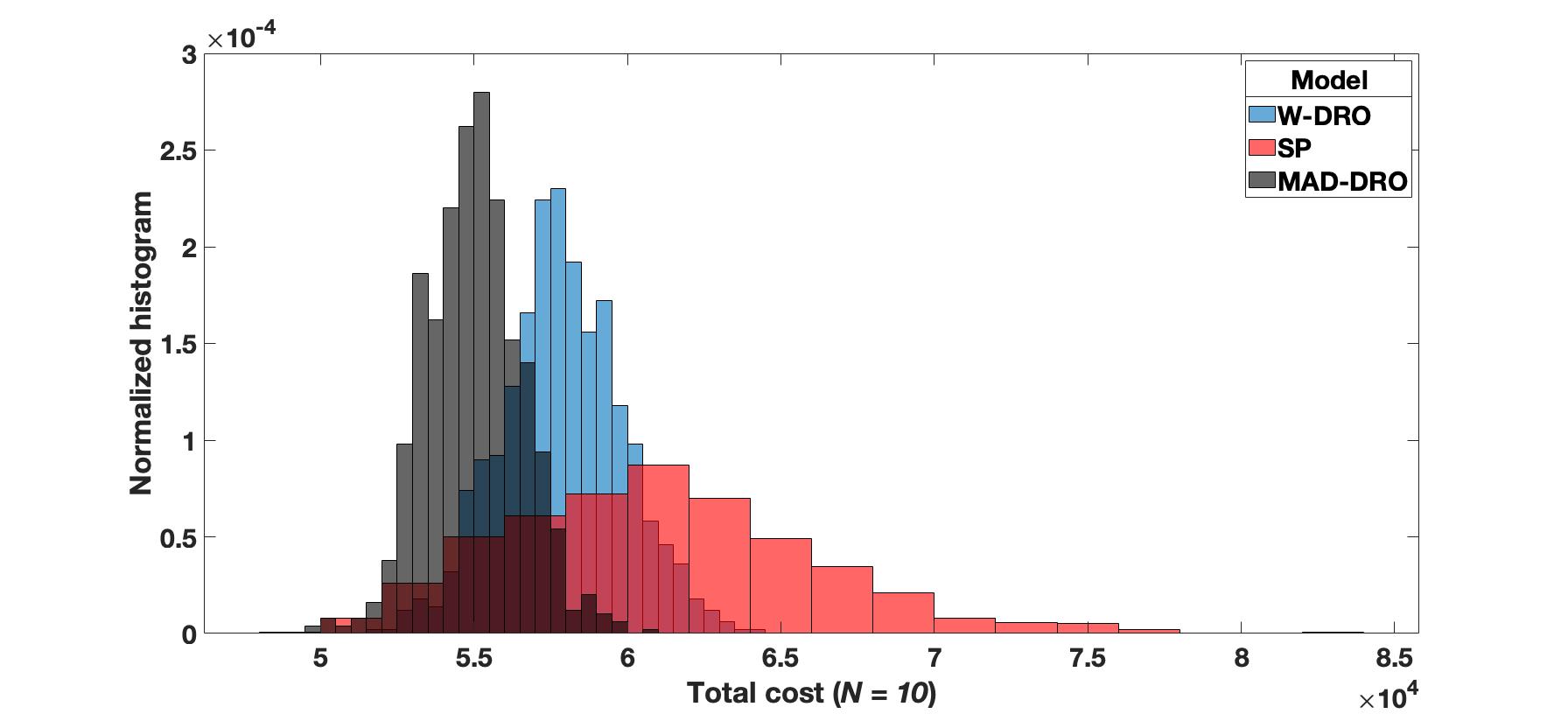}
        \caption{TC (Set 2, $ \Delta=0.25$)}
    \end{subfigure}%
      \begin{subfigure}[b]{0.5\textwidth}
        \centering
        \includegraphics[width=\textwidth]{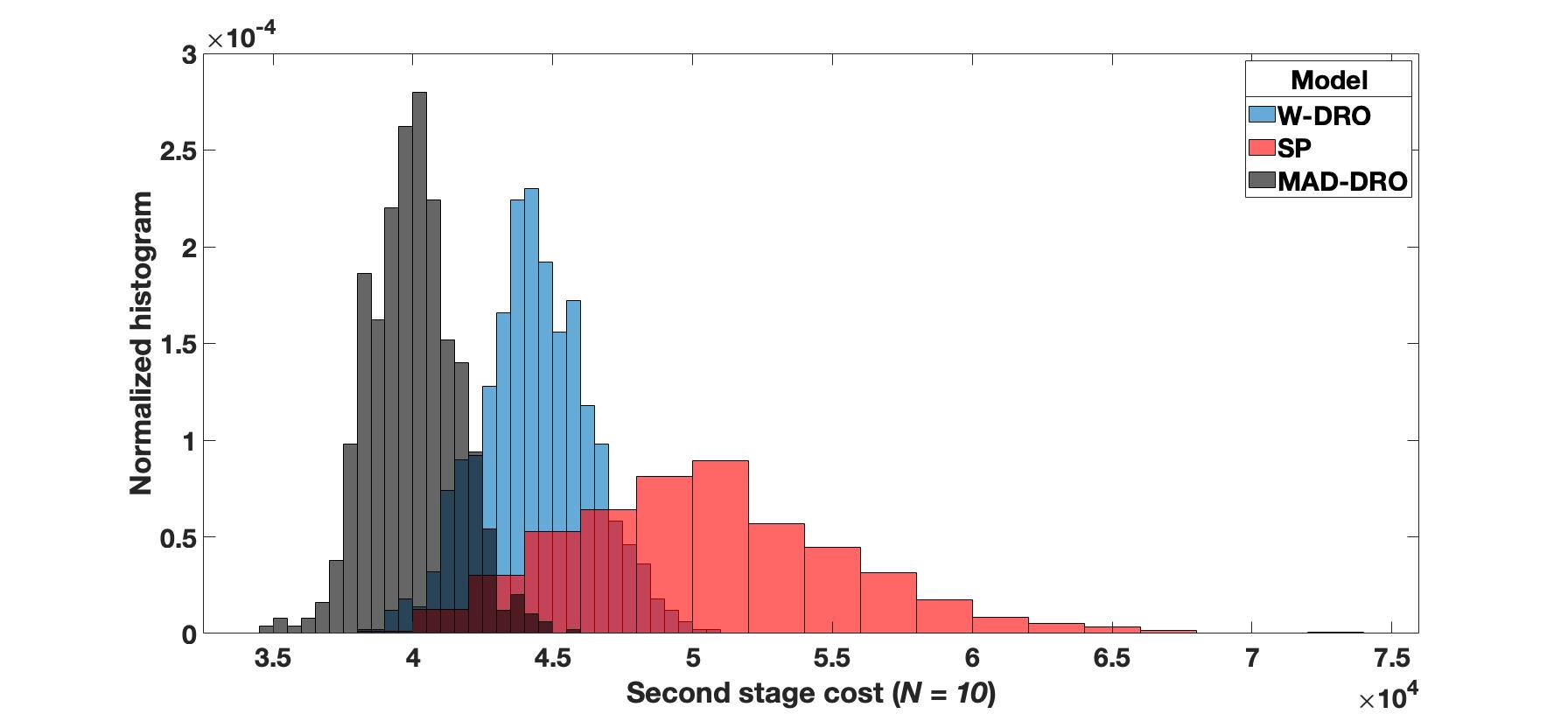}
        \caption{2nd (Set 2, $ \Delta=0.25$)}
    \end{subfigure}%
    
      \begin{subfigure}[b]{0.5\textwidth}
          \centering
        \includegraphics[width=\textwidth]{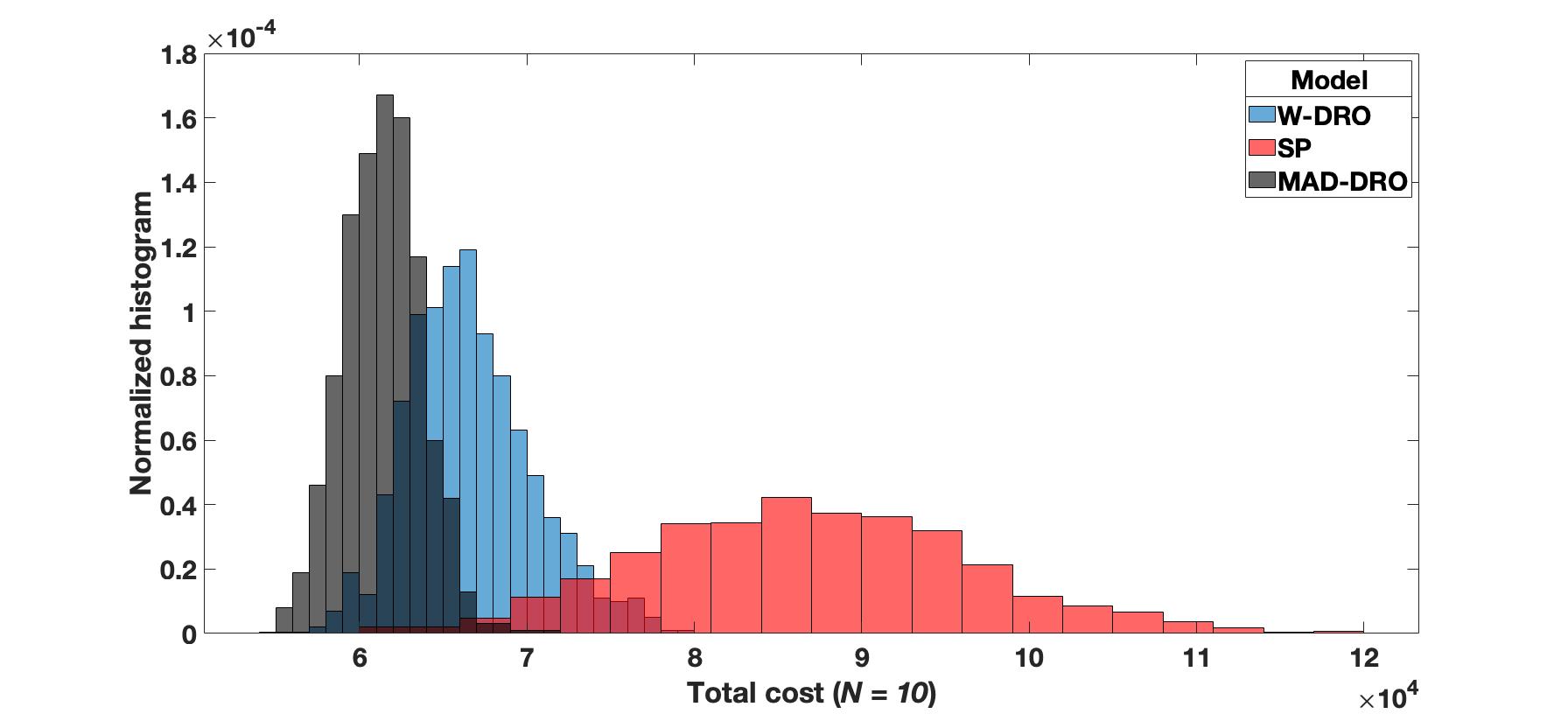}
        \caption{TC (Set 2, $\Delta=0.5$)}
    \end{subfigure}%
    \begin{subfigure}[b]{0.5\textwidth}
          \centering
        \includegraphics[width=\textwidth]{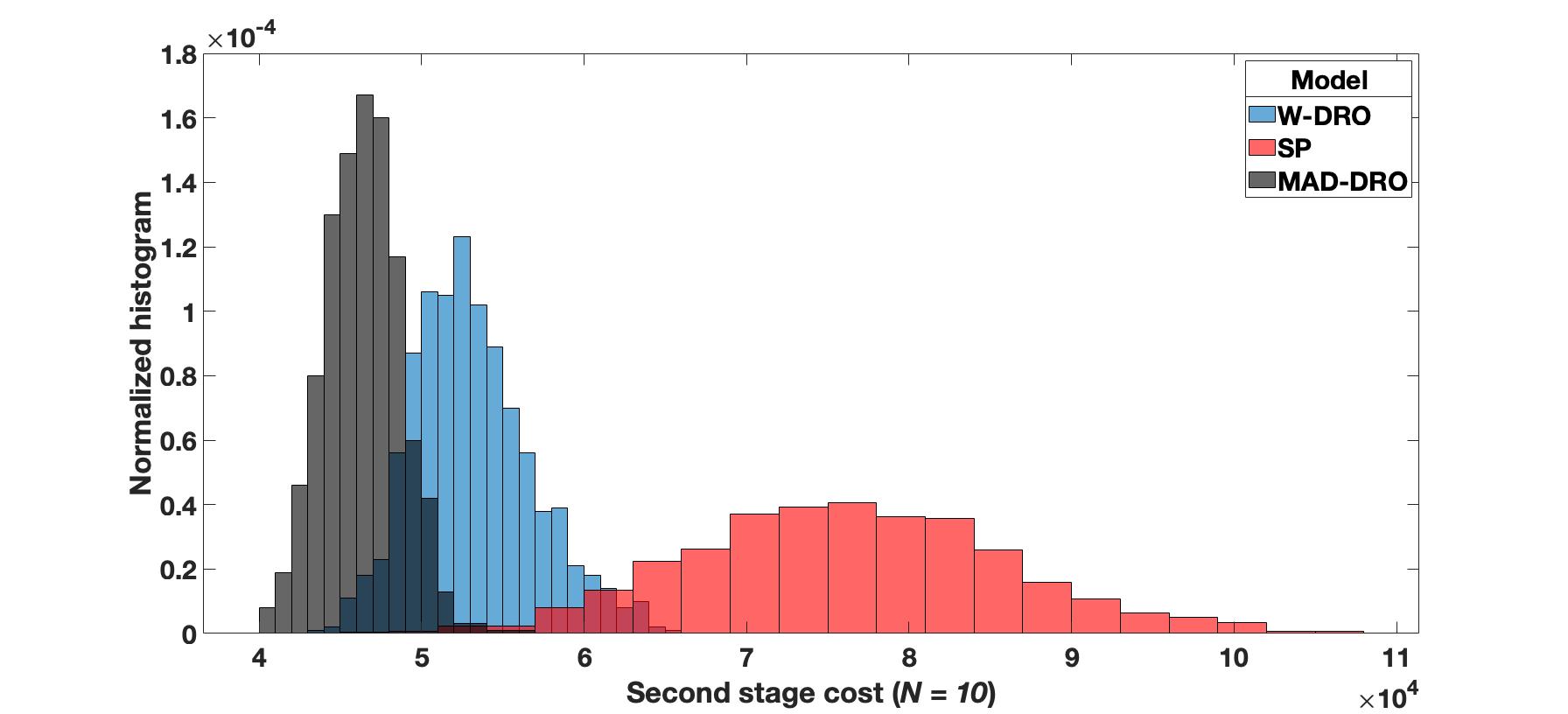}
        \caption{2nd (Set 2, $ \Delta=0.5$)} \label{Inst3_Uni25_2nd}
    \end{subfigure}%
\caption{Normalized histograms of out-of-sample TC and 2nd for Instance 3 ($\Wb \in [20, 60]$, $\pmb{N=10}$) under Set 1 (LogN) and Set 2 (with $\pmb{\Delta \in \{0, 0.25, 0.5\}}$).}\label{Fig3_UniN10_Inst3}
\end{figure}

\begin{figure}[t!]
 \centering
 \begin{subfigure}[b]{0.5\textwidth}
          \centering
        \includegraphics[width=\textwidth]{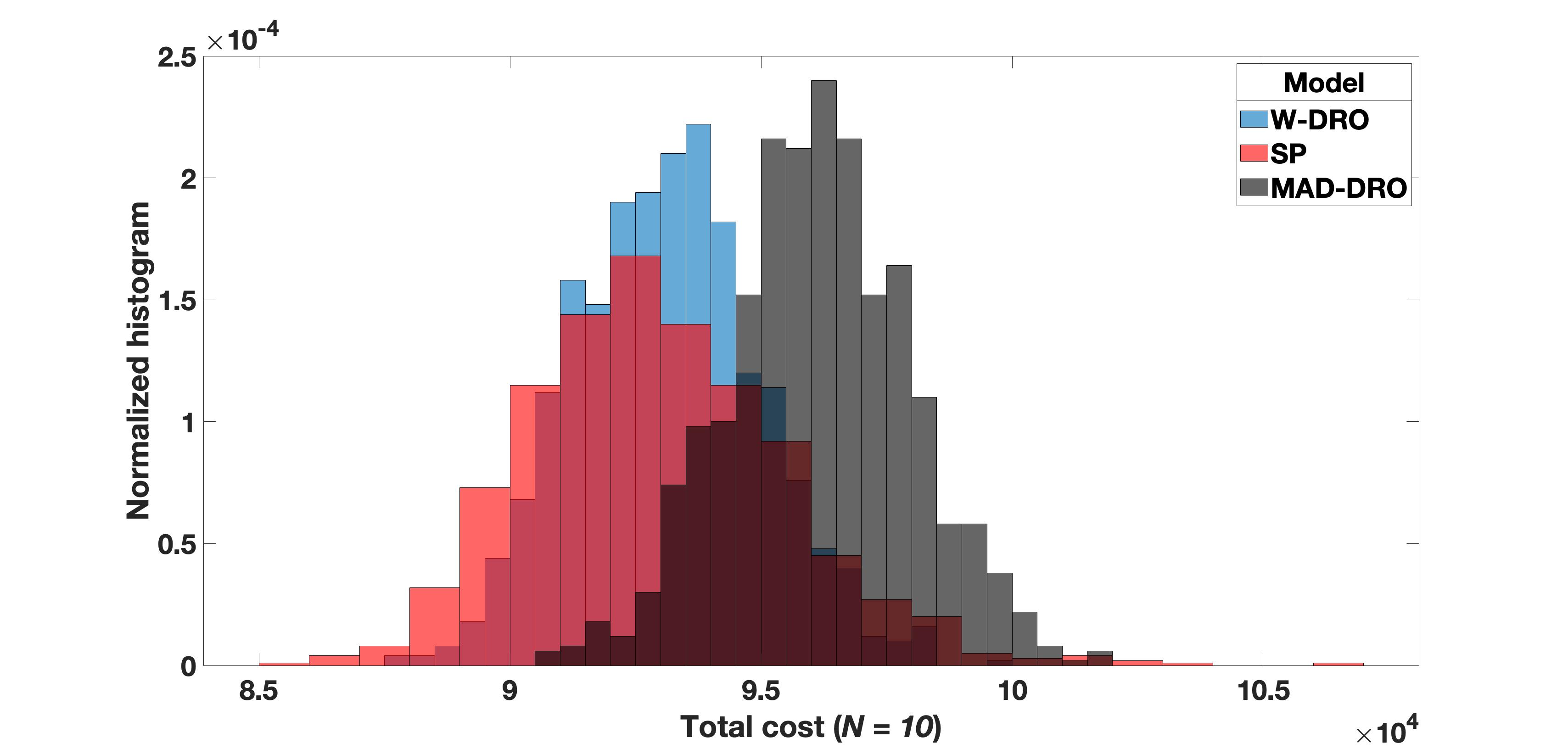}
        \caption{TC (Set 1, LogN)}\label{Inst3_LogNTC_Range2}
    \end{subfigure}%
  \begin{subfigure}[b]{0.5\textwidth}
        \centering
        \includegraphics[width=\textwidth]{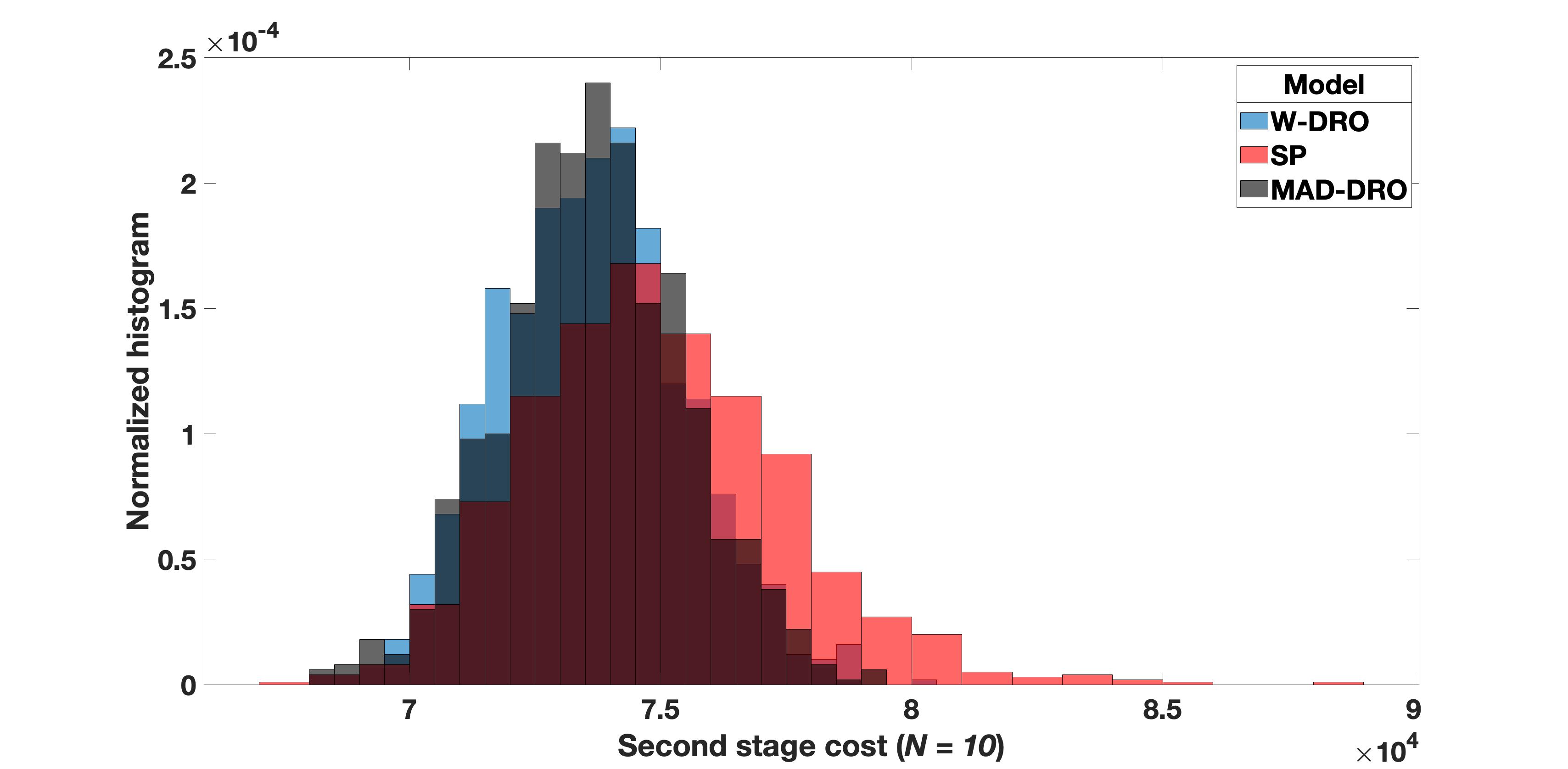}
        \caption{2nd (Set 1, LogN)}\label{Inst3_LogN2nd_Range2}
    \end{subfigure}%

  \begin{subfigure}[b]{0.5\textwidth}
          \centering
        \includegraphics[width=\textwidth]{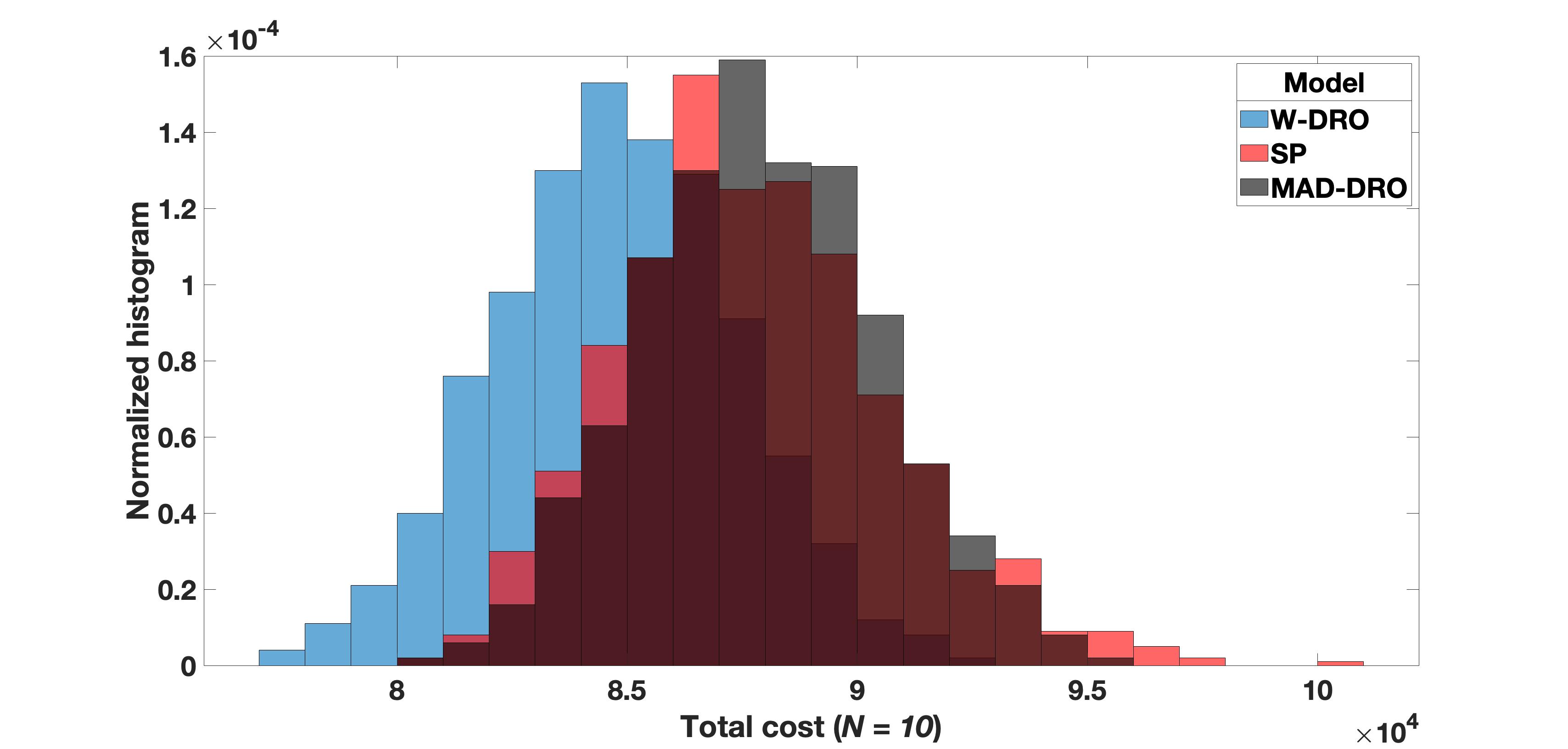}
        \caption{TC (Set 2, $\Delta=0$)}\label{Inst3_Uni0_TC_Range2}
    \end{subfigure}%
      \begin{subfigure}[b]{0.5\textwidth}
          \centering
        \includegraphics[width=\textwidth]{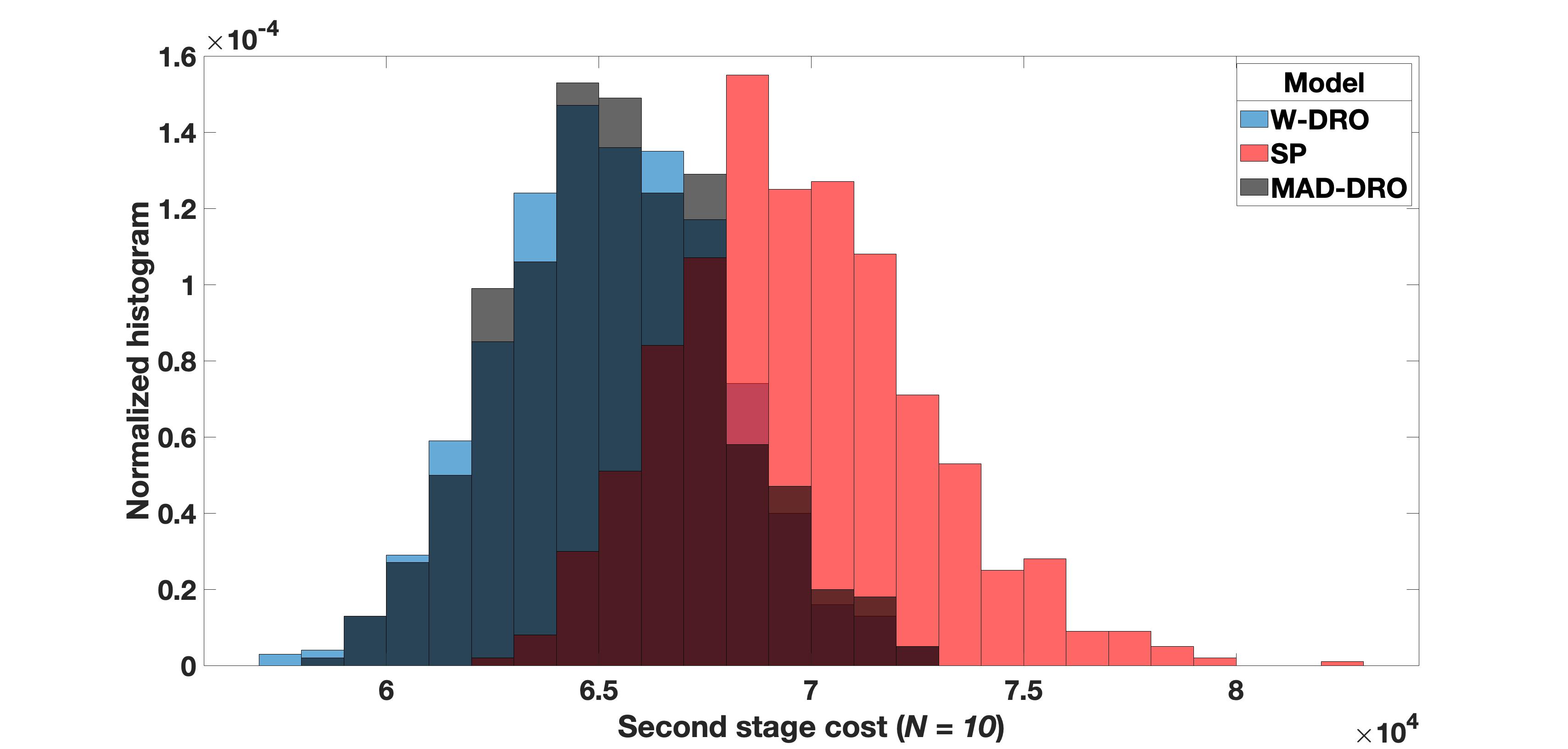} 
        \caption{2nd (Set 2, $\Delta=0$)}\label{Inst3_Uni0_2nd_Range2}
    \end{subfigure}
 

  \begin{subfigure}[b]{0.5\textwidth}
        \centering
        \includegraphics[width=\textwidth]{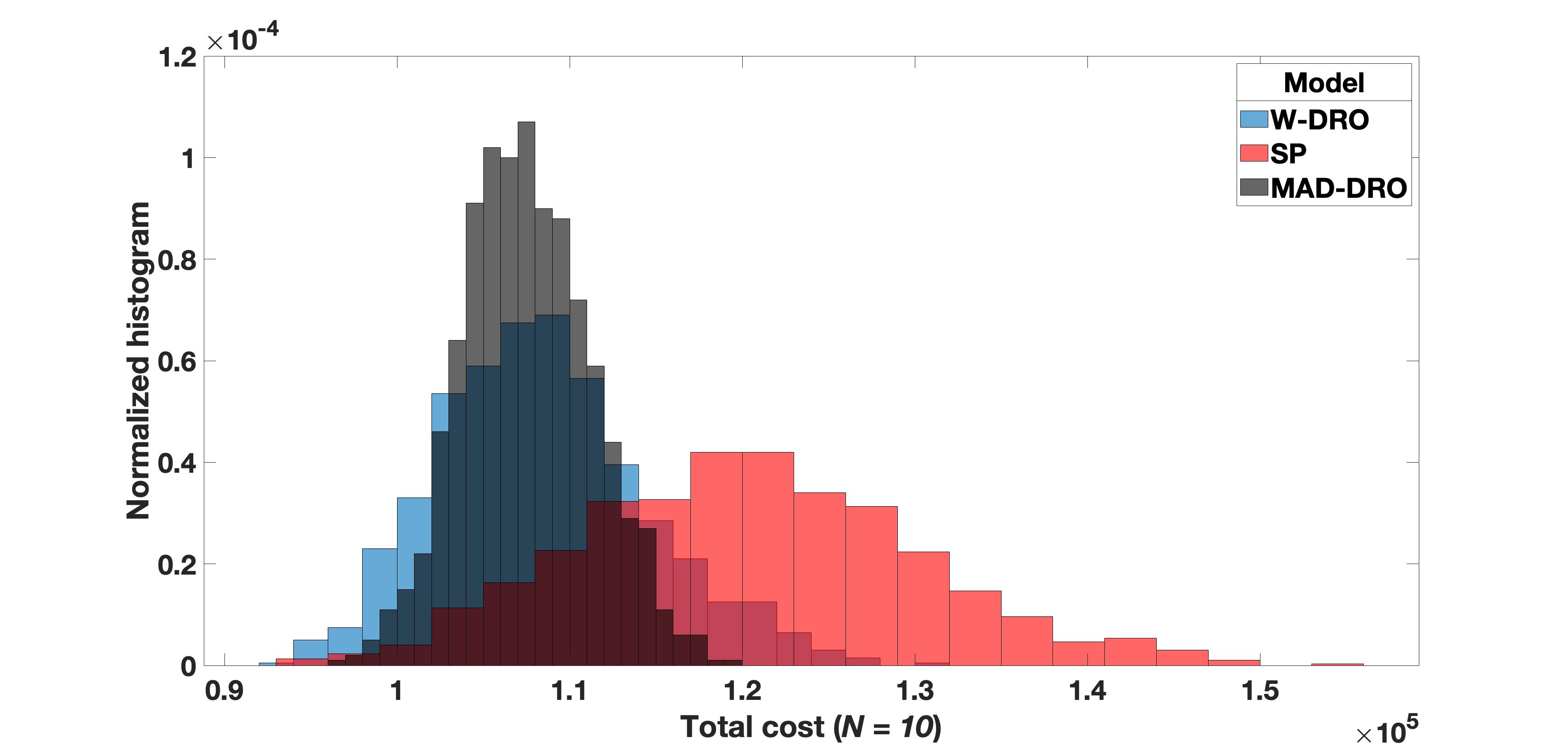}
        \caption{TC (Set 2, $ \Delta=0.25$)}
    \end{subfigure}%
      \begin{subfigure}[b]{0.5\textwidth}
        \centering
        \includegraphics[width=\textwidth]{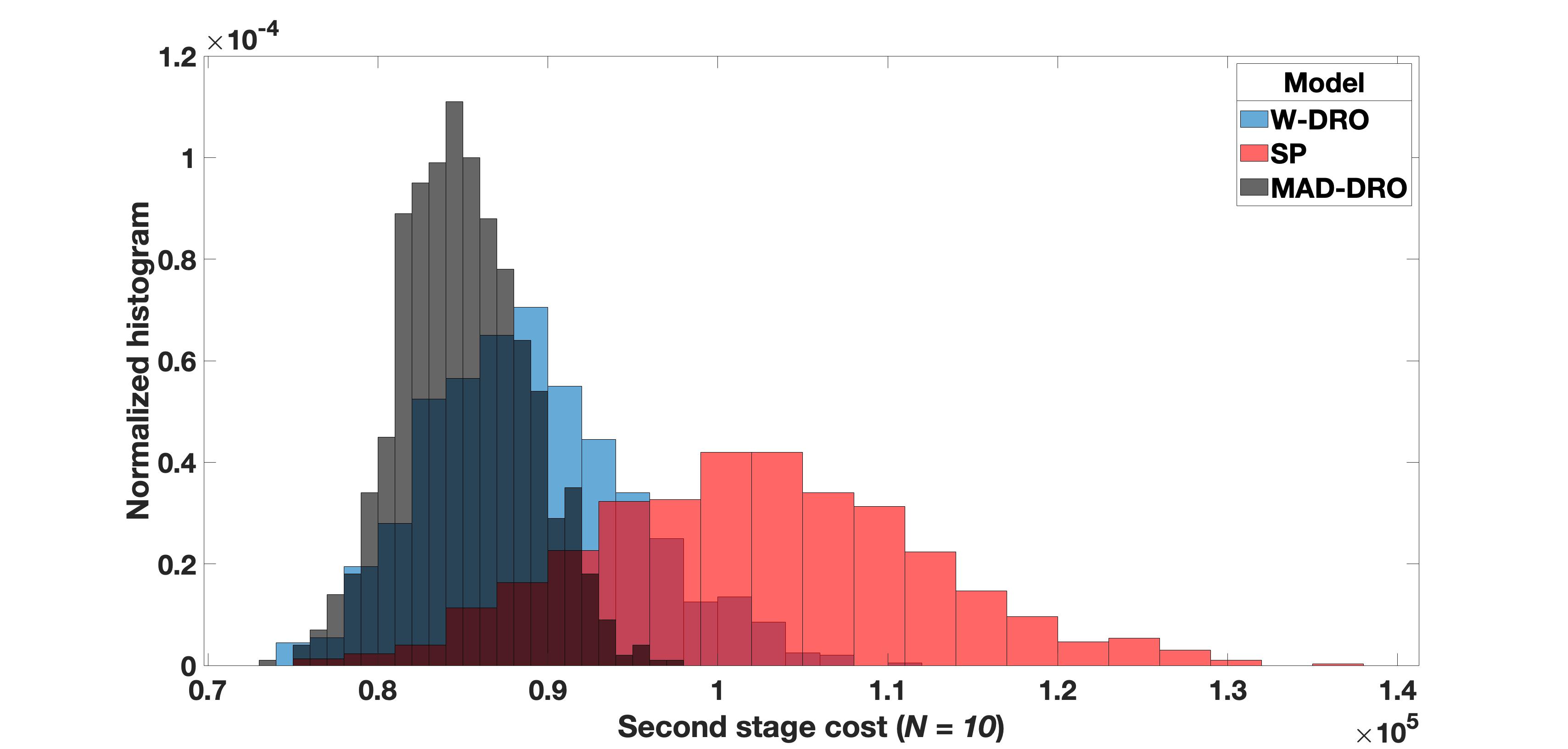}
        \caption{2nd (Set 2, $ \Delta=0.25$)}
    \end{subfigure}%
    
      \begin{subfigure}[b]{0.5\textwidth}
          \centering
        \includegraphics[width=\textwidth]{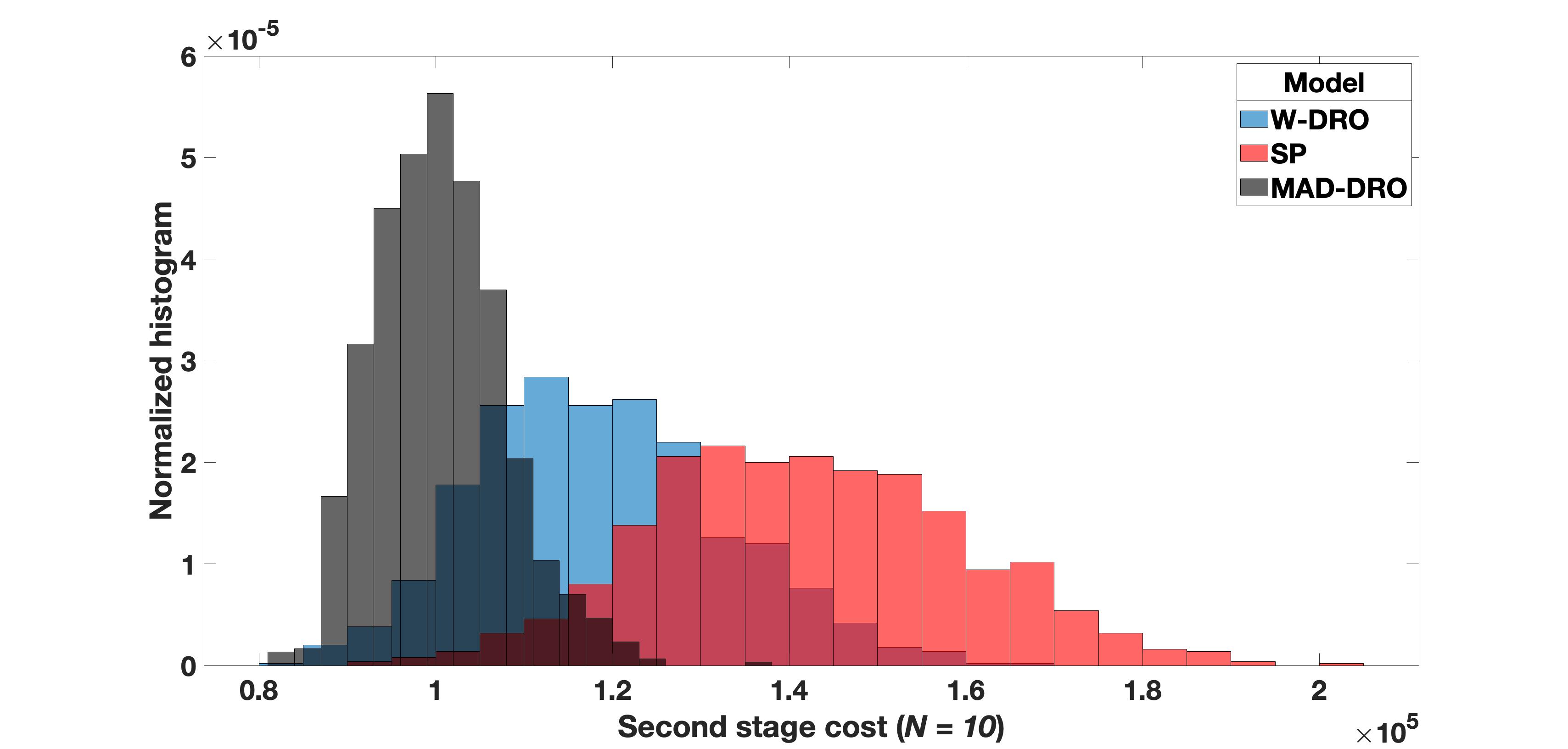}
        \caption{TC (Set 2, $\Delta=0.5$)}
    \end{subfigure}%
    \begin{subfigure}[b]{0.5\textwidth}
          \centering
        \includegraphics[width=\textwidth]{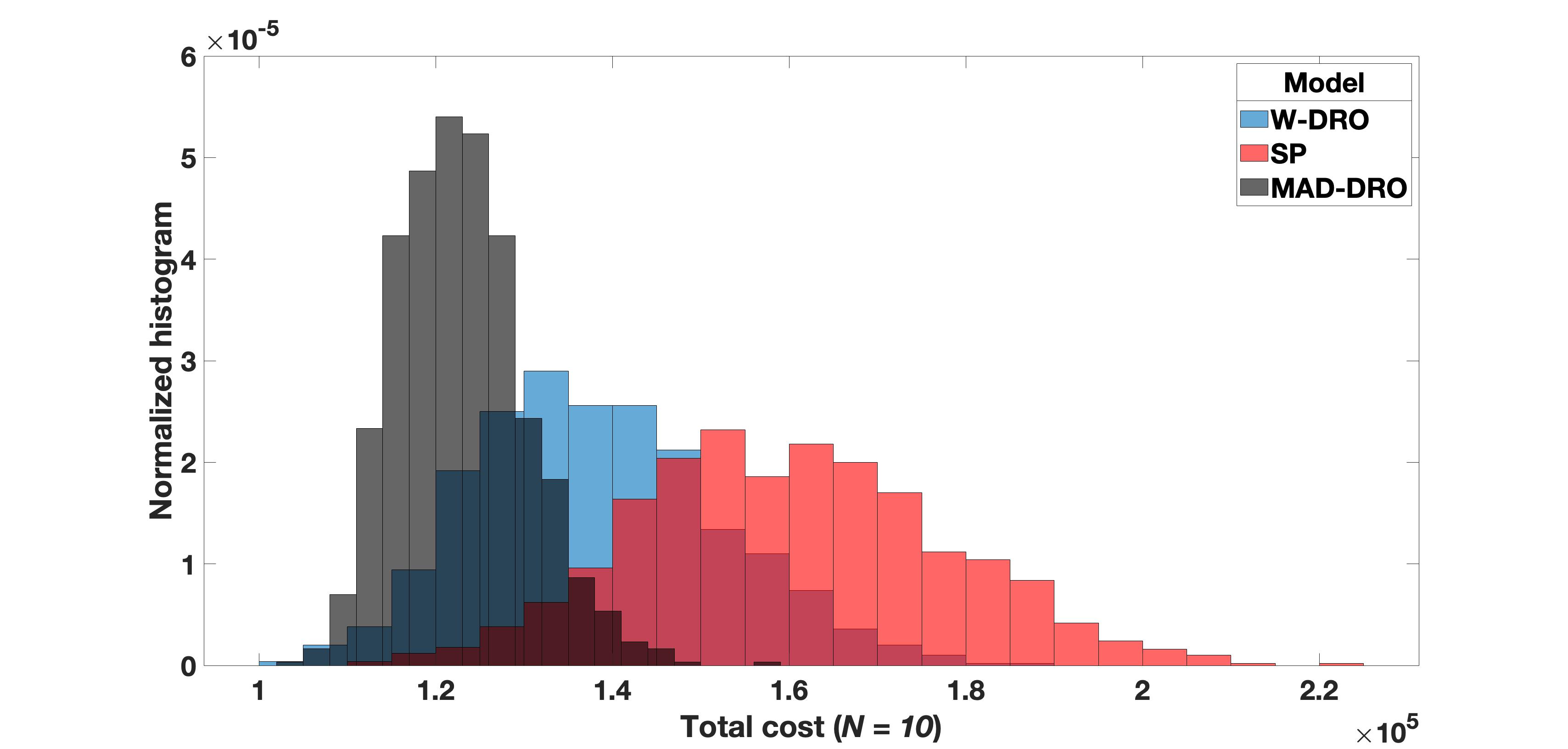}
        \caption{2nd (Set 2, $ \Delta=0.5$)} \label{Inst3_Uni50_2nd_Range2}
    \end{subfigure}%
\caption{Normalized histograms of out-of-sample TC and 2nd for Instance 3 ($\Wb \in [50, 100]$, $\pmb{N=10}$) under Set 1 (LogN) and Set 2 (with $\pmb{\Delta \in \{0, 0.25, 0.5\}}$). }\label{Fig3_UniN10_Inst3_Range2}
\end{figure}

\begin{figure}[t!]
 \centering
   \begin{subfigure}[b]{0.5\textwidth}
          \centering
        \includegraphics[width=\textwidth]{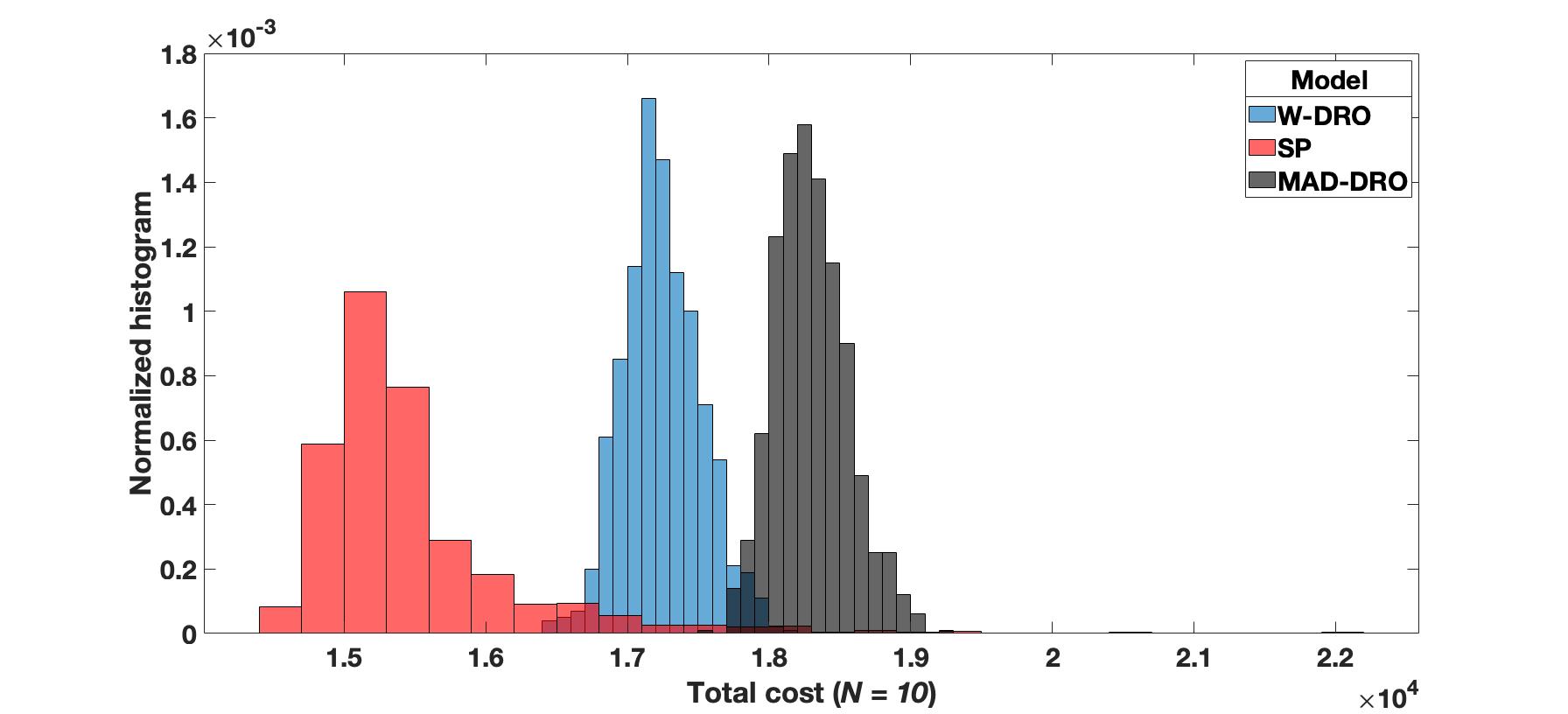}
        \caption{TC (Set 1, LogN)}\label{Lehigh1_LogNTC}
    \end{subfigure}%
      \begin{subfigure}[b]{0.5\textwidth}
          \centering
        \includegraphics[width=\textwidth]{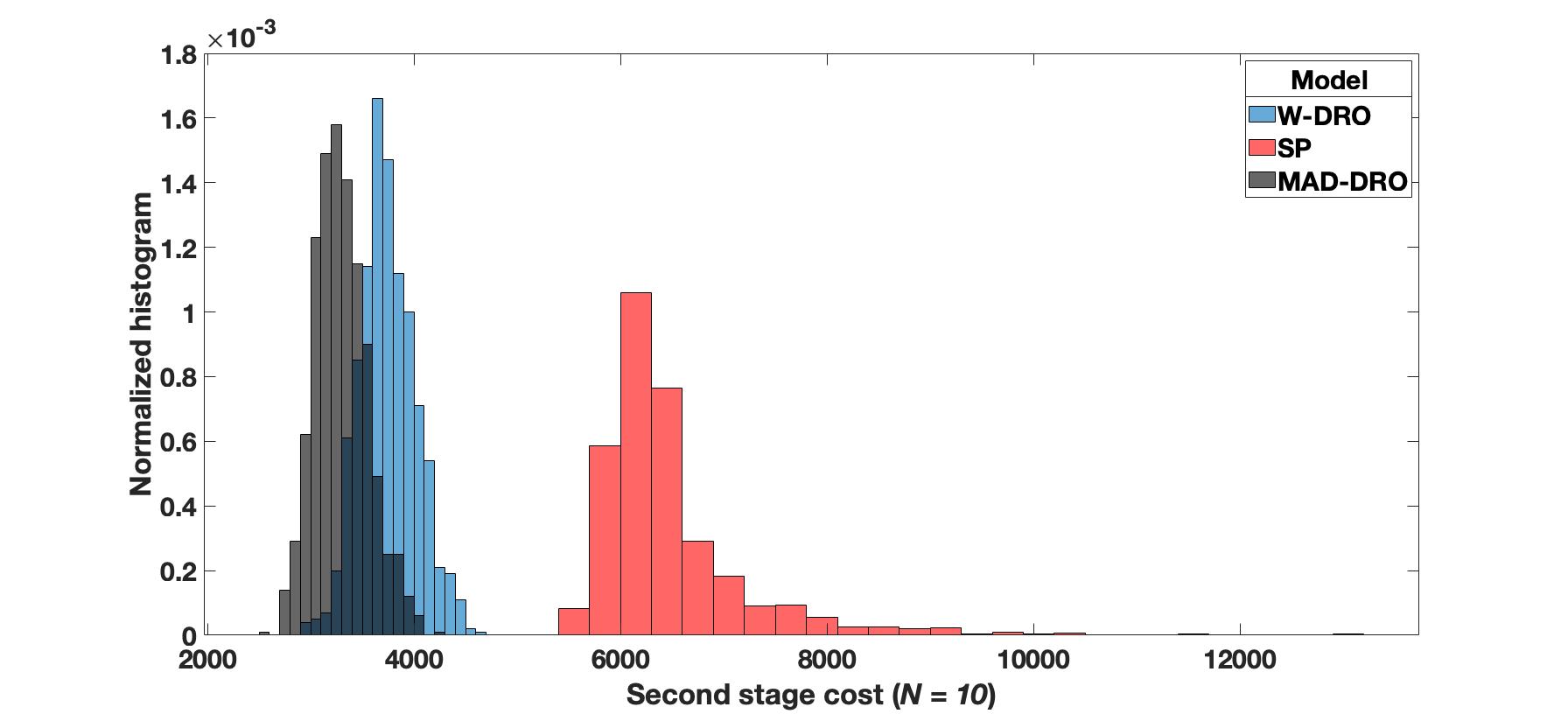}
        \caption{2nd (Set 1, LogN)}\label{Lehigh1_LogN2nd}
    \end{subfigure}%
    
  \begin{subfigure}[b]{0.5\textwidth}
          \centering
        \includegraphics[width=\textwidth]{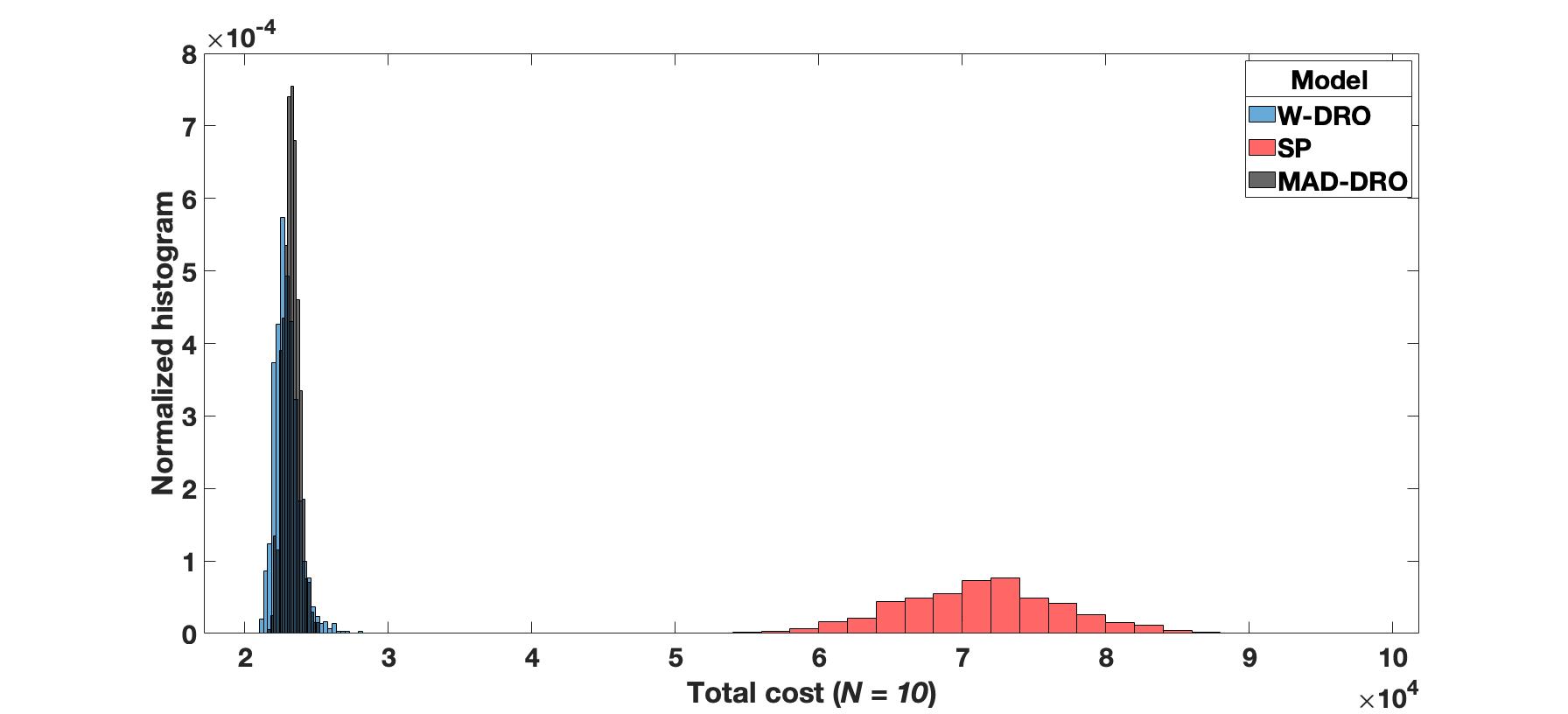}
        \caption{TC (Set 2, $\Delta=0$)}\label{Lehigh1_Uni0TC}
    \end{subfigure}%
      \begin{subfigure}[b]{0.5\textwidth}
          \centering
        \includegraphics[width=\textwidth]{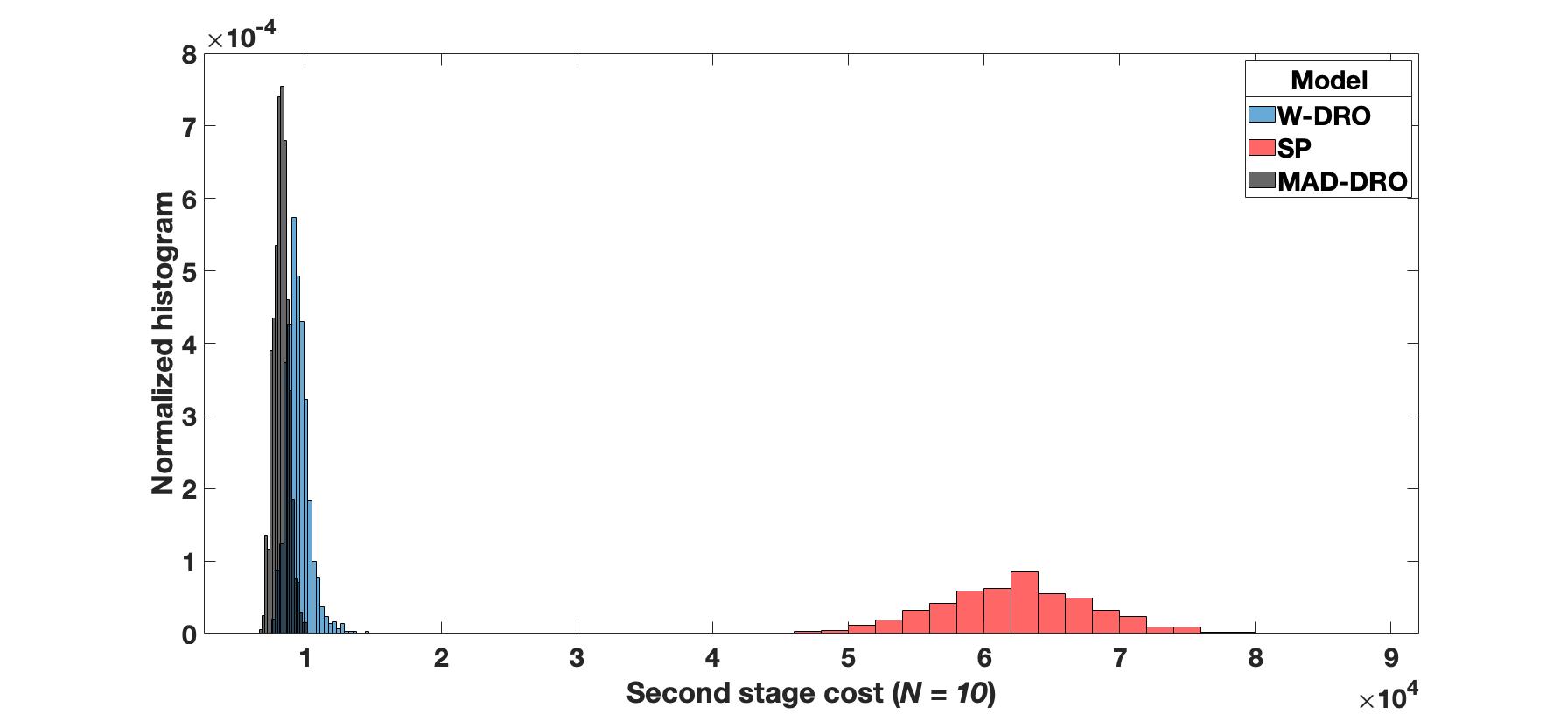}
        \caption{2nd  (Set 2, $\Delta=0$)}
    \end{subfigure}%
    

  \begin{subfigure}[b]{0.5\textwidth}
        \centering
        \includegraphics[width=\textwidth]{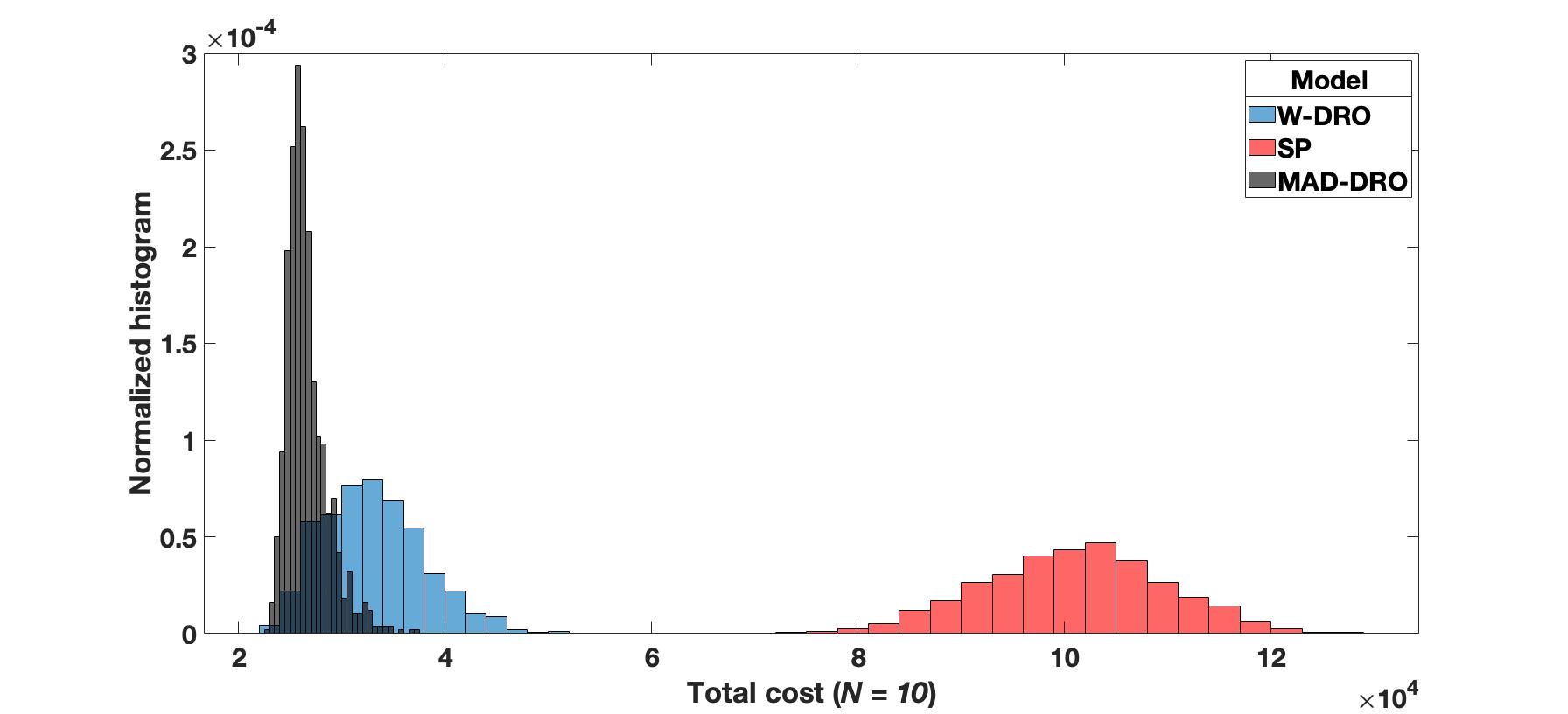}
        \caption{TC  (Set 2, $\Delta=0.25$)}
    \end{subfigure}%
      \begin{subfigure}[b]{0.5\textwidth}
        \centering
        \includegraphics[width=\textwidth]{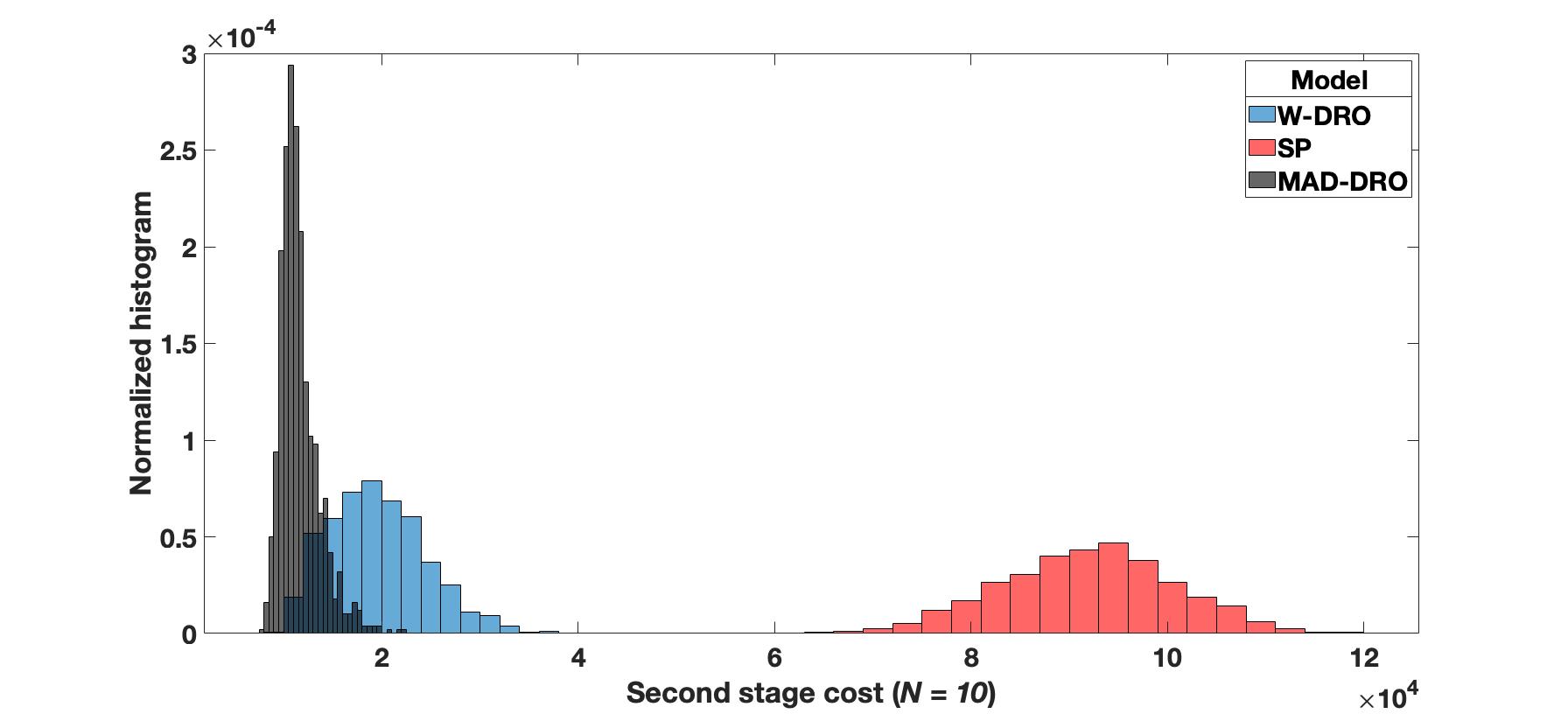}
        \caption{2nd  (Set 2, $ \Delta=0.25$)}
    \end{subfigure}%

      \begin{subfigure}[b]{0.5\textwidth}
          \centering
        \includegraphics[width=\textwidth]{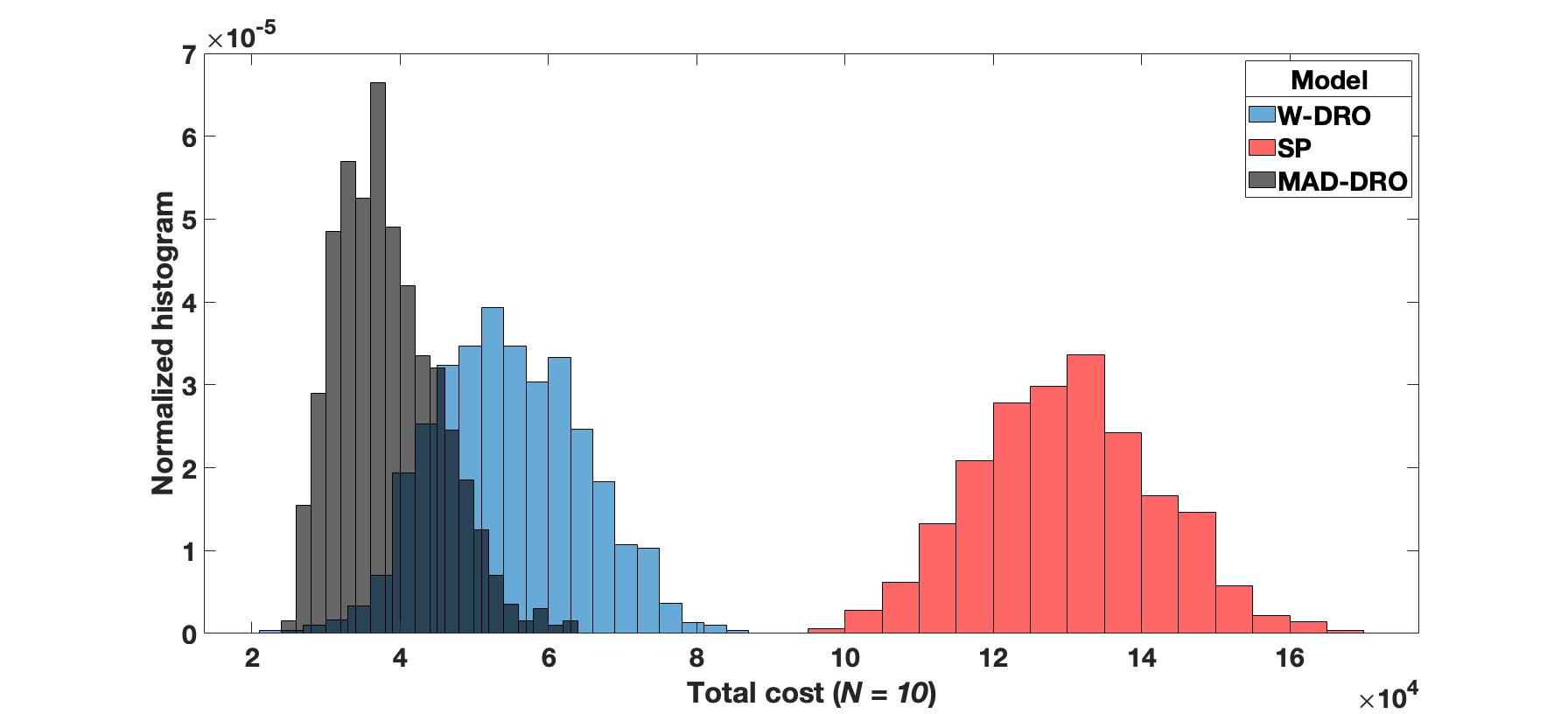}
        \caption{TC  (Set 2, $ \Delta=0.5$)}
    \end{subfigure}%
    \begin{subfigure}[b]{0.5\textwidth}
          \centering
        \includegraphics[width=\textwidth]{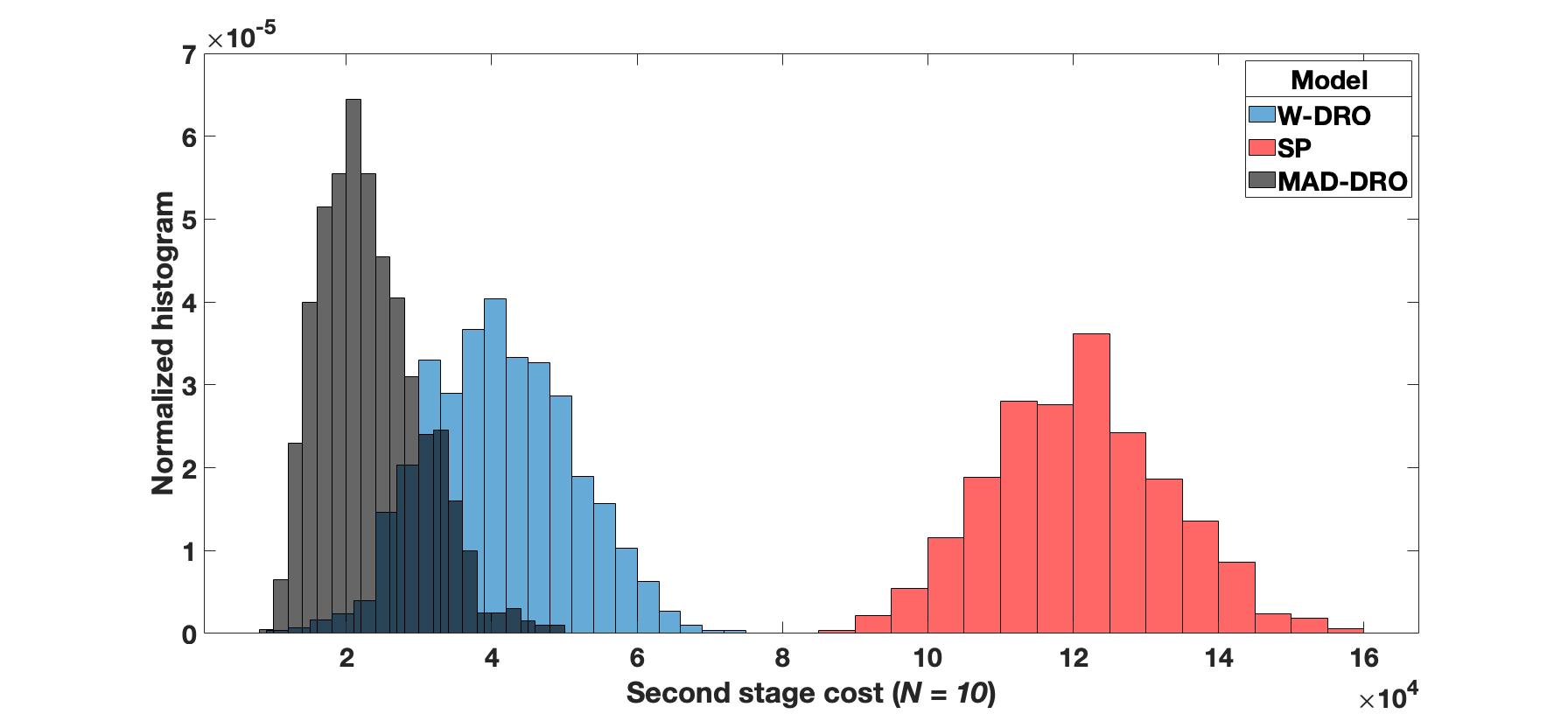}
        \caption{2nd  (Set 2, $ \Delta=0.5$)}\label{Lehigh1_Uni50TC}
    \end{subfigure}%
\caption{Normalized histograms of out-of-sample TC and 2nd for Lehigh 1 ($\pmb{N=10}$) under Set 1 (LogN) and Set 2 (with $\pmb{\Delta \in \{0, 0.25, 0.5\}}$).}\label{Fig3_N10_Lehigh1}
\end{figure}

For brevity, we next discuss simulation results for the solutions obtained with $N=10$. We observe similar results for solutions obtained with $N=50$ (see Appendix~\ref{Appex:AdditionalOut} for these results).  In Figures~\ref{Fig3_UniN10_Inst3}, \ref{Fig3_UniN10_Inst3_Range2}, and \ref{Fig3_N10_Lehigh1}, we present the normalized histograms of out-of-sample total costs (TC) and second-stage costs (2nd) for Instance 3 (with $\Wb \in [20,60]$, $N=10$), Instance 3 (with $\Wb \in [50, 100]$, $N=10$), and Lehigh 1 ($N=10$). We obtained similar results for Lehigh 2 (see Appendix~\ref{Appex:AdditionalOut}). We computed TC as TC$=$first-stage cost+out-of-sample second-stage cost.


Let us first analyze simulation results under Set 1 (i.e., perfect distributional information case) presented in Figures~\ref{Inst3_LogNTC}--\ref{Inst3_LogN2nd},  \ref{Inst3_LogNTC_Range2}--\ref{Inst3_LogN2nd_Range2}, and \ref{Lehigh1_LogNTC}--\ref{Lehigh1_LogN2nd}.  The MAD-DRO model yields a higher TC on average and at upper quantiles than the W-DRO and SP models because it schedules more MFs and thus yields a higher fixed cost (i.e.,  cost of establishing the MF fleet). The W-DRO model yields a slightly higher TC than the SP model because it schedules more MFs (and thus yeild a higher fixed cost). However, the DRO models yield significantly lower second-stage (transportation and unmet demand) costs on average and at all quantiles than the SP model. In addition, the MAD-DRO model yields a lower second-stage cost than the W-DRO model on average and at all quantile, especially for Lehigh 1. Note that a lower second-stage cost indicates a better operational performance (i.e., lower shortage and transportation costs) and thus has a significant practical impact. These results suggest that there are benefits to using the DRO models even when we have perfect distributional information.

We observe the following from simulation results under Set 2 (i.e., misspecified distributional information case) presented in Figures~\ref{Inst3_Uni0_TC}--\ref{Inst3_Uni25_2nd}, \ref{Inst3_Uni0_TC_Range2}--\ref{Inst3_Uni50_2nd_Range2}, and \ref{Lehigh1_Uni0TC}--\ref{Lehigh1_Uni50TC}.  It is clear from these figures that the DRO models consistently outperform the SP model under all levels of variation ($\Delta$) and across the criteria of mean and all quantiles of the the total and second-stage costs.  Interestingly, the DRO models yield substantially lower TC and 2nd than the SP model for Lehigh 1 and Instance 3 (with $\Wb \in [50, 100]$), which have higher demand volume and variations. In addition, the MAD-DRO model yields lower second-stage costs for Instance 3 and substantially lower total and second-stage costs for Lehigh 1.  Finally, the MAD-DRO solutions appear to be more stable with a significantly smaller standard deviation  (i.e., variations) in the total and second-stage costs than the other considered models.  The superior performance of the DRO models reflects the value of modeling uncertainty and distributional ambiguity of the demand.

\color{black}

\color{black}

Next, we investigate the value of distributional robustness from the perspective of out-of-sample disappointment, which measures the extent to which the out-of-sample cost exceeds the model's optimal value \citep{Van-Parys_et_al:2021, wang2020distributionally}. We define OPT and TC as the model's optimal value and the out-of-sample objective value, respectively. That is, OPT and TC can be considered as the estimated and actual costs of implementing the model's optimal solutions in practice, respectively.  Using this notation, we define the out-of-sample disappointment as in \cite{wang2020distributionally} as follows.
\begin{align}
\max \left \{ \frac{\text{TC}-\text{OPT}}{\text{OPT}}, 0 \right \} \times 100\%.
\end{align}
A disappointment of zero indicates that the model's optimal value is equal to or larger than the out-of-sample (actual) cost (i.e., TC$\leq$OPT). This, in turn, indicates that the model is more conservative and avoids underestimating costs. In contrast, a larger disappointment implies a higher level of over-optimism because, in this case, the actual cost (TC) of implementing the optimal solution of a model is larger than the estimated cost (OPT).

\begin{figure}

 \begin{subfigure}[b]{0.5\textwidth}
        \centering
        \includegraphics[width=\textwidth]{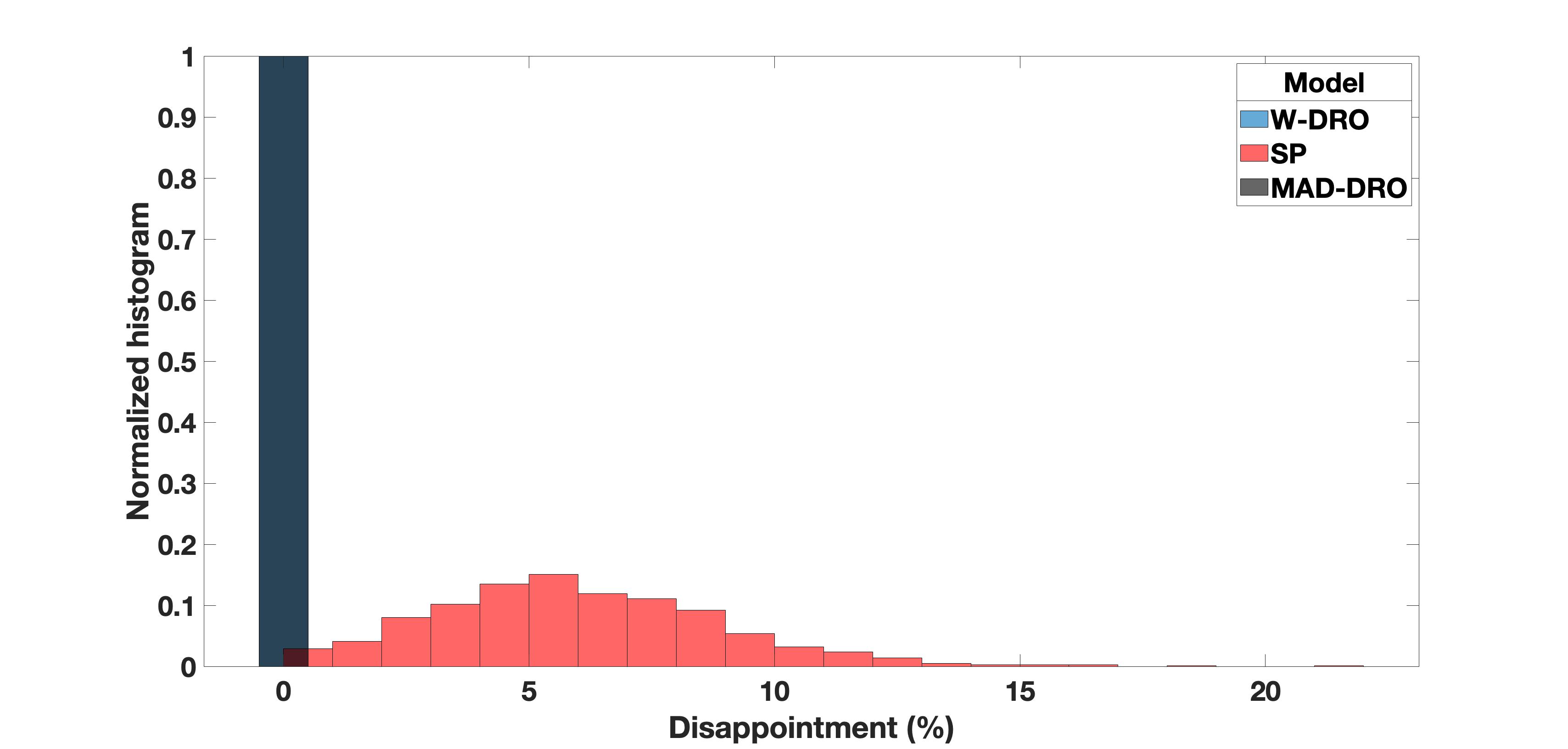}
        \caption{Instance 3, $\Delta=0$}
    \end{subfigure}%
      \begin{subfigure}[b]{0.5\textwidth}
        \centering
        \includegraphics[width=\textwidth]{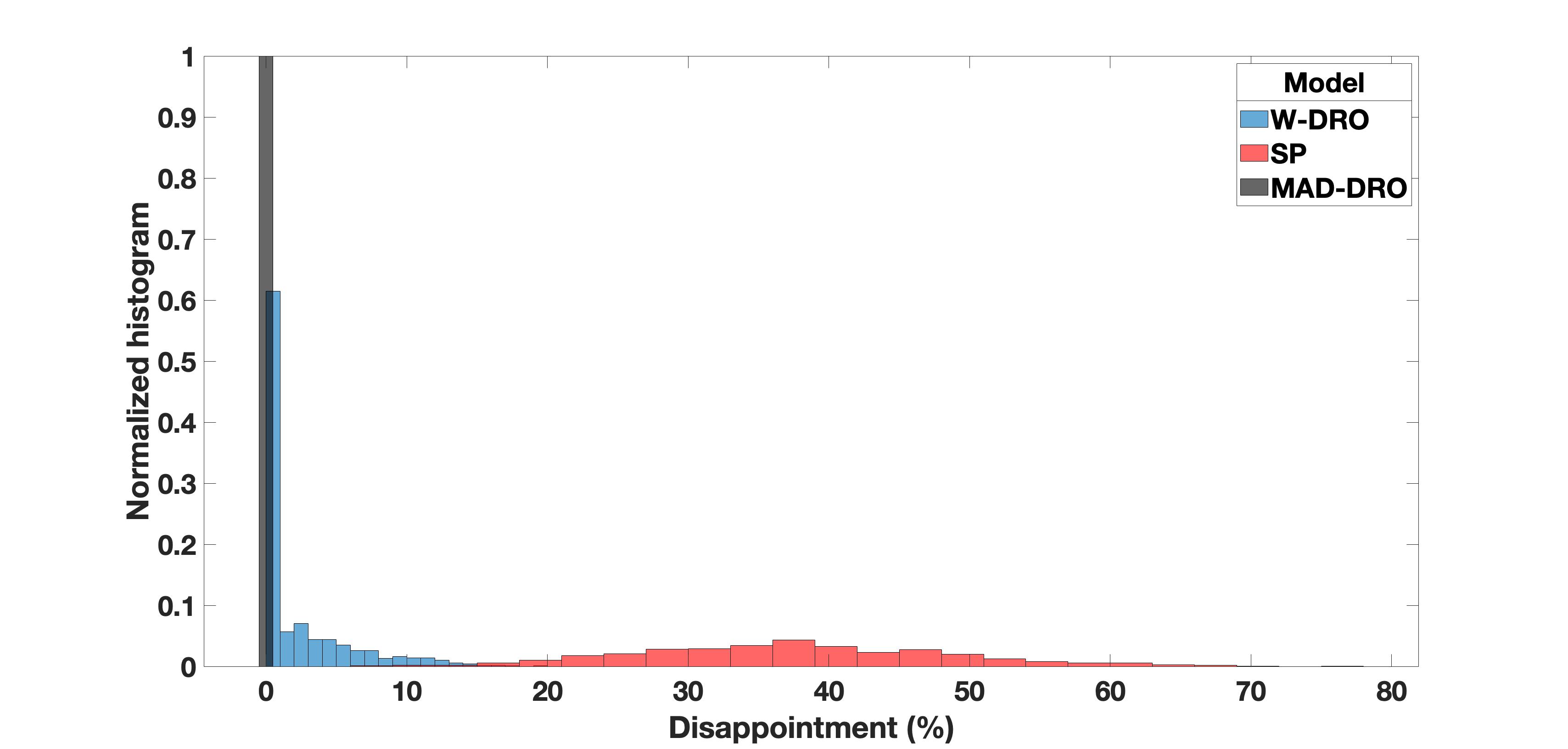}
        \caption{Instance 3, $\Delta=0.25$}
    \end{subfigure}%
    
    
      \begin{subfigure}[b]{0.5\textwidth}
        \centering
        \includegraphics[width=\textwidth]{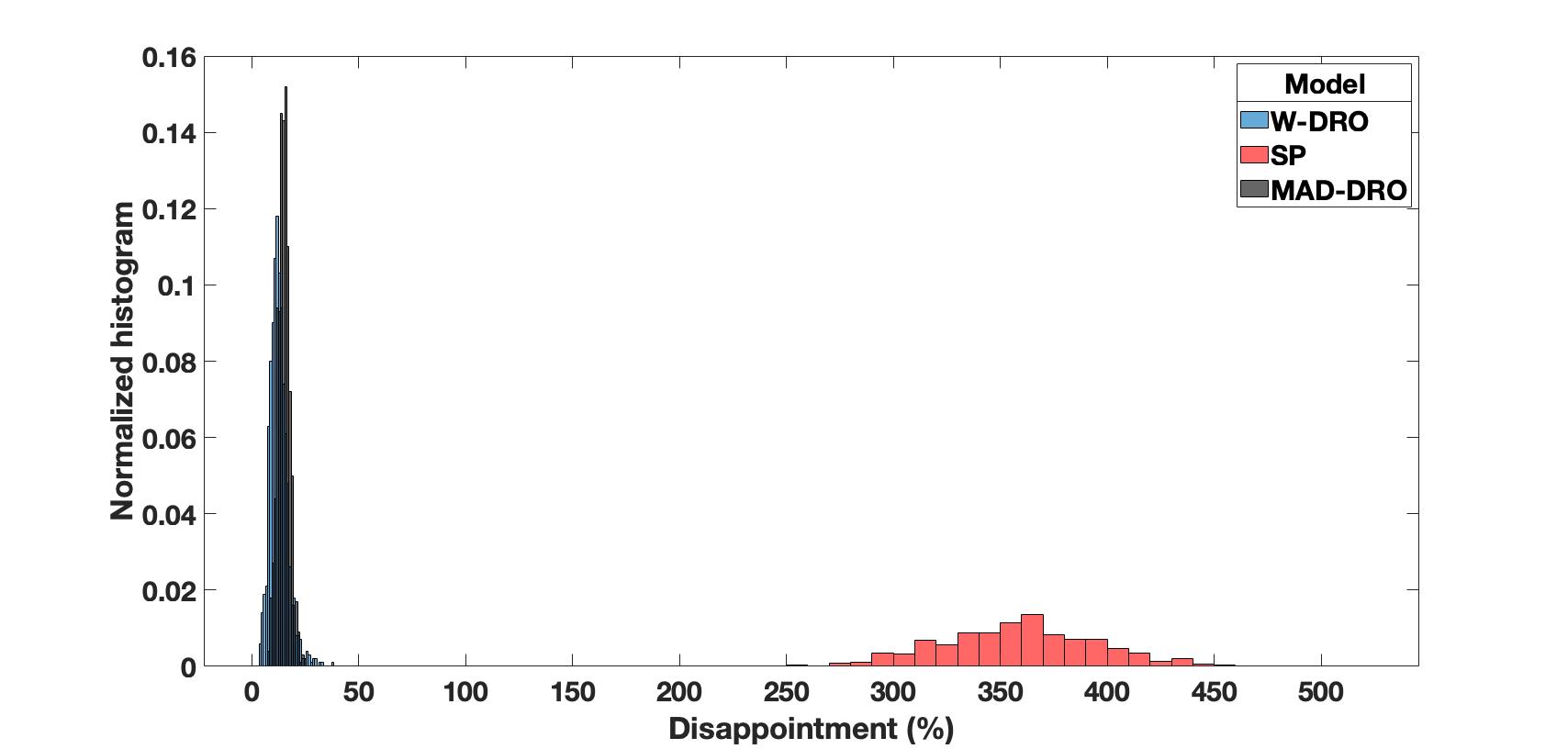}
        \caption{Lehigh 1, $\Delta=0$}
    \end{subfigure}%
      \begin{subfigure}[b]{0.5\textwidth}
        \centering
        \includegraphics[width=\textwidth]{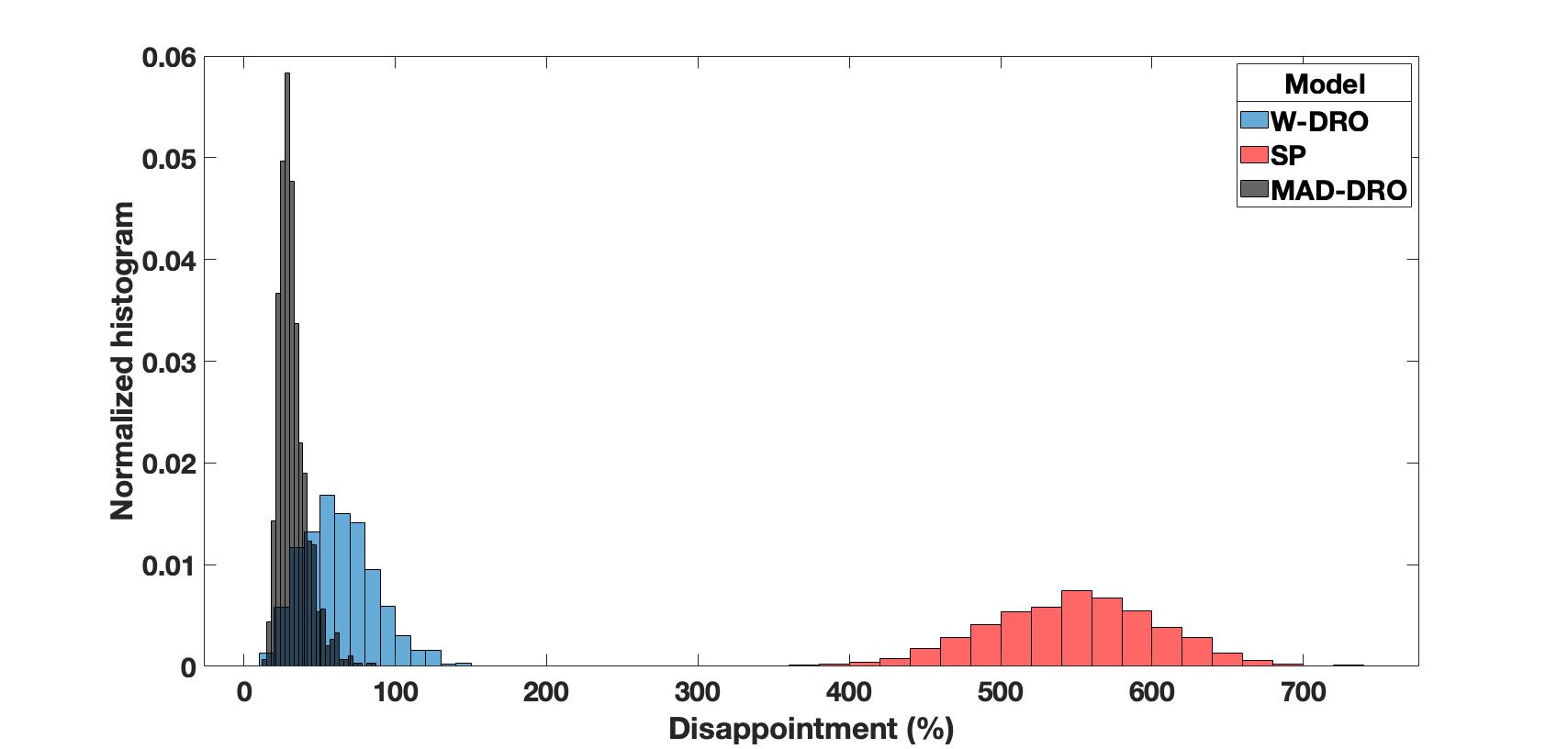}
        \caption{Lehigh 1, $\Delta=0.25$}
    \end{subfigure}%
    
\caption{Normalized histograms of out-of-sample disappointments under Set 2 with $\Delta=0$ and $\Delta=0.25$ for Instance 3 ($\Wb \in[50, 100]$, $N=10$) and Lehigh 1 ($N=10$).}\label{Fig6_Dissappt_Lehigh1}
\end{figure}

Figure~\ref{Fig6_Dissappt_Lehigh1} presents the histograms of the out-of-sample disappointments of the DRO and SP models for Instance 3 ($\Wb \in[50, 100]$, $N=10$) and Lehigh 1 ($N=10$) with $\Delta=0$ and $\Delta=0.25$. Notably, the DRO models yield substantially smaller out-of-sample disappointments on average and at all quantiles.  However, the average and upper quantiles of the disappointments of the MAD-DRO model is smaller than the W-DRO model when $\Delta=0.25$, especially for Lehigh 1. In addition,  it is clear that the average and upper quantiles of the disappointments of the SP model are relatively very large (e.g., exceeding 100\% for Lehigh 1). Finally, we observe that the out-of-sample disappointment of the MAD-DRO model is more stable than the W-DRO and SP models with a smaller standard deviation. We remark that these observations are consistent for the other considered instances, and the results with $\Delta=0.5$ are similar to those with $\Delta=0.25$. This demonstrates that the DRO model provides a more robust estimate of the actual cost that we will incur in practice.

\color{black}

The  results in this section demonstrate that the DRO approaches are effective in an environment where the distribution is hard to estimate (ambiguous), quickly changes, or when there is a small data set on demand variability. Moreover, these results emphasize the value of modeling uncertainty and distributional ambiguity.

\subsection{\textbf{Sensitivity Analysis}}\label{sec5:sensitivity}

\noindent  In this section, we study the sensitivity of DRO models to different parameter settings. Given that we observe similar results for all of the constructed instances, for presentation brevity and illustrative purposes, we present results for Instance 1 ($I$, $J$, $T$)$=$(10, 10, 10) and Instance 5 ($I$, $J$, $T$)$=$(20, 20, 10) as examples of small and relatively large instances.

\begin{figure}[t!]
     \begin{subfigure}[b]{0.5\textwidth}
 \centering
        \includegraphics[width=\textwidth]{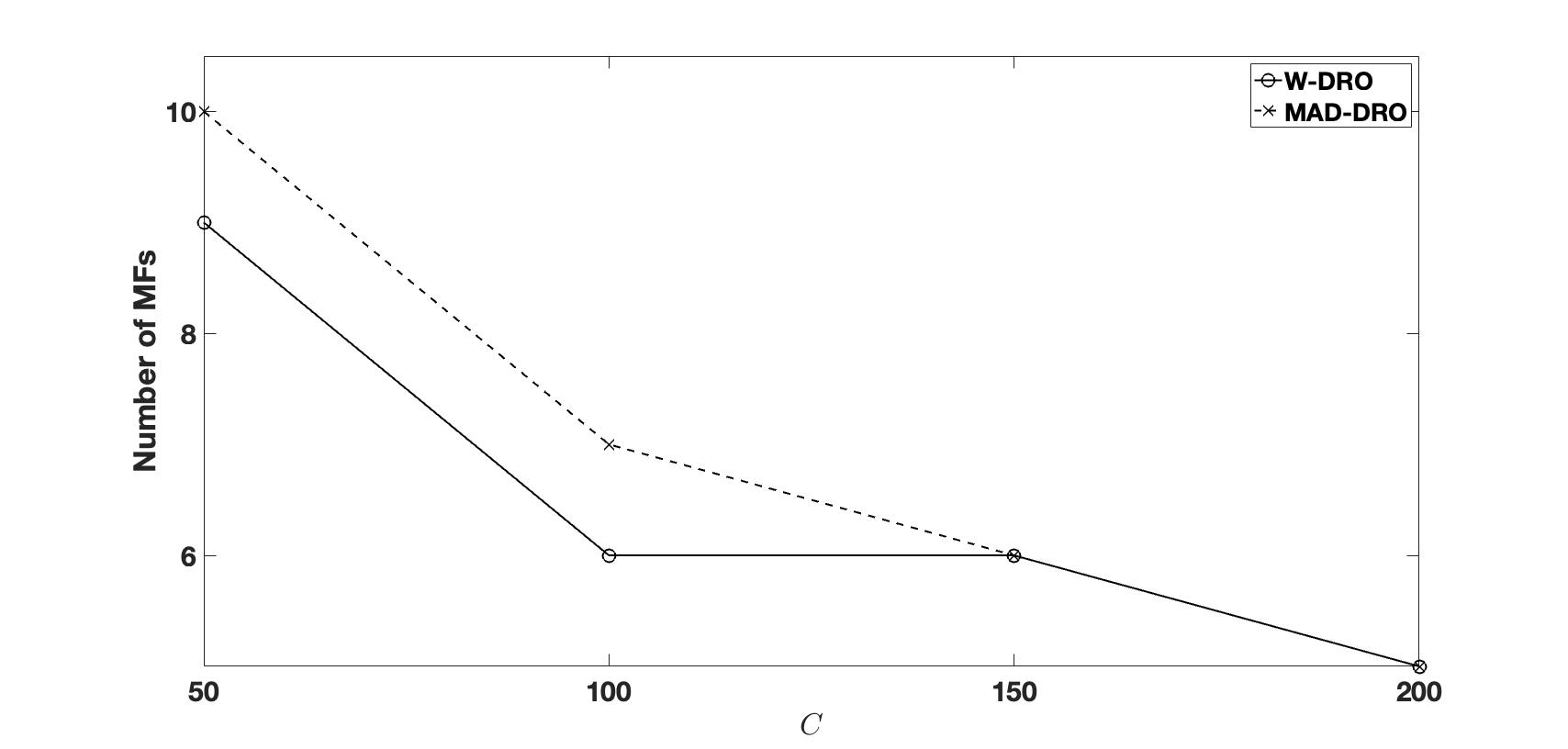}
        \caption{Number of MFs, $f=1,500$}
        \label{Fig6a}
    \end{subfigure}%
    \begin{subfigure}[b]{0.5\textwidth}
            \includegraphics[width=\textwidth]{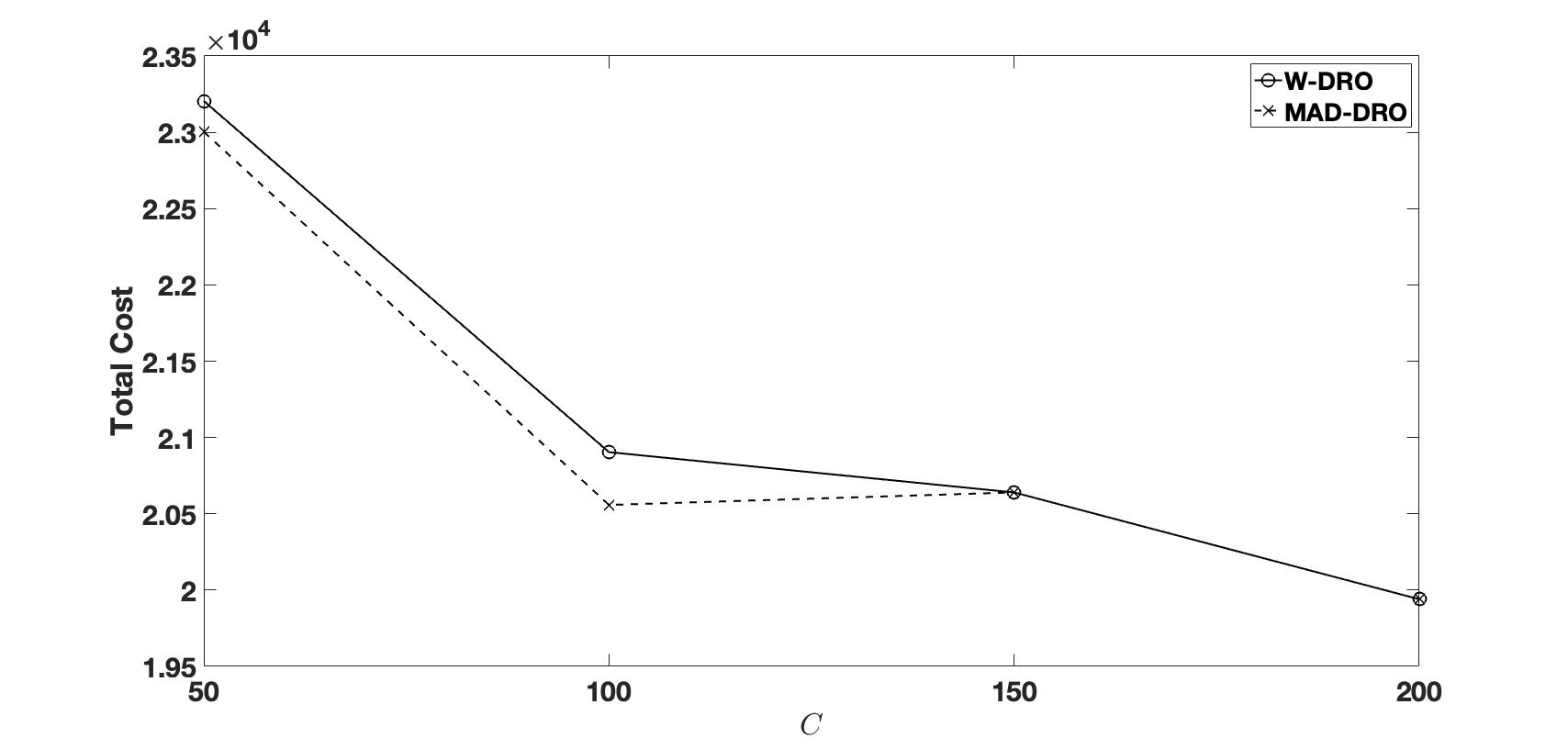}
      \caption{Total cost, $f=1,500$}
      \label{Fig6b}
    \end{subfigure}%
    
        \begin{subfigure}[b]{0.5\textwidth}
 \centering
        \includegraphics[width=\textwidth]{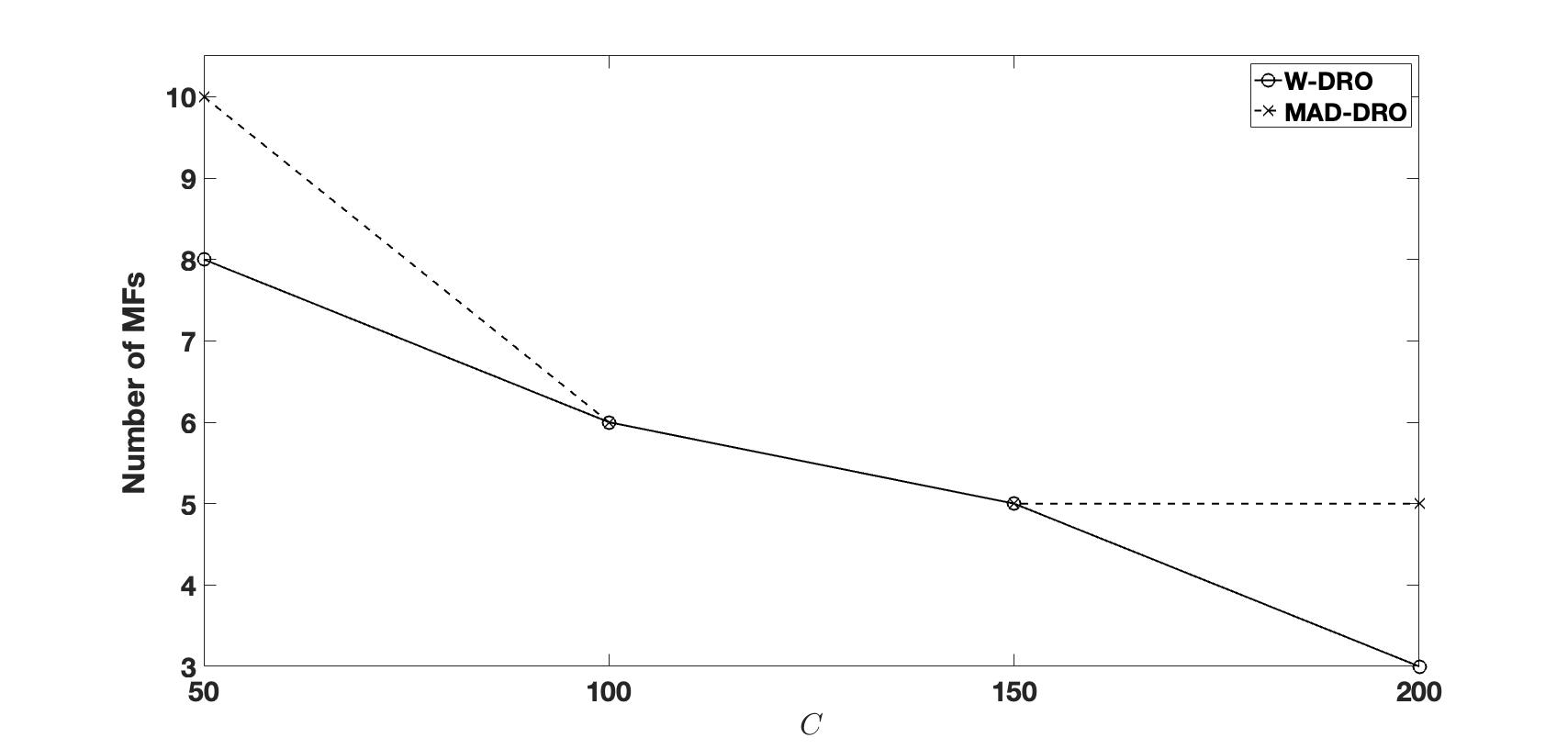}
        \caption{Number of MFs, $f=6,000$}
        \label{Fig6c}
    \end{subfigure}%
    \begin{subfigure}[b]{0.5\textwidth}
            \includegraphics[width=\textwidth]{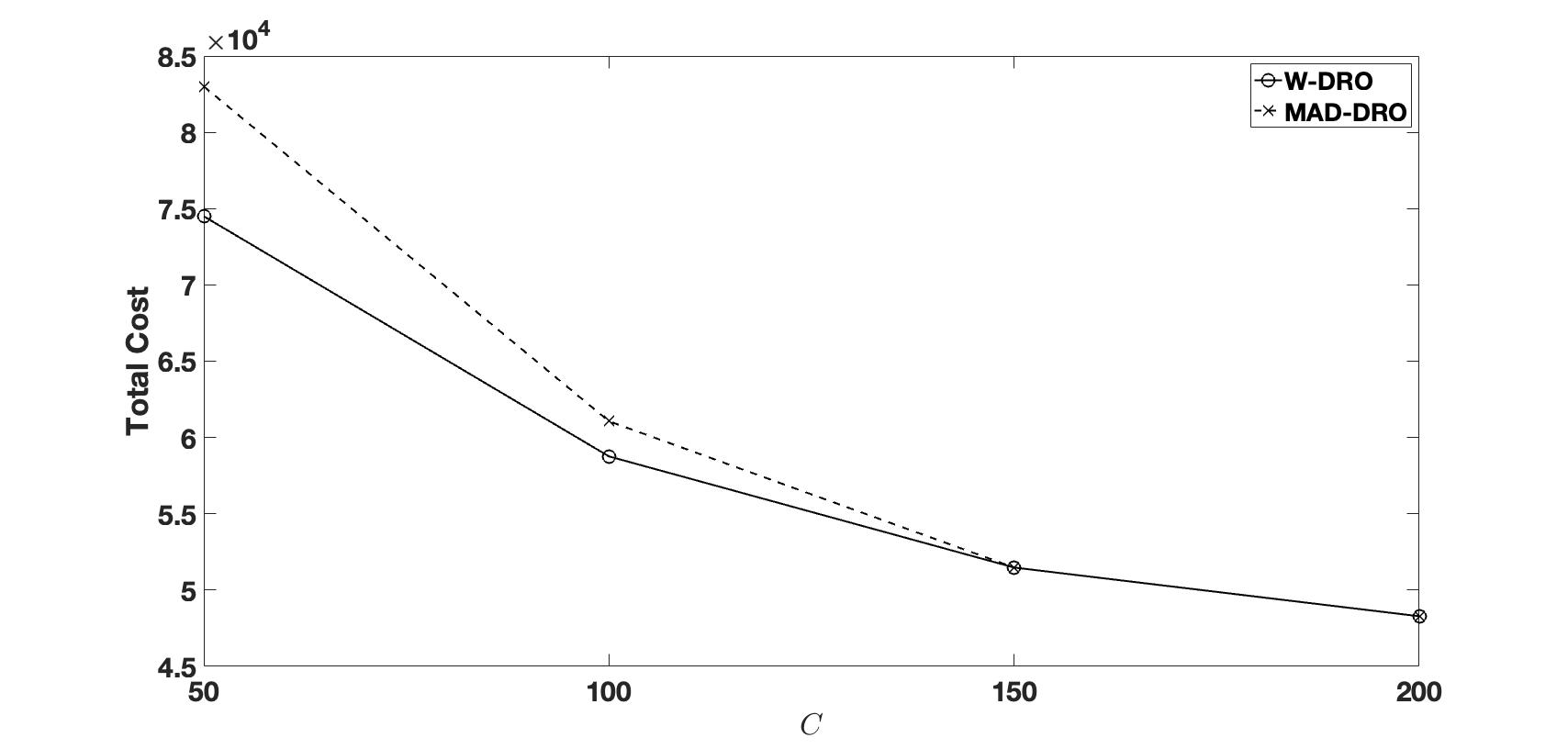}
      \caption{Total cost, $f=6,000$}
      \label{Fig6d}
    \end{subfigure}%
    
            \begin{subfigure}[b]{0.5\textwidth}
 \centering
        \includegraphics[width=\textwidth]{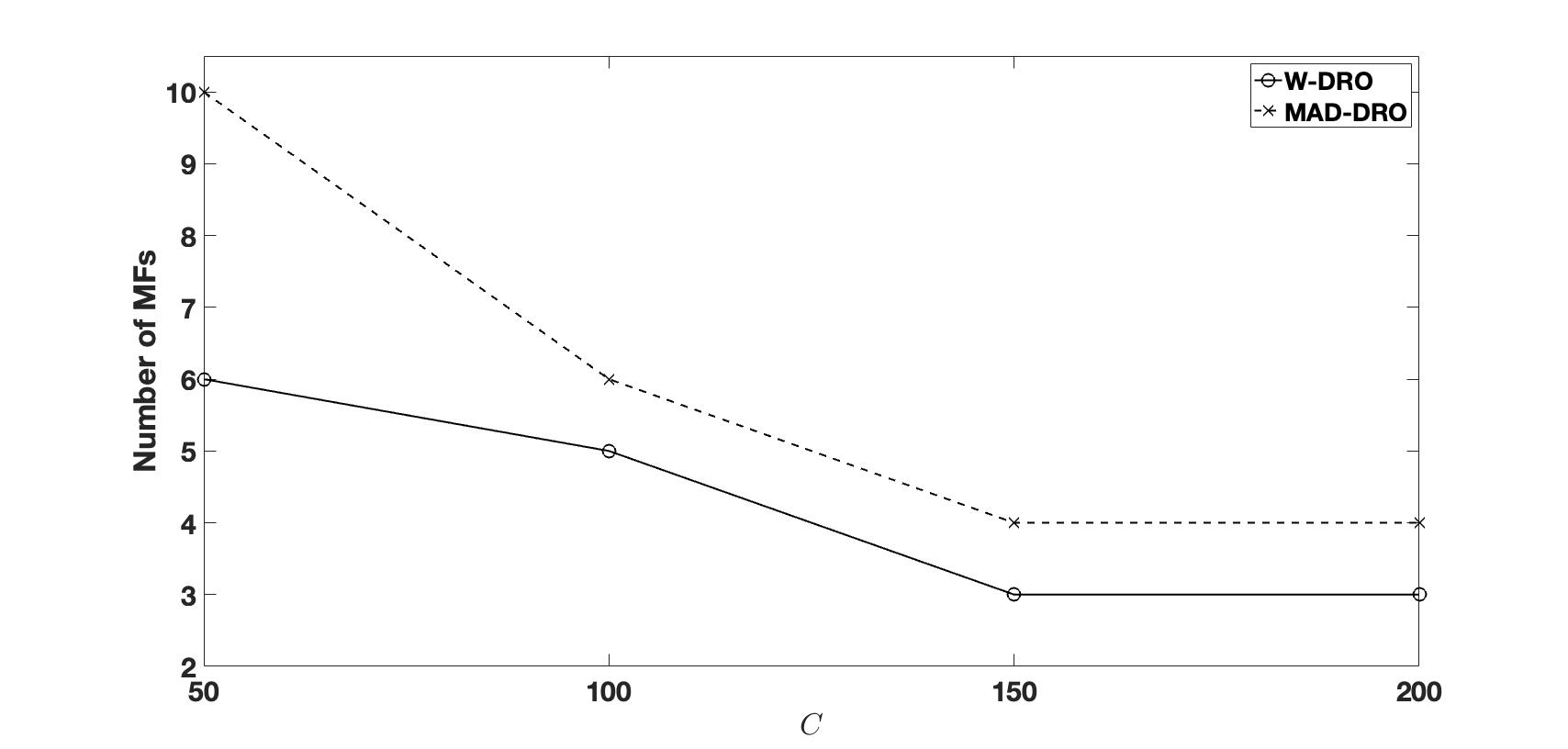}
        \caption{Number of MFs, $f=6,000$}
        \label{Fig6e}
    \end{subfigure}%
    \begin{subfigure}[b]{0.5\textwidth}
            \includegraphics[width=\textwidth]{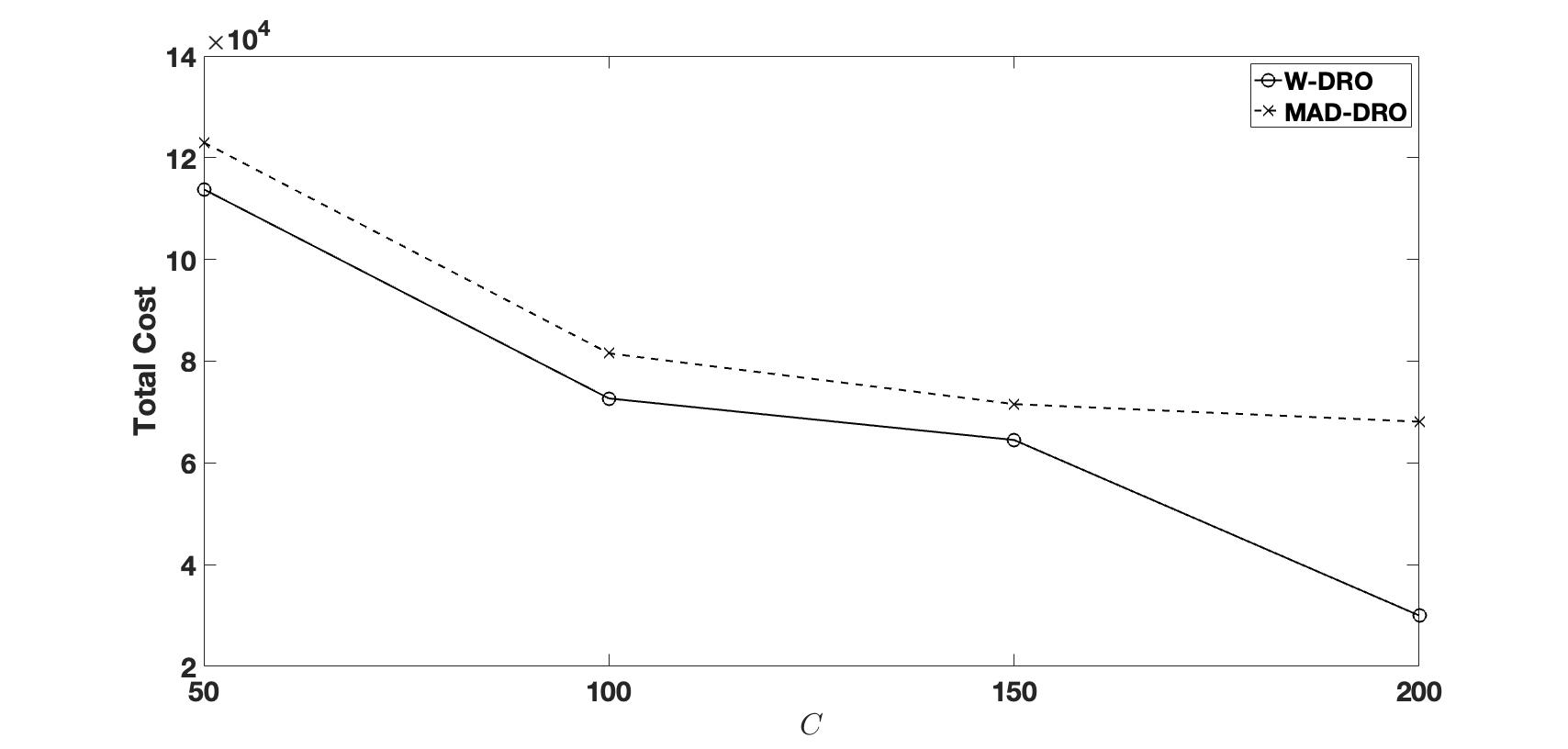}
      \caption{Total cost, $f=10,000$}
      \label{Fig6f}
    \end{subfigure}%
    \caption{Comparison of the results for different values of $C$ and $f$ under $\Wb \in [20, 60]$. Instance 1}\label{Fig6:MF_vs_C_inst1}
\end{figure}

\begin{figure}[t!]
     \begin{subfigure}[b]{0.5\textwidth}
 \centering
        \includegraphics[width=\textwidth]{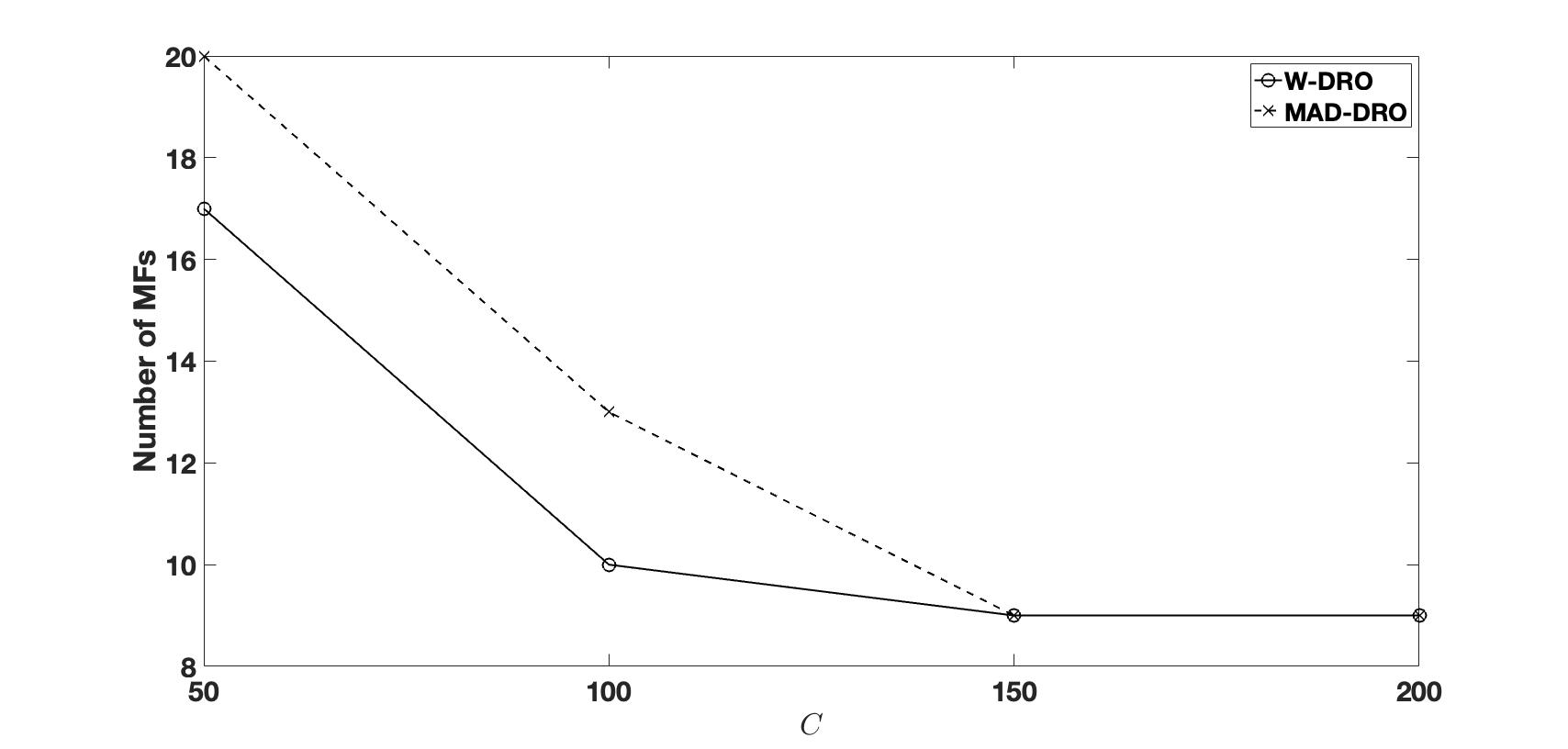}
        \caption{Number of MFs, $f=1,500$}
        \label{Fig9a}
    \end{subfigure}%
    \begin{subfigure}[b]{0.5\textwidth}
            \includegraphics[width=\textwidth]{C_Vs_F_TC_1500.jpg}
      \caption{Total cost, $f=1,500$}
      \label{Fig9b}
    \end{subfigure}%
    
        \begin{subfigure}[b]{0.5\textwidth}
 \centering
        \includegraphics[width=\textwidth]{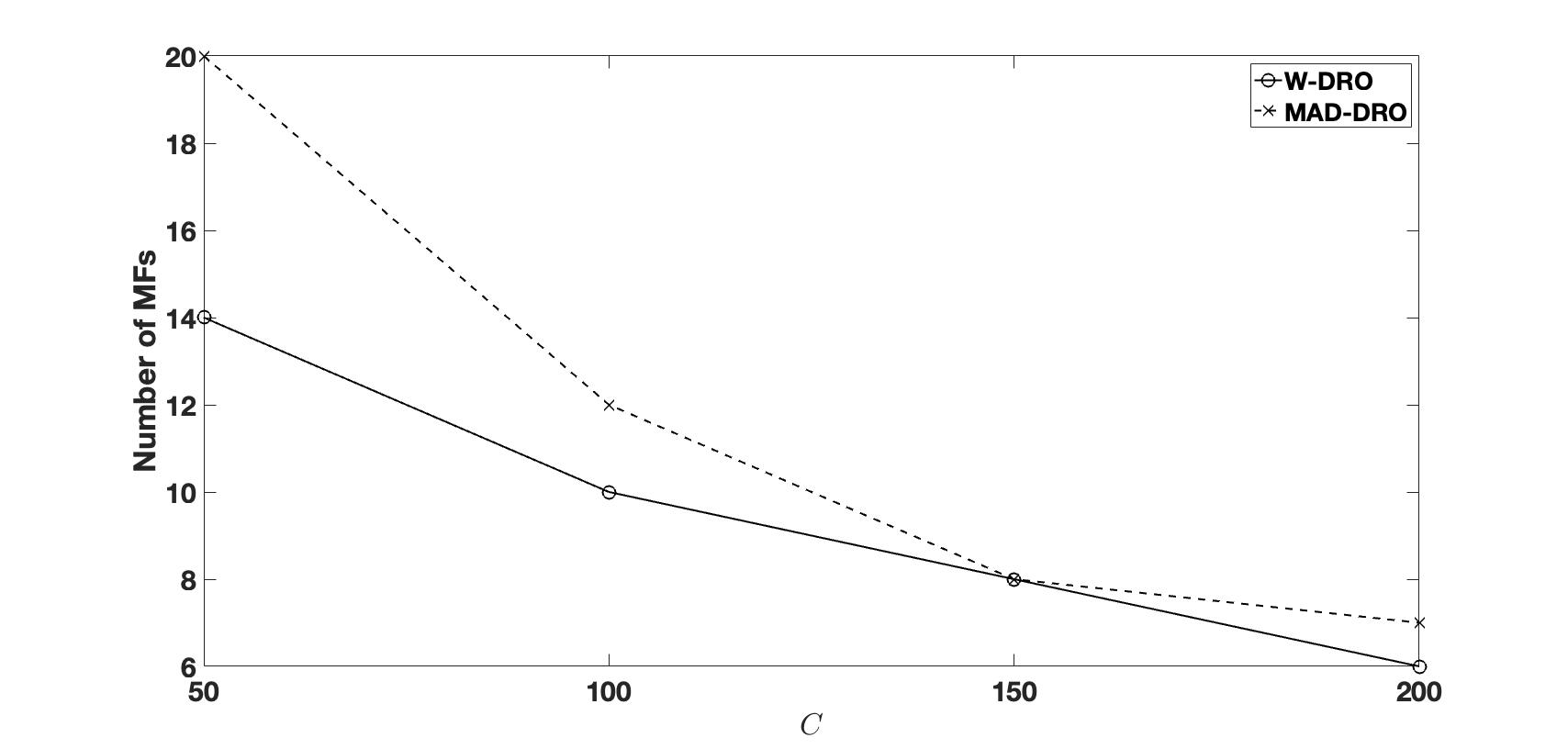}
        \caption{Number of MFs, $f=6,000$}
        \label{Fig9c}
    \end{subfigure}%
    \begin{subfigure}[b]{0.5\textwidth}
            \includegraphics[width=\textwidth]{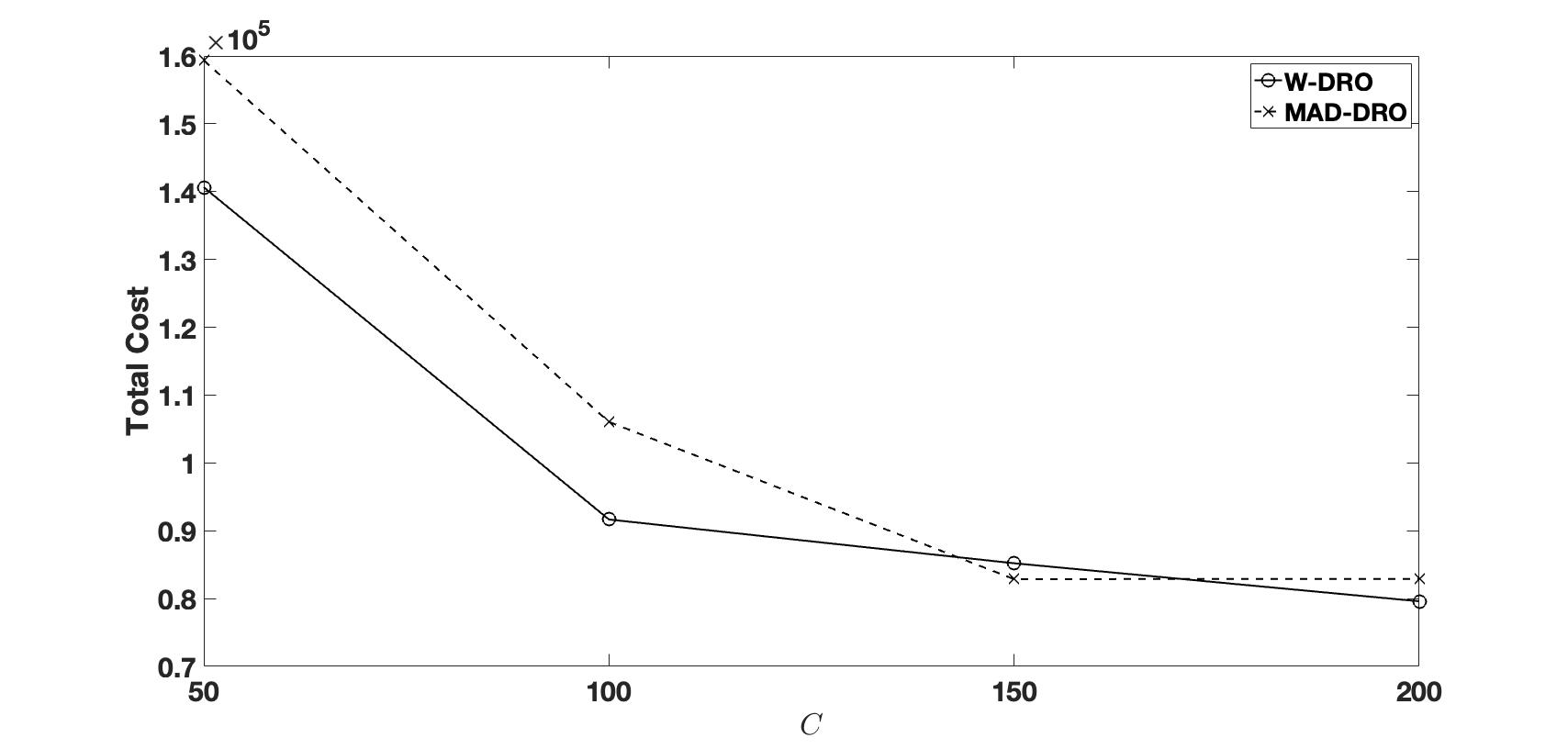}
      \caption{Total cost, $f=6,000$}
      \label{Fig9d}
    \end{subfigure}%
    
            \begin{subfigure}[b]{0.5\textwidth}
 \centering
        \includegraphics[width=\textwidth]{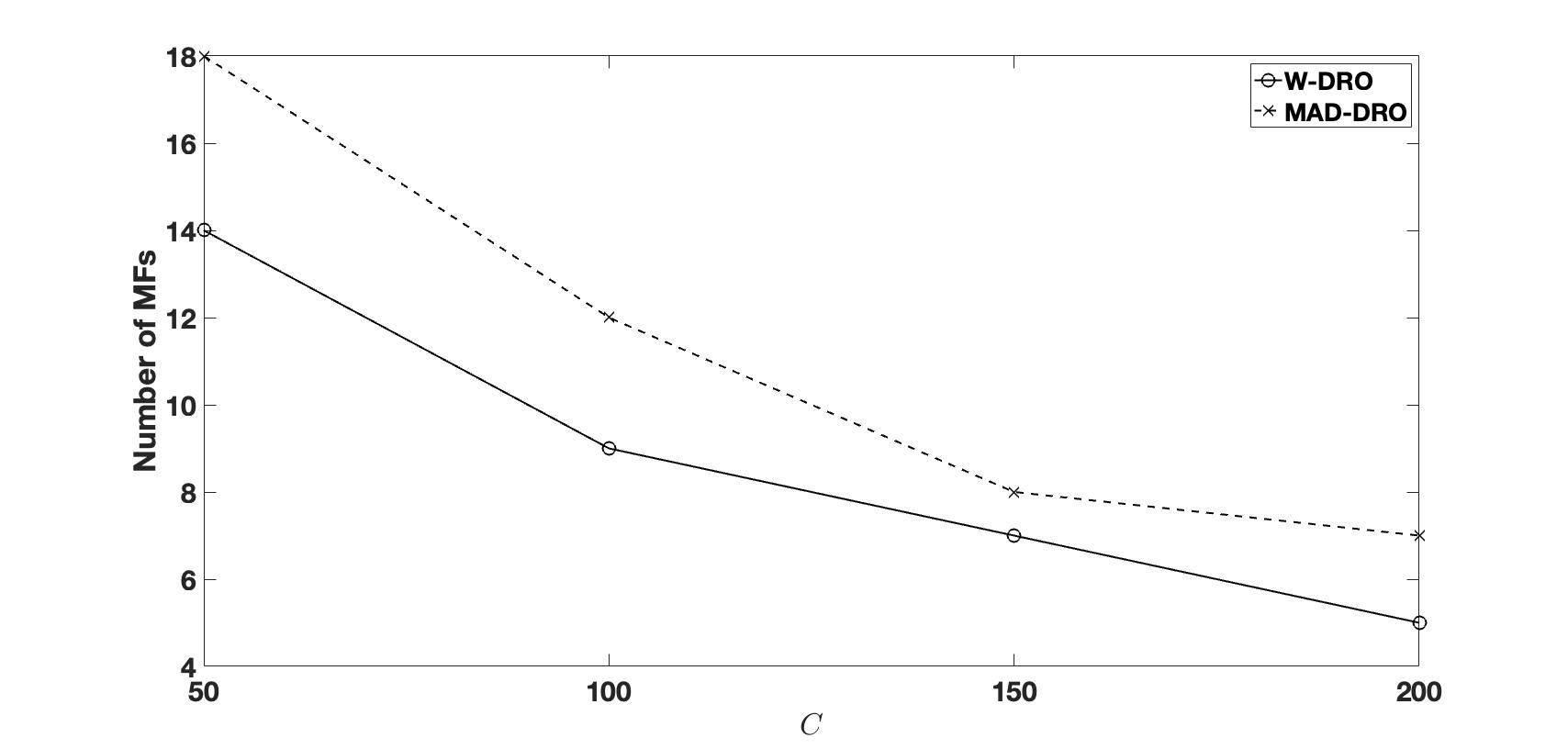}
        \caption{Number of MFs, $f=10,000$}
        \label{Fig9e}
    \end{subfigure}%
    \begin{subfigure}[b]{0.5\textwidth}
            \includegraphics[width=\textwidth]{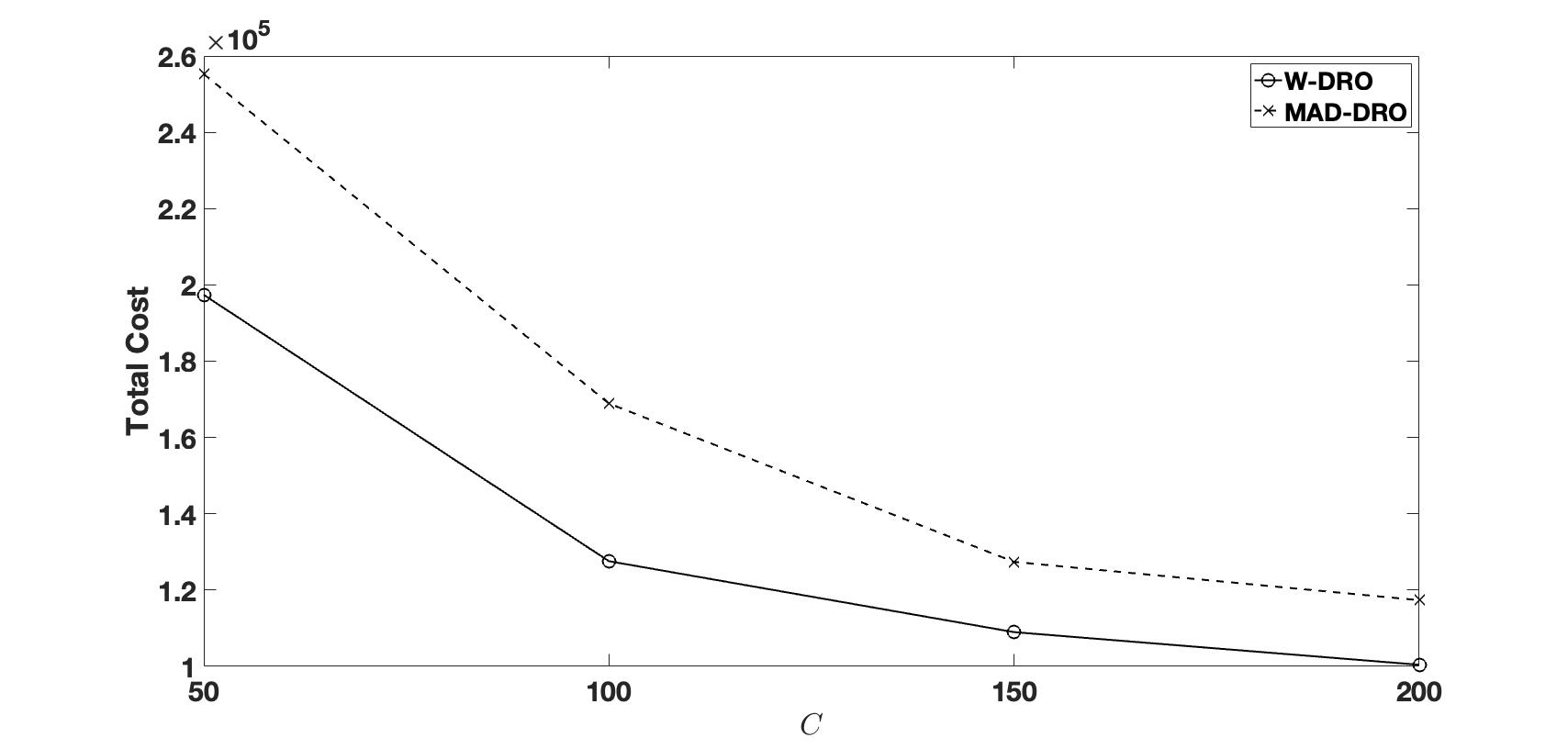}
      \caption{\textcolor{black}{Total cost, $f=10,000$}}
      \label{Fig9f}
    \end{subfigure}%
    \caption{Comparison of the results for different values of $C$ and $f$ under $\Wb \in [20, 60]$. Instance 5}\label{Fig9:MF_vs_C_inst5}
\end{figure}

First, we analyze the optimal number of active MFs as a function of the fixed cost, $f$, MF capacity, $C$, and range of demand.  We fix all parameters as described in Section~\ref{sec5.1:instancegen} and solve the W-DRO and MAD-DRO models with $C\in\{ 50, \ 100, \ 150, \ 200\}$ and $f \in \{1,500 \ (\text{low}), \ 6,000 \ (\text{average}), \ 10,000 \ (\text{high})\} $ under the base range $\Wb \in [20, 60]$ and $\Wb \in [50, 100]$ (a higher volume of the demand). Figures \ref{Fig6:MF_vs_C_inst1} and \ref{Fig9:MF_vs_C_inst5} present the optimal number of active MFs and the associated total cost (under Set 1) for Instance 1 and Instance 5 under $\Wb \in [20, 60]$, respectively.  Figures \ref{Fig7:MF_vs_C_Range2_inst1}--\ref{Fig10:MF_vs_C_Range2_inst5} in Appendix~\ref{Appx:sensitivity} present the results under   $\Wb \in [50, 100]$.

We observe the following from these figures. First, the optimal number of scheduled MFs decreases as $C$ increases irrespective of $f$.  This makes sense because, with a higher capacity, each MF can serve a larger amount of demand in each period. Second, both models schedule more MFs under $\Wb \in [50, 100]$, i.e., a higher volume of the demand. For example, consider Instance 5.  When $f=1,500$ and $C=100$ the (W-DRO, MAD-DRO) models schedule (10, 13) and (18, 19) MFs under $\Wb \in [20, 60]$ and $\Wb \in [50, 100]$, respectively.  Third, the MAD-DRO model always schedules a higher number of MFs, especially when $C$ is tight. As such, the MAD-DRO model often has a slightly higher total cost (due to the higher fixed cost of establishing a larger fleet) and better second-stage cost, i.e., better operational performance  (see Figures \ref{Fig8:MF_vs_C_inst1}-\ref{Fig12:MF_vs_C_range2_inst5}).   For example, consider Instance 1. When $f=6,000$ and $C=50$, the W-DRO and MAD-DRO models schedule 8 and 10 MFs, respectively. The associated (total, second-stage) costs of these solutions are respectively (74,495, 26,495) and (83,004, 23,004). Finally, both models schedule fewer MFs as $f$ increases. 

\begin{figure}[t!]
     \begin{subfigure}[b]{0.5\textwidth}
 \centering
        \includegraphics[width=\textwidth]{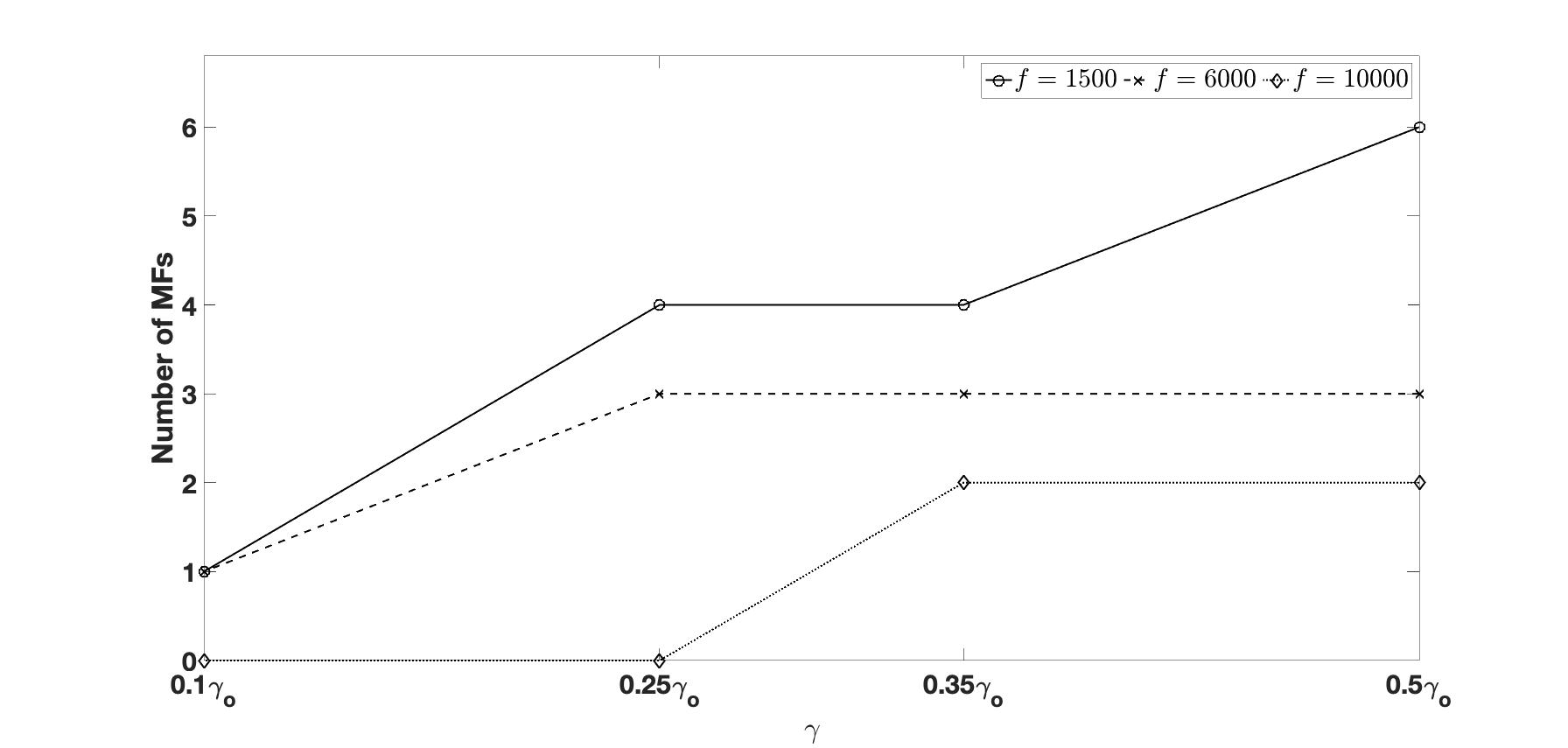}
        \caption{W-DRO, Instance 1}
        \label{Fig12a}
    \end{subfigure}%
    \begin{subfigure}[b]{0.5\textwidth}
            \includegraphics[width=\textwidth]{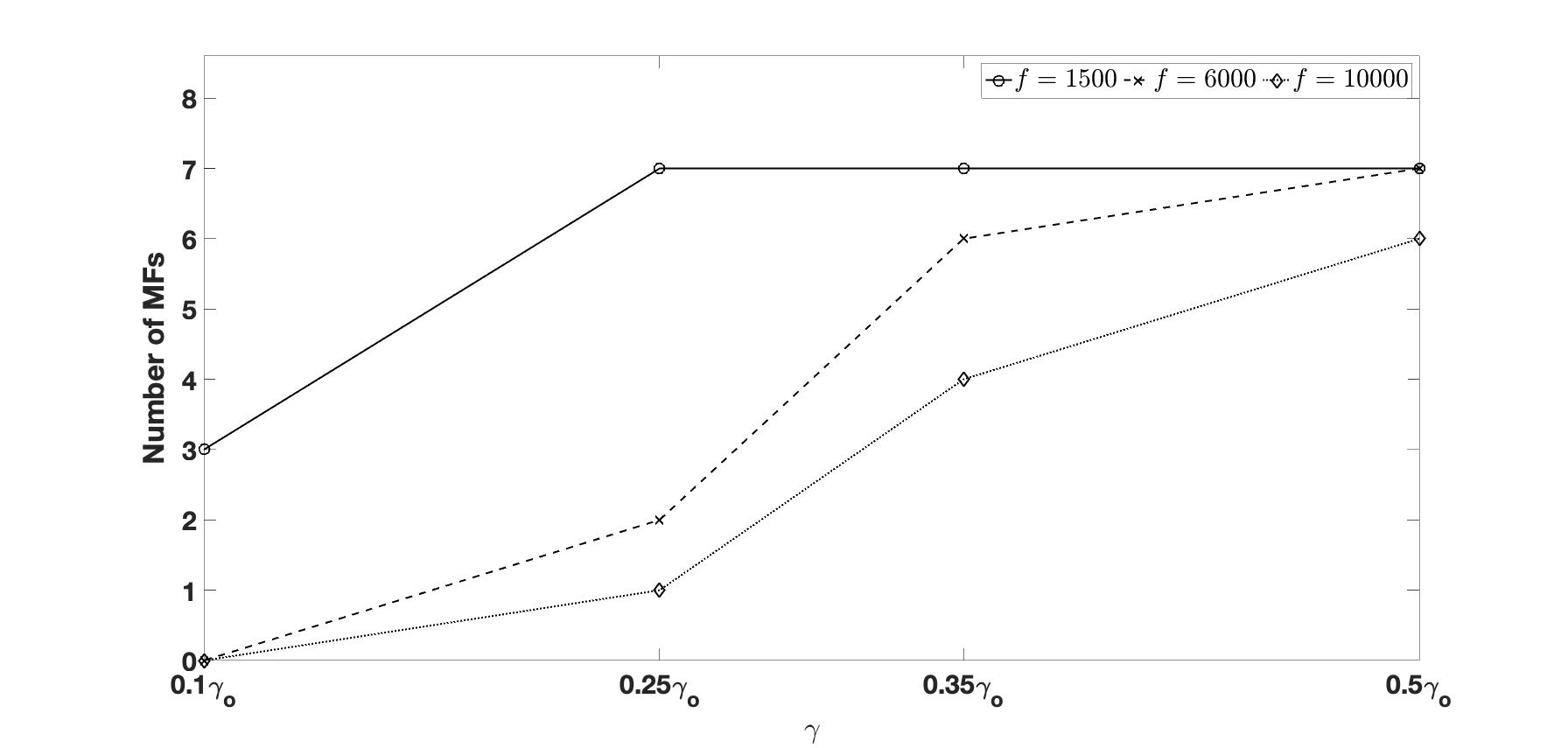}
      \caption{MAD-DRO, Instance 1}
      \label{Fig12b}%
    \end{subfigure}
    
         \begin{subfigure}[b]{0.5\textwidth}
 \centering
        \includegraphics[width=\textwidth]{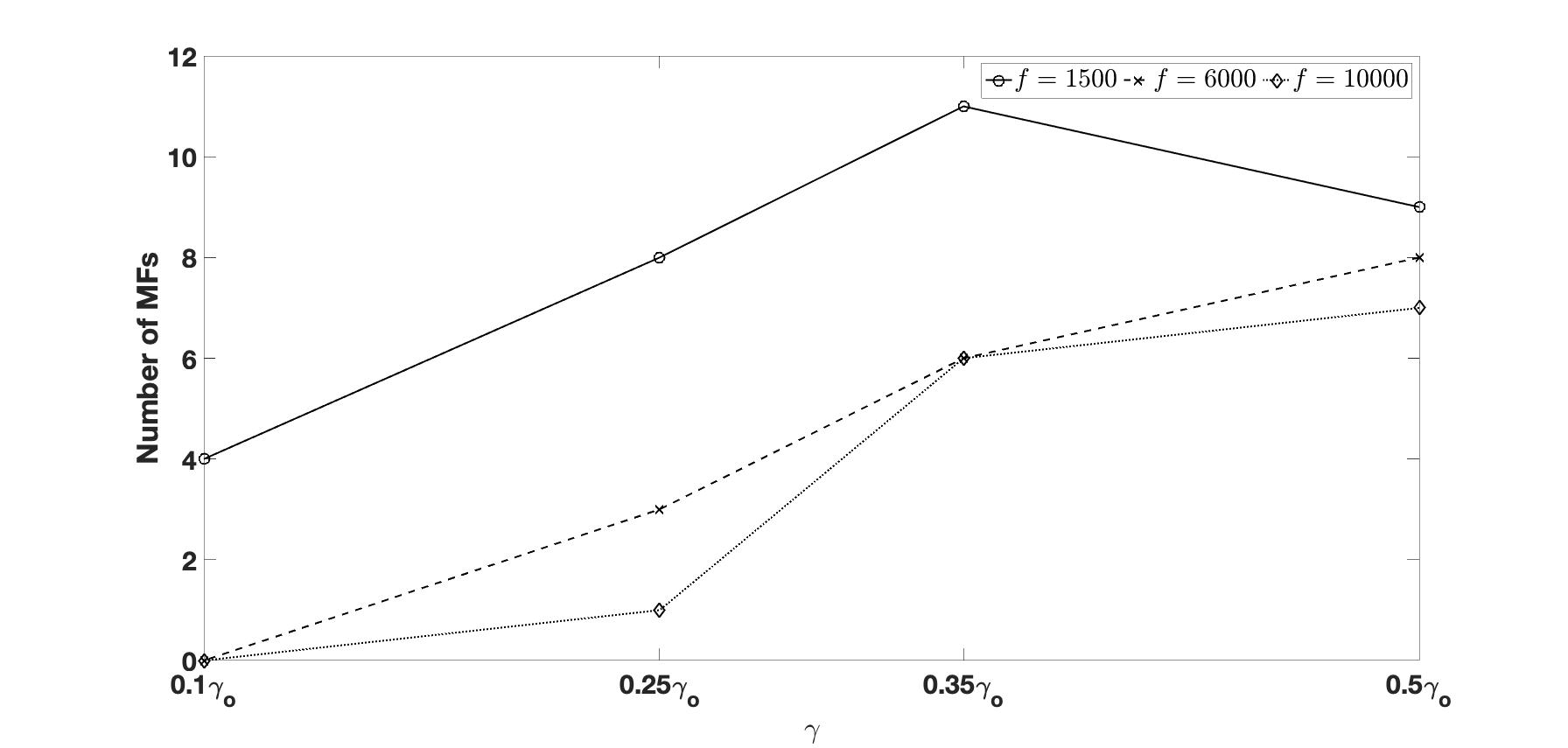}
        \caption{W-DRO, Instance 5}
        \label{Fig12a2}
    \end{subfigure}%
    \begin{subfigure}[b]{0.5\textwidth}
            \includegraphics[width=\textwidth]{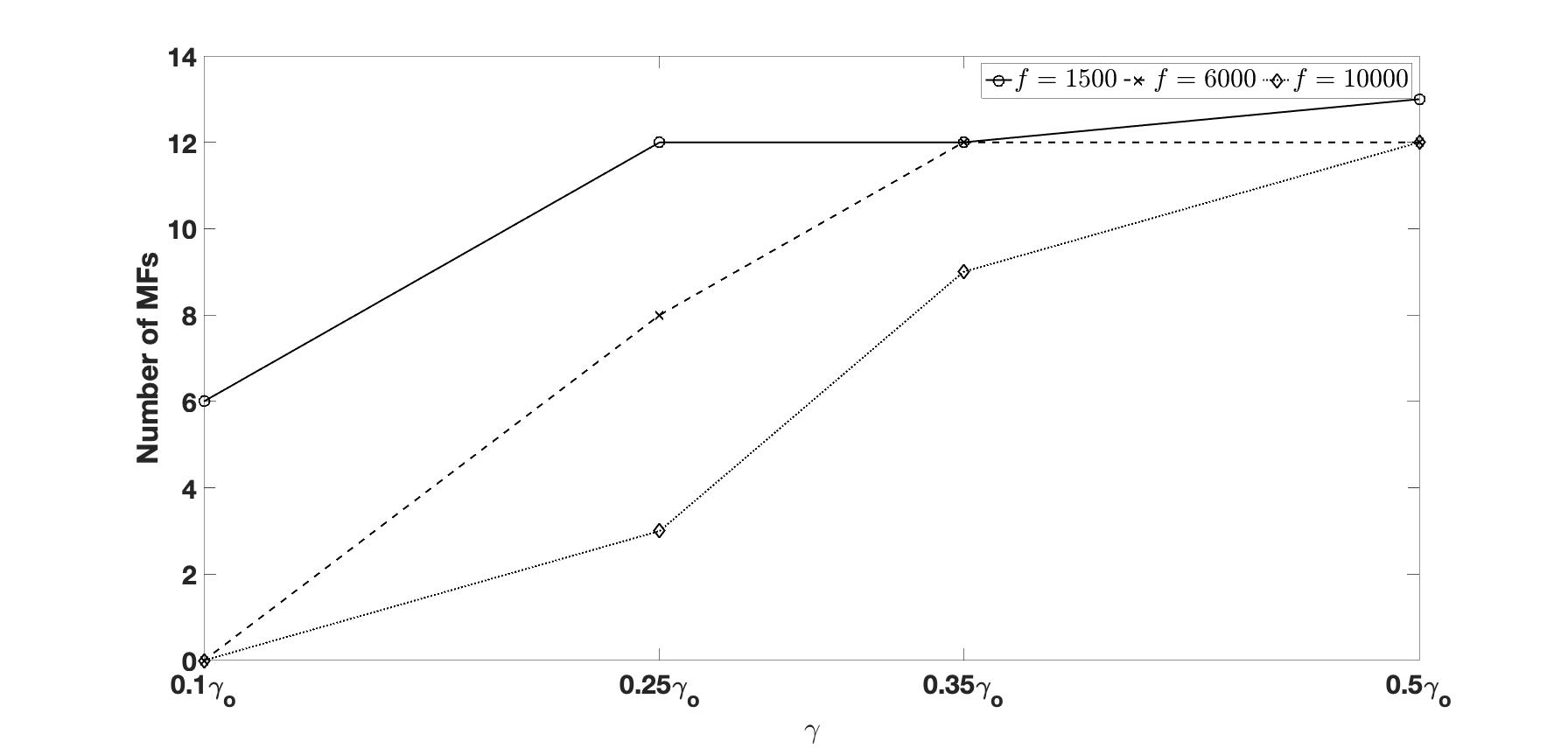}
      \caption{MAD-DRO, Instance 5}
      \label{Fig12b2}
    \end{subfigure}%
    \caption{Comparison of the results for different values of $\gamma$}\label{Fig12:MF_vs_gamma}
\end{figure}

\begin{figure}[t!]
     \begin{subfigure}[b]{0.5\textwidth}
 \centering
        \includegraphics[width=\textwidth]{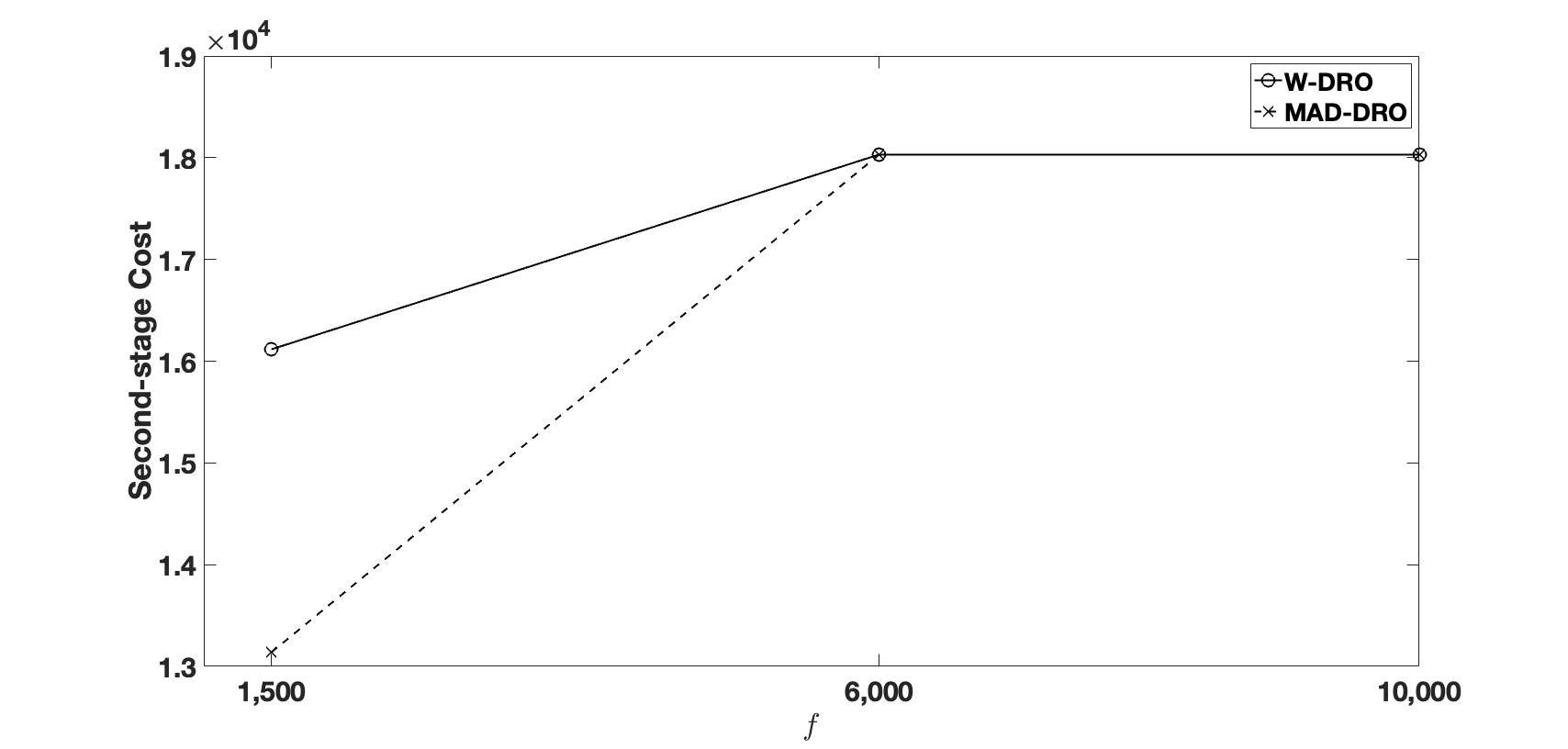}
        \caption{$\gamma=0.1\gamma_o$ Instance 1}
        \label{Fig13a}
    \end{subfigure}%
    \begin{subfigure}[b]{0.5\textwidth}
            \includegraphics[width=\textwidth]{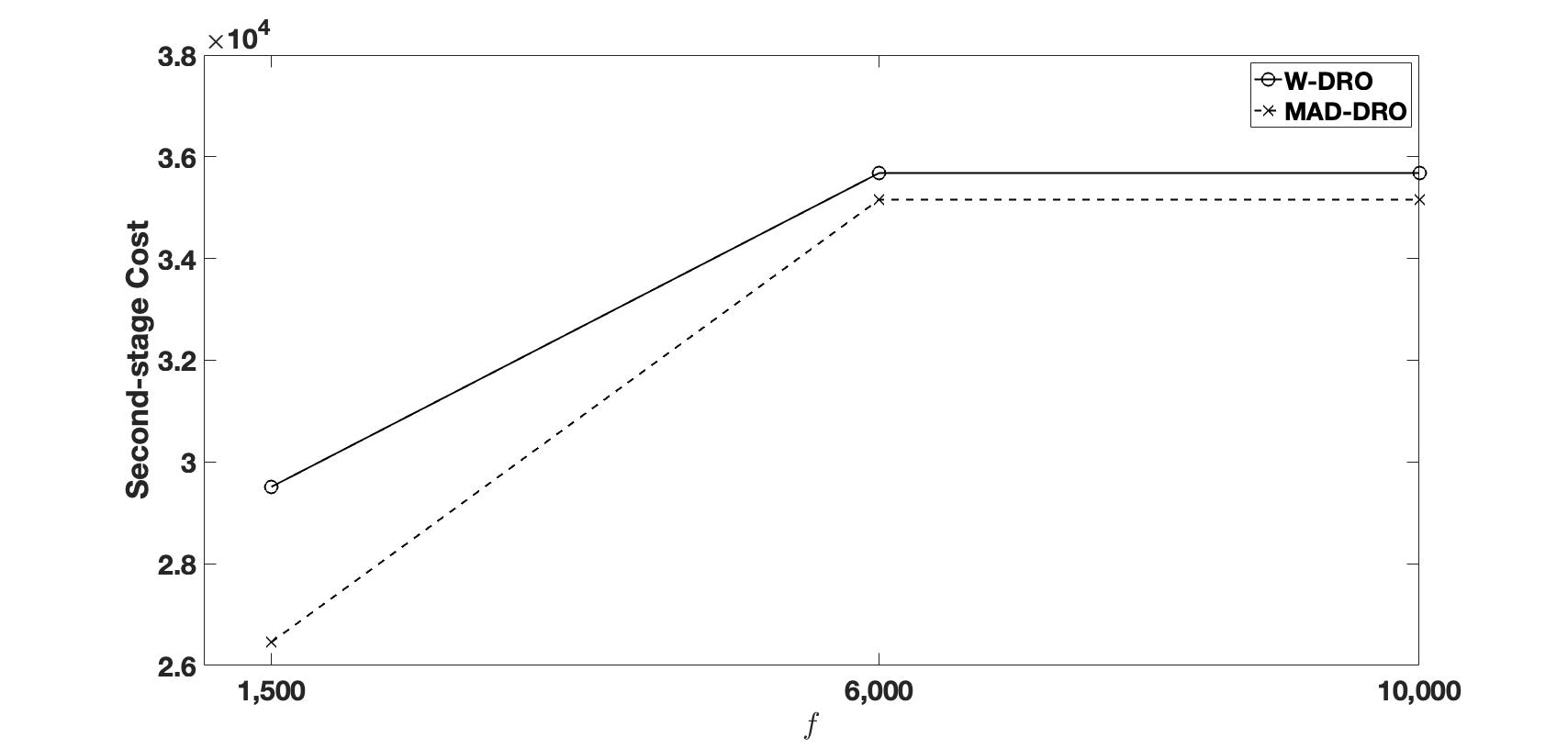}
      \caption{$\gamma=0.1\gamma_o$, Instance 5}
      \label{Fig13b}%
    \end{subfigure}
    
         \begin{subfigure}[b]{0.5\textwidth}
 \centering
        \includegraphics[width=\textwidth]{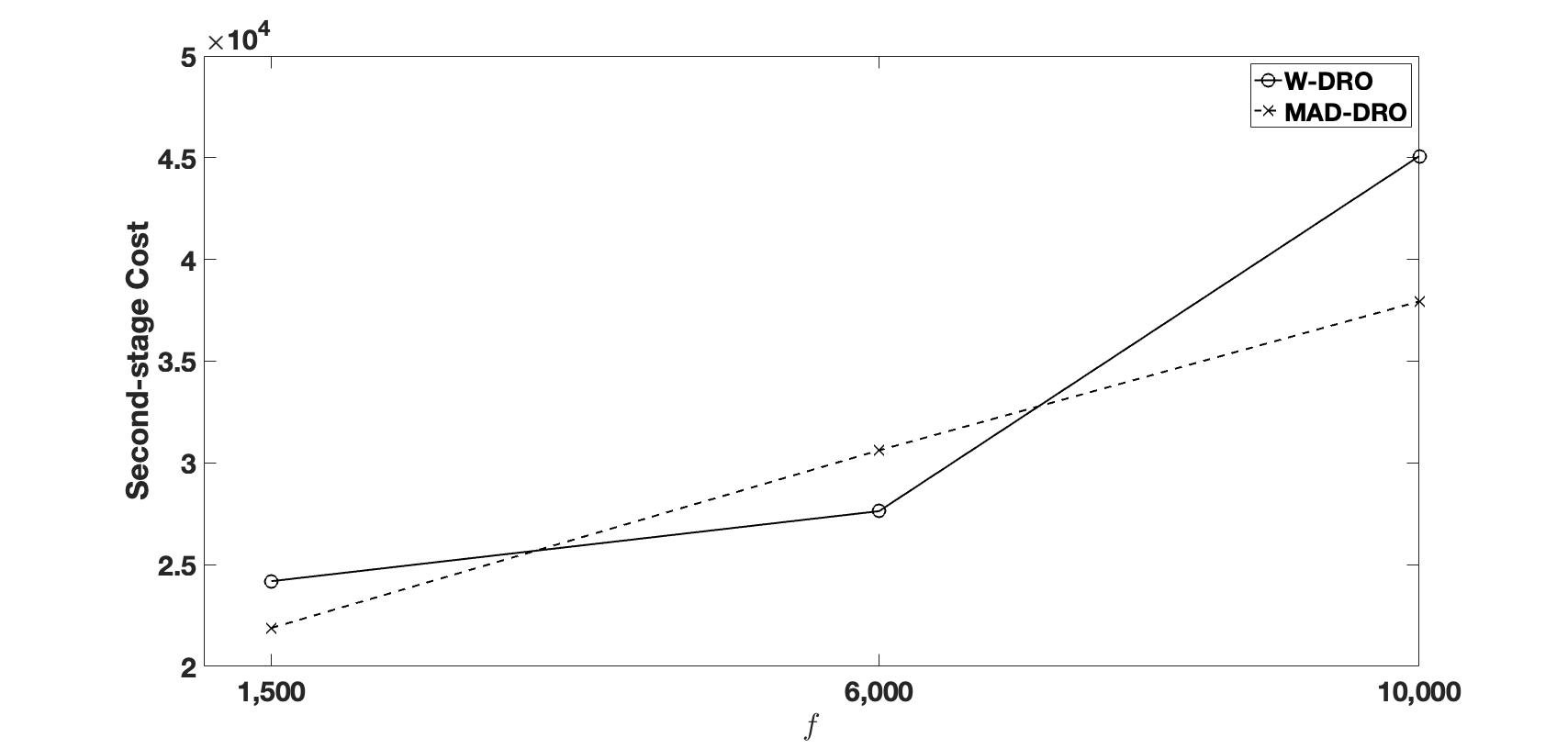}
        \caption{$\gamma=0.25\gamma_o$, Instance 1}
        \label{Fig13a2}
    \end{subfigure}%
    \begin{subfigure}[b]{0.5\textwidth}
            \includegraphics[width=\textwidth]{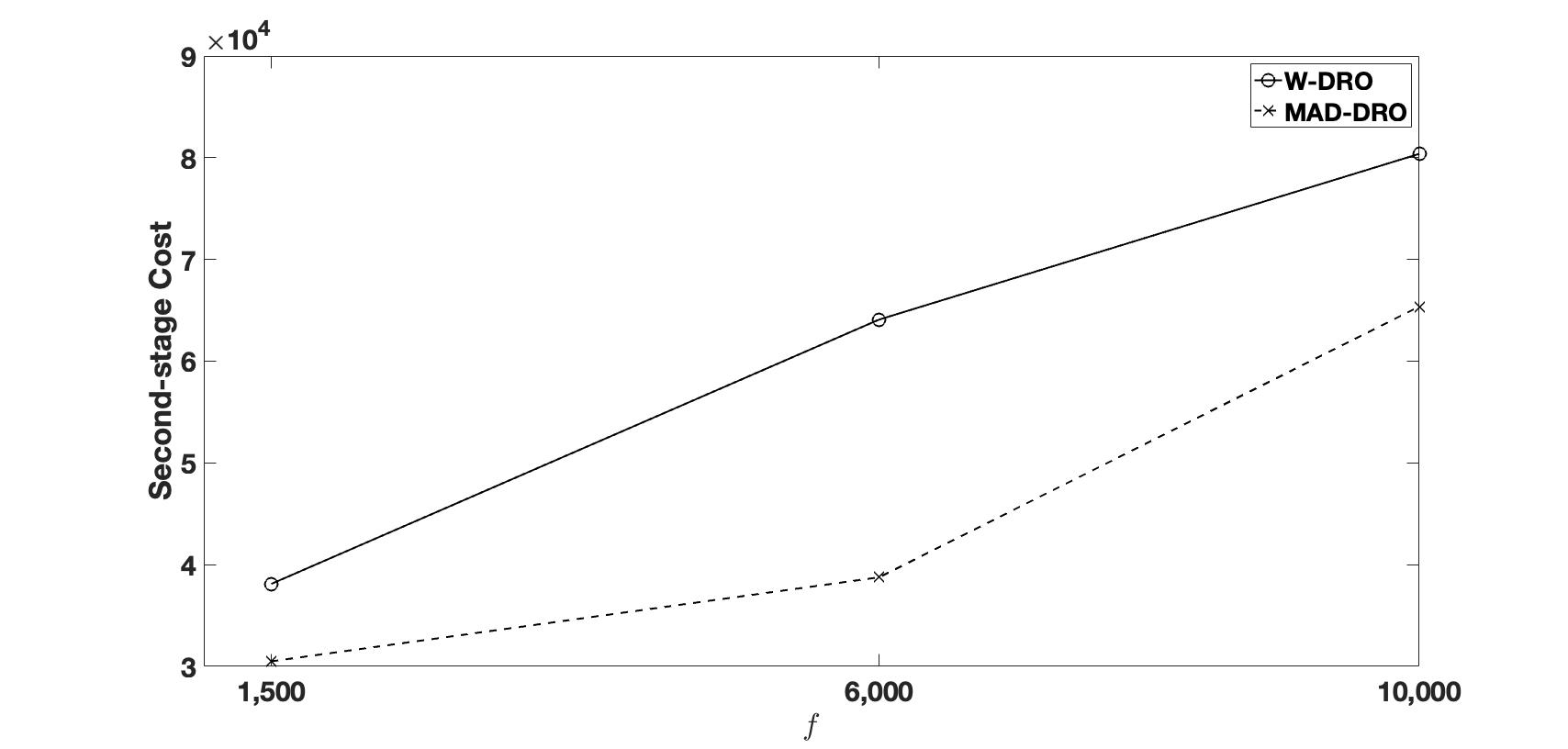}
      \caption{$\gamma=0.25\gamma_o$,Instance 5}
      \label{Fig13b2}
    \end{subfigure}%
    
  \begin{subfigure}[b]{0.5\textwidth}
 \centering
        \includegraphics[width=\textwidth]{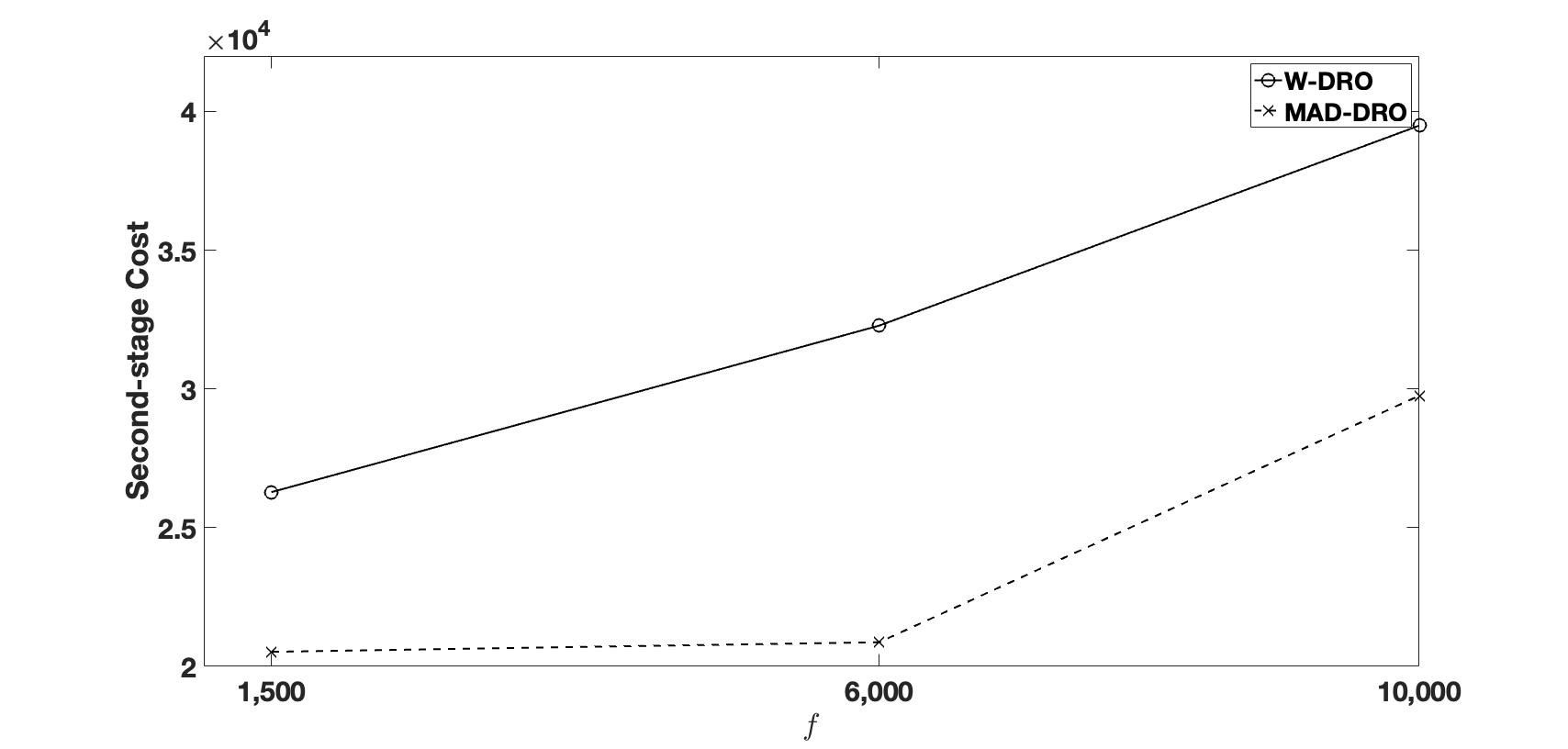}
        \caption{$\gamma=0.35\gamma_o$, Instance 1}
        \label{Fig13a3}
    \end{subfigure}%
    \begin{subfigure}[b]{0.5\textwidth}
            \includegraphics[width=\textwidth]{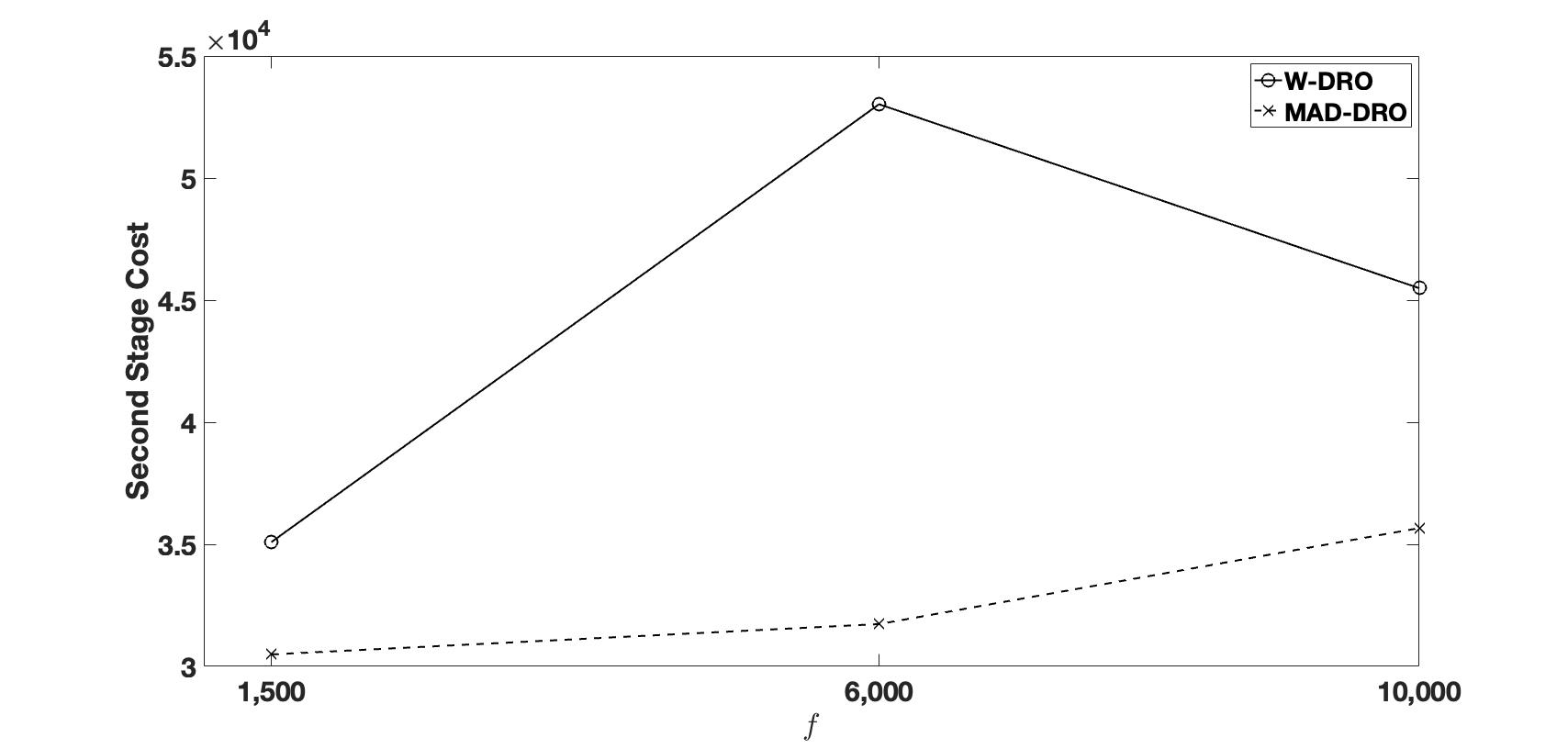}
      \caption{$\gamma=0.35\gamma_o$,Instance 5}
      \label{Fig13b3}
    \end{subfigure}%
    
      \begin{subfigure}[b]{0.5\textwidth}
 \centering
        \includegraphics[width=\textwidth]{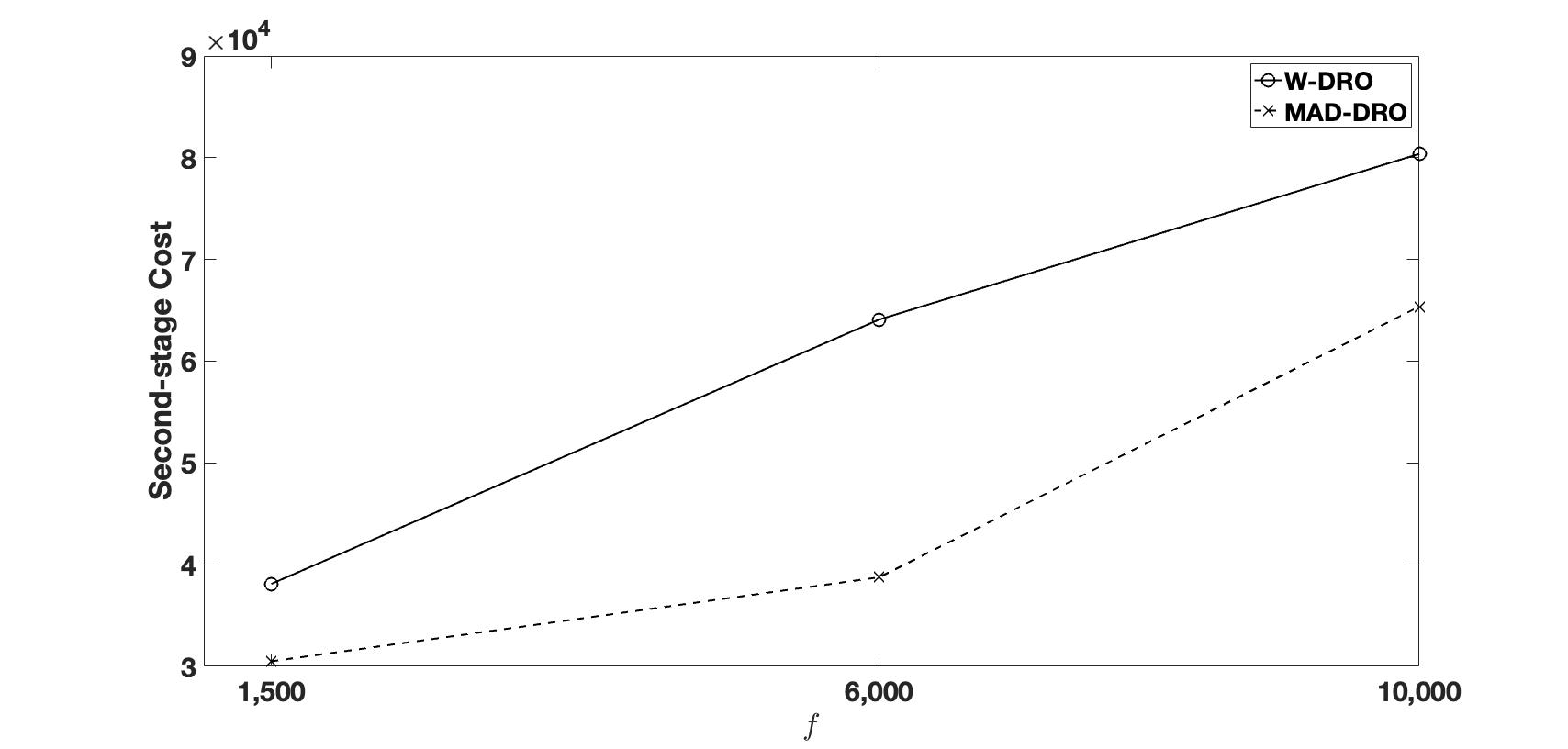}
        \caption{$\gamma=0.50\gamma_o$, Instance 1}
        \label{Fig13a3}
    \end{subfigure}%
    \begin{subfigure}[b]{0.5\textwidth}
            \includegraphics[width=\textwidth]{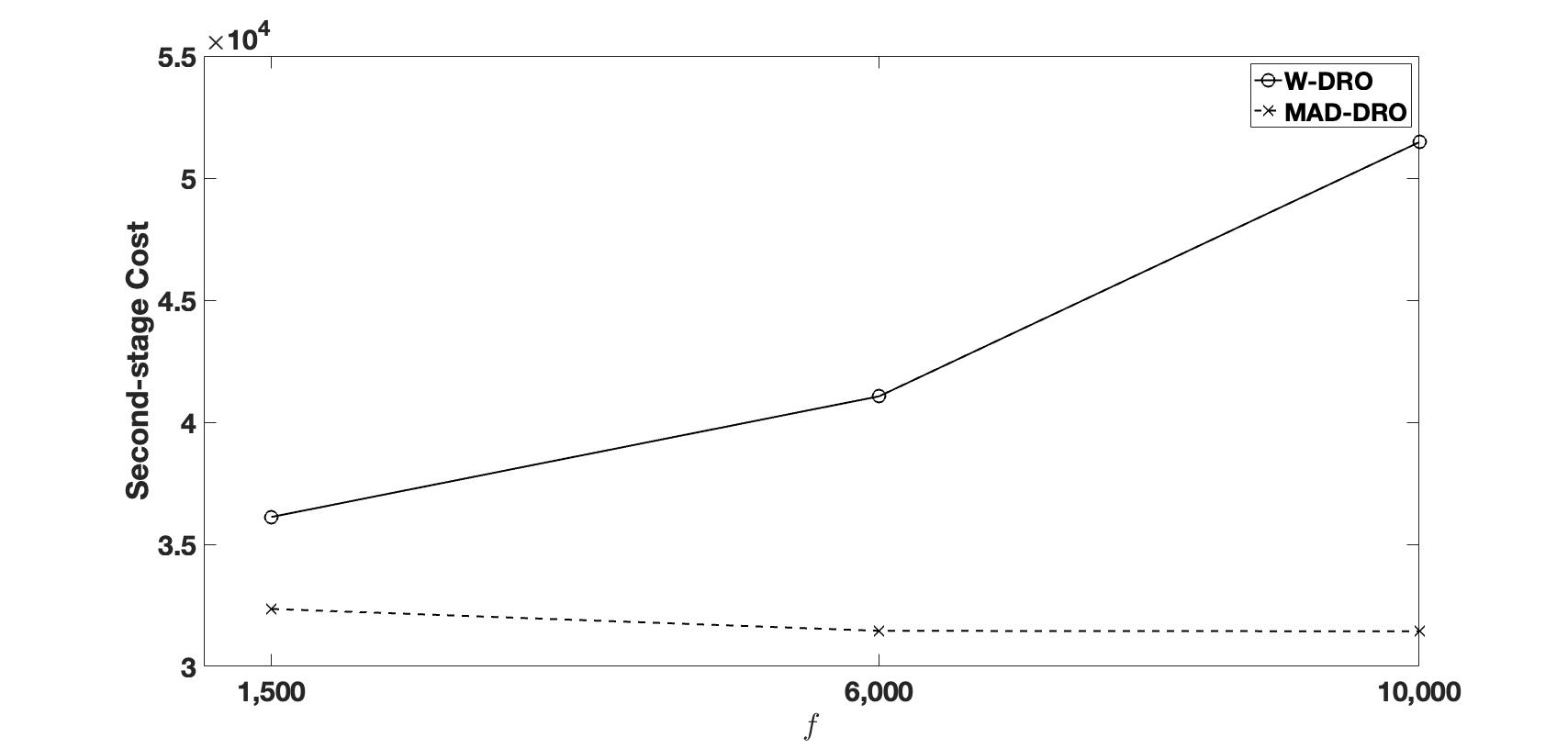}
      \caption{$\gamma=0.50\gamma_o$,Instance 5}
      \label{Fig13b3}
    \end{subfigure}%

    \caption{Comparison of the second-stage cost for different values of $\gamma$.}\label{Fig13:2nd_vs_gamma}
\end{figure}

Second,  we fix $M=20$ and solve the models with unmet demand penalty $\gamma \in \{0.10 \gamma_o, \ 0.25 \gamma_o,$  $0.35\gamma_0, \  0.50\gamma_o\}$ (where $\gamma_o$ is the base case penalty in Section~\ref{sec5.1:instancegen}) and $f \in \{1,500, \ 6,000, \ 10,000\}$. Figure~\ref{Fig12:MF_vs_gamma} presents the number of MFs as a function of $\gamma$ and $f$. Figure~\ref{Fig13:2nd_vs_gamma} presents the associated second-stage cost. It is not surprising that as $\gamma$ increases (i.e., satisfying customer demand becomes more important), the number of scheduled MFs increases. Note that by scheduling a larger number of MFs, we could satisfy a larger amount of demand and reduce the second-stage cost  (see  Figure~\ref{Fig13:2nd_vs_gamma}). However, for fixed $\gamma$, the MAD-DRO model schedules more MFs, and thus yields a lower unmet demand cost (because the MAD-DRO solutions satisfy a  larger amount of demand). For example, consider Instance 1 with $f=1,500$. When $\gamma$ decreases from $0.5\gamma_o$ to $0.1\gamma_o$ the optimal number of scheduled MFs of (W-DRO, MAD-DRO) decreases from (6, 6) to (1, 3) and average unmet demand cost increases from (9, 0) to (16,117,  10,827).

Our experiments in this section provide an example of how decision-makers can use our DRO approaches to generate MFRSP solutions under different parameter settings. Practitioners can use these results to decide whether to adopt the MAD-DRO model (which provides a better operational and computational performance) or the W-DRO model (which provides \textcolor{black}{a lower one-time fixed cost} for establishing the MF fleet).

\subsection{Analysis of the risk-averse solutions}\label{sec:CVaRExp}

\noindent In this section, we analyze the optimal solutions of the mean-CVaR-based models under some critical problem parameters.  Specifically, we solve the models with $f\in$\{1,500, 6,000,  10,000\}, $\gamma \in \{\gamma_o, \ 0.5 \gamma_o, \ 0.25\gamma_o,  \ 0.1 \gamma_o \} $ (where $\gamma_o$ is the base case penalty in Section~\ref{sec5.1:instancegen}), $\kappa=0.95$ (a typical value of $\kappa$; we observe similar results under $\kappa=0.99, \ 0.975$), and $\Theta \in \{0, 0.2, 0.5, 1\}$ (where a smaller $\Theta$ indicate that we are more risk-averse). \textcolor{black}{We use MAD-CVaR (W-CVaR) to denote the mean-CVaR-based DRO model with MAD (1-Wasserstein) ambiguity, and SP-CVaR to denote the mean-CVaR-based SP model.}  Because the SP-CVaR model cannot solve large instances (even with small $N$), for a fair comparison and brevity, we present results for Instance 1. We keep all other parameters as described in Section~\ref{sec5.1:instancegen}.

Table~\ref{tableCVaROptimal} presents the optimal number of scheduled MFs  (i.e., fleet size)  for different $f$, $\gamma$, and $\Theta$. We make the following observations from this table. First, all models schedule fewer MFs as $f$ increases and $\gamma$ decreases irrespective of $\Theta$ (risk-aversion coefficient). This is consistent with our results in Section~\ref{sec5:sensitivity} for the risk neutral models.  Second, the MAD-CVaR model often schedules a higher number of MFs than the W-CVaR model, and the latter model schedule the same or larger number of MFs than the SP-CVaR.  Second, when $f=$\{6,000, 10,000\} (i.e., average and high fixed cost) and $\gamma=0.1\gamma_o$ (very low unmet demand penalty) all models schedule one MF.

\textcolor{black}{Third, the MAD-CVaR model schedules the same number of MFs under all values of $\Theta$ when ($f$, $\gamma$)$=$ (1,500, $\gamma_o$) (and similarly under (6,000, 0.25$\gamma_o$) and (10,000, $\gamma_o$)). Similarly, the W-CVaR model schedules the same number of MFs under all values of $\Theta$ when ($f$, $\gamma$)= (6,000, 0.25$\gamma_o$), and (10,000, 0.5$\gamma_o$). These results indicate that our proposed DRO expectation models with MAD and 1-Wasserstein ambiguity are risk-averse under these settings because they yield the same optimal solutions under all values of the risk-aversion coefficient $\Theta$.}
\begin{table}[t!]
\center 
\caption{Comparison of the optimal number of MFs yielded by each model under different values of $\Theta$, $f$, and $\gamma$. Notation: MAD-CVaR and W-CVaR are respectively the distributionally robust mean-CVaR models with MAD and Wasserstein ambiguity, and SP-CVaR is the SP model based on the  mean-CVaR criterion.}
   \renewcommand{\arraystretch}{0.6}
\begin{tabular}{lllllllllllllllllllllllllllllllllll}
 \hline
  \multicolumn{11}{c}{\textbf{$\pmb{f=1,500}$}}\\
 \hline
&   \multicolumn{4}{c}{$\gamma_o$ }  &&  & \multicolumn{4}{c}{$0.5\gamma_o$ } \\   \cline{2-5} \cline{8-11}
Model	&	0	&	0.2	&	0.5	&	1	&	&	Model	 &	0	&	0.2	&	0.5	&	1	\\
\hline
MAD-CVaR	&	8	&	8	&	8	&	8	&	&	MAD-CVaR	&	8	&	8	&	8	&	7	\\
W-CVaR	&	8	&	8	&	8	&	6	&	&	W-CVaR	&	7	&	7	&	6	&	6	\\
SP-CVaR	&	6	&	6	&	6	&	5	&	&	SP-CVaR	&	6	&	6	&	6	&	5	\\
\\
\hline
&   \multicolumn{4}{c}{$0.25\gamma_o$ }  &&  &  \multicolumn{4}{c}{$0.1\gamma_o$ } \\   \cline{2-5} \cline{8-11}
Model	&	0	&	0.2	&	0.5	&	1	&	&	Model	&	0	&	0.2	&	0.5	&	1	\\
\hline
MAD-CVaR	&	6	&	6	&	6	&	5	&	&	MAD-CVaR	&	3	&	3	&	3	&	1	\\
W-CVaR	&	\textcolor{red}{6}	&	5	&	5	&	5	&	&	W-CVaR	&	2	&	2	&	2	&	1	\\
SP-CVaR	&	4	&	4	&	4	&	4	&	&	SP-CVaR	&	2	&	2	&	2	&	1	\\
\hline
  \multicolumn{11}{c}{\textbf{$\pmb{f=6,000}$}}\\
  \hline
&   \multicolumn{4}{c}{$\gamma_o$ }  &&  &  \multicolumn{4}{c}{$0.5\gamma_o$ } \\   \cline{2-5} \cline{8-11}
Model	&	0	&	0.2	&	0.5	&	1	&	&	Model	&	0	&	0.2	&	0.5	&	1	\\
\hline
MAD-CVaR	&	\textcolor{black}{7}	&	\textcolor{black}{7}		&	7	&	6	&	&	MAD-CVaR	&	7	&	7	&	6	&	5	\\
W-CVaR	&	6	&	6	&	5	&	5	&	&	W-CVaR	&	6	&	5	&	4	&	4	\\
SP-CVaR	&	5	&	5	&	5	&	5	&	&	SP-CVaR	&	5	&	5	&	4	&	4	\\
\hline 
&   \multicolumn{4}{c}{$0.25\gamma_o$ }  &&   & \multicolumn{4}{c}{$0.1\gamma_o$ } \\   \cline{2-5} \cline{8-11}
Model	&	0	&	0.2	&	0.5	&	1	&	&	Model	&	0	&	0.2	&	0.5	&	1\\
\hline
MAD-CVaR	&	3	&	3	&	3	&	3	&	&	MAD-CVaR	&	1	&	1	&	1	&	1	\\
W-CVaR	&	2	&	2	&	2	&	2	&	&	W-CVaR	&	1	&	1	&	1	&	1	\\
SP-CVaR	&	2	&	2	&	2	&	2	&	&	SP-CVaR	&	1	&	1	&	1	&	1	\\
\hline 
  \multicolumn{11}{c}{\textbf{$\pmb{f=10,000}$}}\\
  \hline
&   \multicolumn{4}{c}{$\gamma_o$ }  &&   & \multicolumn{4}{c}{$0.5\gamma_o$ } \\   \cline{2-5} \cline{8-11}
Model	&	0	&	0.2	&	0.5	&	1	&	&	Model	&	0	&	0.2	&	0.5	&	1	\\
\hline
MAD-CVaR	&	6	&	6	&	6	&	6	&	&	MAD-CVaR	&	6	&	5	&	5	&	4	\\
W-CVaR	&	6	&	6	&	6	&	5	&	&	W-CVaR	&	4	&	4	&	4	&	4	\\
SP-CVaR	&	5	&	4	&	4	&	4	&	&	SP-CVaR	&	4	&	4	&	4	&	3	\\
\hline 
&   \multicolumn{4}{c}{$0.25\gamma_o$ }  &&  &   \multicolumn{4}{c}{$0.1\gamma_o$ } \\   \cline{2-5} \cline{8-11}
Model	&	0	&	0.2	&	0.5	&	1	&	&	Model	&	0	&	0.2	&	0.5	&	1	\\
\hline 
MAD-CVaR	&	1	&	1	&	1	&	1	&	&	MAD-CVaR	&	1	&	1	&	1	&	1	\\
W-CVaR	&	1	&	1	&	1	&	1	&	&	W-CVaR	&	1	&	1	&	1	&	1	\\
SP-CVaR	&	1	&	1	&	1	&	1	&	&	SP-CVaR	&	1	&	1	&	1	&	1	\\
\hline 
\end{tabular}\label{tableCVaROptimal}
\end{table}


Fourth, all models schedule more MFs under a smaller $\Theta$, especially when $\gamma=(0.25\gamma_o, 0.1\gamma_o)$ with $f=1,500$ (i.e., a low cost and lower unmet demand penalty), $\gamma=(0.5\gamma_o,0.25\gamma_o)$ with $f=6,000$, and $\gamma=(\gamma_o, 0.5\gamma_o)$ with $f=10,000$. In particular, more MFs are scheduled when $\Theta=0$ (risk-averse models with CVaR criterion) as compared to $\Theta=1$ (risk-neutral models).  This makes sense because a risk-averse decision-maker may schedule more MFs to avoid high operational cost and, in particular, excessive shortages.

\color{black}

\begin{figure}
      \begin{subfigure}[b]{0.5\textwidth}
    \centering
        \includegraphics[width=\textwidth]{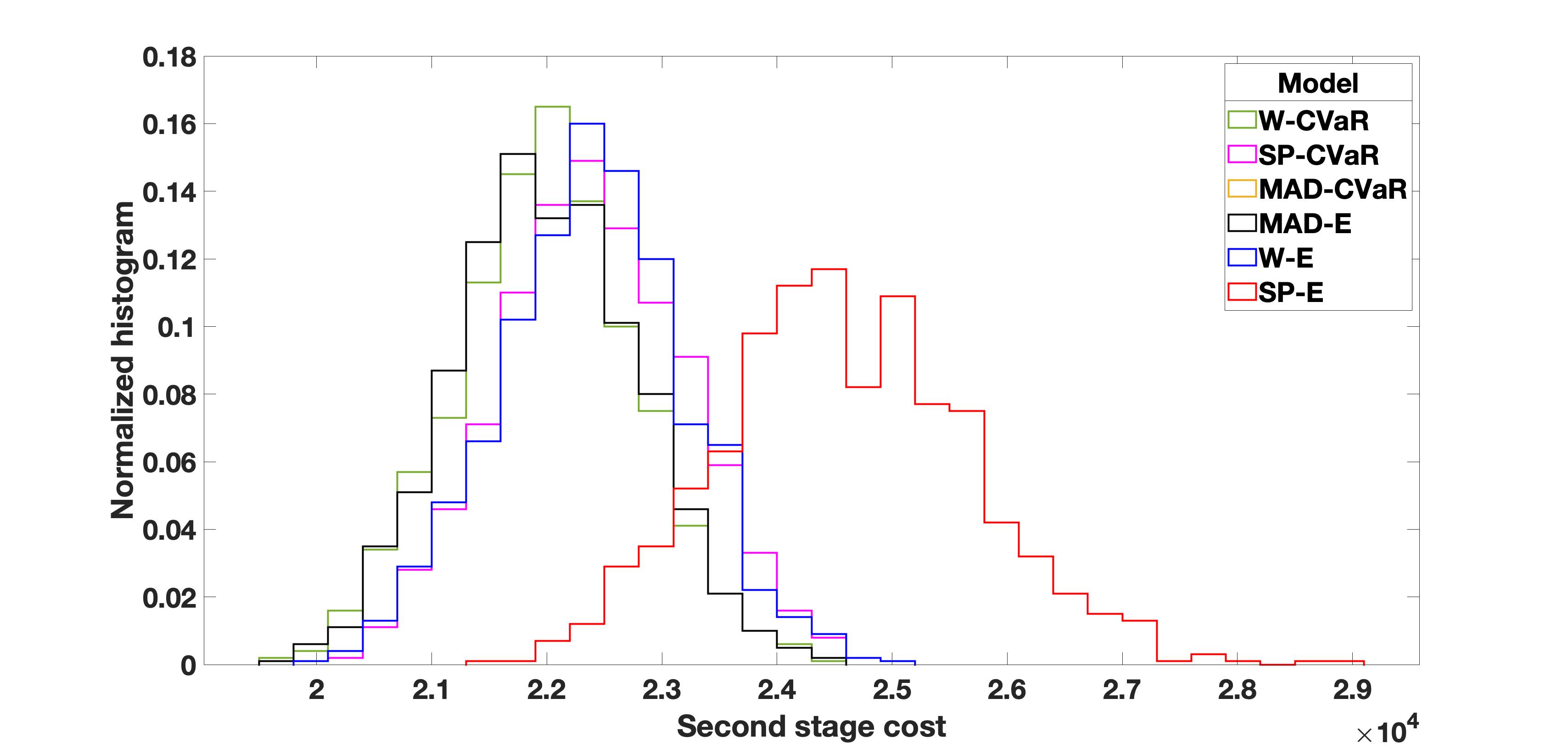}
        \caption{2nd, $\Delta=0$}
    \end{subfigure}%
      \begin{subfigure}[b]{0.5\textwidth}
        \centering
        \includegraphics[width=\textwidth]{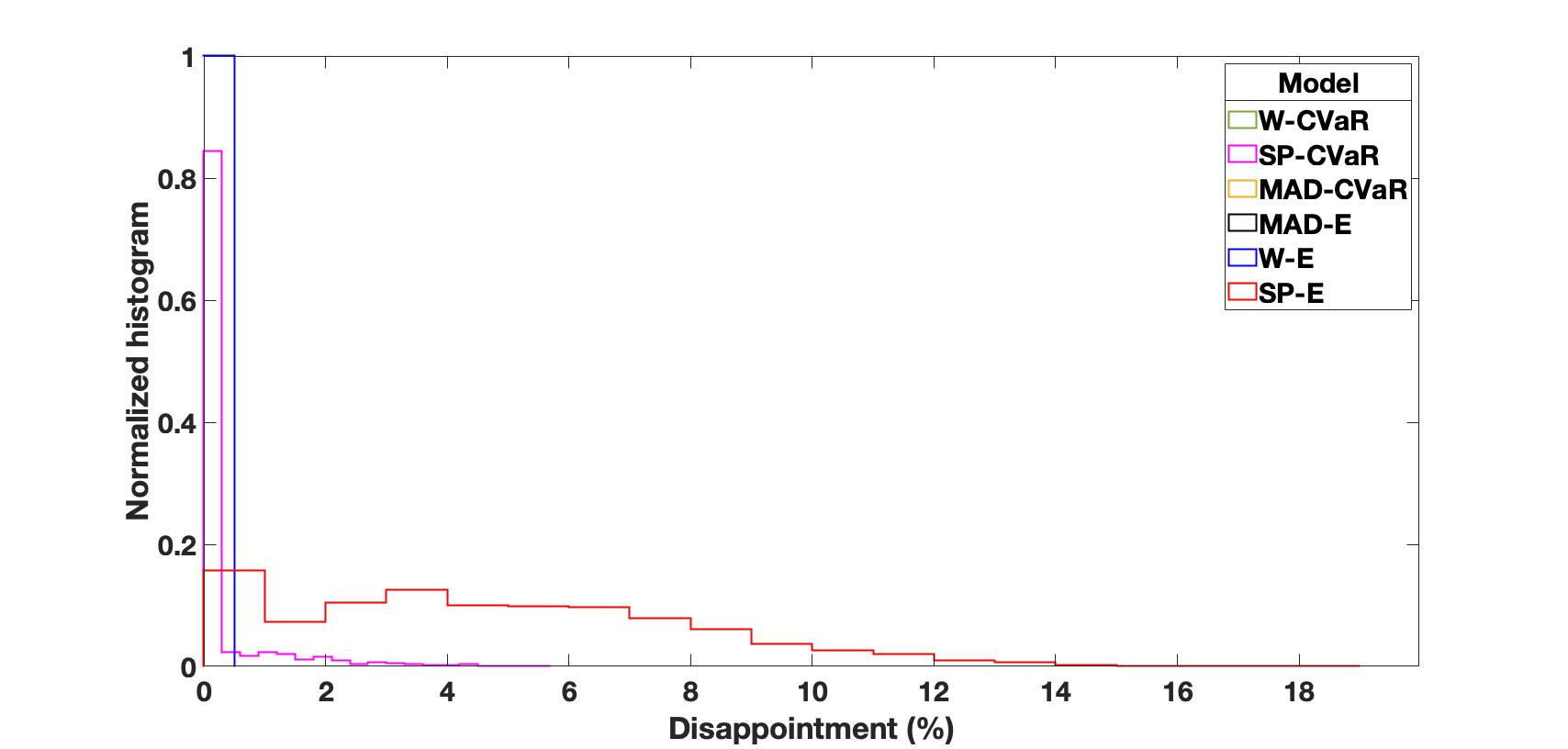}
        \caption{Disappointment, $\Delta=0$}
    \end{subfigure}%
    
            \centering
      \begin{subfigure}[b]{0.5\textwidth}
          \centering
        \includegraphics[width=\textwidth]{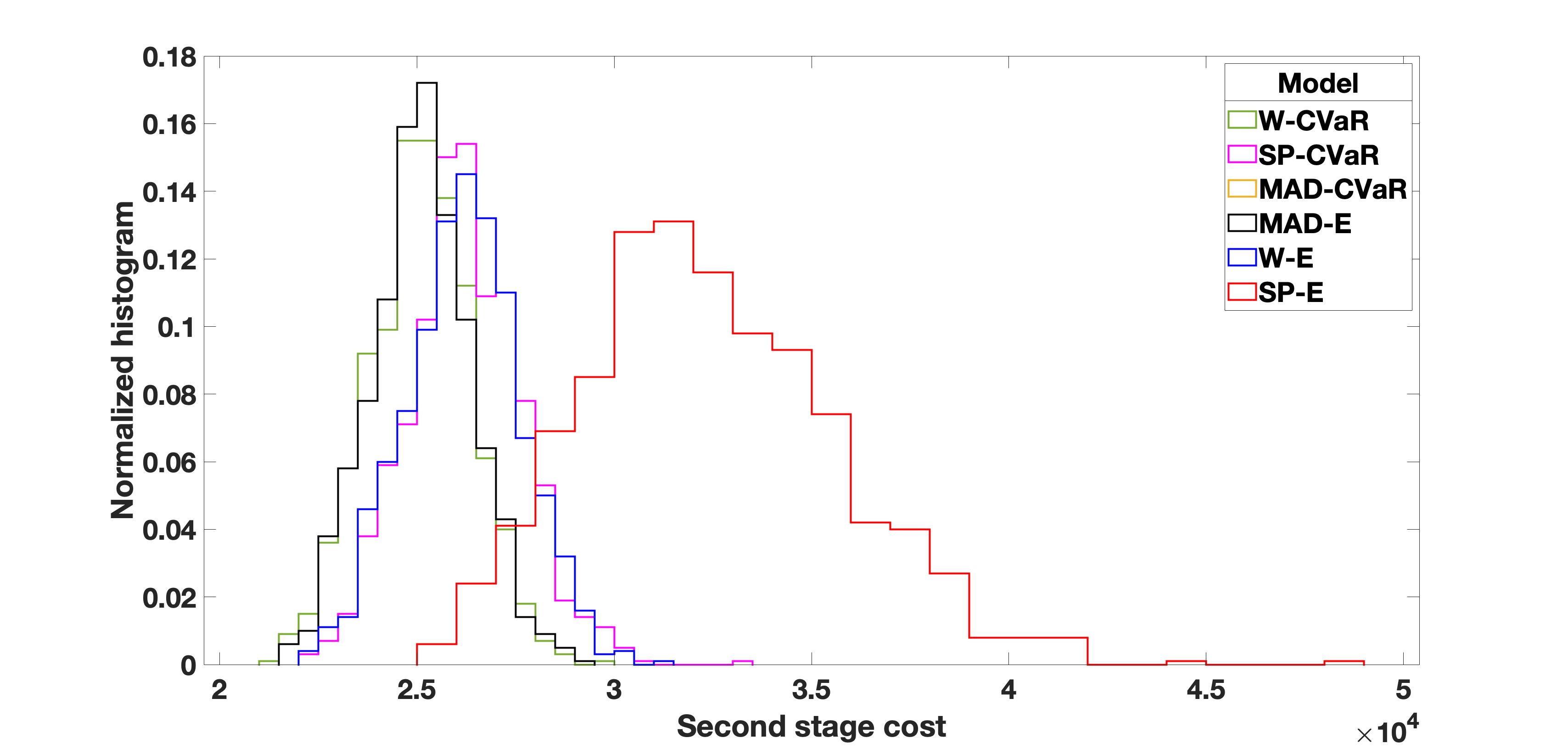}
        \caption{2nd, $ \Delta=0.25$}
    \end{subfigure}%
          \begin{subfigure}[b]{0.5\textwidth}
        \centering
        \includegraphics[width=\textwidth]{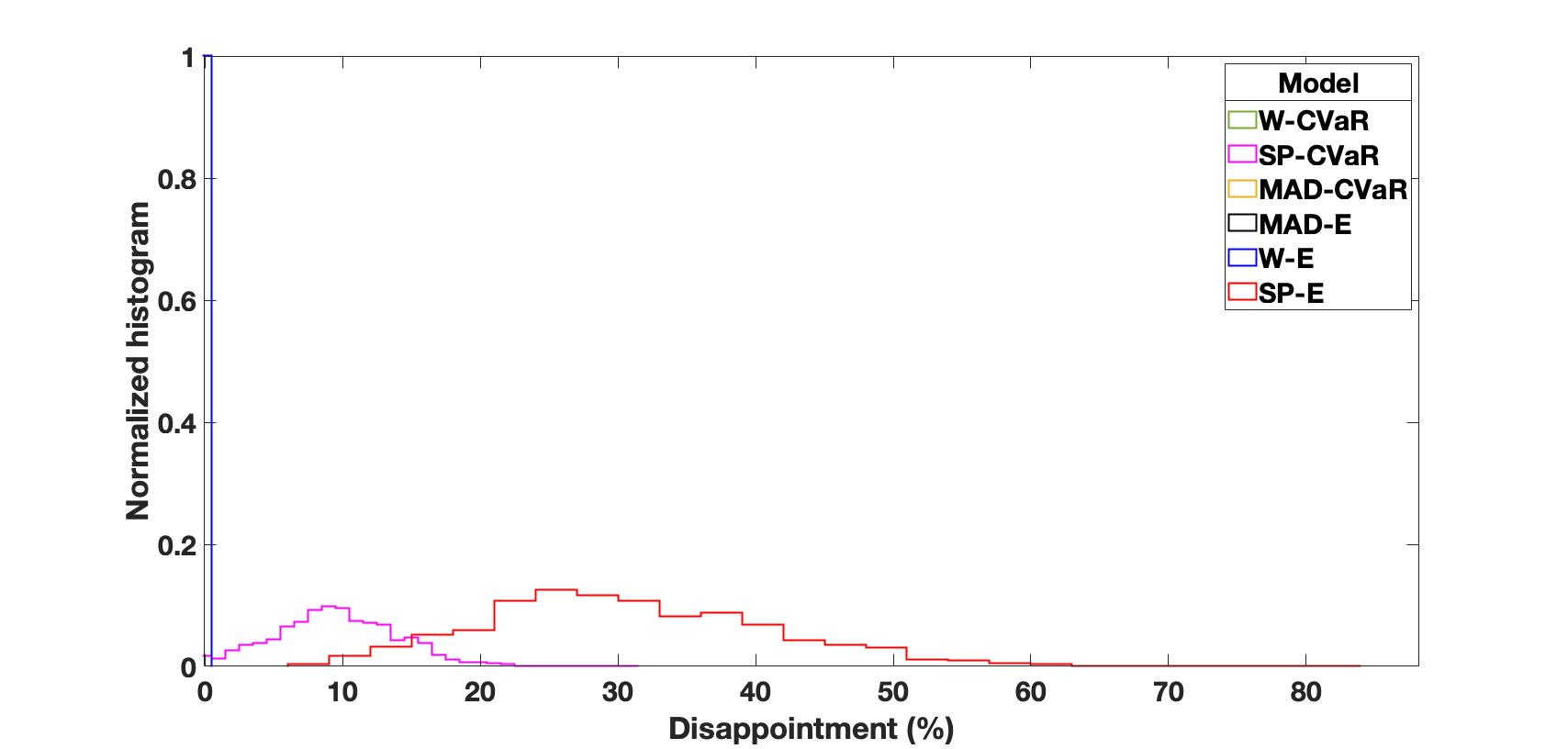}
        \caption{Disappointment, $\Delta=0.25$}
    \end{subfigure}%

     \begin{subfigure}[b]{0.5\textwidth}
          \centering
        \includegraphics[width=\textwidth]{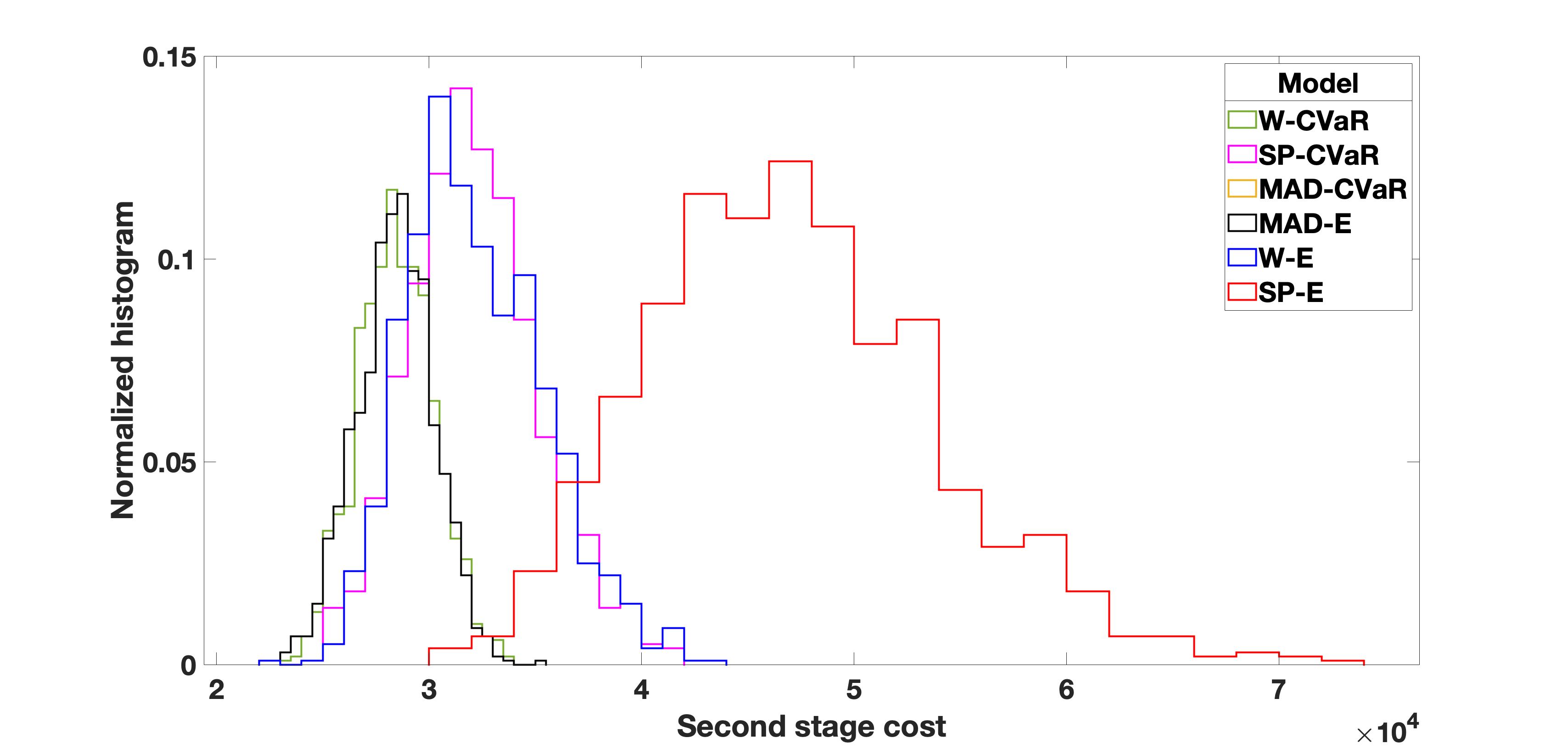}
        \caption{2nd, $ \Delta=0.5$}
    \end{subfigure}%
          \begin{subfigure}[b]{0.5\textwidth}
        \centering
        \includegraphics[width=\textwidth]{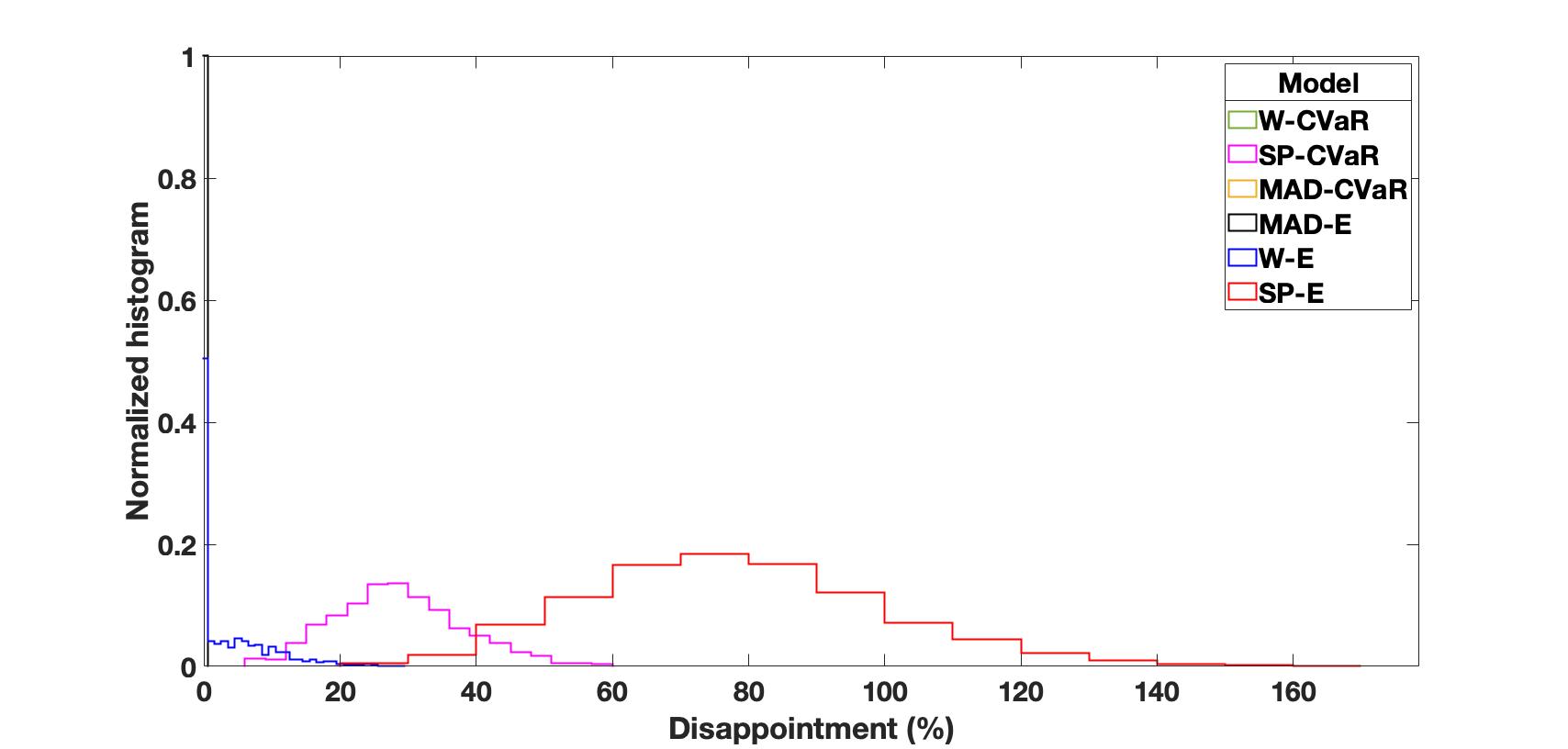}
        \caption{Disappointment, $\Delta=0.5$}
    \end{subfigure}%

\caption{Normalized histograms of second-stage cost (2nd) and out-of-sample disappointments under Set 2 with $\gamma=\gamma_o$, $f=1,500$, and $\pmb{\Delta \in \{0, 0.25, 0.5\}}$.}\label{CVaR_gamma}
\end{figure}

\begin{figure}
      \begin{subfigure}[b]{0.5\textwidth}
    \centering
        \includegraphics[width=\textwidth]{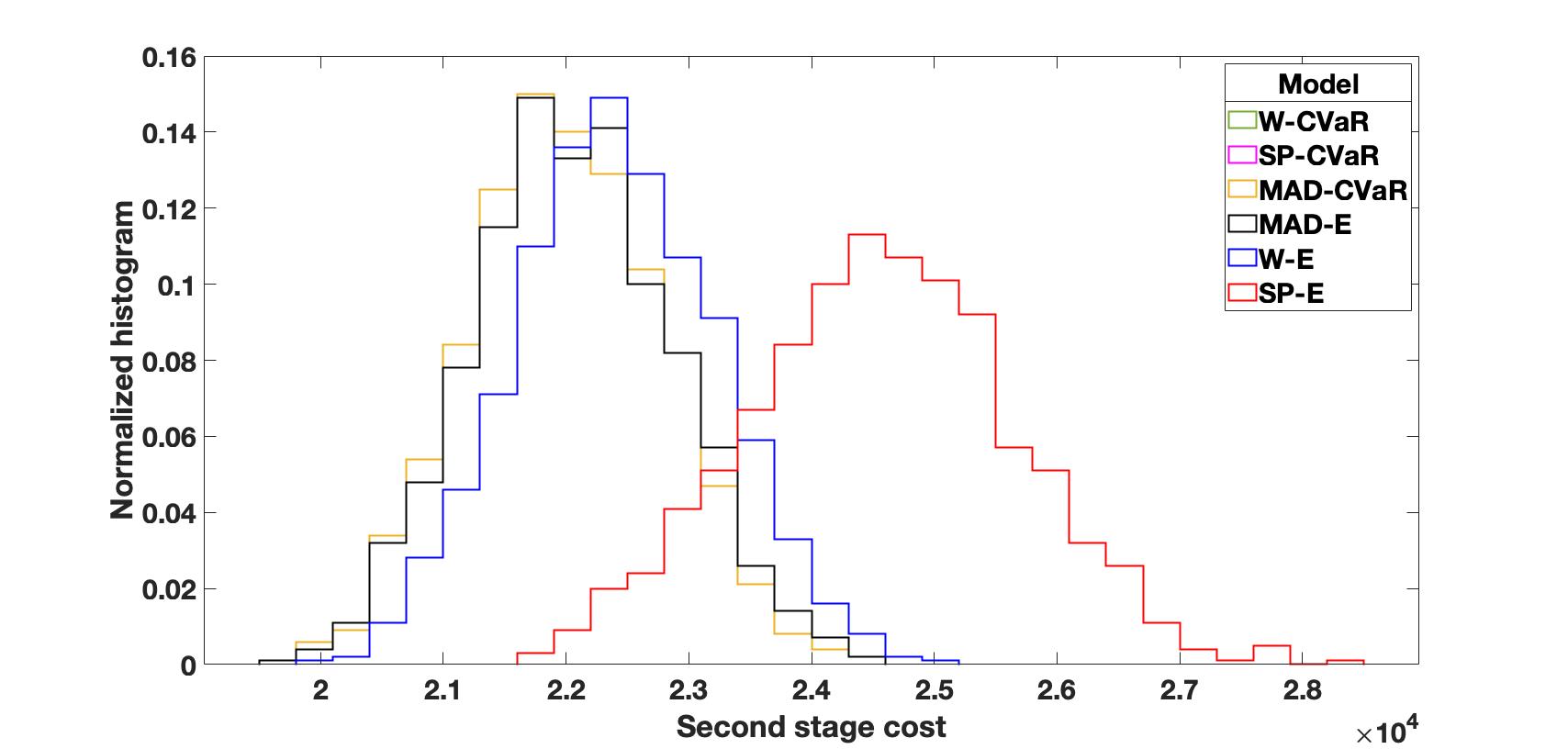}
        \caption{2nd, $\Delta=0$}
    \end{subfigure}%
      \begin{subfigure}[b]{0.5\textwidth}
        \centering
        \includegraphics[width=\textwidth]{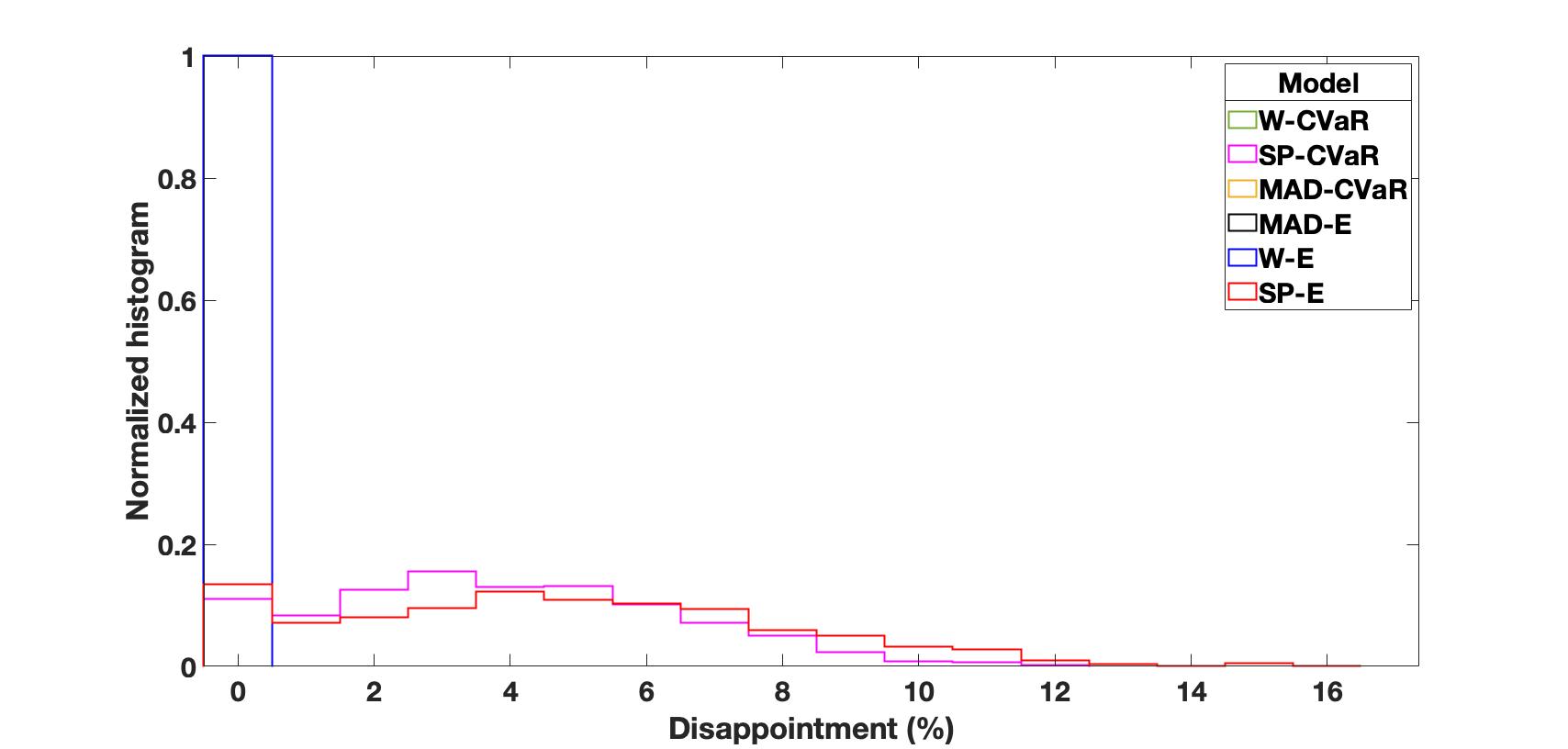}
        \caption{Disappointment, $\Delta=0$}
    \end{subfigure}%
    
            \centering
      \begin{subfigure}[b]{0.5\textwidth}
          \centering
        \includegraphics[width=\textwidth]{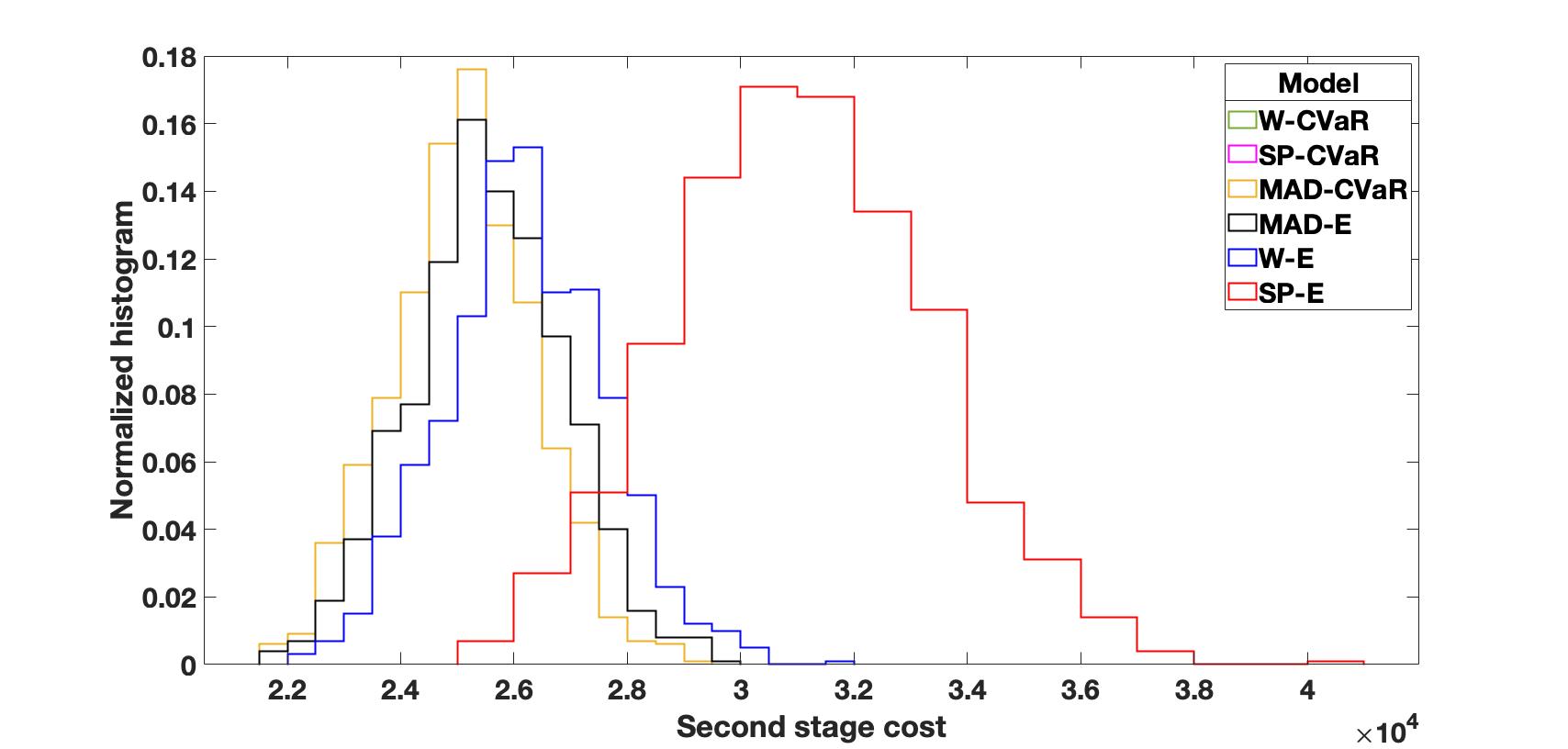}
        \caption{2nd, $ \Delta=0.25$}
    \end{subfigure}%
          \begin{subfigure}[b]{0.5\textwidth}
        \centering
        \includegraphics[width=\textwidth]{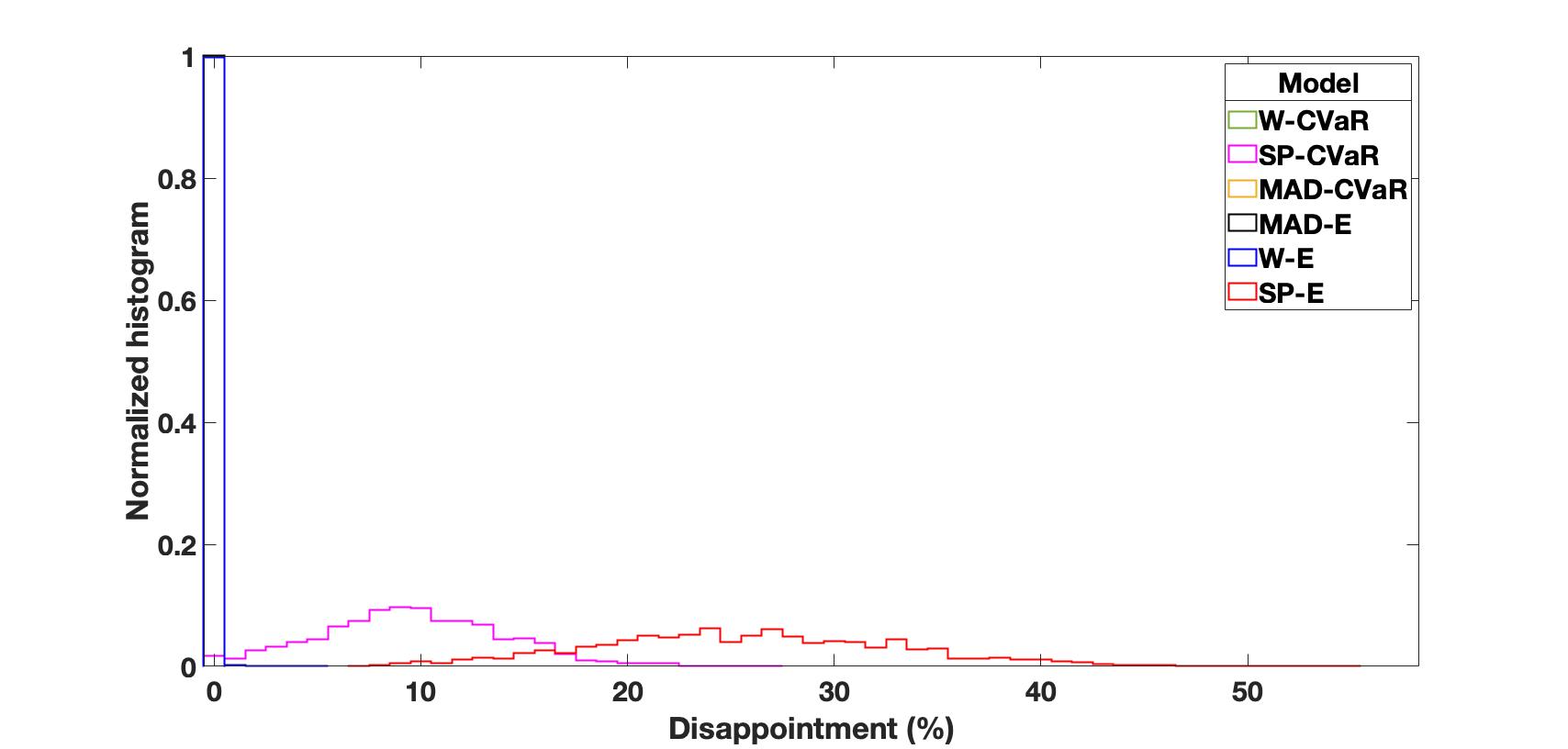}
        \caption{Disappointment, $\Delta=0.25$}
    \end{subfigure}%

     \begin{subfigure}[b]{0.5\textwidth}
          \centering
        \includegraphics[width=\textwidth]{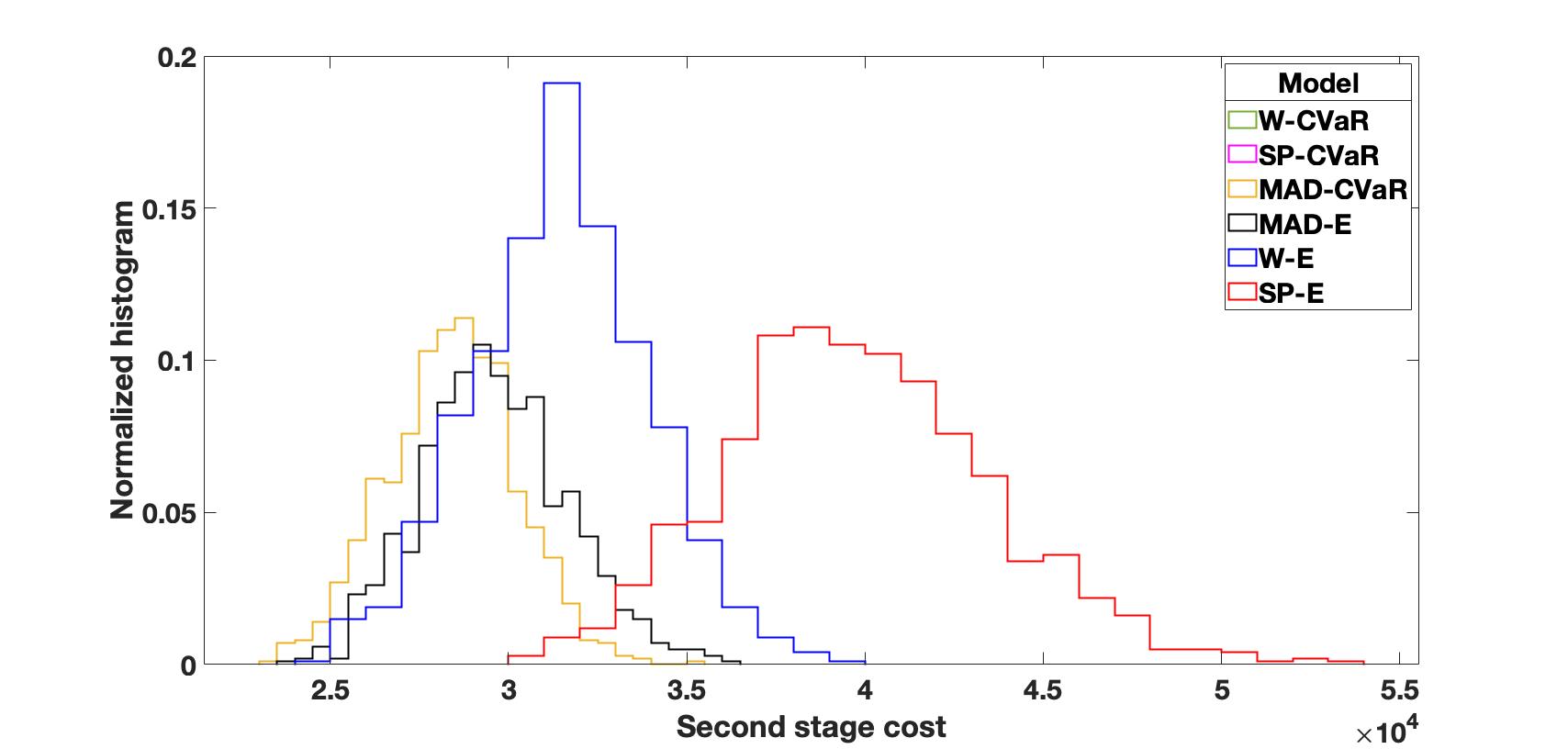}
        \caption{2nd, $ \Delta=0.5$}
    \end{subfigure}%
          \begin{subfigure}[b]{0.5\textwidth}
        \centering
        \includegraphics[width=\textwidth]{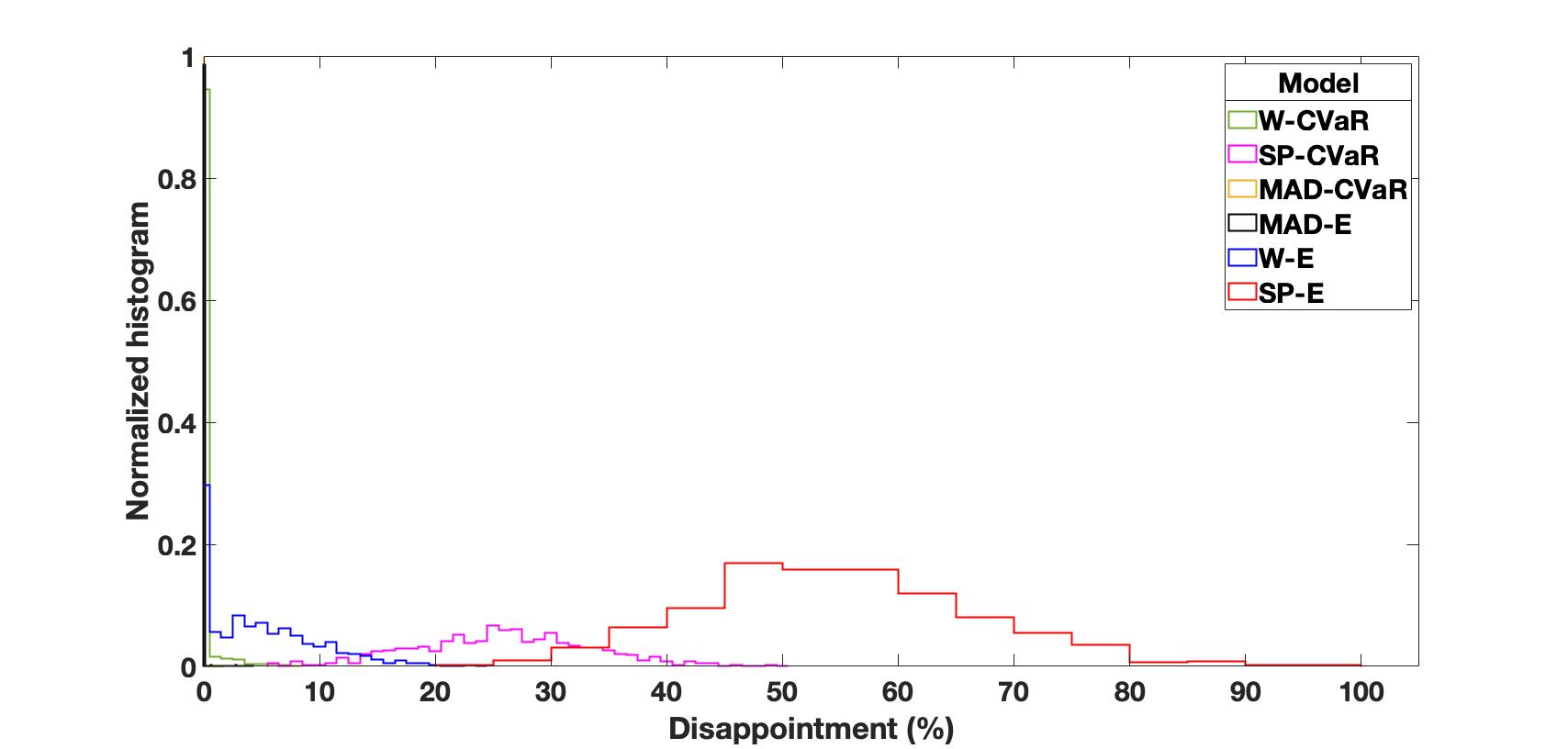}
        \caption{Disappointment, $\Delta=0.5$}
    \end{subfigure}%

\caption{Normalized histograms of second-stage cost (2nd) and out-of-sample disappointments under Set 2 with $\gamma=0.5\gamma_o$, $f=1,500$, and $\pmb{\Delta \in \{0, 0.25, 0.5\}}$.}\label{CVaR_50gamma}
\end{figure}

\begin{figure}
      \begin{subfigure}[b]{0.5\textwidth}
    \centering
        \includegraphics[width=\textwidth]{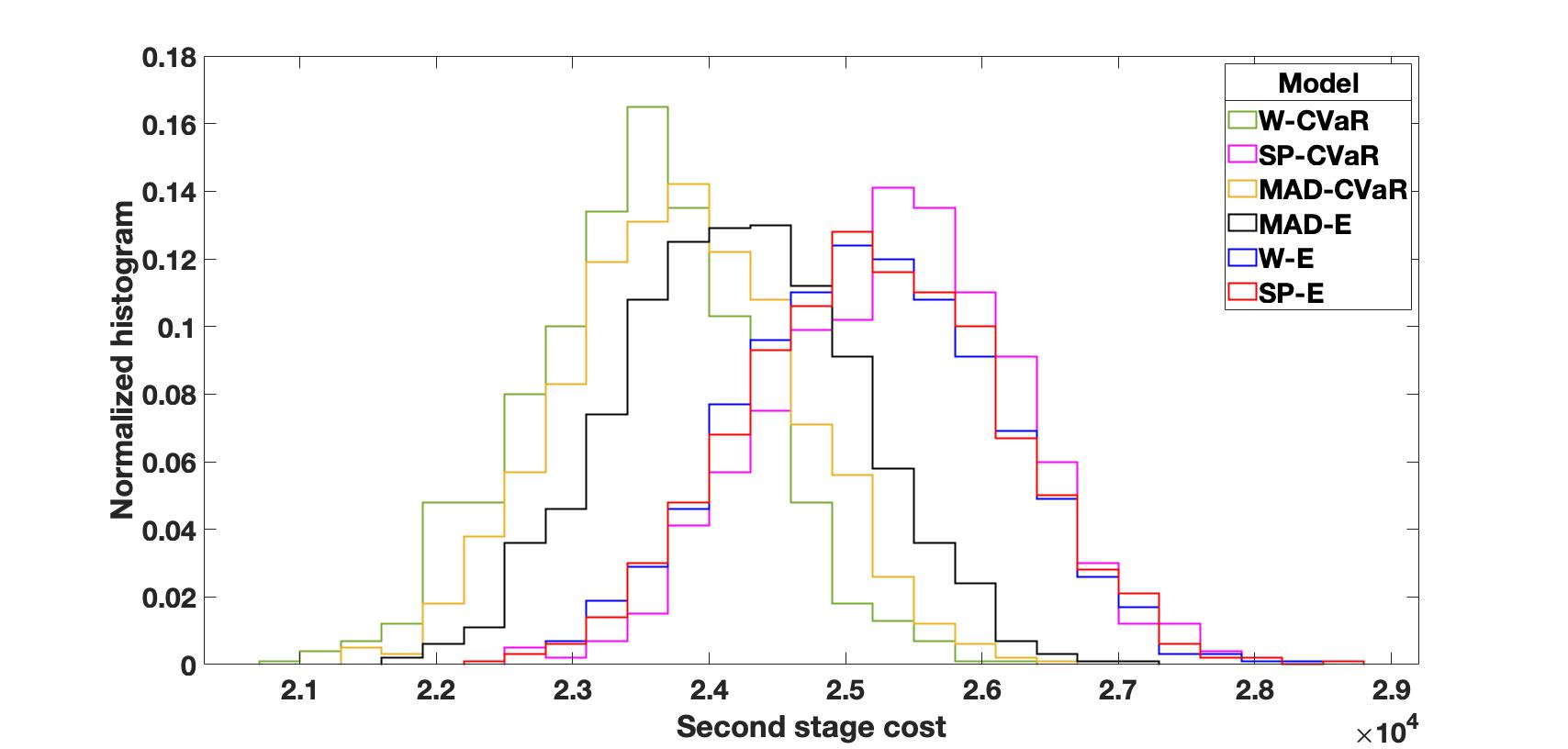}
        \caption{2nd, $\Delta=0$}
    \end{subfigure}%
      \begin{subfigure}[b]{0.5\textwidth}
        \centering
        \includegraphics[width=\textwidth]{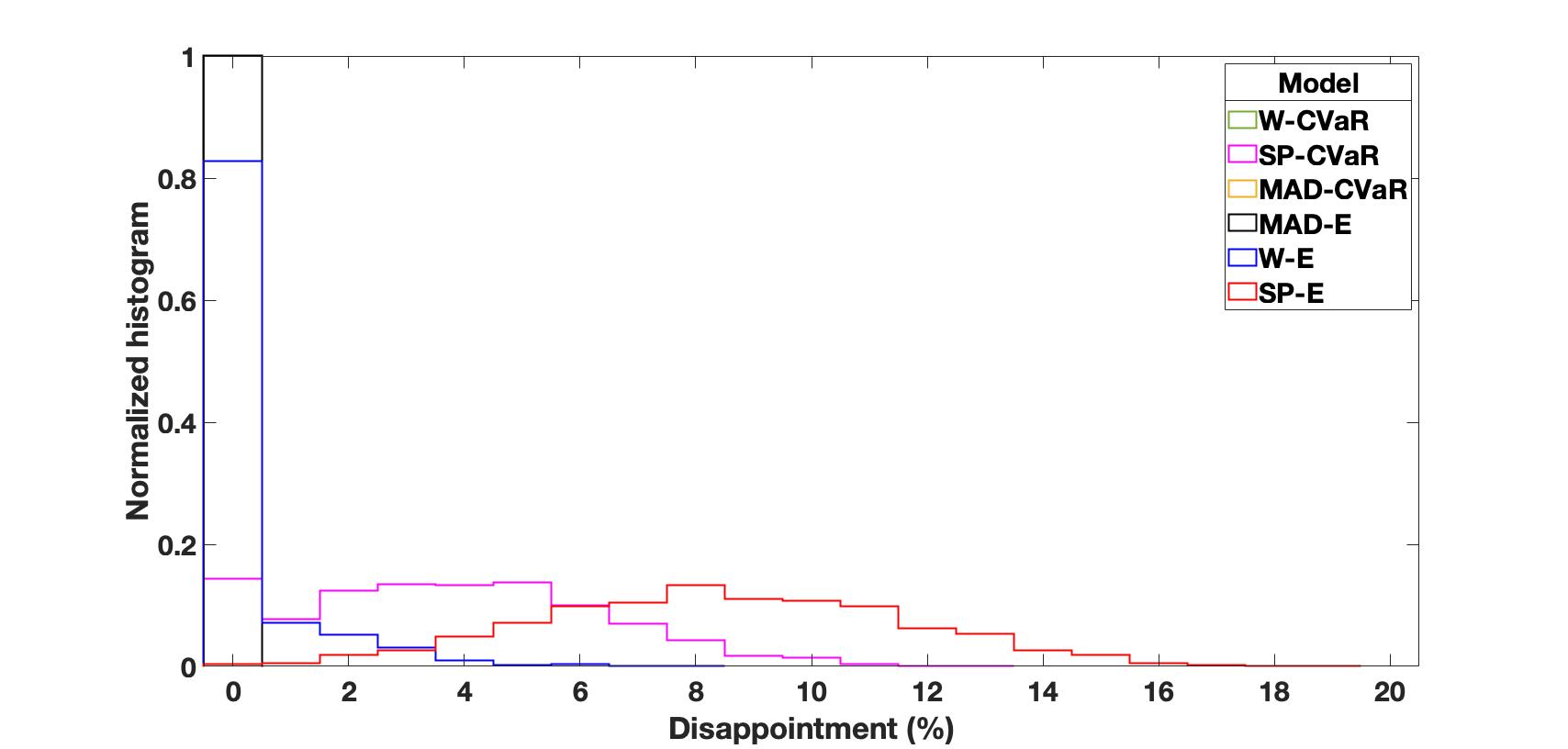}
        \caption{Disappointment, $\Delta=0$}
    \end{subfigure}%
    
            \centering
      \begin{subfigure}[b]{0.5\textwidth}
          \centering
        \includegraphics[width=\textwidth]{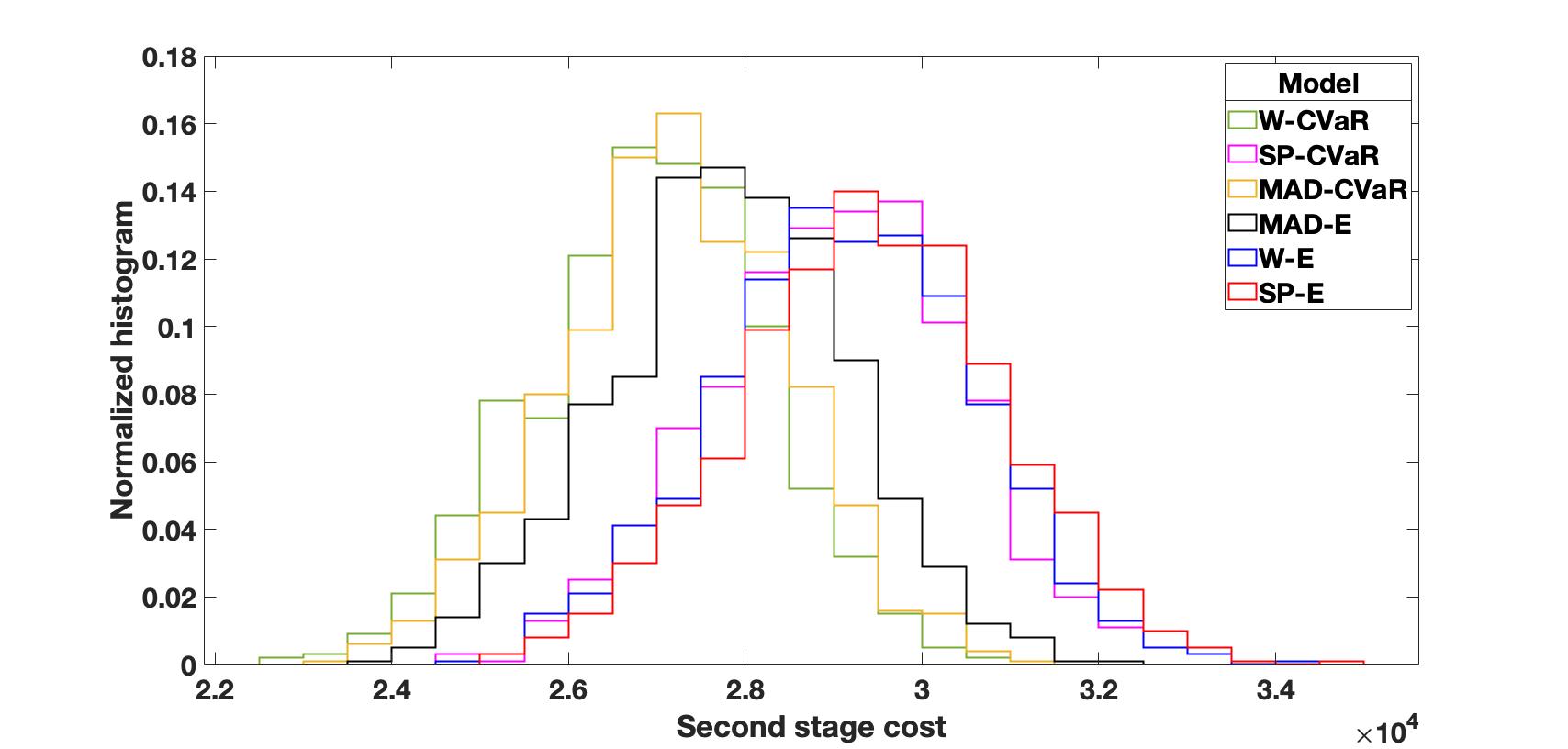}
        \caption{2nd, $ \Delta=0.25$}
    \end{subfigure}%
          \begin{subfigure}[b]{0.5\textwidth}
        \centering
        \includegraphics[width=\textwidth]{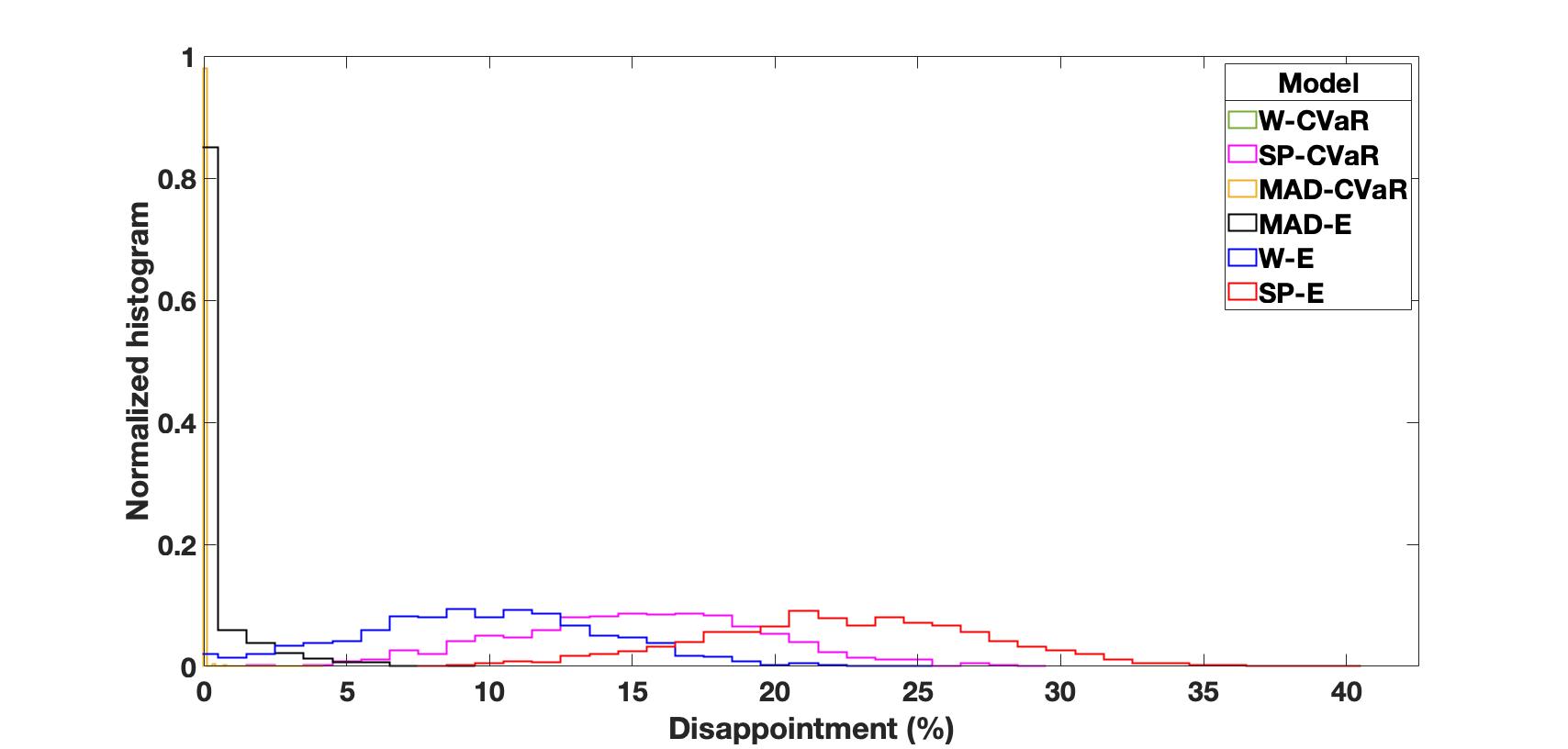}
        \caption{Disappointment, $\Delta=0.25$}
    \end{subfigure}%

     \begin{subfigure}[b]{0.5\textwidth}
          \centering
        \includegraphics[width=\textwidth]{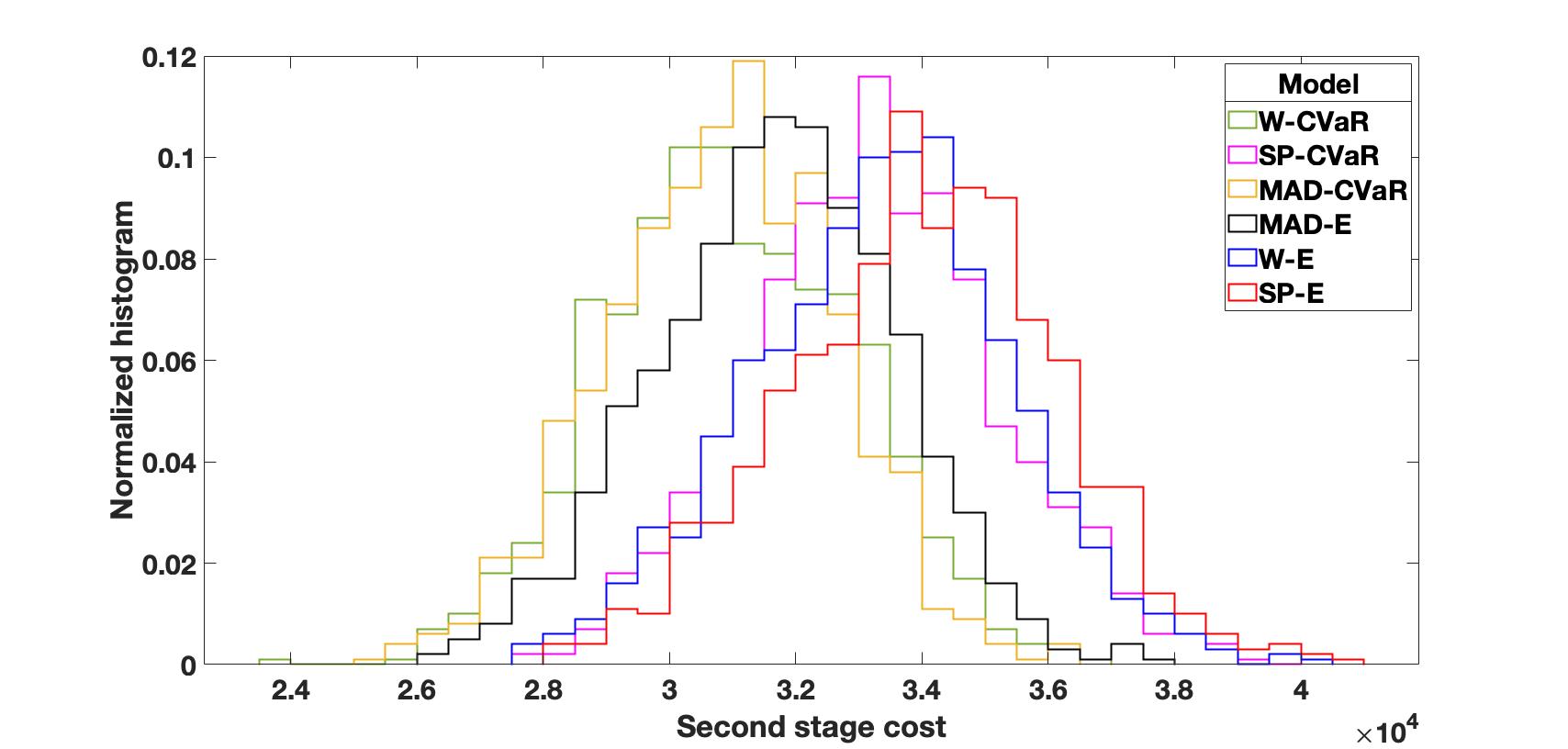}
        \caption{2nd, $ \Delta=0.5$}
    \end{subfigure}%
          \begin{subfigure}[b]{0.5\textwidth}
        \centering
        \includegraphics[width=\textwidth]{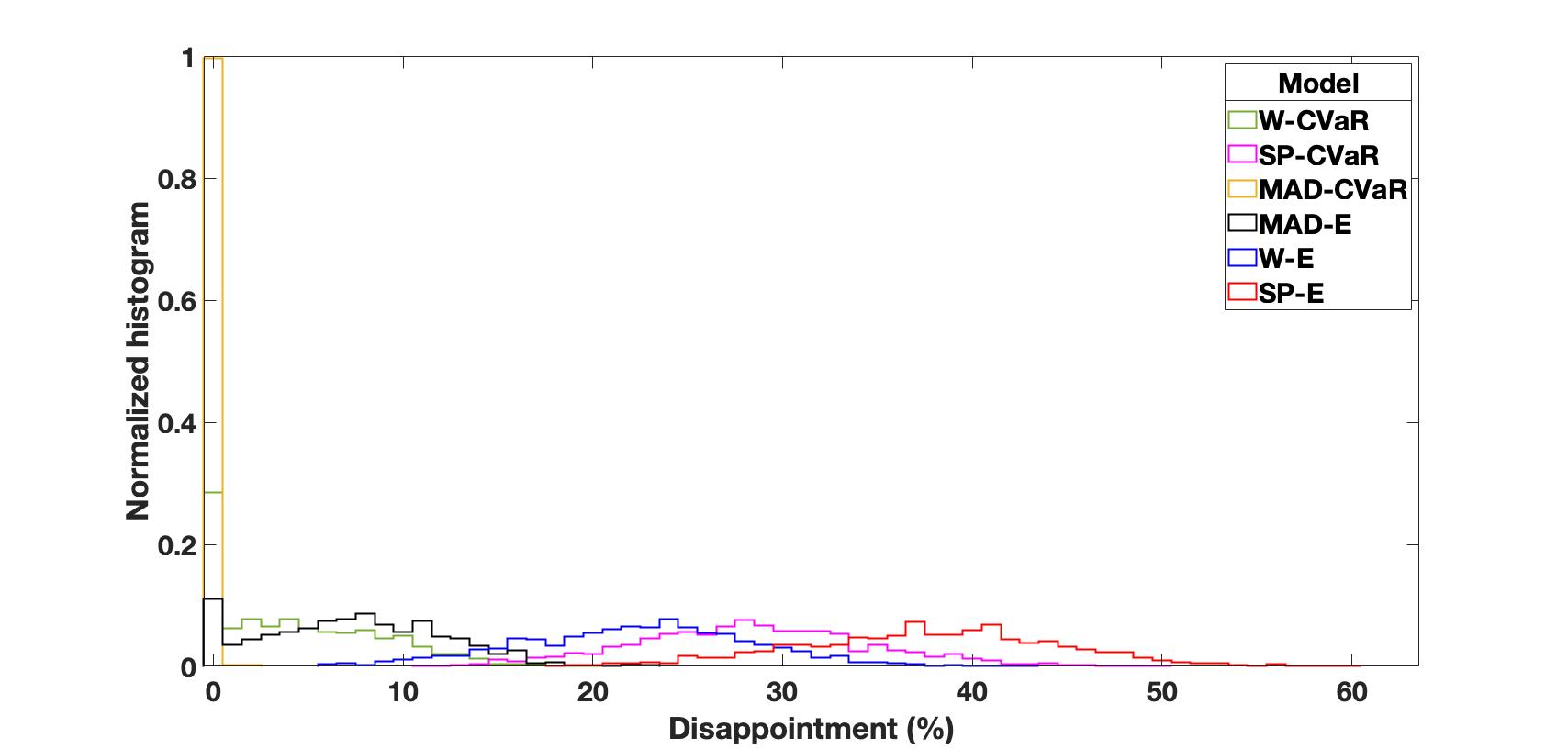}
        \caption{Disappointment, $\Delta=0.5$}
    \end{subfigure}%

\caption{Normalized histograms of second-stage cost (2nd) and out-of-sample disappointments under Set 2 with $\gamma=0.25\gamma_o$, $f=1,500$, and $\pmb{\Delta \in \{0, 0.25, 0.5\}}$.}\label{CVaR_25gamma}
\end{figure}

Next, we compare the out-of-sample operational performance (i.e., second-stage cost) and disappointment of the risk-neutral (i.e., expectation models) and risk-averse models (i.e., mean-CVaR-based models with $\Theta=0$, or equivalently, risk-averse models with CVaR criterion).  We use MAD-E (W-E) to denote the risk-neutral DRO model with MAD (1-Wasserstein) ambiguity  presented in Section~\ref{sec:MAD-DRO} (Section~\ref{sec:WDMFRS_model}). In addition, we use SP-E to denote the risk-neutral SP model.  In figures~\ref{CVaR_gamma}, \ref{CVaR_50gamma}, and \ref{CVaR_25gamma}, we present histograms of the out-of-sample second-stage costs and disappointments under Set 2 with $f=$1,500 and $\gamma=\gamma_o$, $\gamma=0.5\gamma_o$, and $\gamma=0.25\gamma_o$, respectively. We obtained similar observations for the other considered values of $f$.

Let us first compare the performance of the risk-neutral and risk-averse SP models. Notably, the SP-E solutions have the worst performance with significantly higher second-stage costs and larger positive disappointments under all values of $\gamma$ and $\Delta$ than the other considered models. In contrast, the SP-CVaR solutions yield smaller second-stage costs and disappointments than the SP-E model.  This makes sense because the SP-CVaR model schedules a larger numbers of MFs. In addition, when $\gamma=\gamma_o$ and $0.5\gamma_o$, the SP-CVaR model yields the same second-stage costs as the W-E model because both models schedule 6 MFs under these settings (this is why we do not see histograms for the second-stage cost of the SP-CVaR model). However, the W-E model controls the disappointments in a smaller range, while the SP-CVaR model yield larger positive disappointments than the W-E model and the other DRO models. 


Let us now compare the performance of the risk-neutral and risk-averse DRO models.  First, when $\gamma= \gamma_o$ (i.e., the largest unmet demand penalty), the MAD-CVaR, MAD-E, and W-CVaR models  have the same and best performance because they schedule a larger fleet of 8 MFs than the other considered models. Second, when $\gamma=$0.5$\gamma_o$, the MAD-CVaR  model has the lowest second-stage costs and zero disappointments under all values of $\Delta$. This makes sense because the MAD-CVaR model schedules a larger number of MFs than the other considered models when $\gamma=$0.5$\gamma_o$.  Third, when $\gamma=0.25\gamma_o$, the MAD-CVaR and W-CVaR model yield the lowest second-stage costs and disappointments because they schedule a larger fleet (6 MFs) than the other considered models. However, the W-CVaR model yields slightly higher disappointments than the MAD-CVaR when $\Delta=0.5$.

Fourth, the second-stage costs and disappointments of the W-CVaR model are smaller than those of the W-E model because the W-CVaR model schedules a larger number of MFs. Fifth, the W-CVaR and MAD-E models yield the same second-stage costs when $\gamma=0.5 \gamma_o$ because they schedule 7 MFs  (this is why we can only see black histograms for the MAD-E model).  Finally, the MAD-E model have lower second-stage costs and disappointments than the W-E model under all values of $\gamma$ and $\Delta$, which is consistent with our results in  Section~\ref{sec5:OutSample}.



Our results in this section demonstrate that the distributionally robust CVaR models tend to hedge against uncertainty, ambiguity, and risk by scheduling more MFs. Our results also indicate that the proposed DRO expectation models may be risk-averse under some parameter settings (e.g., high unmet demand penalty and low cost). 


\color{black}


\section{\textcolor{black}{Conclusion}}\label{Sec:Conclusion}
\color{black}
\noindent In this paper, we propose two DRO models for the MFRSP. Specifically, given a set of MFs, a planning horizon, and a service region, our models aim to find the number of MFs to use within the planning horizon and a route and schedule for each MF in the fleet. The objective is to minimize the fixed cost of establishing the MF fleet plus a risk measure (expectation or mean-CVaR) of the operational cost over all demand distributions defined by an ambiguity set. In the first model (MAD-DRO), we use an ambiguity set based on the demand's mean, support, and mean absolute deviation. In the second model (W-DRO), we use an ambiguity set that incorporates all distributions within a 1-Wasserstein distance from a reference distribution. To solve the proposed DRO models, we propose a decomposition-based algorithm.  We also derive lower bound inequalities and two families of symmetry breaking  constraints to improve the solvability of the proposed models.

\color{black}


Our computational results demonstrate (1) how the DRO approaches have superior operational performance in terms of satisfying customers demand as compared to the SP approach, (2) the MAD-DRO model is more computationally efficient than the W-DRO model, (3) the MAD-DRO model yield more conservative decisions than the W-DRO model, which often have a higher fixed cost but significantly lower operational cost, (4) how mobile facilities can move from one location to another to accommodate the change in demand over time and location, (5) efficiency of the proposed symmetry breaking constraints and lower bound inequalities, (6) the trade-off between cost, number of MFs, MF capacity, and operational performance, and (7) the trade-off between the risk-neutral and risk-averse approaches. Most importantly, our results show the value of modeling uncertainty and distributional ambiguity.

Note that we have used benchmark instances from the literature in our computational experiments, which may be a limitation of our results. However, in the sensitivity analysis section, we tested the proposed approaches under different parameters settings, demonstrating how decision-makers can use our approaches to generate MFRSP solutions under different parameter settings relevant to their specific application. Moreover, these benchmark instances represent a wide range of potential service regions, which we can efficiently solve.  If we account for the scale of the problem in the sense of static facility location problem, we have demonstrated that we can solve instances of $30 \times 20= 600$ customers (instance 10), which are relatively large for many practical applications.

We suggest the following areas for future research. First, we aim to extend our models to optimize the capacity and size of the MF fleet. Second, we want to extend our approach by incorporating multi-modal probability distributions and more complex relationships between random parameters (e.g., correlation). Third, we aim to extend our approach to more comprehensive MF planning models, which consider all relevant organizational and technical constraints and various sources of uncertainties (e.g., travel time) with a particular focus on real-life healthcare settings. Although conceptually and theoretically advanced, stochastic optimization approaches such as SP and DRO are not intuitive or transparent to decision-makers who often do not have optimization expertise. Thus, future efforts should also focus on closing the gap between theory and practice.

\vspace{2mm}

\ACKNOWLEDGMENT{%
 We want to thank all colleagues who have contributed significantly to the related literature. We are grateful to the anonymous reviewers for their insightful comments and suggestions that allowed us to improve the paper. Special thanks to Mr. Man Yiu Tsang (a Ph.D. student at the Department of Industrial and Systems Engineering, Lehigh University) for helping with Figure~\ref{LehighMap} and proofreading the paper. Dr. Karmel S. Shehadeh dedicates her effort in this paper to every little dreamer in the whole world who has a dream so big and so exciting. Believe in your dreams and do whatever it takes to achieve them--the best is yet to come for you.}
\vspace{1mm}

\newpage
 \begin{APPENDICES}
\section{Derivation of feasible region $\mathcal{X}$ in \eqref{eq:RegionX}} \label{Appx:FirstStageDec}

\noindent In this Appendix, we provide additional details on the derivation of the constraints defining the feasible region $\calX$ of variables ($\xbb,\yb$). As described in \cite{lei2014multicut}, we can enforce the requirement that an MF can only be in service when it is stationary using the following constraints:
\begin{align}\label{eq1}
& x_{j,m}^t+x_{j',m}^t \leq 1, && \forall t, m, j, \ j\neq j^\prime, \ t^\prime \in \lbrace t, \ldots, \min \lbrace t+t_{j,j^\prime}, T\rbrace \rbrace.
\end{align}

\noindent If $x_{j,m}^t=1$ (i.e., MF $m$ is stationary at some location $j$ in period $t$), it can only be available at location $j'\neq j$ after a certain period, depending on the time it takes to travel from location $j$ to location $j'$. It follows by \eqref{eq1} that $ x_{j',m}^{t'}=0$ for all $j\neq j^\prime, \ t^\prime \in \lbrace t, \ldots, \min \lbrace t+t_{j,j^\prime}, T\rbrace \rbrace$. As pointed out by \cite{lei2014multicut}, this indicates that \textit{an earlier decision of deploying an MF at one candidate location would directly affect future decisions both temporally and spatially. In fact, this correlation is a major source of complexity for optimizing the MFRSP.}  

Since the MF has to be in an active condition before providing
service, we have to include the following constraints: 
\begin{align}\label{eq2}
x_{j,m}^t \leq y_m, && \forall j,m,t.
\end{align}

It is straightforward to verify that constraint sets \eqref{eq1} and \eqref{eq2} can be combined into the following compact form:
 \begin{align}
\mathcal{X}&=\left\{ (\xbb, \yb) :  \begin{array}{l} x_{j,m}^t+x_{j^\prime,m}^{t^\prime} \leq y_m, \ \ \forall t, m, j, \ j\neq j^\prime, \ t^\prime \in \lbrace t, \ldots, \min \lbrace t+t_{j,j^\prime}, T\rbrace \rbrace \\
  x_{j,m}^t \in \lbrace 0, 1\rbrace, \ y_m \in \lbrace 0, 1\rbrace,  \ \forall j, m , t \end{array} \right\} .
\end{align} 

\section{Proof of Proposition 1 }\label{Appx:ProofofProp1}

\noindent \textit{Proof}. For a fixed $x \in \mathcal{X}$,  we can explicitly write the inner problem $\sup [ \cdot ]$ in \eqref{MAD-DMFRS} as the following functional linear optimization problem. 
\begin{subequations}
\begin{align}
 \max_{\Prob\geq 0} & \int_{\calS} Q(\xbb,\Wb)  \ d \Prob \\
\ \text{s.t.} & \    \int_{\calS}  W_{i,t} \ d\Prob= \mu_{i,t} \quad \quad \forall i \in I, \ t \in T, \label{ConInner:W}\\
& \int_{\calS} |W_{i,t}-\mu_{i,t}| \ d\Prob\leq \eta_{i,t} \quad \quad \forall i \in I, \ t \in T, \label{ConInner:K}\\
& \  \int_{\calS}   d\Prob= 1. \label{ConInner:Distribution}
\end{align} \label{InnerMax2}
\end{subequations}

Letting  $\rho_{i,t}, \psi_{i,t}$ and $\theta$ be the dual variables associated with constraints \eqref{ConInner:W}, \eqref{ConInner:K}, \eqref{ConInner:Distribution}, respectively, we present problem \eqref{InnerMax2} (problem (9) in the main manuscript) in its dual form:
\begin{subequations}
\begin{align}
& \min_{\pmb{\rho, \theta, \psi \geq 0}} \ \sum \limits_{t \in T} \sum \limits_{i \in I}( \mu_{i,t} \rho_{i,t}+\eta_{i,t} \psi_{i,t})+ \theta \label{DualInner:Obj} \\
& \ \  \text{s.t.} \  \sum \limits_{t \in T} \sum \limits_{i \in I} (W_{i,t} \rho_{i,t}+  |W_{i,t}-\mu_{i,t}| \psi_{i,t})+ \theta \geq Q(\xbb,\Wb)  && \forall \pmb{W}\in \calS, \label{DualInner:PrimalVariabl}
\end{align} \label{DualInnerMax}
\end{subequations}
 where $\pmb{\rho}$ and $\theta$ are unrestricted in sign, $\pmb \psi \geq 0$, and constraint \eqref{DualInner:PrimalVariabl} is associated with the primal variable $\Prob$. Note that for fixed ($\pmb{\rho, \ \psi,} \  \theta$), constraint \eqref{DualInner:PrimalVariabl} is equivalent to
$$\theta \geq \max \limits_{\Wb \in \calS } \Big \lbrace Q(\xbb,\Wb) +  \sum \limits_{t \in T} \sum \limits_{i \in I}  -(W_{i,t} \rho_{i,t}+ |W_{i,t}-\mu_{i,t}| \psi_{i,t}) \Big\rbrace. $$

  Since we are minimizing $\theta$ in \eqref{DualInnerMax}, the dual formulation of \eqref{InnerMax2} is equivalent to:
\begin{align*} 
& \min_{\pmb{\rho}, \pmb{\psi} \geq 0}  \ \left \lbrace  \sum \limits_{t \in T} \sum \limits_{i \in I}( \mu_{i,t} \rho_{i,t}+\eta_{i,t} \psi_{i,t})+\max \limits_{\Wb \in \calS} \Big\lbrace Q(\xbb,\Wb) +   \sum \limits_{t \in T} \sum \limits_{i \in I} -(W_{i,t} \rho_{i,t}+ |W_{i,t}-\mu_{i,t}| \psi_{i,t})\Big \rbrace  \right  \rbrace .
\end{align*}

\section{Proof of Proposition 2}\label{Appx:Prop2}

\noindent First, note that the feasible region $\Omega=\{\eqref{Const1:QudalX}-\eqref{Const3:Qdualv}, \eqref{McCormick1}-\eqref{McCormick2}, \eqref{Const2_finalinner}\}$ and $\calS$ are both independent of $\xbb$, $\pmb \rho$, and $\pmb \psi$ and bounded. In addition, the MFRSP has a complete recourse (i.e., the recourse problem is feasible for any feasible $(\xbb,\yb) \in \mathcal{X}$). Therefore,   $\max\limits_{\pmb{\lambda, v, W, \pi, k} }\Big[ \sum \limits_{t \in T} \sum\limits_{i \in I} \pi_{i,t}+ \sum\limits_{t \in T} \sum\limits_{j \in J} \sum \limits_{m \in M} Cx_{j,m}^t v_{j,m}^t +\sum \limits_{t \in T} \sum \limits_{i \in I} -(W_{i,t} \rho_{i,t}+k_{i,t}\psi_{i,t}) \Big]< +\infty$. Second, for any fixed $\pmb{\pi, v, W, k}$, $\Big [\sum \limits_{t \in T} \sum\limits_{i \in I} \pi_{i,t}+ \sum\limits_{t \in T} \sum\limits_{j \in J} \sum \limits_{m \in M} Cx_{j,m}^t v_{j,m}^t +\sum \limits_{t \in T} \sum \limits_{i \in I} -(W_{i,t} \rho_{i,t}+k_{i,t}\psi_{i,t}) \Big] $ is a linear function of $\xbb$, $\pmb \rho$, and $\pmb \psi$. It follows that $\max\limits_{\pmb{\lambda, v, W, \pi, k} }\Big[ \sum \limits_{t \in T} \sum\limits_{i \in I} \pi_{i,t}+ \sum\limits_{t \in T} \sum\limits_{j \in J} \sum \limits_{m \in M} Cx_{j,m}^t v_{j,m}^t +\sum \limits_{t \in T} \sum \limits_{i \in I} -(W_{i,t} \rho_{i,t}+k_{i,t}\psi_{i,t}) \Big]$ is the maximum of  linear functions of $\xbb$, $\pmb \rho$, and $\pmb \psi$, and hence convex and piecewise linear. Finally, it is easy to see that each linear piece of this function is  associated with one distinct extreme point of $\Omega$ and $\calS$. Given that each of these polyhedra has a finite number of extreme points, the number of pieces of this function is finite. This completes the proof.

\section{Number of Points in The Worst-Case Distribution of MAD-DRO}\label{Appx:Points_MAD}
\color{black}

The results of \cite{longsupermodularity} indicate that if the second-stage optimal value is supermodular in the realization of uncertainties under the MAD ambiguity set, the worst-case is a distribution supported on  ($2n+1$) points, where $n$ is the dimension of the random vector. In this Appendix, we show that even if our recourse is supermodular in demand realization, the number of points in the worst-case distribution of the demand is large, which renders our two-stage MAD-DRO model computationally challenging to solve using \cite{longsupermodularity}' approach

 Recall that the demand is indexed by $i$ and $t$, i.e., $W_{i,t}$, for all $i \in I$ and $t \in T$. Thus, the dimension of our random vector is $|I | |T|$. Accordingly, assuming that the second-stage optimal value is supermodular, then the resullts of \cite{longsupermodularity} suggest that the worst-case distribution in MAD-DRO has ($2|I| |T|$+1) points or scenarios. Note that $|I|$ and $|T|$ and thus $2|I||T|$ is large for most instances of our problem (See Example 1--2 below). Since our computational results and prior literature indicate that solving the scenario-based model using a small set of scenarios is challenging, solving a reformulation of MAD-DRO using the ($2|I||T|$+1) points is expected to be computationally challenging.  
 

%
%
%

\noindent \textbf{Example 1.} Instance 6 ($I=20$ and $T=20$). The worst-case distribution has 801 points.

\noindent \textbf{Example 3.} Instance 10 ($I=30$ and $T=20$). The worst-case distribution has 1201 points.

\color{black}

\newpage

\section{Proof of Proposition~\ref{Prop1}}\label{Proof_Prop1}

Recall that $\hat{\Prob}^N=\frac{1}{N}\sum_{n=1}^n \delta_{\hat{\Wb}^n}$. The definition  of Wasserstein distance indicates that there exist a joint distribution $\Pi$  of $(\Wb, \hat{\Wb}$) such that $\E_{\Pi} [||\Wb-\hat{\Wb}||] \leq \epsilon$. In other words, for any $\Prob\in\calP(\calS)$, we can rewrite any joint distribution $\Pi\in\calP(\Prob,\hat{\Prob}^N)$ by the conditional distribution of $\Wb$ given $\hat{\Wb}=\hat{\Wb}^n$ for $n=1,\dots,N$, denoted as $\F^n$. That is, $\Pi=\frac{1}{N}\sum_{n=1}^N \F^n \times \delta_{\hat{\xi}_n}$. Notice that if we find one joint distribution $\Pi\in\calP(\Prob,\hat{\Prob}^N)$ such that $\int ||\Wb-\hat{\Wb}||  \ d\Pi \leq \epsilon$, then $\text{dist}(\Prob,\hat{\Prob}^N)\leq\epsilon$. Hence, we can drop the infimum operator in Wasserstein distance and arrive at the following equivalent problem.
\begin{subequations}\label{AppxInner}
\begin{align}
& \sup_{\F^n\in\calP(\calS), n\in[N]} \frac{1}{N}\sum_{n=1}^N\int_\calS Q(\xbb,\Wb) \ d\F^n\\
& \text{s.t.}  \ \ \ \ \ \ \ \  \frac{1}{N} \sum_{n=1}^N  \int_\calS ||\Wb-\hat{\Wb}^n||  \ d\F^n \leq \epsilon.
\end{align}
\end{subequations}
Using a standard strong duality argument and letting $\rho\geq 0$ be the dual multiplier, we can reformulate problem \eqref{AppxInner} by its dual, i.e.
\begin{align}
&\quad\,\inf_{\rho\geq 0} \sup_{\F^n\in\calP(\calS),n\in[N]}\left\{\frac{1}{N}\sum_{n=1}^N\int_\calS Q(\xbb,\Wb) d\F^n+\rho\left[\epsilon-\frac{1}{N} \sum_{n=1}^N  \int_\calS ||\Wb-\hat{\Wb}^n||\ d\F^n \right] \right\}\nonumber\\
&=\inf_{\rho\geq 0}  \Big\{ \epsilon \rho+ \frac{1}{N}\sum_{n=1}^N \sup_{\F^n\in\calP(\calS)}\int_\calS \left[Q(\xbb,\Wb)-\rho||\Wb-\hat{\Wb}^n||\right] d\F^n\Big\} \nonumber\\
&=\inf_{\rho\geq 0} \Big\{ \epsilon \rho + \frac{1}{N} \sum_{n=1}^N \sup_{\Wb \in \calS}  \{ Q(\xbb,\Wb)-\rho || \Wb -\hat{\Wb}^n ||  \} \Big\}.
\end{align}

\newpage

\section{DRO with mean-CVaR as a Risk Measure}\label{Sec:meanCVAR}

\noindent Both the MAD-DRO  and W-DRO model presented in Section \ref{sec:DRO_Models} assume that the decision-maker is risk-neutral (i.e., adopt the expected  value of the recourse as a risk measure). In some applications of the MFRSP, however, decision-makers might be risk-averse. Therefore, as one of our reviewers suggested, in this section, we present a distributionally robust risk-averse model for the MFRSP. 

To model the decision maker's risk aversion, most studies adopt the CVaR, i.e., set $\varrho(\cdot)=\CVAR_\kappa(\cdot)$, where $\kappa\in(0,1)$.  CVaR is the conditional expectation of $Q(\cdot)$ above the value-at-risk VaR (informally, VaR is the $\kappa$ quantile of the distribution of $Q(\cdot)$, see \cite{pacc2014robust, rockafellar2002conditional, sarykalin2008value, van2015distributionally}). CVaR  is a popular coherent risk measure widely used to avoid solutions influenced by a bad scenario with a low probability. However, as pointed out by \cite{wang2021two} and a reviwer of this paper, neither expected value nor CVaR can capture the variability of uncertainty in a comprehensive manner. Alternatively, we consider minimizing the mean-CVaR, which balances the cost on average and avoids high-risk levels. As pointed out by \cite{lim2011conditional}, \cite{wang2021two}, and our reviewer, the traditional CVaR criterion is sensitive to the misspecification of the underlying loss distribution and lacks robustness.  Therefore, we propose a distributionally robust mean-CVaR model to remedy such fragility, reflecting both risk-averse and ambiguity-averse attitudes. For brevity, we use the MAD ambiguity set to formulate and analyze this model because similar formulation and reformulation steps can be used to derive a solvable mean-CVaR-based model based on the 1-Wasserstein ambiguity.

Let us now introduce our distributionally robust mean-CVaR-based model with MAD-ambiguity (MAD-CVaR).  First, following \cite{rockafellar2000optimization, rockafellar2002conditional}, and \cite{ van2015distributionally}, we formally define CVaR as
\begin{align}
\CVAR_\kappa (Q(\xbb,\Wb))= \inf \limits_{\zeta \in \R} \Big \{ \zeta+ \frac{1}{1-\kappa}  \E_\mathbb{P} [Q(\xbb,\Wb)-\zeta ]^+ \Big\}, \label{CVAR_Def_main}
\end{align}

\noindent where $[c]^+:=\max \{c,0\}$ for $c \in \R$.  Parameter $\kappa$ measures a wide range of risk preferences, where $\kappa=0$ corresponds to the risk-neutral formulation. In contrast, when $\kappa \rightarrow 1$, the decision-makers become more risk-averse. Using \eqref{CVAR_Def_main}, we formulate the following MAD-CVaR model (see, e.g., \cite{wang2021two} for a recent application in facility location):
\begin{align}
 \min \limits_{(\xbb, \yb)\in \mathcal{X}}  & \Bigg \{ \sum_{m \in M} f y_{m} -\sum_{t \in T} \sum_{j \in J} \sum_{m \in M } \alpha x_{j,m}^t + \Bigg[ \Theta \sup \limits_{\Prob \in \calF(\calS, \mub, \etab)} \E_\Prob [Q(\xbb,\Wb)] + (1-\Theta)  \sup \limits_{\Prob \in \calF(\calS, \mub, \etab)} \CVAR (Q(\xbb,\Wb))   \Bigg]  \Bigg \}, \label{CVaR-DMFRS}
\end{align}

\noindent where  $0 \leq \Theta \leq 1$  is the risk-aversion coefficient, which represents a trade-off between the risk-neutral (i.e., $\E[\cdot]$) and risk-averse (i.e., $\CVAR(\cdot)$) objectives. A larger $\Theta$ implies less aversion to risk, and vice verse. In extreme cases, when $\Theta =1$, the decision maker is risk-neutral and \eqref{CVaR-DMFRS} reduces to the MAD-DRO expectation model in \eqref{MAD-DMFRS}. When $\Theta =0$,  the decision maker is risk and ambiguity averse.

Next, we derive a solvable reformulation of \eqref{CVaR-DMFRS}. Let us first consider the inner maximization problem $\sup \limits_{\Prob \in \calF(\calS, \mub, \etab)} \CVAR (Q(\xbb,\Wb))  $ in \eqref{CVaR-DMFRS}. It is easy to verify that
\begin{subequations}\label{SupCVaR}
\begin{align}
 \sup \limits_{\Prob \in \calF(\calS, \mub, \etab)} \CVAR (Q(\xbb,\Wb))&=  \sup \limits_{\Prob \in \calF(\calS, \mub, \etab)}  \inf \limits_{\zeta \in \R} \Bigg \{  \zeta+ \frac{1}{1-\kappa}  \E_\mathbb{P} [Q(\xbb,\Wb)-\zeta ]^+ \Bigg\} \label{CVAR_alter1} \\
 &= \inf \limits_{\zeta \in \R}  \sup \limits_{\Prob \in \calF(\calS, \mub, \etab)} \Bigg \{ \zeta+ \frac{1}{1-\kappa} \E_\mathbb{P} [Q(\xbb,\Wb)-\zeta ]^+  \Bigg\}  \label{CVAR_alter2} \\
 &=  \inf \limits_{\zeta \in \R} \Bigg \{ \zeta+ \frac{1}{1-\kappa}  \sup \limits_{\Prob \in \calF(\calS, \mub, \etab)}   \E_\mathbb{P} [Q(\xbb,\Wb)-\zeta ]^+ \Bigg\}. \label{CVAR_alter3} \
\end{align}
\end{subequations}

\noindent Interchanging the order of $\sup \limits_{\Prob \in \calF(\calS, \mub, \etab)}$ and $\inf \limits_{\zeta \in \R} $ follows from Sion's minimax theorem \citep{sion1958general}  because $\zeta+ \frac{1}{1-\kappa}  \E_\mathbb{P} [Q(\xbb,\Wb)-\zeta ]^+$ is convex in $\zeta$ and concave in $\Prob$. Next, we apply the same techniques in Section~\ref{sec3:reform} and Appendix~\ref{Appx:ProofofProp1} to reformulate the inner maximization problem  $\sup \limits_{\Prob \in \calF(\calS, \mub, \etab)}   \E_\mathbb{P} [Q(\xbb,\Wb)-\zeta ]^+$ in \eqref{CVAR_alter3} as a minimization problem and combine it with the outer minimization problem to obtain
\begin{subequations}
\begin{align}
 \inf \limits _{\zeta, \pmb{\rho}, \theta, \pmb{\psi} \geq 0} &   \ \    \Bigg \{\zeta+ \frac{1}{1-\kappa} \Bigg [\sum \limits_{t \in T} \sum \limits_{i \in I}( \mu_{i,t} \rho_{i,t}+\eta_{i,t} \psi_{i,t})+ \theta \Bigg] \Bigg\} \label{CVAR_Obj} \\
   \text{s.t.} &   \ \  \sum \limits_{t \in T} \sum \limits_{i \in I} (W_{i,t} \rho_{i,t}+  |W_{i,t}-\mu_{i,t}| \psi_{i,t})+ \theta \geq Q(\xbb,\Wb)-\zeta,  && \forall \pmb{W}\in \calS , \label{CVAR_C1} \\
&  \ \    \sum \limits_{t \in T} \sum \limits_{i \in I} (W_{i,t} \rho_{i,t}+  |W_{i,t}-\mu_{i,t}| \psi_{i,t})+ \theta \geq 0, && \forall \pmb{W}\in \calS,  \label{CVAR_C2}
\end{align}
\end{subequations}
 \noindent where \eqref{CVAR_C1} and \eqref{CVAR_C2} follows from the definition of $[\cdot]^+$. Accordingly, problem \eqref{SupCVaR} (equivalently last term in formulation in \eqref{CVaR-DMFRS}) is equivalent to
 \begin{subequations}\label{CVARmodel_objective}
 \begin{align}
 \inf \limits_{\zeta, \pmb{\rho}, \theta, \pmb{\psi} \geq 0} &   \ \    \Bigg \{\zeta+ \frac{1}{1-\kappa} \Bigg [\sum \limits_{t \in T} \sum \limits_{i \in I}( \mu_{i,t} \rho_{i,t}+\eta_{i,t} \psi_{i,t})+ \theta \Bigg] \Bigg\} \label{CVAR_Obj2} \\
   \text{s.t.} &  \ \  \zeta   \geq \max_{\Wb \in \calS} \Big \{ Q(\xbb,\Wb)-  \sum \limits_{t \in T} \sum \limits_{i \in I} (W_{i,t} \rho_{i,t}+  |W_{i,t}-\mu_{i,t}| \psi_{i,t}) \Big\}-\theta, \label{CVAR_equiv2}\\
& \ \   \min_{\Wb \in \calS} \Big \{ \sum \limits_{t \in T} \sum \limits_{i \in I} (W_{i,t} \rho_{i,t}+  |W_{i,t}-\mu_{i,t}| \psi_{i,t}) \Big \} + \theta \geq 0. \label{CVAR_equiv3}
 \end{align}
 \end{subequations}

 \textcolor{black}{Since we are minimizing $\zeta$ in \eqref{CVAR_Obj2}, we can equivalently re-write  \eqref{CVAR_Obj2} as}
 \color{black}
 \begin{align}
 \inf \limits_{\pmb{\rho}, \theta, \pmb{\psi} \geq 0} &  \    \Bigg \{ \max_{\Wb \in \calS} \Big \{ Q(\xbb,\Wb)-  \sum \limits_{t \in T} \sum \limits_{i \in I} (W_{i,t} \rho_{i,t}+  |W_{i,t}-\mu_{i,t}| \psi_{i,t}) \Big\}+ \big (\frac{\kappa}{1-\kappa} \big) \theta + \frac{1}{1-\kappa} \sum \limits_{t \in T} \sum \limits_{i \in I}( \mu_{i,t} \rho_{i,t}+\eta_{i,t} \psi_{i,t}) \Bigg\}  \nonumber \\
 \equiv \inf \limits_{\pmb{\rho}, \theta, \pmb{\psi} \geq 0}  &  \ \    \Bigg \{ h(\xbb, \pmb \rho, \pmb \psi)   + \big (\frac{\kappa}{1-\kappa} \big) \theta + \frac{1}{1-\kappa} \sum \limits_{t \in T} \sum \limits_{i \in I}( \mu_{i,t} \rho_{i,t}+\eta_{i,t} \psi_{i,t}) \Bigg\},  \label{CVAR_Obj3} 
 \end{align}
\noindent where \eqref{CVAR_Obj3} follows from the derivation of $h(\xbb, \pmb \rho, \pmb \psi) \equiv \max\limits_{\Wb \in \calS} \big \{ Q(\xbb,\Wb)-  \sum \limits_{t \in T} \sum \limits_{i \in I} (W_{i,t} \rho_{i,t}+  |W_{i,t}-\mu_{i,t}| \psi_{i,t}) \big \}$ in Section~\ref{sec3:reform}. Accordingly, problem \eqref{CVARmodel_objective} is equivalent to
 \begin{subequations}\label{CVARmodel_objective4}
 \begin{align}
\inf \limits_{\pmb{\rho}, \theta, \pmb{\psi} \geq 0}  &  \ \    \Bigg \{ h(\xbb, \pmb \rho, \pmb \psi)   + \big (\frac{\kappa}{1-\kappa} \big) \theta + \frac{1}{1-\kappa} \sum \limits_{t \in T} \sum \limits_{i \in I}( \mu_{i,t} \rho_{i,t}+\eta_{i,t} \psi_{i,t}) \Bigg\}  \label{CVAR_Obj4} \\
  \text{s.t.} &  \ \  \min_{\Wb \in \calS} \Big \{ \sum \limits_{t \in T} \sum \limits_{i \in I} (W_{i,t} \rho_{i,t}+  |W_{i,t}-\mu_{i,t}| \psi_{i,t}) \Big \} + \theta \geq 0. \label{CVAR_equiv4}
 \end{align}
 \end{subequations}

 Next, derive equivalent linear constraints of the embedded minimization problem in constraint \eqref{CVAR_equiv4}. For fixed $\pmb \rho$ and $\pmb \psi$, we re-write constraint \eqref{CVAR_equiv4} as
\begin{subequations}\label{CVAR_equiv2_reform}
\begin{align}
 \theta + & \min_{\Wb} \ \   \Big \{ \sum \limits_{t \in T} \sum \limits_{i \in I} (W_{i,t} \rho_{i,t}+  k_{i,t} \psi_{i,t}) \Big \}  \geq 0, \label{CVAR_equiv2_reform_C1}\\
 & \text{s.t.} \ \   \  W_{i,t}  \leq \WU_{i,t}, && \forall i, t,  \label{CVAR_equiv2_reform_C2}\\
& \ \ \ \ \ \  \   W_{i,t} \geq  \WL_{i,t},  && \forall i, t,  \label{CVAR_equiv2_reform_C22}\\
& \ \ \ \ \ \  \ k_{i,t} \geq W_{i,t}-\mu_{i,t},  &&   \forall i, t,  \label{CVAR_equiv2_reform_C3}\\
& \ \ \ \ \ \ \   k_{i,t} \geq \mu_{i,t}-W_{i,t},  &&  \ \forall i, t.\label{CVAR_equiv2_reform_C4}
\end{align}
\end{subequations}

\noindent Letting $a_{i,t}$, $b_{i,t}$, $g_{i,t}$, and $o_{i,t}$ be the dual variables associated with constraints \eqref{CVAR_equiv2_reform_C2}, \eqref{CVAR_equiv2_reform_C22}, \eqref{CVAR_equiv2_reform_C3}, and \eqref{CVAR_equiv2_reform_C4}, respectively, we present the linear program in \eqref{CVAR_equiv2_reform_C1}-\eqref{CVAR_equiv2_reform_C4} in its dual form as
 \begin{subequations}
 \begin{align}
& \theta + \sum_{i \in I} \sum_{t \in T} \Big [ \WU_{i,t}a_{i,t}+  \WL_{i,t} b_{i,t} -\mu_{i,t}g_{i,t}+\mu_{i,t}o_{i,t}   \Big]   \geq 0, \label{dual_equiv2_C1} \\
 &a_{i,t}+b_{i,t}-g_{i,t}+o_{i,t} \leq \rho_{i,t}, && \forall i, t , \label{dual_equiv2_C2} \\
&g_{i,t}+o_{i,t}  \leq \psi_{i,t}, && \forall i, t, \label{dual_equiv2_C3} \\
&  (b_{i,t}, g_{i,t}, \  o_{i,t}) \geq 0,  \ a_{i,t}  \leq 0, && \forall i, t. \label{dual_equiv2_C4}
 \end{align}
 \end{subequations}

\noindent Replacing \eqref{CVAR_equiv4} in \eqref{CVARmodel_objective4}  with \eqref{dual_equiv2_C1}--\eqref{dual_equiv2_C4}, we derive the following equivalent reformulation of problem \eqref{CVARmodel_objective4} (equivalently, problem $\sup \limits_{\Prob \in \calF(\calS, \mub, \etab)} \CVAR (Q(\xbb,\Wb))$ in \eqref{CVaR-DMFRS}):
 \begin{subequations}\label{CVARmodel_obj5}
 \begin{align}
 \inf \limits_{\pmb{\rho}, \theta, \pmb{\psi} \geq 0}  &  \ \    \Bigg \{ h(\xbb, \pmb \rho, \pmb \psi)   + \big (\frac{\kappa}{1-\kappa} \big) \theta + \frac{1}{1-\kappa} \sum \limits_{t \in T} \sum \limits_{i \in I}( \mu_{i,t} \rho_{i,t}+\eta_{i,t} \psi_{i,t}) \Bigg\} \label{CVAR_Obj5} \\
   \text{s.t.} &  \ \  \eqref{dual_equiv2_C1}-\eqref{dual_equiv2_C3}.
 \end{align}
 \end{subequations}

\noindent Combining the inner problem  $\sup \limits_{\Prob \in \calF(\calS, \mub, \etab)} \CVAR (Q(\xbb,\Wb))$ in the form of \eqref{CVARmodel_obj5} and problem  $\sup \limits_{\Prob \in \calF(\calS, \mub, \etab)} \E_\Prob [Q(\xbb,\Wb)]$  in the form of  \eqref{FinalDR} with the outer minimization problem in \eqref{CVaR-DMFRS}, we derive the following equivalent reformulation of the  MAD-CVaR  model in \eqref{CVaR-DMFRS}. 
\begin{subequations}\label{Final_CVaR}
\begin{align}
 \min  & \Bigg \{ \sum_{m \in M} f y_{m} -\sum_{t \in T} \sum_{j \in J} \sum_{m \in M } \alpha x_{j,m}^t+ \Theta \Bigg [ \sum \limits_{t \in T} \sum \limits_{i \in I} \Big(\mu_{i,t} \rho_{i,t}+\eta_{i,t} \psi_{i,t}\Big) + h(\xbb, \pmb \rho, \pmb \psi)   \Bigg] \nonumber \\
 & \ \ +  (1-\Theta)    \Bigg [ h(\xbb, \pmb \rho, \pmb \psi)   + \big (\frac{\kappa}{1-\kappa} \big) \theta + \frac{1}{1-\kappa} \sum \limits_{t \in T} \sum \limits_{i \in I}( \mu_{i,t} \rho_{i,t}+\eta_{i,t} \psi_{i,t})  \Bigg] \Bigg \} \label{CVaR-DMFRS3} \\
 \text{s.t.} & \ \  \ (\xbb, \yb)\in \mathcal{X},   \ \pmb \psi \geq 0, \ \eqref{dual_equiv2_C1}-\eqref{dual_equiv2_C3}. 
 \end{align}
\end{subequations}
  
\noindent It is easy to verify that problem \eqref{Final_CVaR} is equivalent to 

\begin{subequations}\label{Final_CVaR2}
\begin{align}
 \min  & \Bigg \{ \sum_{m \in M} f y_{m} -\sum_{t \in T} \sum_{j \in J} \sum_{m \in M } \alpha x_{j,m}^t+ \Theta \Bigg [ \sum \limits_{t \in T} \sum \limits_{i \in I} \Big(\mu_{i,t} \rho_{i,t}+\eta_{i,t} \psi_{i,t}\Big) \Bigg] \nonumber \\
 & \ \ + \delta + (1-\Theta)    \Bigg [  \big (\frac{\kappa}{1-\kappa} \big) \theta + \frac{1}{1-\kappa} \sum \limits_{t \in T} \sum \limits_{i \in I}( \mu_{i,t} \rho_{i,t}+\eta_{i,t} \psi_{i,t})  \Bigg] \Bigg \} \label{CVaR-DMFRS4_1} \\
 \text{s.t.} & \ \  \ (\xbb, \yb)\in \mathcal{X},   \ \pmb \psi \geq 0, \ \eqref{dual_equiv2_C1}-\eqref{dual_equiv2_C3}, \\
 &  \ \  \  \delta \geq h(\xbb, \pmb \rho, \pmb \psi), \label{CVaR-DMFRS4_2}
\end{align}
\end{subequations}
  
where  $\textcolor{black}{h(\xbb, \pmb \rho, \pmb \psi)}=\max \limits_{\pmb{\lambda, v, W, \pi, k} } \ \Big\{ \sum \limits_{t \in T} \sum\limits_{i \in I} \pi_{i,t}+ \sum\limits_{t \in T} \sum\limits_{j \in J} \sum \limits_{m \in M} Cx_{j,m}^t v_{j,m}^t +\sum \limits_{t \in T} \sum \limits_{i \in I} -(W_{i,t} \rho_{i,t}+k_{i,t}\psi_{i,t}):  \eqref{Const1:QudalX}-\eqref{Const3:Qdualv}, \eqref{McCormick1}-\eqref{McCormick2}, \eqref{Const2_finalinner} \Big\} $ from Section~\ref{sec3:reform}.  Finally, we observe that the right-hand side (RHS) of constraints \eqref{CVaR-DMFRS4_2} is equivalent to the RHS of constraints \eqref{const1:FinalDR} in the equivalent reformulation of the risk-neutral MAD-DRO model in \eqref{FinalDR}. Therefore, we can easily adapt Algorithm 1 to solve \eqref{Final_CVaR2} (see Algorithm~\ref{Alg2:Decomp2}).

\begin{algorithm}[t!]
\color{black}
\small
   \renewcommand{\arraystretch}{0.3}
\caption{Decomposition algorithm for the MAD-CVaR Model.}
\label{Alg2:Decomp2}
\noindent \textbf{1. Input.}  Feasible region $\mathcal{X}$; support $\calS$; set of cuts $ \lbrace L (\xbb, \delta)\geq 0 \rbrace=\emptyset $; $LB=-\infty$ and $UB=\infty.$

\vspace{2mm}

\noindent \textbf{2. Master Problem.} Solve the following master problem
\begin{subequations}\label{MastCVAR}
\begin{align}
\min \quad &  \Bigg\{\sum_{m \in M} f y_{m}-\sum_{t \in T} \sum_{j \in J} \sum_{m \in M} \alpha x_{j, m}^{t}+\Theta\left[\sum_{t \in T} \sum_{i \in I}\left(\mu_{i, t} \rho_{i, t}+\eta_{i, t} \psi_{i, t}\right)\right]\nonumber \\
& \left.+(1-\Theta)\left[\left(\frac{\kappa}{1-\kappa}\right) \theta+\frac{1}{1-\kappa} \sum_{t \in T} \sum_{i \in I}\left(\mu_{i, t} \rho_{i, t}+\eta_{i, t} \psi_{i, t}\right)\right] +\delta  \right\} \\
\text{s.t.} &  \qquad  (\xbb,\yb ) \in \mathcal{X},  \ \ \pmb{\psi} \geq 0, \ \eqref{dual_equiv2_C1}-\eqref{dual_equiv2_C3}, \ \{ L (\xbb, \delta)\geq 0 \},
\end{align}%
\end{subequations}
$\ \ \ $  and record an optimal solution $(\xbb^	*, \pmb{\rho}^*, \pmb{\psi}^*,  \delta^*)$ and optimal value $Z^*$. Set $LB=Z^*$.

\noindent \textbf{3. Sub-problem.} 
\begin{itemize}
\item[3.1.] With $(\xbb, \pmb{\rho}, \pmb{\psi})$ fixed to $(\xbb^*, \pmb{\rho}^*, \pmb{\psi}^*)$, solve the following problem 
\begin{subequations}\label{MILPSep2}
\begin{align}
h(\xbb, \pmb \rho, \pmb \psi)=  &\max \limits_{\pmb{\lambda, v, W, \pi, k} } \ \Big\{ \sum \limits_{t \in T} \sum\limits_{i \in I} \pi_{i,t}+ \sum\limits_{t \in T} \sum\limits_{j \in J} \sum \limits_{m \in M} Cx_{j,m}^t v_{j,m}^t +\sum \limits_{t \in T} \sum \limits_{i \in I} -(W_{i,t} \rho_{i,t}+k_{i,t}\psi_{i,t})  \Big\}\\
& \ \ \ \text{s.t. } \  \eqref{Const1:QudalX}-\eqref{Const3:Qdualv}, \eqref{McCormick1}-\eqref{McCormick2}, \eqref{Const2_finalinner},
\end{align}
\end{subequations}

$\qquad \ $ and record optimal solution $(\pmb{\pi^*, \lambda^*, W^*, v^*, k^*})$ and $h(\xbb, \pmb \rho, \pmb \psi)^*$. 
\item[3.2.] Set $UB=\min \{ UB, \ h(\xbb, \pmb \rho, \pmb \psi)^*+ (LB-\delta^*) \}$.
\end{itemize}

\noindent \textbf{4. if} $\delta^* \geq  \sum \limits_{t \in T} \sum\limits_{i \in I} \pi_{i,t}^*+ \sum\limits_{t \in T} \sum\limits_{j \in J} \sum \limits_{m \in M} Cx_{j,m}^{t*} v_{j,m}^{t*} +\sum \limits_{t \in T} \sum \limits_{i \in I} -(W_{i,t}^* \rho_{i,t}^*+k_{i,t}^*\psi_{i,t}^*) $ \textbf{then}

$\qquad \ \ $ stop and return $\xbb^*$ and $\yb^*$ as the optimal solution to problem \eqref{MastCVAR} (equivalently, \eqref{Final_CVaR2}).

\noindent $\ \ $ \textbf{else} add the cut $\delta \geq   \sum \limits_{t \in T} \sum \limits_{i \in I} \pi_{i,t}^*+ \sum\limits_{t \in T} \sum \limits_{j \in J} \sum\limits_{m \in M} Cx_{j,m} v_{j,m}^{t*}+\sum \limits_{t \in T} \sum \limits_{i \in I} -(W_{i,t}^* \rho_{i,t}+k_{i,t}^*\psi_{i,t} )$  to the set of

$\qquad \quad$  cuts   $ \{L (\xbb, \delta) \geq 0 \}$ and go to step 2.

\noindent $\ \ $ \textbf{end if}

\end{algorithm}

\color{black}
\newpage
\section{Two-stage SP Model with Mean-CVaR Objective}\label{Appex:SP_CVAR}
\noindent The mean-CVaR-based SP model (denoted as SP-CVaR) is as follow:

\begin{align}
 \min \limits_{(\xbb, \yb)\in \mathcal{X}}  & \Bigg \{ \sum_{m \in M} f y_{m} -\sum_{t \in T} \sum_{j \in J} \sum_{m \in M } \alpha x_{j,m}^t + \Theta \E_\mathbb{P} [Q(\xbb,\Wb))] + (1-\Theta)   \CVAR (Q(\xbb,\Wb))     \Bigg \} \nonumber\\
  \equiv  &\nonumber \\
\min \limits_{(\xbb, \yb)\in \mathcal{X}, \zeta \in \R^+}  &\Big \{ \sum_{m \in M} f y_{m} -\sum_{t \in T} \sum_{j \in J} \sum_{m \in M } \alpha x_{j,m}^t + \Theta \E_\mathbb{P}  [Q(\xbb,\Wb))]+ (1-\Theta) \Big ( \zeta+ \frac{1}{1-\kappa}   \E_\mathbb{P} [Q(\xbb,\Wb)-\zeta ]^+ \Big) \Bigg\}. \label{CVaR-SP}
\end{align}
\noindent The sample average deterministic equivalent of \eqref{CVaR-SP} based on $N$ scenarios, $\Wb^1,\ldots, \Wb^N$, is as follows: 
\begin{subequations}
\begin{align}
\min \limits_{(\xbb, \yb)\in \mathcal{X}, \zeta \in \R^+}   &  \Bigg \{ \sum_{m \in M} f y_{m} -\sum_{t \in T} \sum_{j \in J} \sum_{m \in M } \alpha x_{j,m}^t + \Theta  \sum_{n=1}^N \frac{1}{N} Q(\xbb,\Wb^n) + (1-\Theta) \Big ( \zeta+  \frac{1}{1-\kappa}   \sum_{n=1}^N \frac{1}{N} \Phi^n \Big) \Bigg\} \\
\text{s.t. } & \ \Phi^n \geq  Q(\xbb,\Wb^n) -\zeta, \qquad \forall n \in [N],\\
 & \  \Phi^n \geq 0, \qquad \forall n \in [N],
\end{align}
\end{subequations}

\noindent where for each $n \in [N]$, $Q(\xbb, \Wb^n)$ is the recourse problem defined in \eqref{2ndstage}, and $\Theta \in [0, 1]$/

\section{Proof of Proposition~\ref{Prop3:LowerB1}}\label{Appx:Prop3}

\noindent Recall from the definition of the support set that the lowest demand of each customer $i$ in period $t$ equals to the integer parameter $\underline{W}_{i,t}$. Now, if we treat the MFs as uncapacitated facilities, then we can fully satisfy $\underline{W}_{i,t}$ at the lowest assignment cost from the nearest location $j \in J^\prime$, where $J^\prime:= \lbrace j: x_{j,m}^t=1 \rbrace$.  Note that $J^\prime \subseteq J$. Thus, the lowest assignment cost  must be at least equal to or larger than $\sum_{i \in I} \min \limits_{j \in J} \lbrace  \beta d_{i,j} \rbrace \underline{W}_{i,t} $. If  $\gamma < \min \limits_{j \in J} \lbrace  \beta d_{i,j} \rbrace$,  then the recourse must be at least equal to or larger than $ \sum_{i \in I} \gamma  \underline{W}_{i,t}$. Accordingly, $\sum \limits_{i \in I} \min \lbrace \gamma,  \min \limits_{j \in J} \lbrace \beta d_{i,j} \rbrace \rbrace \underline{W}_{i,t}$ is a valid lower bound on recourse $Q(\xbb, \Wb)$ for each $t \in T$. 

\section{Sample Average Approximation}\label{Appx:SAA}
\begin{subequations}\label{SPModel}
\begin{align}
& \min  \Big [  \sum_{m \in M} f y_{m} -\sum_{t \in T} \sum_{j \in J} \sum_{m \in M } \alpha x_{j,m}^t+ \frac{1}{N} \sum_{n=1}^N\Big( \sum_{j \in J} \sum_{i \in I} \sum_{m \in M} \sum_{t \in T} \beta d_{i,j} z_{i,j,m}^{t,n}+ \gamma \sum_{t \in T} \sum_{i \in I} u_{i,t}^n  \Big)  \Big]\label{ObjSP}\\
& \ \text{s.t.} \ \ (\xbb, \yb) \in \mathcal{X}, \\
&\quad  \   \quad  \ \ \sum_{j \in J} \sum_{m \in M} z_{i,j,m}^{t,n}+u_{i,t}^n=  W_{i,t}^n, \qquad  \forall i \in I, \ t \in T , \ n \in [N], \\
&\quad  \   \quad  \quad  \sum_{i \in I} z_{i,j,m}^{t,n}\leq C x_{j,m}^t \qquad  \forall j \in J, \ m \in M, \ t \in T, \ n \in [N], \\
&\quad  \   \quad  \quad u_{i,t}^n \geq 0, \ z_{i,j,m}^{t,n} \geq 0, \qquad   \forall i \in I, \ j \in J, m \in M, \ t \in [T, \  n \in [N].
\end{align}
\end{subequations}

%
%
%
%
%

\newpage

\section{\textcolor{black}{Details of Lehigh County Instances}}\label{AppexLehigh}

\begin{table}[h!]
 \center 
 \footnotesize
   \renewcommand{\arraystretch}{0.7}
  \caption{A subset of 20 nodes in Lehigh County and their population based on the 2010 census of Lehigh County (from \cite{wikiLehighCounty}).} 
\begin{tabular}{lll}
\hline
\textbf{City/Town/etc.} &	Population	\\
\hline 
Allentown	&	118,032	\\
Bethlehem	&	74,982	\\
Emmaus	&	11,211	\\
Ancient Oaks	&	6,661	\\
Catasauqua	&	6,436	\\
Wescosville	&	5,872	\\
Fountain Hill	&	4,597	\\
Dorneyville	&	4,406	\\
Slatington	&	4,232	\\
Breinigsville	&	4,138	\\
Coplay	&	3,192	\\
Macungie	&	3,074	\\
Schnecksville	&	2,935	\\
Coopersburg	&	2,386	\\
Alburtis	&	2,361	\\
Cetronia	&	2,115	\\
Trexlertown	&	1,988	\\
Laurys Station	&	1,243	\\
New Tripoli	&	898	\\
Slatedale	&	455	\\
			\hline														
\end{tabular} 
\label{table:LehighInst}
\end{table}


\begin{table}[h!]
 \center 
 \footnotesize
   \renewcommand{\arraystretch}{0.7}
  \caption{Average demand of each node in Lehigh 1 and Lehigh 2. } 
\begin{tabular}{lll}
\hline
\textbf{Node} &	\textbf{Lehigh 1 }& \textbf{Lehigh 2}	\\
\hline 
Allentown	&		40 & 60\\
Bethlehem	&	40 & 60	\\
Emmaus	&	40 & 43 	\\
Ancient Oaks	&	30 & 25	\\
Catasauqua	&	30 &25	\\
Wescosville	&	30 & 22	\\
Fountain Hill	&	20 & 18	\\
Dorneyville	&	20 &17	\\
Slatington	&	20	& 16\\
Breinigsville	&	20 & 16	\\
Coplay	&	20	& 12\\
Macungie	&	20 & 12	\\
Schnecksville	&	20 & 11	\\
Coopersburg	&	20 &9	\\
Alburtis	&	20 &9	\\
Cetronia	&	20 & 8	 \\
Trexlertown	&	20 & 8\\
Laurys Station	&20	& 5\\
New Tripoli	&	15 &3	\\
Slatedale	&	15 &3	\\
			\hline														
\end{tabular} 
\label{table:AvgDemand_LehighInst}
\end{table}

\begin{table}[h!]
 \center 
 \footnotesize
   \renewcommand{\arraystretch}{0.7}
  \caption{Average demand of each node in period 1 and period 2.} 
\begin{tabular}{lll}
\hline
\textbf{Node} &	\textbf{Period 1}& \textbf{Period 2}	\\
\hline 
Allentown	&	60	&	9	\\
Bethlehem	&	60	&	8	\\
Emmaus	&	43	&	8	\\
AncientOaks	&	25	&	5	\\
Catasauqua	&	25	&	3	\\
Wescosville	&	22	&	3	\\
FountainHill	&	18	&	18	\\
Dorneyville	&	17	&	17	\\
Slatington	&	16	&	16	\\
Breinigsville	&	16	&	16	\\
Coplay	&	12	&	12	\\
Macungie	&	12	&	12	\\
Schnecksville	&	11	&	11	\\
Coopersburg	&	9	&	9	\\
Alburtis	&	9	&	60	\\
Cetronia	&	8	&	60	\\
Trexlertown	&	8	&	43	\\
LaurysStation	&	5	&	25	\\
NewTripoli	&	3	&	25	\\
Slatedale	&	3	&	22	\\
			\hline														
\end{tabular} 
\label{table:period}
\end{table}

\begin{table}[t!]
 \center 
 \footnotesize
   \renewcommand{\arraystretch}{0.3}
  \caption{Optimal MFs locations in period 1 and period 2.}
\begin{tabular}{llllllllllllll}
\hline
 \textbf{MF} & \multicolumn{3}{c}{\textbf{MAD-DRO}} & \multicolumn{3}{c}{\textbf{W-DRO}}  & &\multicolumn{3}{c}{\textbf{SP}}  \\ \cline{2-3}  \cline{5-6} \cline{8-9}
	&	Period 1	&	Period 2	&	&	Period 1	&	Period 2	&	&	Period 1	&	Period 2	\\ \cline{1-9}
\hline
1	&	Allentown	&	New Tripoli	&	&	Allentown	&	Cetronia	&	&	Allentown	&	Fountain Hill	\\
2	&	Bethlehem	&	Schnecksville	&	&	Bethlehem	&	New Tripoli	&	&	Bethlehem	&	Slatedale	\\
3	&	Emmaus	&	Trexlertown	&	&	Emmaus	&	Ancient Oaks	&	&	Dorneyville	&	Cetronia	\\
4	&	Wescosville	&	Alburtis	&	&	Breinigsville	&	Trexlertown	&	&	Ancient Oaks	&	Alburtis	\\ \cline{1-9}
\hline																	
\end{tabular} 
\label{table:OptimalLocations}
\end{table}

\newpage

 \section{Calibrating the Wasserstein Radius in the W-DRO Model}\label{WDRO_Radius}

The Wasserstein ball’s radius $\epsilon$ in the W-DRO model is an input parameter, where a larger $\epsilon$ implies that we seek more distributionally robust solutions. For each training data, different values of $\epsilon$ may result in robust solutions $\xbb (\epsilon, N)$ with very different out-of-sample performance $\hat{\E} [Q(\xbb (\epsilon, N),\Wb)]$. On the one hand, the radius should not be too small. \textcolor{black}{Otherwise,}  the problem may behave like sample average approximation and hence, losing the purpose of robustification. In particular, if we set  the radius to zero,  the ambiguity set shrinks to a singleton that contains only the nominal distribution, in which case the DRO problem reduces to an ambiguity-free SP \citep{esfahani2018data}. But, on the other hand, the radius should not be too large to avoid conservative solutions, which is one of the major criticism faced by traditional RO methods. Given that the true distribution $\Prob$ is possibly unknown, it is impossible to compute $\epsilon$  that minimizes $\hat{\E} [Q(\xbb (\epsilon, N),\Wb)]$. Thus, as detailed in \cite{esfahani2018data}, the best we can hope for is to approximate  $\epsilon^{\mbox{\tiny opt}}$ using the training data set. 

 As pointed out by \cite{esfahani2018data} and \cite{gao2020finite}, practically, the radius is often selected via cross-validation. We employ the following widely used cross-validation method to estimate $\epsilon^{\mbox{\tiny opt}}$  as in \cite{jiang2019data} and \cite{esfahani2018data}.  First,  for each $N\in \{10, \ 50, \ 100\} $ and each $\epsilon \in \{ 0.01, \ 0.02,\ldots, \ 0.09, \ 0.1, \ldots, \ 0.9,1, \ldots, \ 10 \}$ (i.e., a log-scale interval as in \cite{esfahani2018data}, \cite{jian2017integer}, \cite{tsang2021distributionally}), we randomly partition the data into a training ($N'$) and testing set ($N''$).  Using the training set, we solve W-DRO to obtain the optimal first-stage solution $\xbb (\epsilon,N)$ for each $\epsilon$ and $N$. Then, we use the \textcolor{black}{testing} data to evaluate these solutions by computing $\hat{\E}_{\Prob_{N''}}[Q(\xbb (\epsilon, N),\Wb)]$ (where $\Prob_{N''}$ is the empirical distribution based on the testing data $N''$) via sample average approximation. That is, we solve the second-stage with $\xbb$ fixed to $\xbb (\epsilon,N')$ and $N''$ data and compute the corresponding second-stage cost $\hat{\E}_{\Prob_{N''}}[Q(\xbb (\epsilon, N''),\Wb)]$. Finally, we set $\epsilon^{\mbox{\tiny best}}_N$ to any $\epsilon$ that minimizes $\hat{\E}_{\Prob_{N''}}[Q(\xbb (\epsilon, N''),\Wb)]$. We repeat this procedure 30 times for each $N$ and set $\epsilon$ to the average of the $\epsilon^{\mbox{\tiny best}}_N$ across these 30 replications. 
 
 We found that $\epsilon^{\mbox{\tiny best}}_N$ equals (7, 5, 2) when $N$=(10, 50, 100) for most instances. It is expected that $\epsilon$  decreases with $N$   (see, e.g., \cite{esfahani2018data}, \cite{jiang2019data}, \cite{tsang2021distributionally}).  Intuitively, a small sample does not have sufficient distributional information, and thus a larger $\epsilon$ produces distributionally robust solutions that better hedge against ambiguity. In contrast, with a larger sample, we may have more information from the data, and thus we can make a less conservative decision using a smaller $\epsilon$ value.

While Wasserstein ambiguity sets offer powerful out-of-sample performance guarantees and enable practitioners to control the model's conservativeness by choosing $\epsilon$, moment-based ambiguity sets often display better tractability properties. In fact, various studies provided evidence that DRO models with moment ambiguity sets are more tractable than the corresponding SP because the intractable high-dimensional integrals in the objective function are replaced with tractable (generalized) moment problems (see, e.g., \cite{esfahani2018data, delage2010distributionally, goh2010distributionally, wiesemann2014distributionally}).  In contrast, DRO models with Wasserstein ambiguity sets tend to be more computationally challenging than some moment-based DRO model and their SP counterparts, especially when $N$ is large. In this paper, we obtained similar observations. For a detailed discussion we refer to \cite{esfahani2018data} and references therein.

\newpage

\section{Additional CPU Time Results}\label{Appx:CPU2}

\color{black}

\begin{table}[h!]
\color{black}
\center 
      \small 
\caption{Computational details of solving the MAD-CVaR model ($ \Wb \in [20, 60]$).}
   \renewcommand{\arraystretch}{0.6}
\begin{tabular}{lllllllllllllllllllllllllllllllllll}
 \hline
&  &  & \multicolumn{6}{c}{$C=60$ }  &&   \multicolumn{6}{c}{$C=100$ } \\  \cline{5-10} \cline{12-17}
Inst & $I$  &  $T$  & &  \multicolumn{3}{c}{CPU time} &   \multicolumn{3}{c}{iteration}  &&  \multicolumn{3}{c}{CPU time}   & \multicolumn{3}{c}{iteration} \\ \cline{5-10} \cline{12-17}
& &  & &  Min & Avg & Max & Min & Avg & Max & & Min & Avg & Max & Min & Avg & Max\\
\hline
1	&	10	&	10	&	&	1	&	3	&	6	&	4	&	12	&	25	&	&	2	&	7	&	15	&	4	&	12	&	25	\\
2	&	10	&	20	&	&	5	&	11	&	22	&	11	&	27	&	50	&	&	18	&	21	&	29	&	35	&	44	&	61	\\
3	&	15	&	10	&	&	3	&	8	&	19	&	5	&	16	&	36	&	&	5	&	13	&	32	&	9	&	24	&	58	\\
4	&	15	&	20	&	&	28	&	63	&	107	&	23	&	56	&	93	&	&	35	&	60	&	85	&	29	&	86	&	121	\\
5	&	20	&	10	&	&	13	&	25	&	37	&	12	&	23	&	33	&	&	19	&	75	&	107	&	17	&	47	&	68	\\
6	&	20	&	20	&	&	207	&	407	&	617	&	59	&	106	&	165	&	&	433	&	969	&	1576	&	103	&	198	&	303	\\
7	&	25	&	10	&	&	111	&	242	&	494	&	25	&	43	&	67	&	&	147	&	454	&	794	&	42	&	75	&	125	\\
8	&	25	&	20	&	&	2083	&	2944	&	3903	&	11	&	94	&	169	&	&	2979	&	3274	&	3600	&	164	&	237	&	316	\\
\hline 
\end{tabular}\label{table:MADCVAR_2060}
\end{table}

\begin{table}[h!]
\color{black}
\center 
      \small 
\caption{Computational details of solving the MAD-CVaR model ($ \Wb \in [50, 100]$).}
   \renewcommand{\arraystretch}{0.6}
\begin{tabular}{lllllllllllllllllllllllllllllllllll}
 \hline
&  &  & \multicolumn{6}{c}{$C=60$ }  &&   \multicolumn{6}{c}{$C=100$ } \\  \cline{5-10} \cline{12-17}
Inst & $I$  &  $T$  & &  \multicolumn{3}{c}{CPU time} &   \multicolumn{3}{c}{iteration}  &&  \multicolumn{3}{c}{CPU time}   & \multicolumn{3}{c}{iteration} \\ \cline{5-10} \cline{12-17}
& &  & &  Min & Avg & Max & Min & Avg & Max & & Min & Avg & Max & Min & Avg & Max\\
\hline
1	&	10	&	10	&	&	2	&	6	&	10	&	5	&	17	&	30	&	&	2	&	5	&	8	&	4	&	14	&	26	\\
2	&	10	&	20	&	&	14	&	22	&	29	&	22	&	35	&	46	&	&	18	&	23	&	29	&	19	&	31	&	44	\\
3	&	15	&	10	&	&	18	&	31	&	57	&	9	&	21	&	49	&	&	16	&	34	&	45	&	8	&	22	&	47	\\
4	&	15	&	20	&	&	67	&	106	&	156	&	23	&	68	&	139	&	&	64	&	132	&	255	&	27	&	62	&	106	\\
5	&	20	&	10	&	&	52	&	78	&	103	&	18	&	33	&	51	&	&	38	&	67	&	117	&	16	&	28	&	56	\\
6	&	20	&	20	&	&	549	&	715	&	935	&	96	&	124	&	172	&	&	757	&	1459	&	2802	&	91	&	124	&	184	\\
7	&	25	&	10	&	&	584	&	685	&	760	&	28	&	48	&	83	&	&	743	&	2277	&	2552	&	11	&	42	&	90	\\
\hline 
\end{tabular}\label{table:MADCVAR__50100}
\end{table}

\begin{table}[h!]
\color{black}
\center 
      \small 
\caption{Computational details of solving the W-CVaR model ($ \Wb \in [20, 60]$, $N=10$).}
   \renewcommand{\arraystretch}{0.6}
\begin{tabular}{lllllllllllllllllllllllllllllllllll}
 \hline
&  &  & \multicolumn{6}{c}{$C=60$ }  &&   \multicolumn{6}{c}{$C=100$ } \\  \cline{5-10} \cline{12-17}
Inst & $I$  &  $T$  & &  \multicolumn{3}{c}{CPU time} &   \multicolumn{3}{c}{iteration}  &&  \multicolumn{3}{c}{CPU time}   & \multicolumn{3}{c}{iteration} \\ \cline{5-10} \cline{12-17}
& &  & &  Min & Avg & Max & Min & Avg & Max & & Min & Avg & Max & Min & Avg & Max\\
\hline
1	&	10	&	10	&	&	7	&	19	&	17	&	5	&	5	&	6	&	&	16	&	20	&	24	&	12	&	15	&	19	\\
2	&	10	&	20	&	&	12	&	13	&	16	&	6	&	6	&	8	&	&	37	&	44	&	57	&	17	&	18	&	22	\\
3	&	15	&	10	&	&	27	&	44	&	70	&	4	&	5	&	5	&	&	22	&	28	&	38	&	10	&	12	&	13	\\
4	&	15	&	20	&	&	9	&	12	&	17	&	4	&	5	&	7	&	&	15	&	33	&	60	&	6	&	11	&	17	\\
5	&	20	&	10	&	&	34	&	35	&	39	&	4	&	5	&	6	&	&	72	&	88	&	98	&	14	&	16	&	20	\\
6	&	20	&	20	&	&	20	&	31	&	46	&	4	&	6	&	9	&	&	240	&	574	&	827	&	6	&	13	&	21	\\
\hline 
\end{tabular}\label{table:Wass_CVAR_2060_N10}
\end{table}

\begin{table}[h!]
\color{black}
\center 
      \small 
\caption{Computational details of solving the W-CVaR model ($ \Wb \in [50, 100]$, $N=10$).}
   \renewcommand{\arraystretch}{0.6}
\begin{tabular}{lllllllllllllllllllllllllllllllllll}
 \hline
&  &  & \multicolumn{6}{c}{$C=60$ }  &&   \multicolumn{6}{c}{$C=100$ } \\  \cline{5-10} \cline{12-17}
Inst & $I$  &  $T$  & &  \multicolumn{3}{c}{CPU time} &   \multicolumn{3}{c}{iteration}  &&  \multicolumn{3}{c}{CPU time}   & \multicolumn{3}{c}{iteration} \\ \cline{5-10} \cline{12-17}
& &  & &  Min & Avg & Max & Min & Avg & Max & & Min & Avg & Max & Min & Avg & Max\\
\hline
1	&	10	&	10	&	&	13	&	25	&	53	&	3	&	4	&	4	&	&	11	&	16	&	25	&	5	&	6	&	7	\\
2	&	10	&	20	&	&	6	&	8	&	14	&	3	&	3	&	4	&	&	10	&	12	&	16	&	5	&	6	&	8	\\
3	&	15	&	10	&	&	16	&	26	&	57	&	3	&	3	&	3	&	&	29	&	45	&	75	&	5	&	5	&	6	\\
4	&	15	&	20	&	&	38	&	45	&	61	&	3	&	3	&	3	&	&	63	&	82	&	136	&	5	&	6	&	7	\\
5	&	20	&	10	&	&	45	&	49	&	55	&	3	&	3	&	3	&	&	55	&	70	&	79	&	7	&	7	&	8	\\
\hline 
\end{tabular}\label{table:Wass_CVAR_50100_N10}
\end{table}

\newpage
\begin{table}[h!]
\color{black}
\center 
      \small 
\caption{Computational details of solving the W-CVaR model ($ \Wb \in [20, 60]$, $N=50$).}
   \renewcommand{\arraystretch}{0.6}
\begin{tabular}{lllllllllllllllllllllllllllllllllll}
 \hline
&  &  & \multicolumn{6}{c}{$C=60$ }  &&   \multicolumn{6}{c}{$C=100$ } \\  \cline{5-10} \cline{12-17}
Inst & $I$  &  $T$  & &  \multicolumn{3}{c}{CPU time} &   \multicolumn{3}{c}{iteration}  &&  \multicolumn{3}{c}{CPU time}   & \multicolumn{3}{c}{iteration} \\ \cline{5-10} \cline{12-17}
& &  & &  Min & Avg & Max & Min & Avg & Max & & Min & Avg & Max & Min & Avg & Max\\
\hline
1	&	10	&	10	&	&	28	&	35	&	40	&	5	&	6	&	7	&	&	104	&	117	&	129	&	15	&	19	&	22\\
2	&	10	&	20	&	&	33	&	42	&	51	&	6	&	7	&	9	&	&	84	&	112	&	149	&	15	&	19	&	24\\
3	&	15	&	10	&	&	28	&	34	&	43	&	4	&	5	&	5	&	&	109	&	136	&	148	&	14	&	17	&	19\\
4	&	15	&	20	&	&	34	&	59	&	102	&	4	&	5	&	7	&	&	189	&	1471	&	3730	&	17	&	20	&	23\\
5	&	20	&	10	&	&	68	&	88	&	125	&	4	&	6	&	7	&	&	455	&	1,084	&	1,556	&	24	&	30	&	34\\
\hline 
\end{tabular}\label{table:Wass_CVAR_2060_N50}
\end{table}

\begin{table}[h!]
\color{black}
\center 
      \small 
\caption{Computational details of solving the W-CVaR model ($ \Wb \in [50, 100]$, $N=50$).}
   \renewcommand{\arraystretch}{0.6}
\begin{tabular}{lllllllllllllllllllllllllllllllllll}
 \hline
&  &  & \multicolumn{6}{c}{$C=60$ }  &&   \multicolumn{6}{c}{$C=100$ } \\  \cline{5-10} \cline{12-17}
Inst & $I$  &  $T$  & &  \multicolumn{3}{c}{CPU time} &   \multicolumn{3}{c}{iteration}  &&  \multicolumn{3}{c}{CPU time}   & \multicolumn{3}{c}{iteration} \\ \cline{5-10} \cline{12-17}
& &  & &  Min & Avg & Max & Min & Avg & Max & & Min & Avg & Max & Min & Avg & Max\\
\hline
1	&	10	&	10	&	&	24	&	25	&	26	&	4	&	4	&	4	&	&	37	&	46	&	51	&	6	&	7	&	8	\\
2	&	10	&	20	&	&	29	&	30	&	32	&	4	&	4	&	4	&	&	52	&	73	&	124	&	7	&	8	&	9	\\
3	&	15	&	10	&	&	33	&	50	&	75	&	3	&	3	&	3	&	&	90	&	100	&	111	&	6	&	6	&	7	\\
4	&	15	&	20	&	&	50	&	59	&	68	&	3	&	3	&	3	&	&	130	&	310	&	743	&	5	&	6	&	9	\\
5	&	20	&	10	&	&	107	&	130	&	170	&	3	&	3	&	3	&	&	118	&	182	&	321	&	5	&	6	&	8	\\
\hline 
\end{tabular}\label{table:Wass_CVAR_50100_N50}
\end{table}

\textcolor{white}{sabjhasvghavwejfce}

\clearpage
\newpage

\section{Additional Out-of-Sample Results}\label{Appex:AdditionalOut}


\begin{figure}[h!]
 \centering
 \begin{subfigure}[b]{0.5\textwidth}
          \centering
        \includegraphics[width=\textwidth]{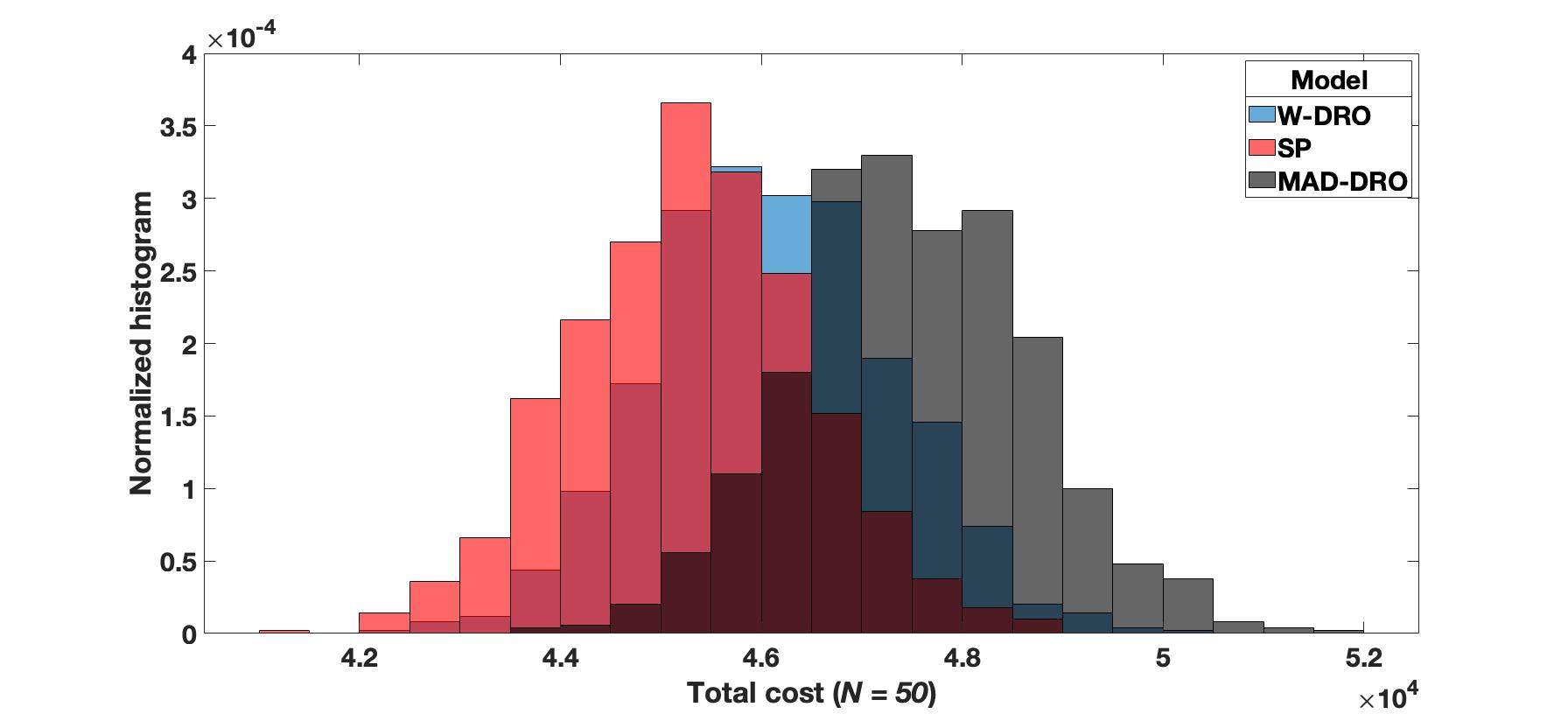}
        \caption{TC (Set 1, LogN)}
    \end{subfigure}%
  \begin{subfigure}[b]{0.5\textwidth}
        \centering
        \includegraphics[width=\textwidth]{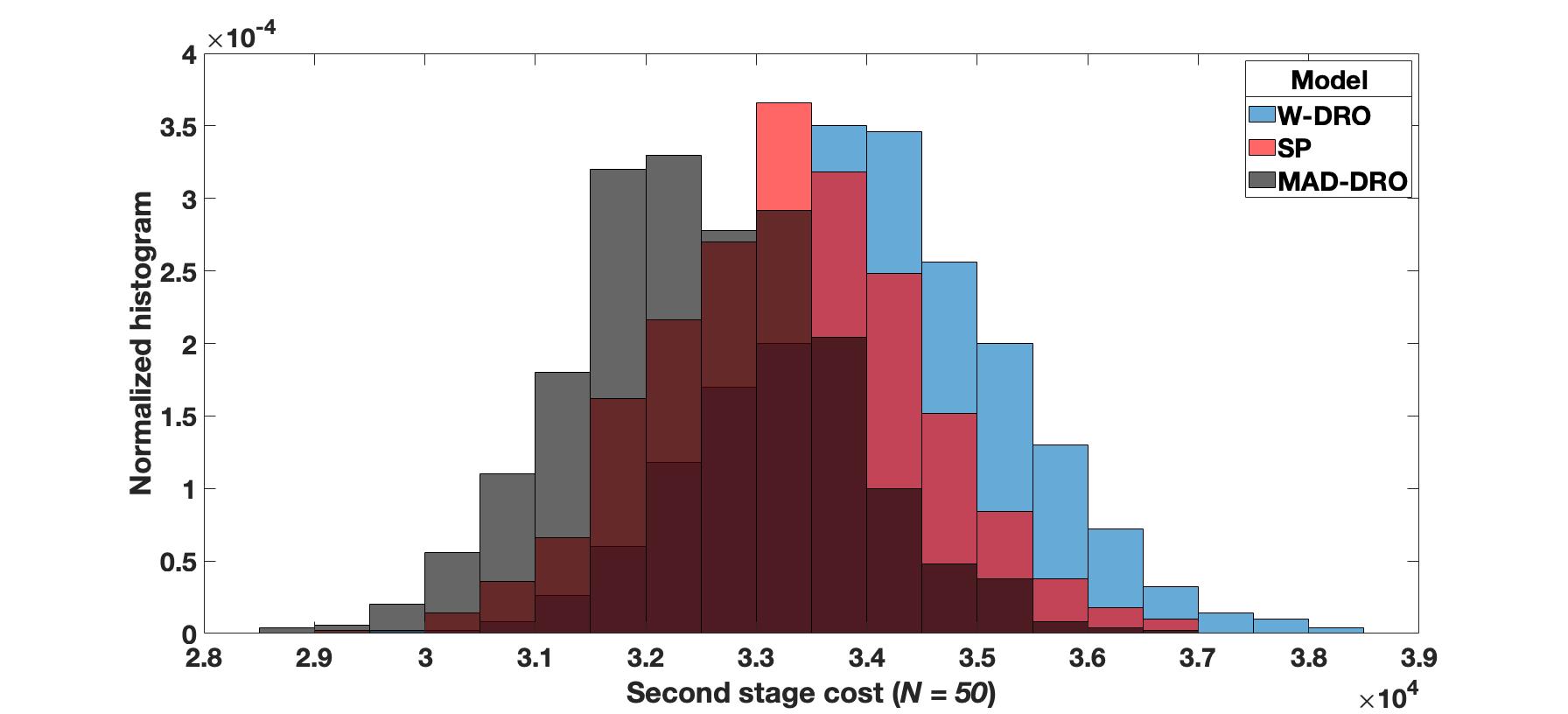}
        \caption{2nd (Set 1, LogN)}
    \end{subfigure}%

  \begin{subfigure}[b]{0.5\textwidth}
          \centering
        \includegraphics[width=\textwidth]{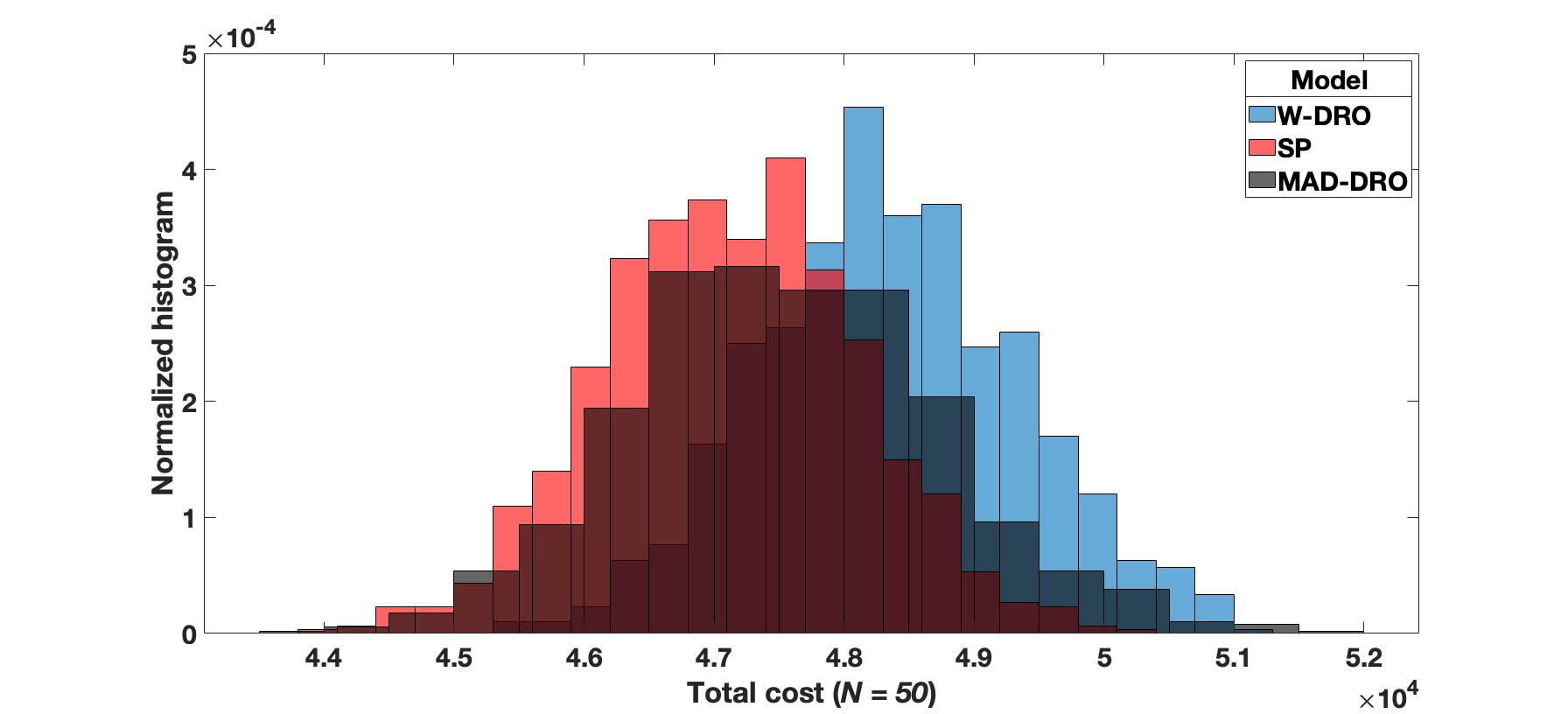}
        \caption{TC (Set 2, $\Delta=0$)}
    \end{subfigure}%
      \begin{subfigure}[b]{0.5\textwidth}
          \centering
        \includegraphics[width=\textwidth]{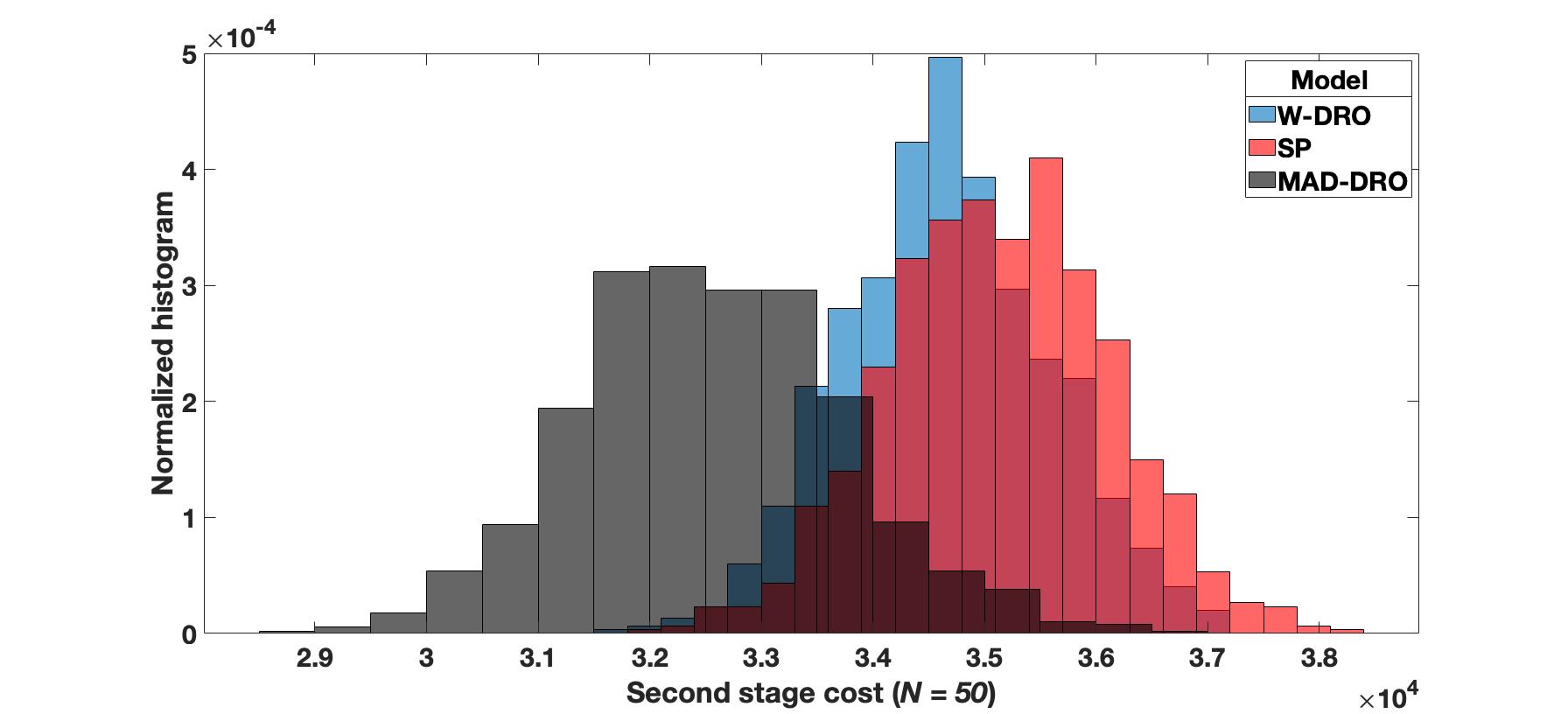}
        \caption{2nd (Set 2, $\Delta=0$)}
    \end{subfigure}%

  \begin{subfigure}[b]{0.5\textwidth}
        \centering
        \includegraphics[width=\textwidth]{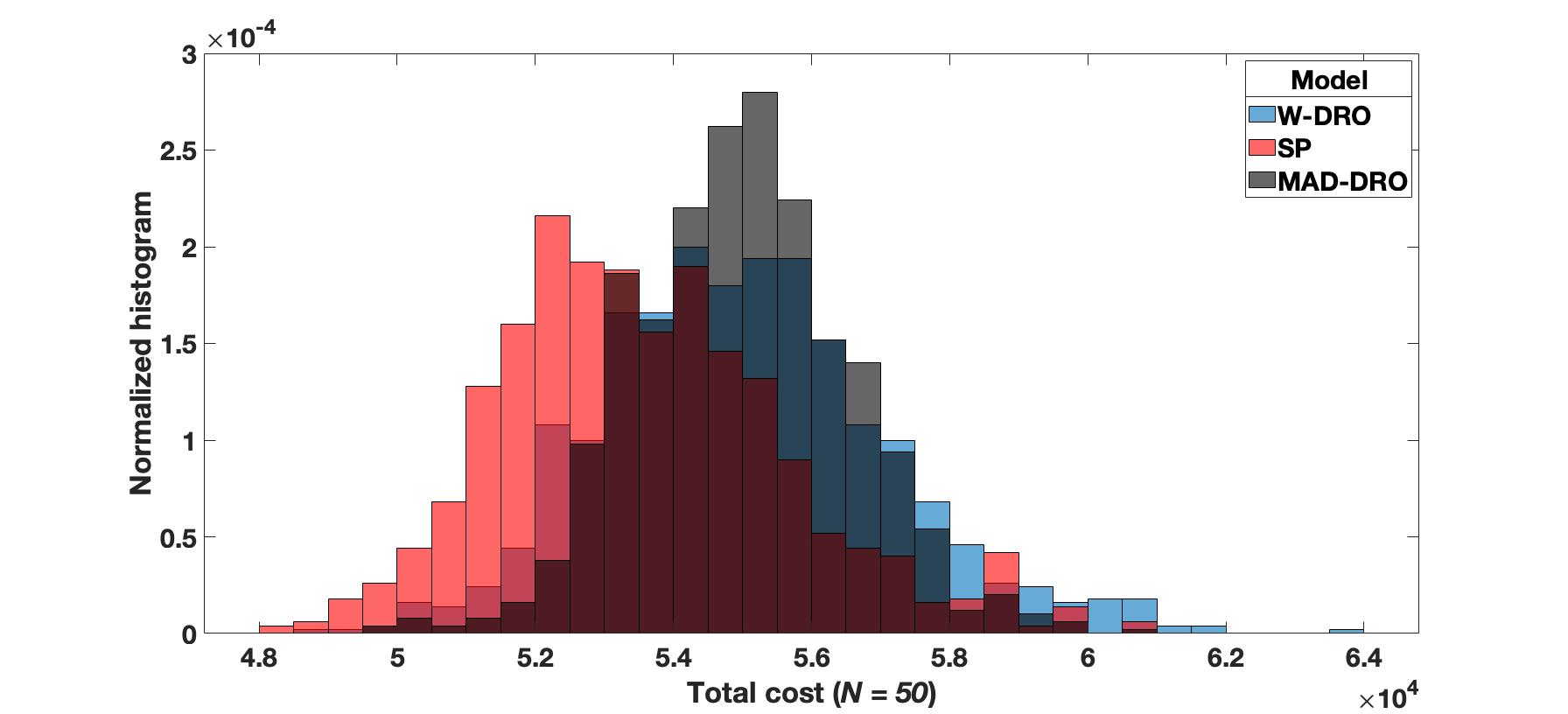}
        \caption{TC (Set 2, $ \Delta=0.25$)}
    \end{subfigure}%
      \begin{subfigure}[b]{0.5\textwidth}
        \centering
        \includegraphics[width=\textwidth]{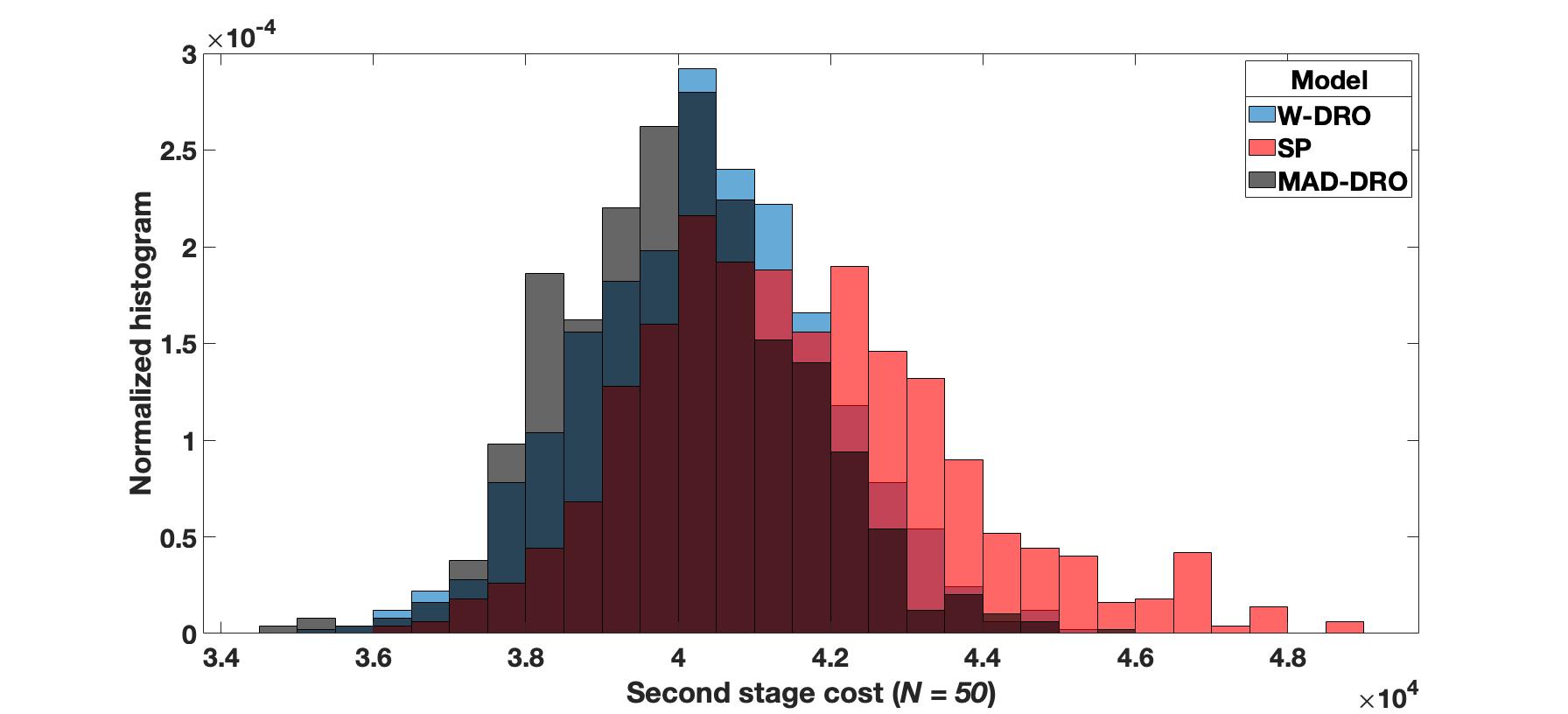}
        \caption{2nd (Set 2, $ \Delta=0.25$)}
    \end{subfigure}%
    
      \begin{subfigure}[b]{0.5\textwidth}
          \centering
        \includegraphics[width=\textwidth]{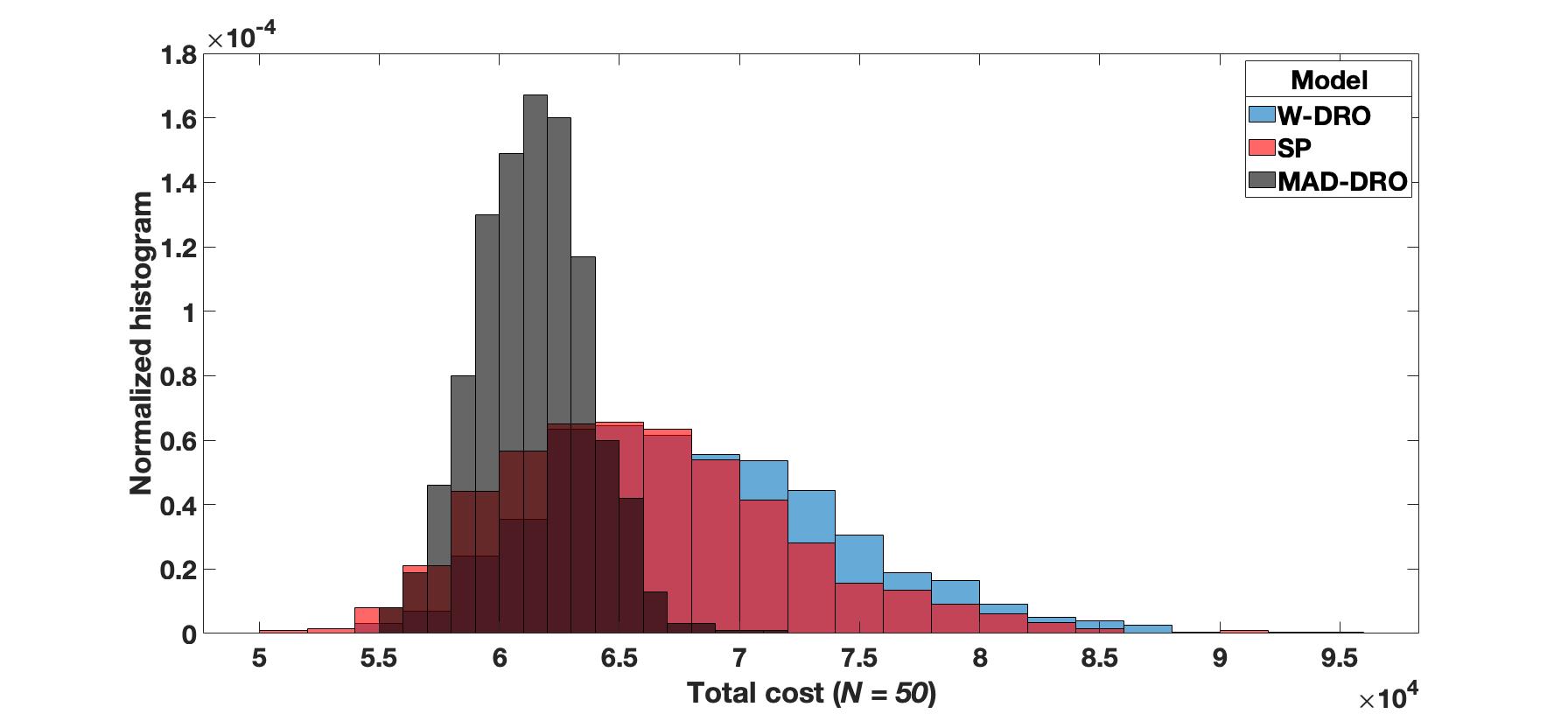}
        \caption{TC (Set 2, $\Delta=0.5$)}
    \end{subfigure}%
    \begin{subfigure}[b]{0.5\textwidth}
          \centering
        \includegraphics[width=\textwidth]{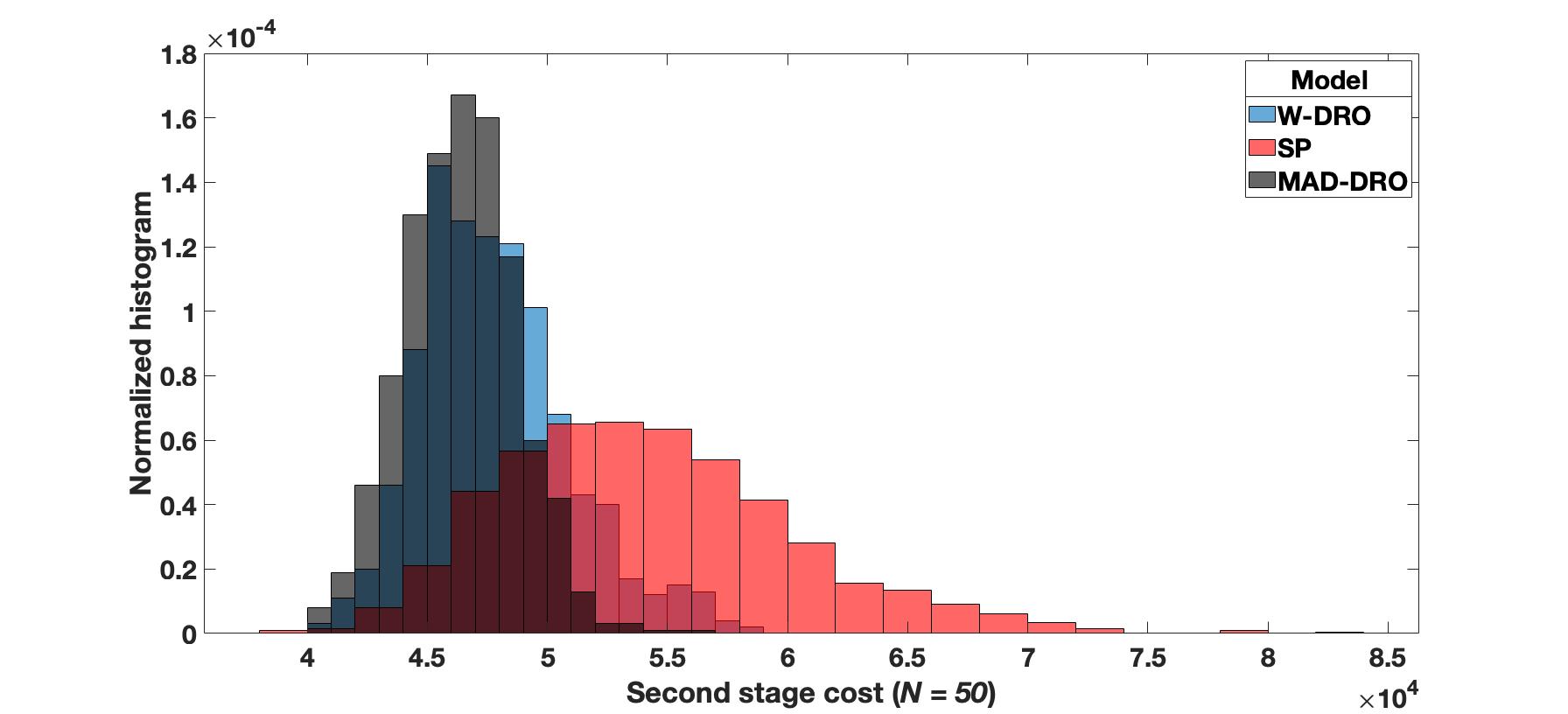}
        \caption{2nd (Set 2, $ \Delta=0.5$)}
    \end{subfigure}%
\caption{Normalized histograms of out-of-sample for Instance 3  ($\Wb \in [20, 60]$, $\pmb{N=50}$) under Set 1 (LogN) and Set 2 (with $\pmb{\Delta \in \{0, 0.25, 0.5\}}$).}\label{Fig3_UniN50_Inst3}
\end{figure}


\begin{figure}[t!]
 \centering
 \begin{subfigure}[b]{0.5\textwidth}
          \centering
        \includegraphics[width=\textwidth]{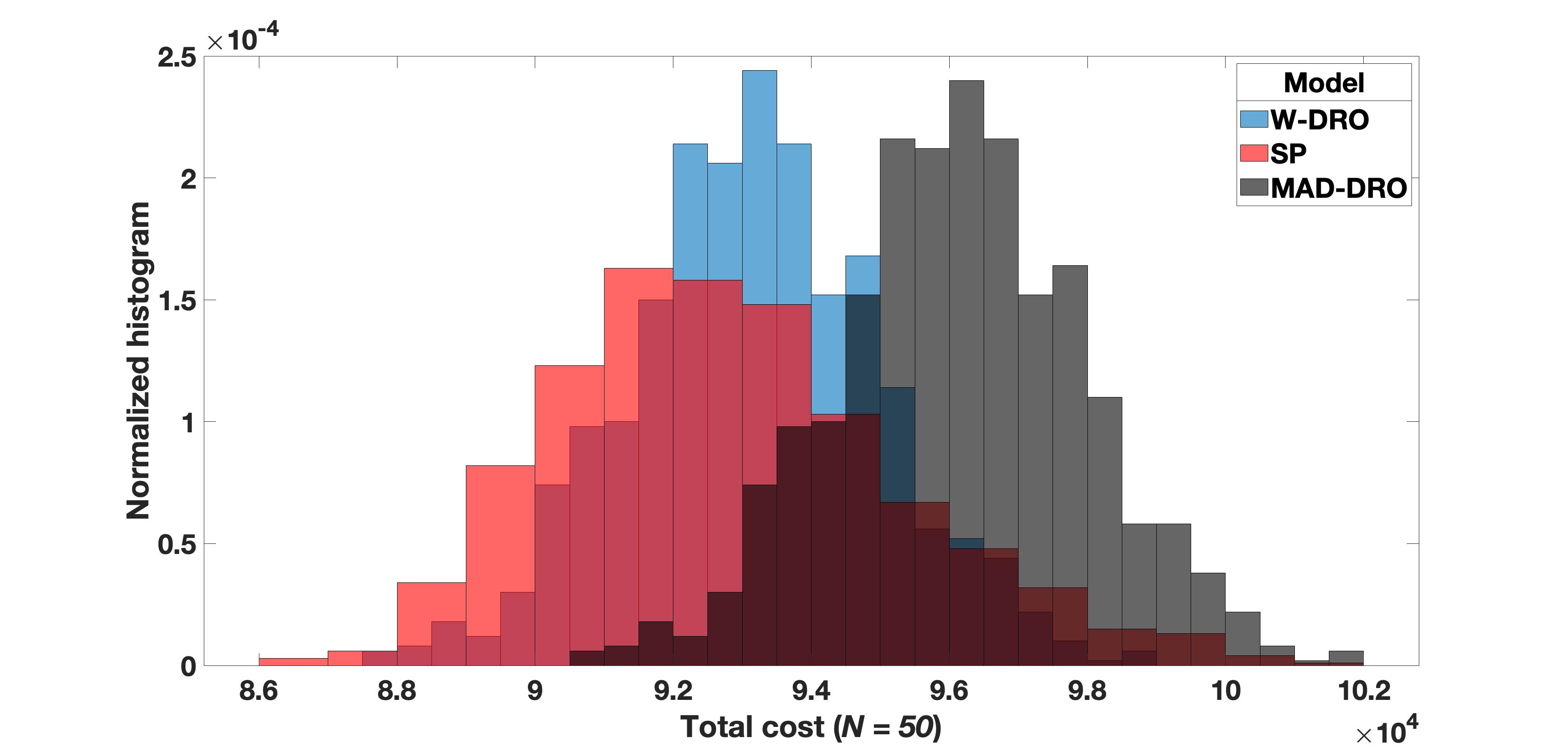}
        \caption{TC (Set 1, LogN)}\label{Inst3_LogNTC_Range2N50}
    \end{subfigure}%
  \begin{subfigure}[b]{0.5\textwidth}
        \centering
        \includegraphics[width=\textwidth]{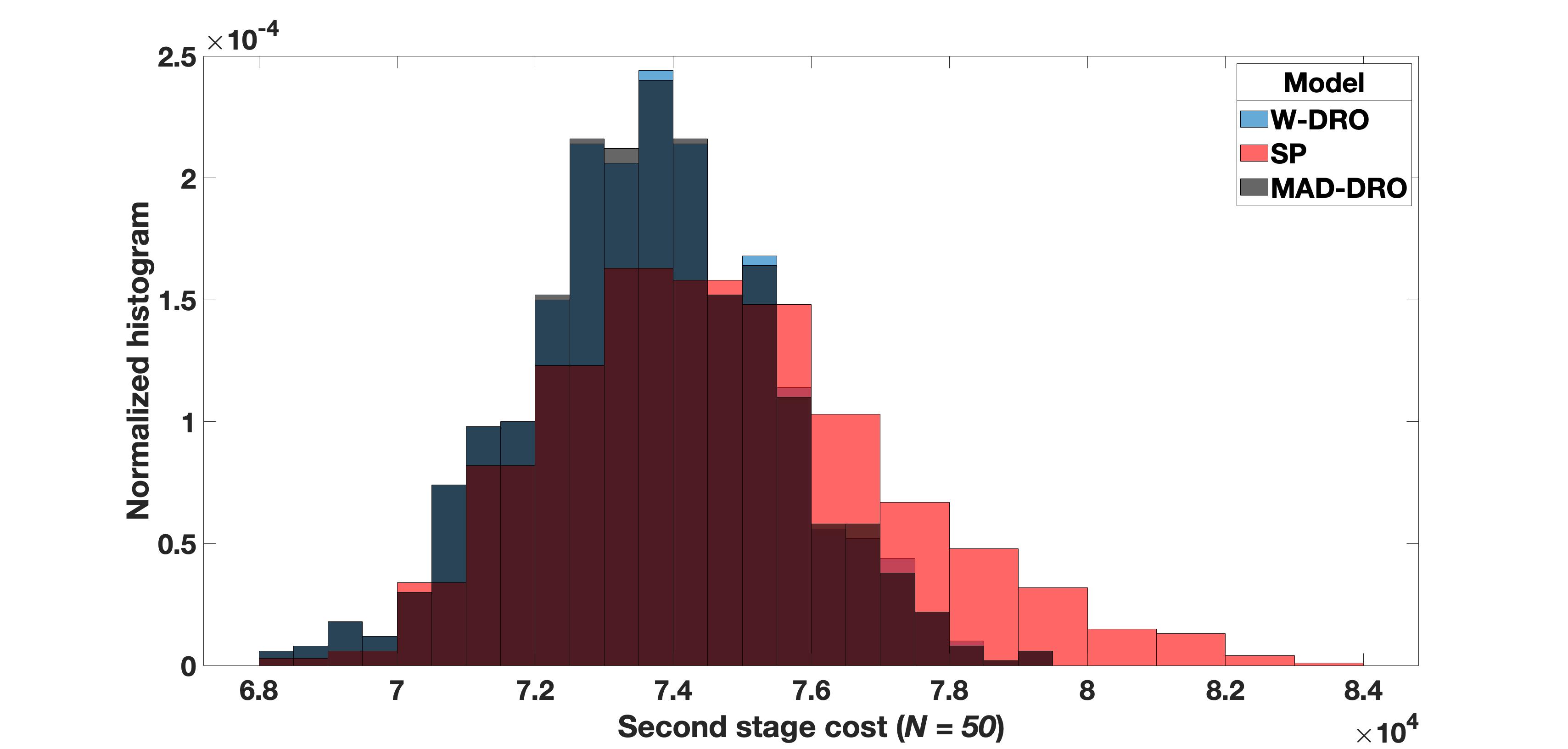}
        \caption{2nd (Set 1, LogN)}\label{Inst3_LogN2nd_Range2N50}
    \end{subfigure}%

  \begin{subfigure}[b]{0.5\textwidth}
          \centering
        \includegraphics[width=\textwidth]{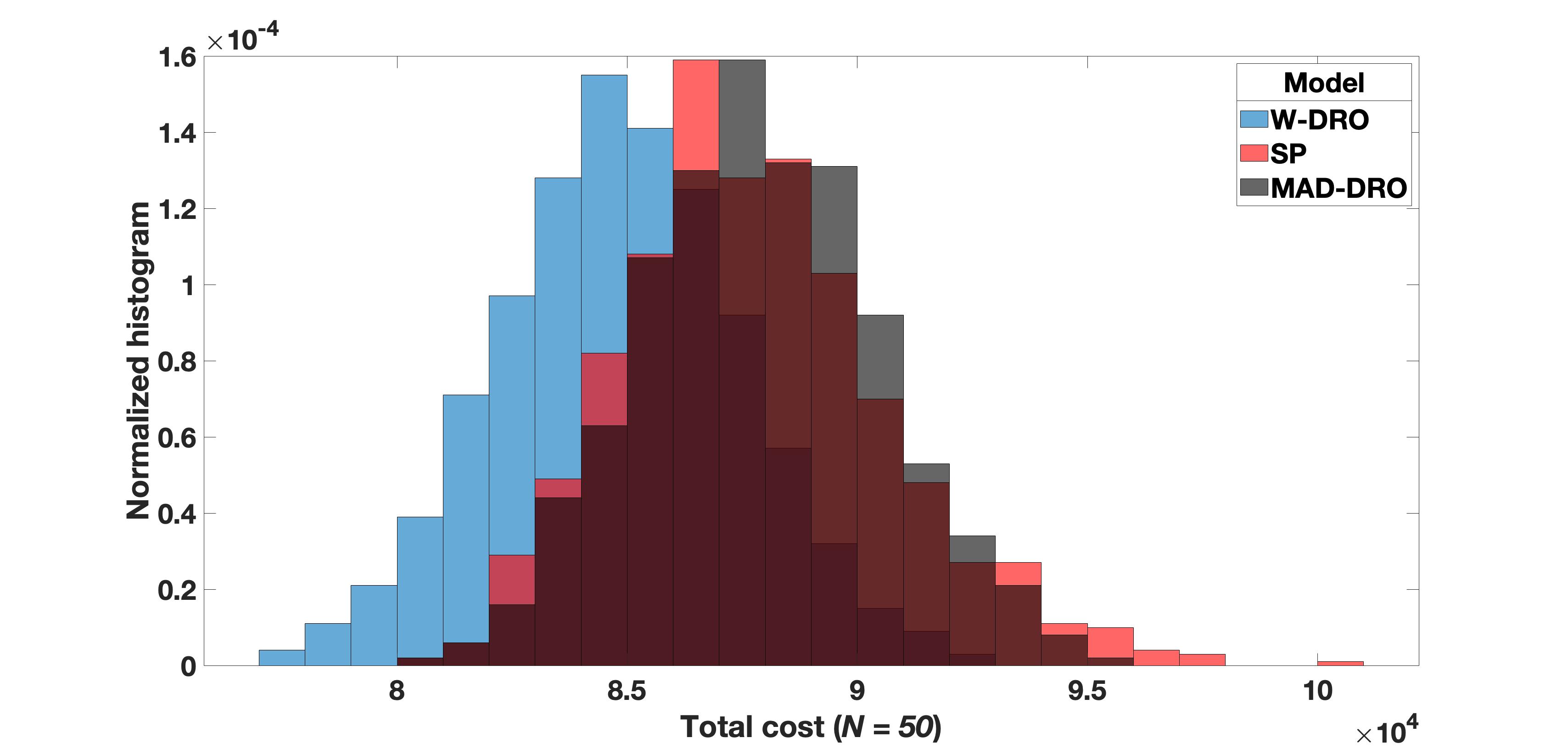}
        \caption{TC (Set 2, $\Delta=0$)}\label{Inst3_Uni0_TC_Range2N50}
    \end{subfigure}%
      \begin{subfigure}[b]{0.5\textwidth}
          \centering
        \includegraphics[width=\textwidth]{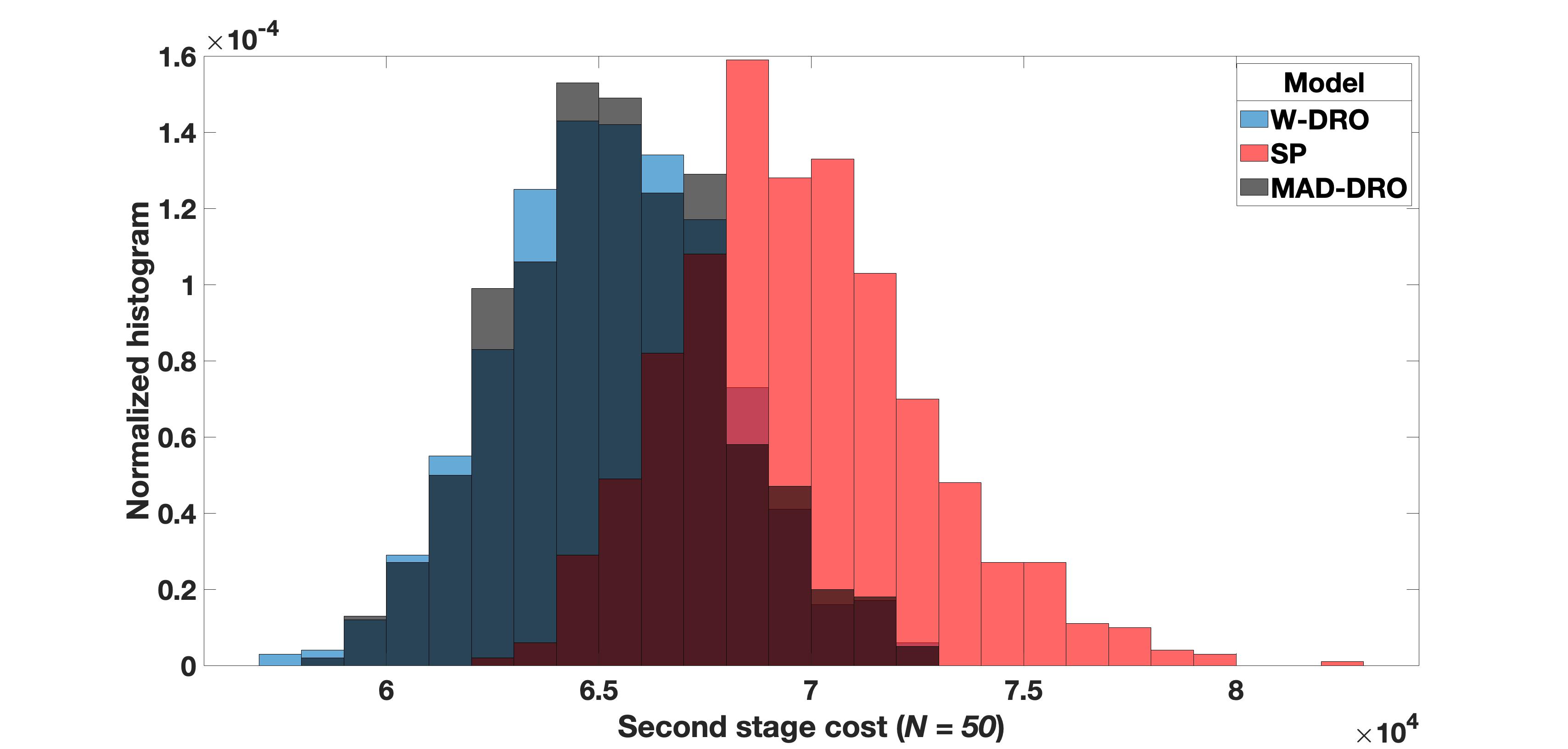} 
        \caption{2nd (Set 2, $\Delta=0$)}\label{Inst3_Uni0_2nd_Range2N50}
    \end{subfigure}
 

  \begin{subfigure}[b]{0.5\textwidth}
        \centering
        \includegraphics[width=\textwidth]{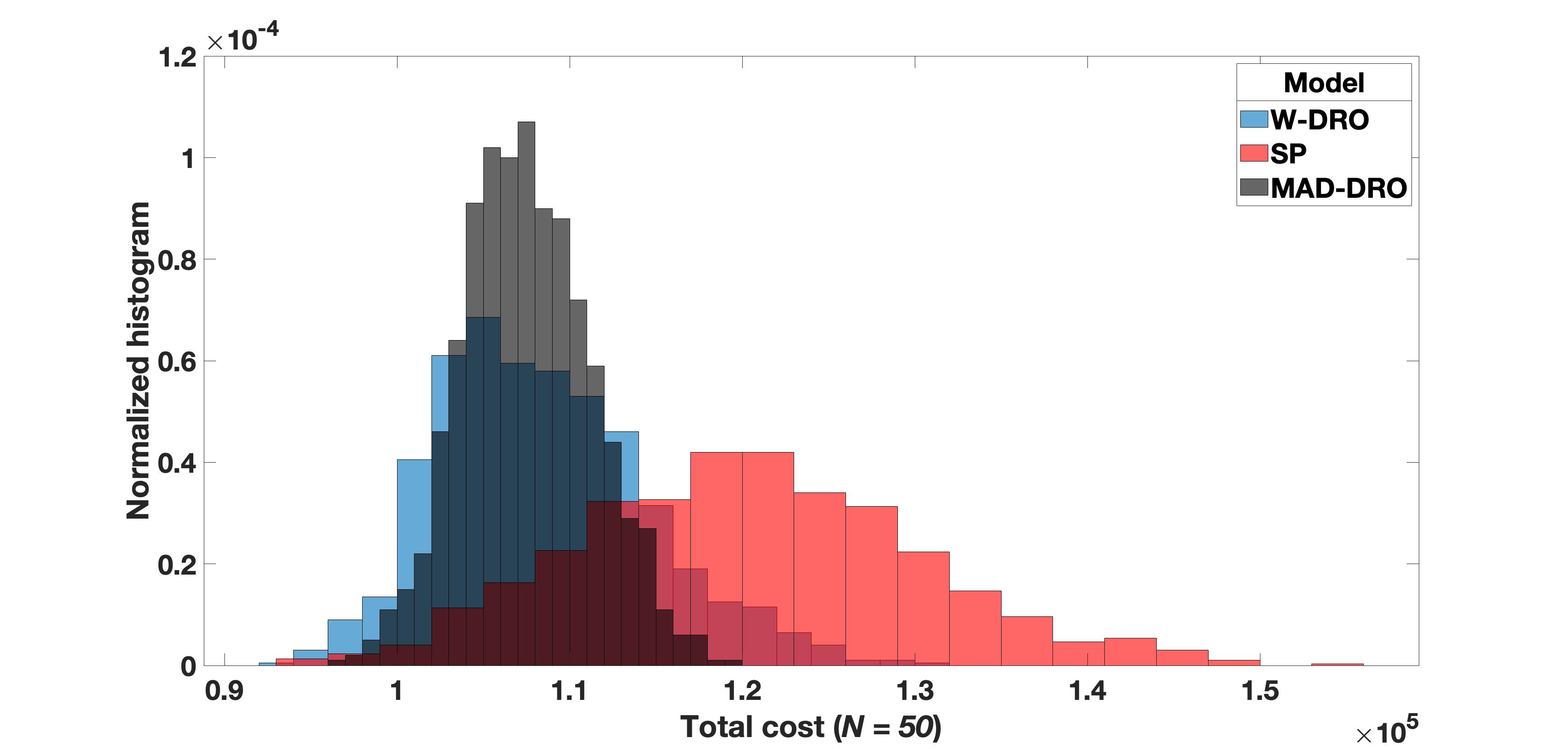}
        \caption{TC (Set 2, $ \Delta=0.25$)}
    \end{subfigure}%
      \begin{subfigure}[b]{0.5\textwidth}
        \centering
        \includegraphics[width=\textwidth]{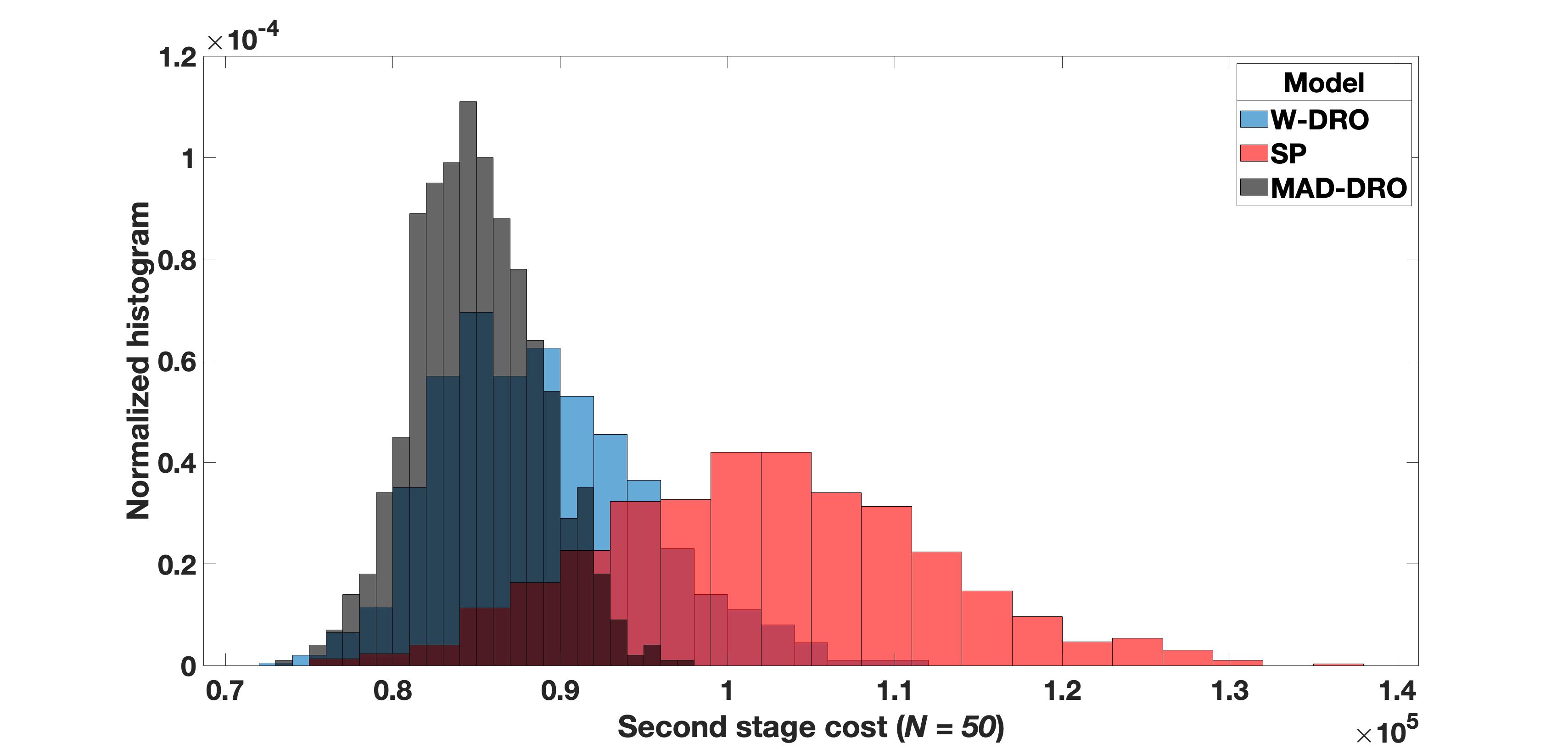}
        \caption{2nd (Set 2, $ \Delta=0.25$)}
    \end{subfigure}%
    
      \begin{subfigure}[b]{0.5\textwidth}
          \centering
        \includegraphics[width=\textwidth]{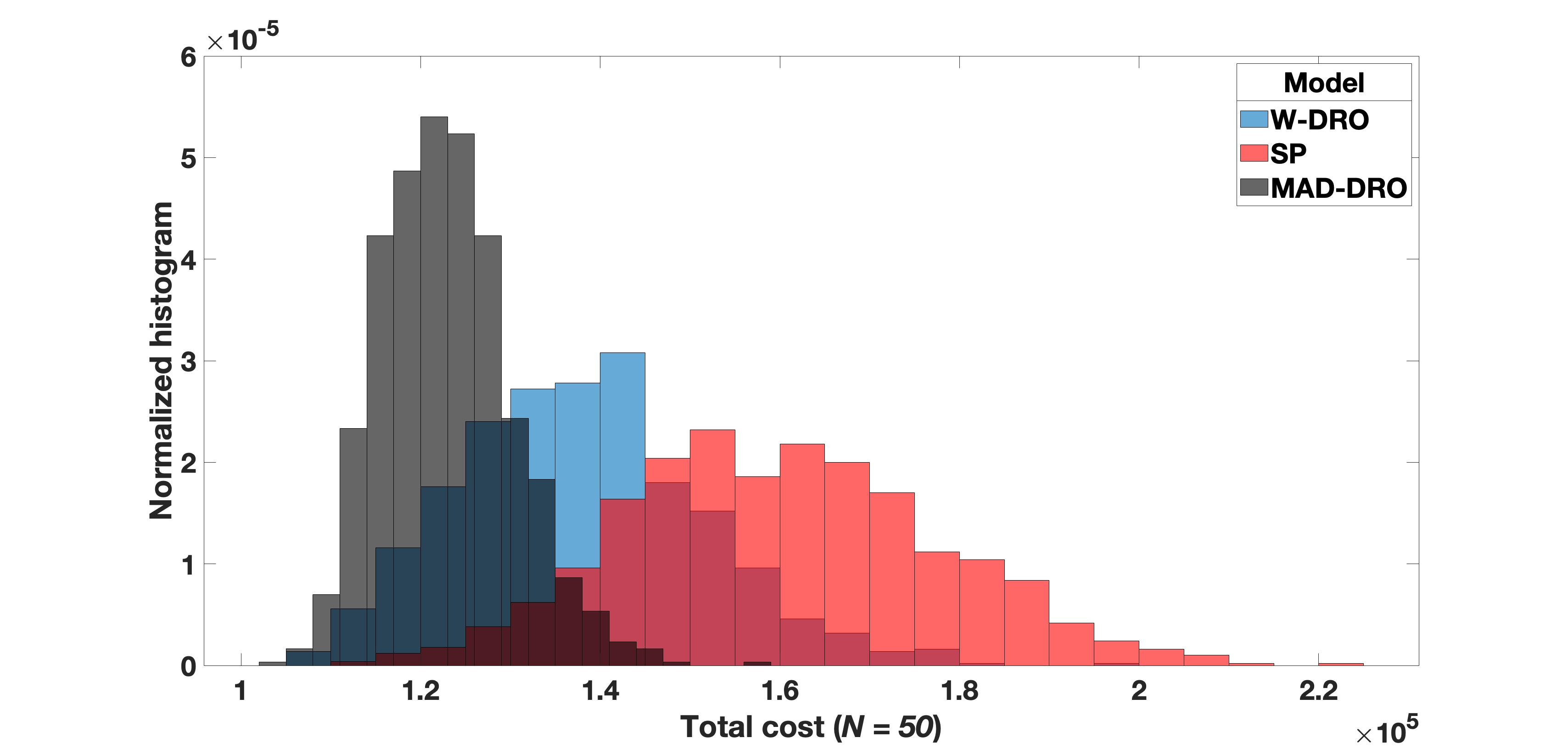}
        \caption{TC (Set 2, $\Delta=0.5$)}
    \end{subfigure}%
    \begin{subfigure}[b]{0.5\textwidth}
          \centering
        \includegraphics[width=\textwidth]{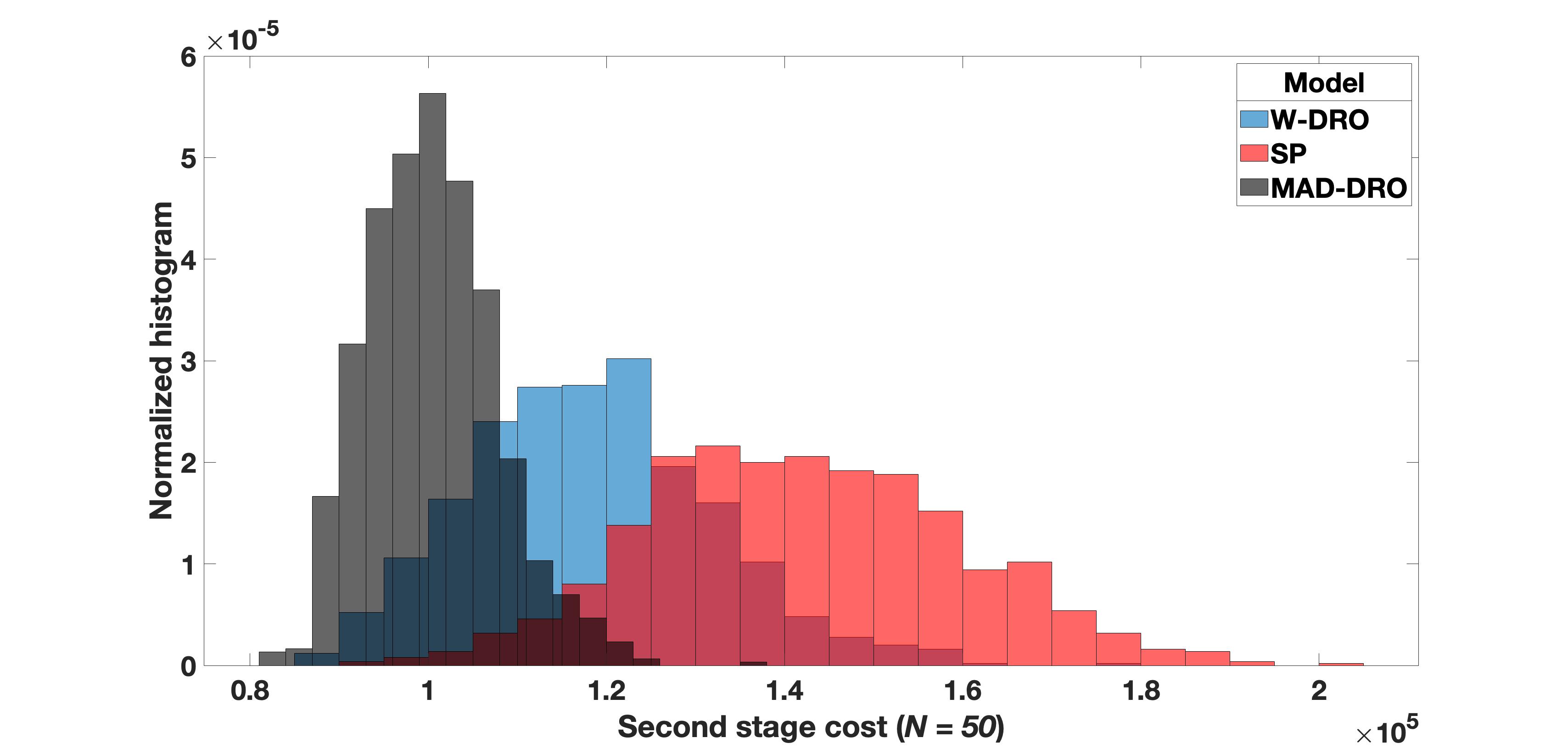}
        \caption{2nd (Set 2, $ \Delta=0.5$)} \label{Inst3_Uni25_2nd_Range2N50}
    \end{subfigure}%
\caption{Normalized histograms of out-of-sample TC and 2nd for Instance 3 ($\Wb \in [50, 100]$, $\pmb{N=50}$) under Set 1 (LogN) and Set 2 (with $\pmb{\Delta \in \{0, 0.25, 0.5\}}$). }\label{Fig3_UniN50_Inst3_Range2N50}
\end{figure}

\begin{figure}[t!]
 \centering
   \begin{subfigure}[b]{0.5\textwidth}
          \centering
        \includegraphics[width=\textwidth]{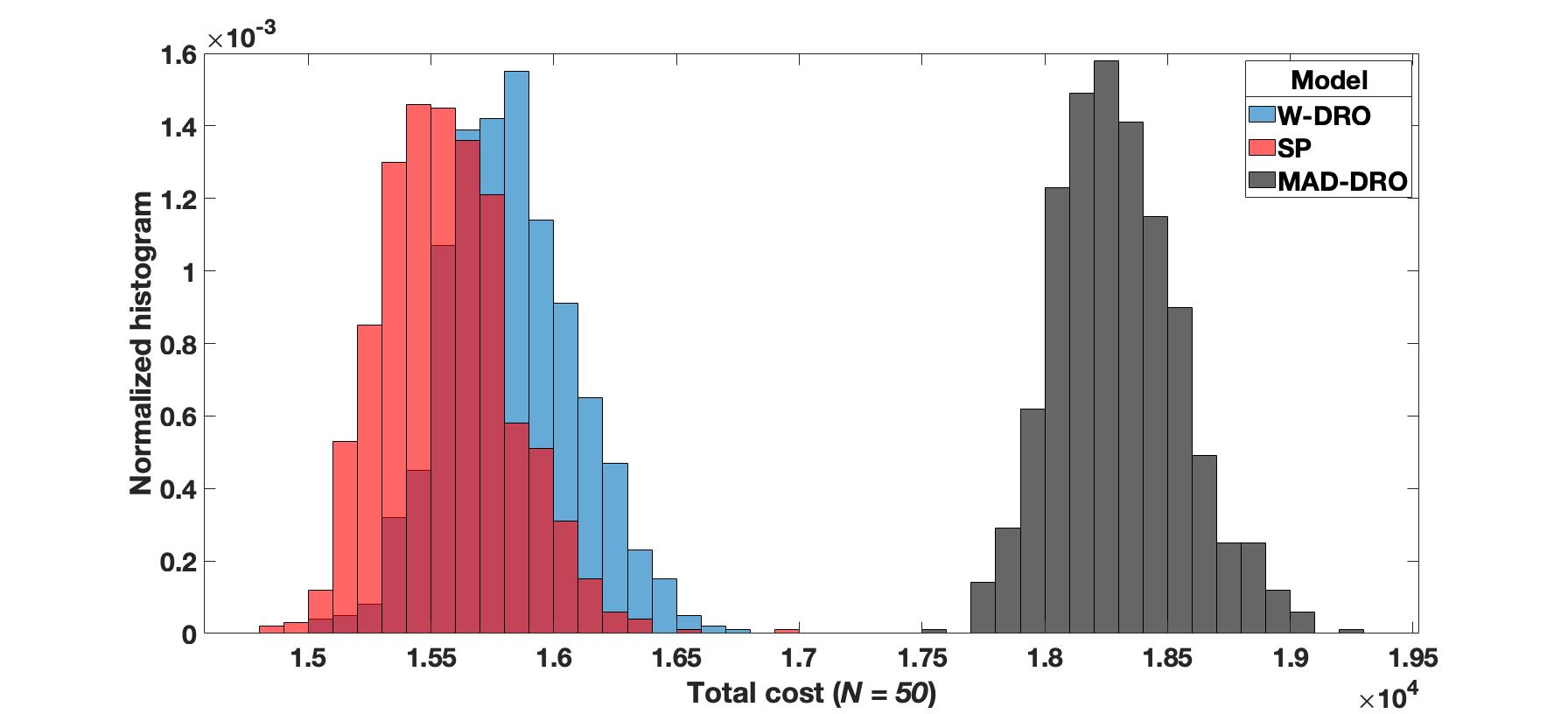}
        \caption{TC (Set 1, LogN)}
    \end{subfigure}%
      \begin{subfigure}[b]{0.5\textwidth}
          \centering
        \includegraphics[width=\textwidth]{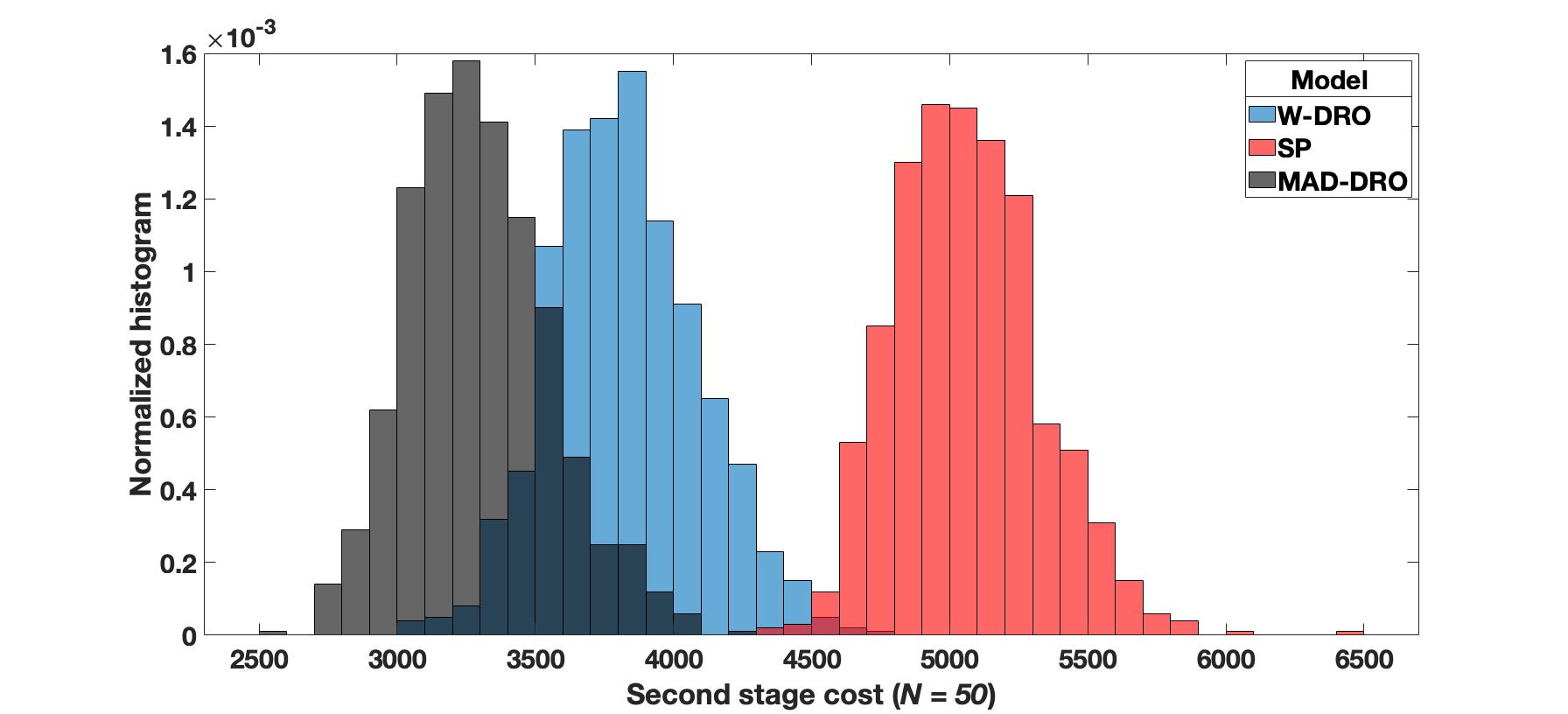}
        \caption{2nd (Set 1, LogN)}
    \end{subfigure}%
    
  \begin{subfigure}[b]{0.5\textwidth}
          \centering
        \includegraphics[width=\textwidth]{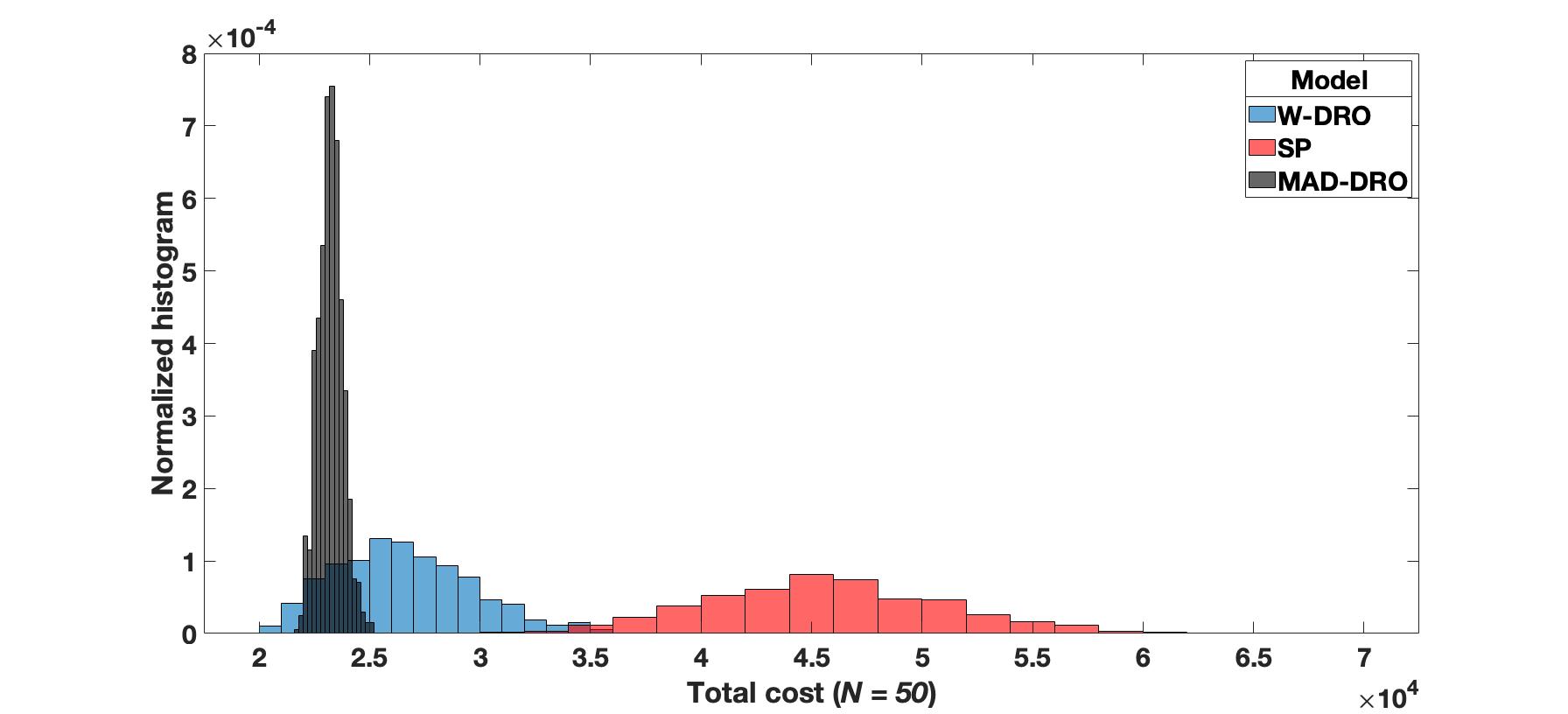}
        \caption{TC (Set 2, $\Delta=0$)}
    \end{subfigure}%
      \begin{subfigure}[b]{0.5\textwidth}
          \centering
        \includegraphics[width=\textwidth]{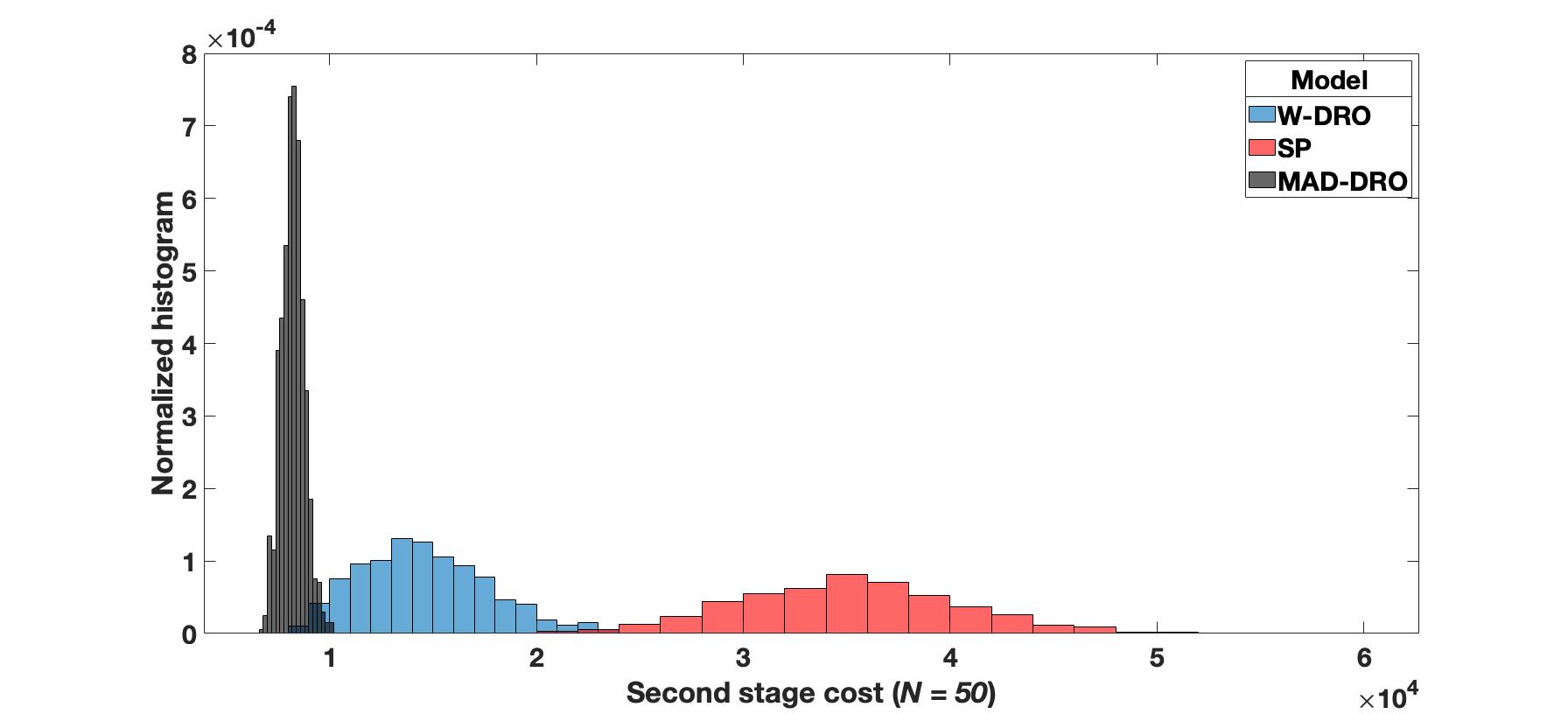}
        \caption{2nd  (Set 2, $\Delta=0$)}
    \end{subfigure}%

  \begin{subfigure}[b]{0.5\textwidth}
        \centering
        \includegraphics[width=\textwidth]{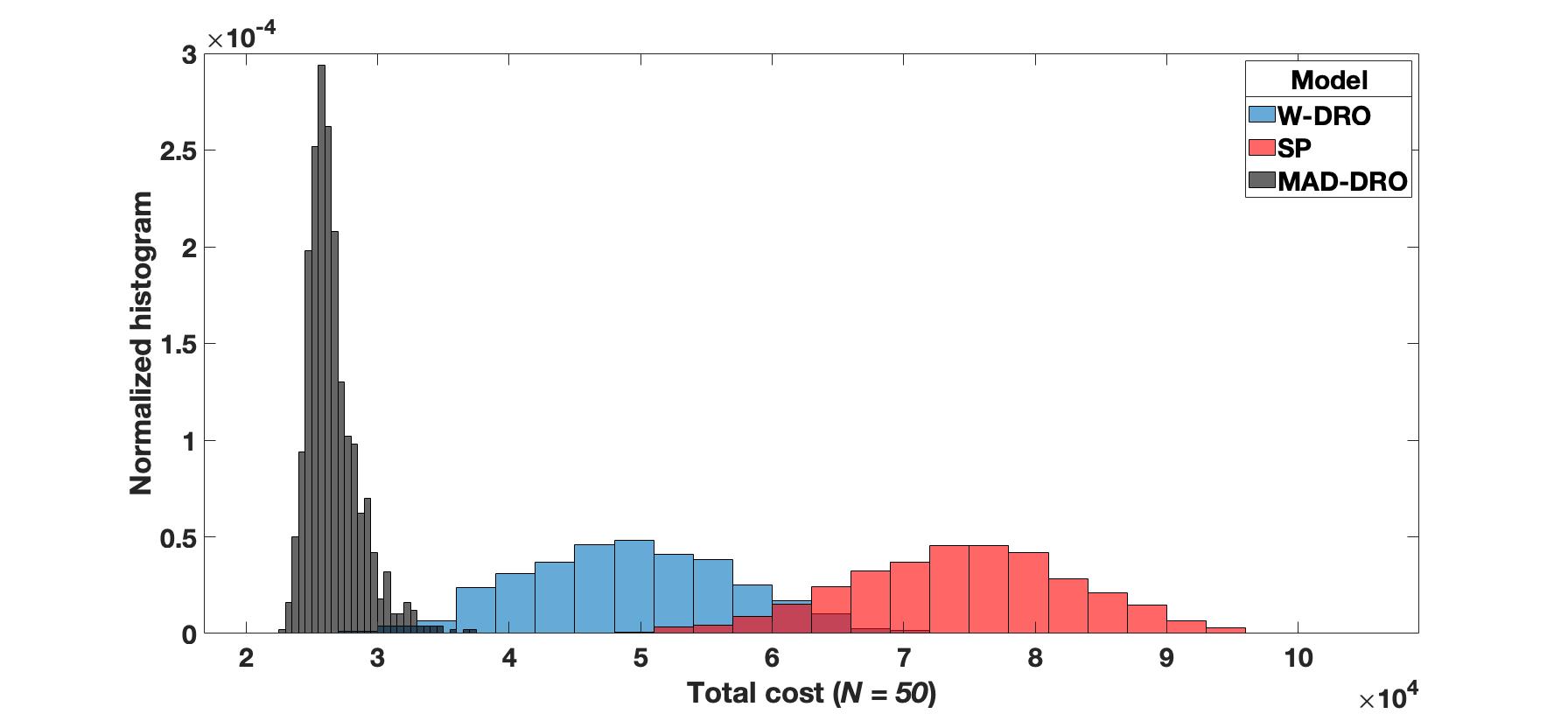}
        \caption{TC  (Set 2, $\Delta=0.25$)}
    \end{subfigure}%
      \begin{subfigure}[b]{0.5\textwidth}
        \centering
        \includegraphics[width=\textwidth]{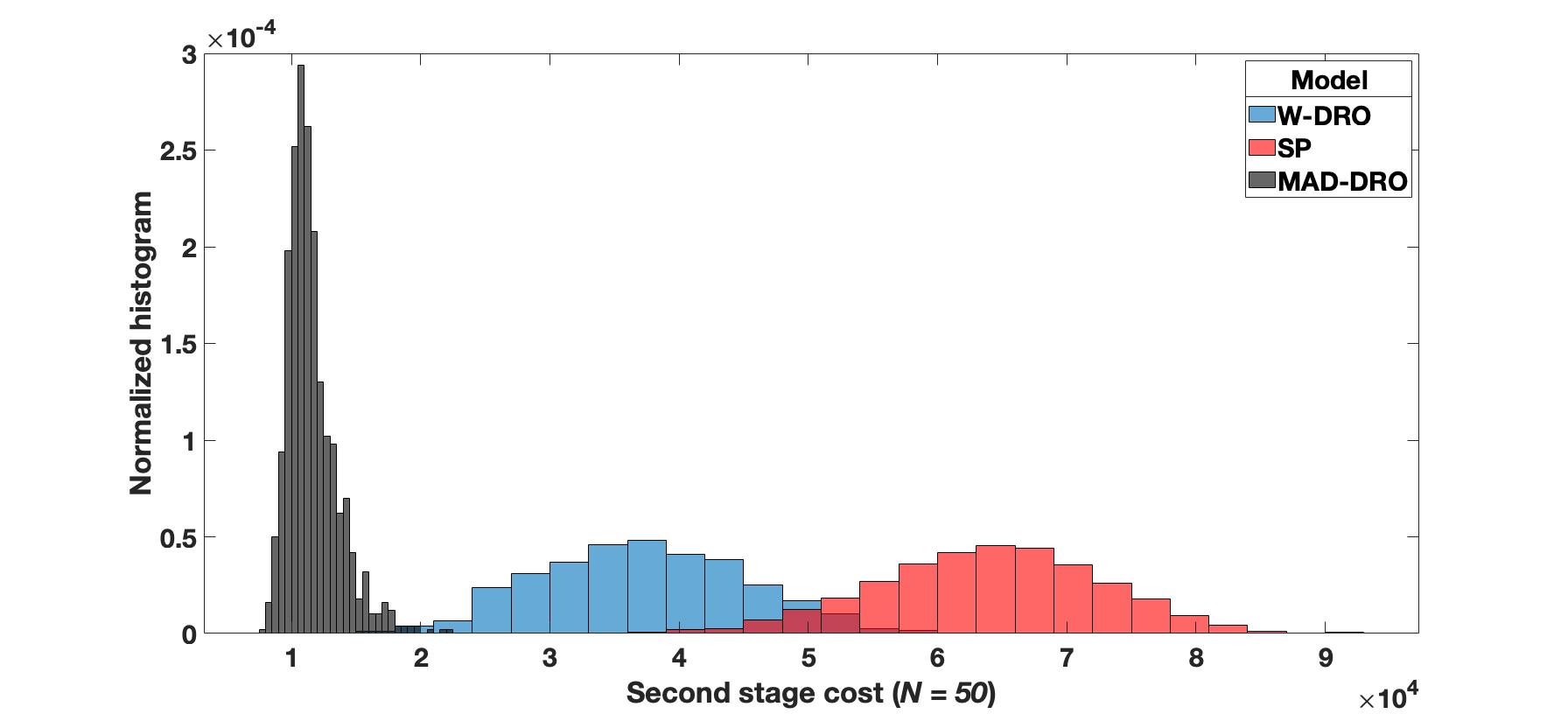}
        \caption{2nd  (Set 2, $ \Delta=0.25$)}
    \end{subfigure}%

      \begin{subfigure}[b]{0.5\textwidth}
          \centering
        \includegraphics[width=\textwidth]{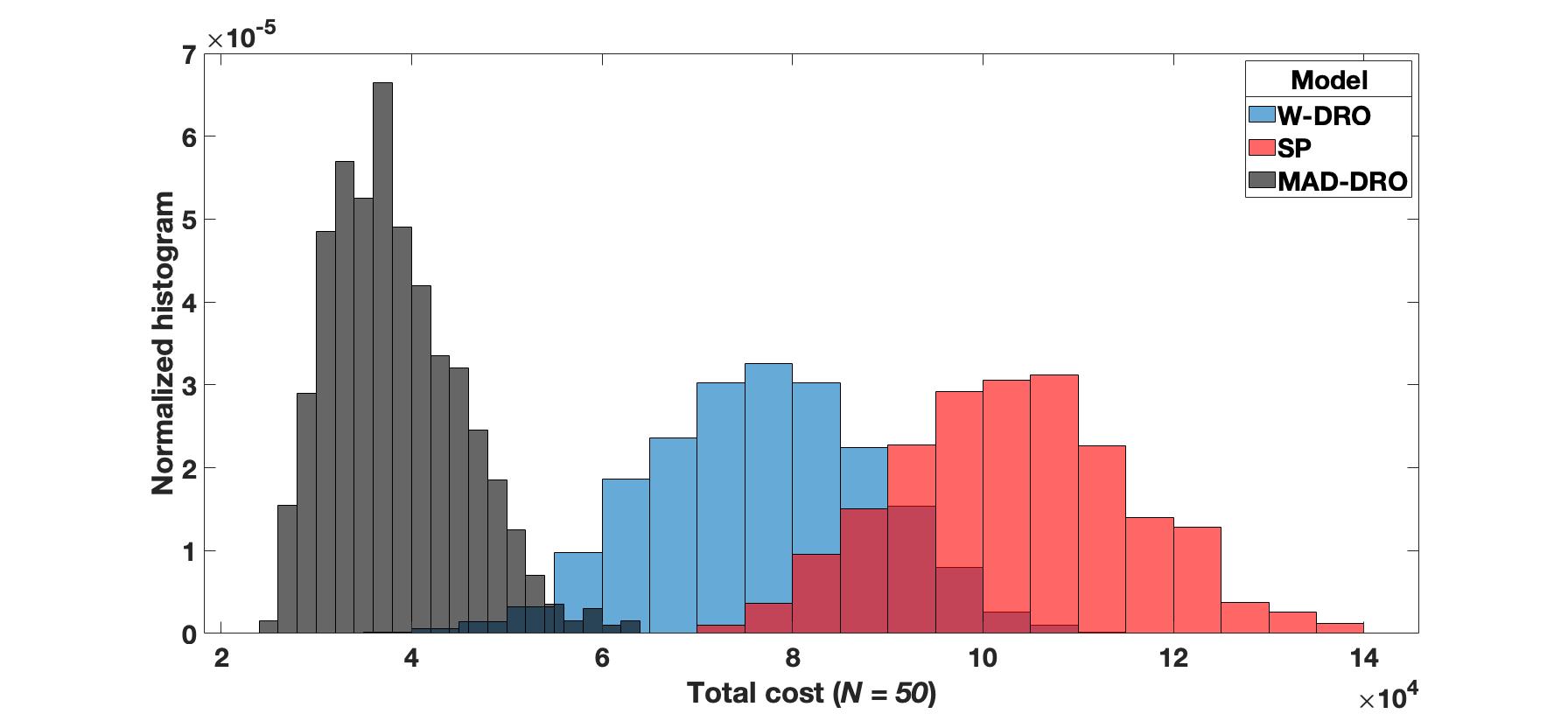}
        \caption{TC  (Set 2, $ \Delta=0.5$)}
    \end{subfigure}%
    \begin{subfigure}[b]{0.5\textwidth}
          \centering
        \includegraphics[width=\textwidth]{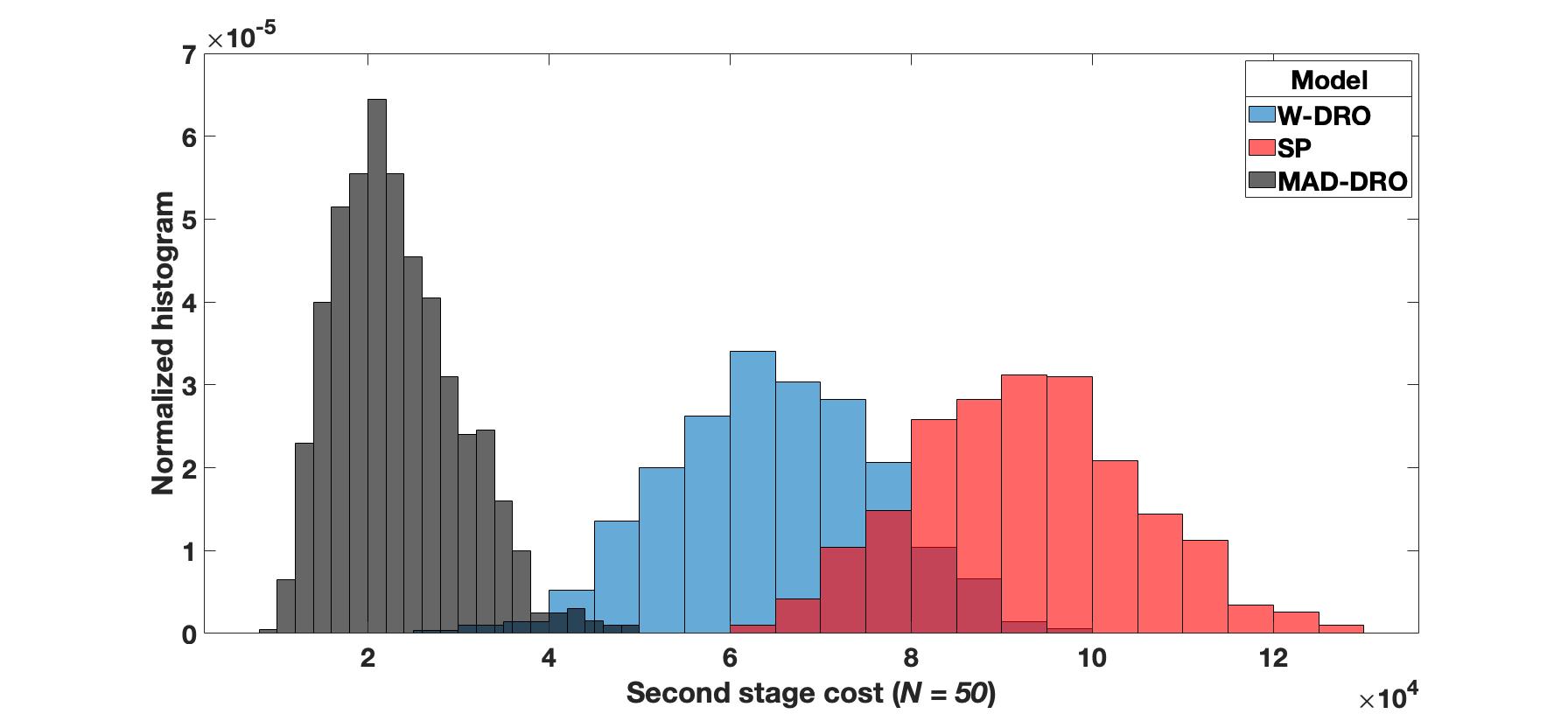}
        \caption{2nd  (Set 2, $ \Delta=0.5$)}
    \end{subfigure}%
\caption{Normalized histograms of out-of-sample TC and 2nd for Lehigh 1 ($\pmb{N=50}$) under Set 1 (LogN) and Set 2 (with $\pmb{\Delta \in \{0, 0.25, 0.5\}}$). }\label{Fig3_N50_Lehigh1}
\end{figure}

\begin{figure}[t!]
 \centering

  \begin{subfigure}[b]{0.5\textwidth}
          \centering
        \includegraphics[width=\textwidth]{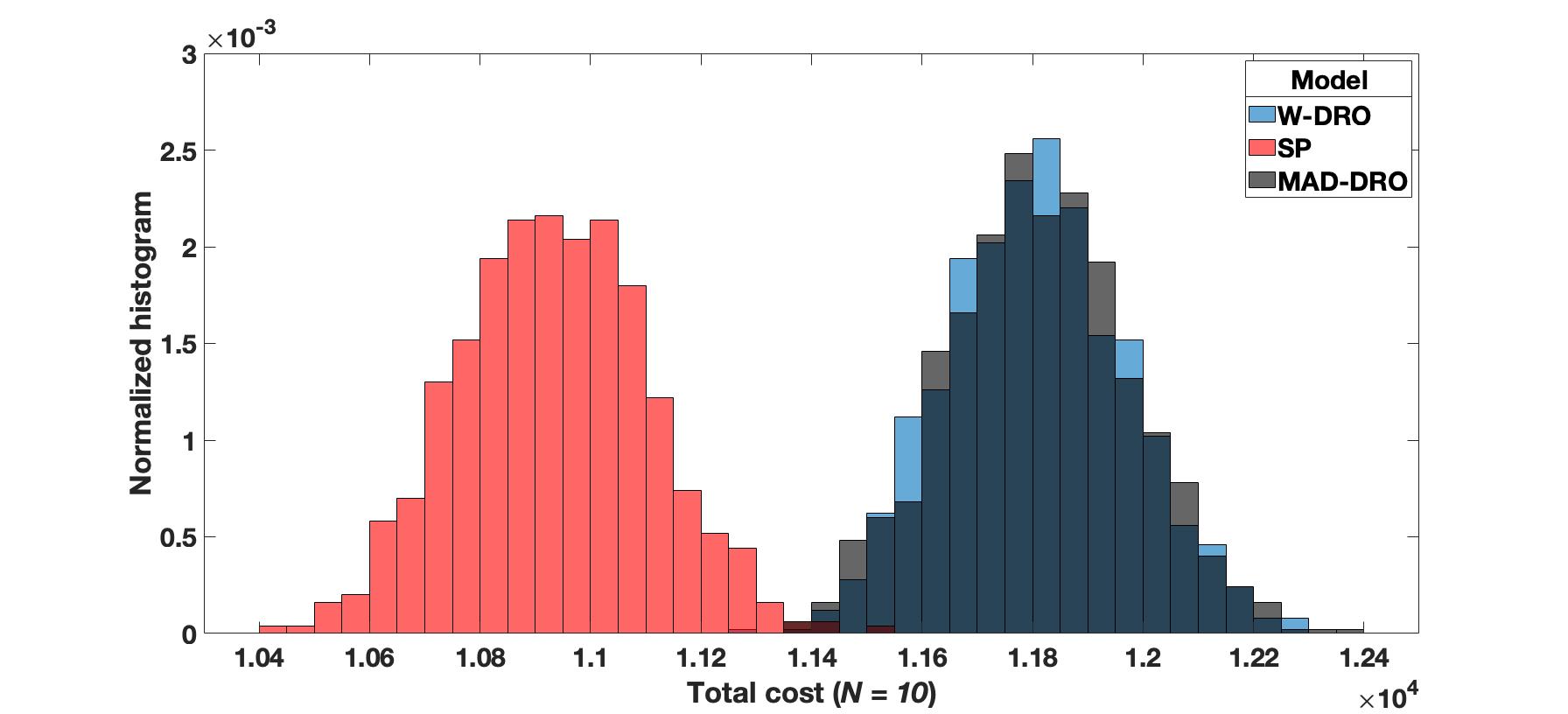}
        \caption{TC (Set 1, LogN)}
    \end{subfigure}%
      \begin{subfigure}[b]{0.5\textwidth}
          \centering
        \includegraphics[width=\textwidth]{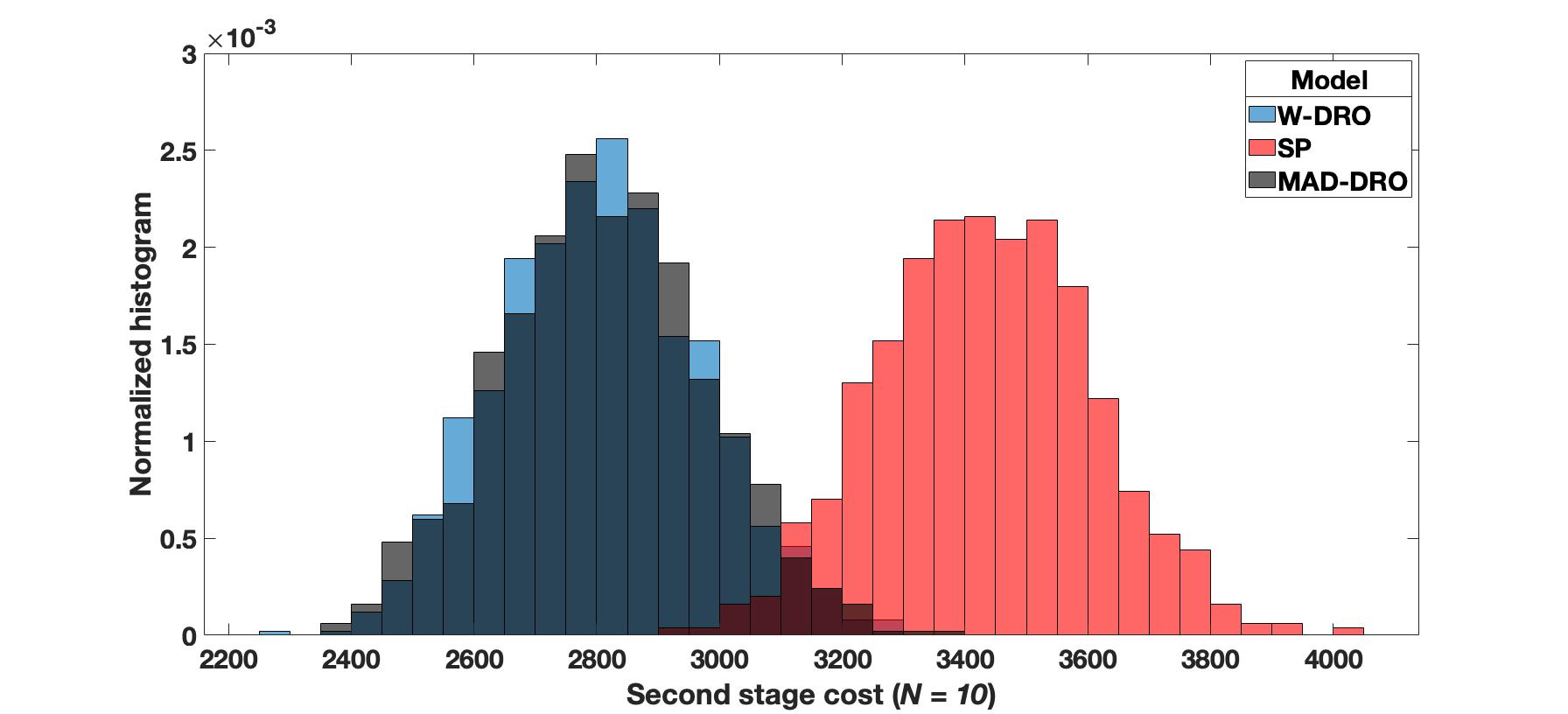}
        \caption{2nd (Set 1, LogN))}
    \end{subfigure}%
 
  \begin{subfigure}[b]{0.5\textwidth}
          \centering
        \includegraphics[width=\textwidth]{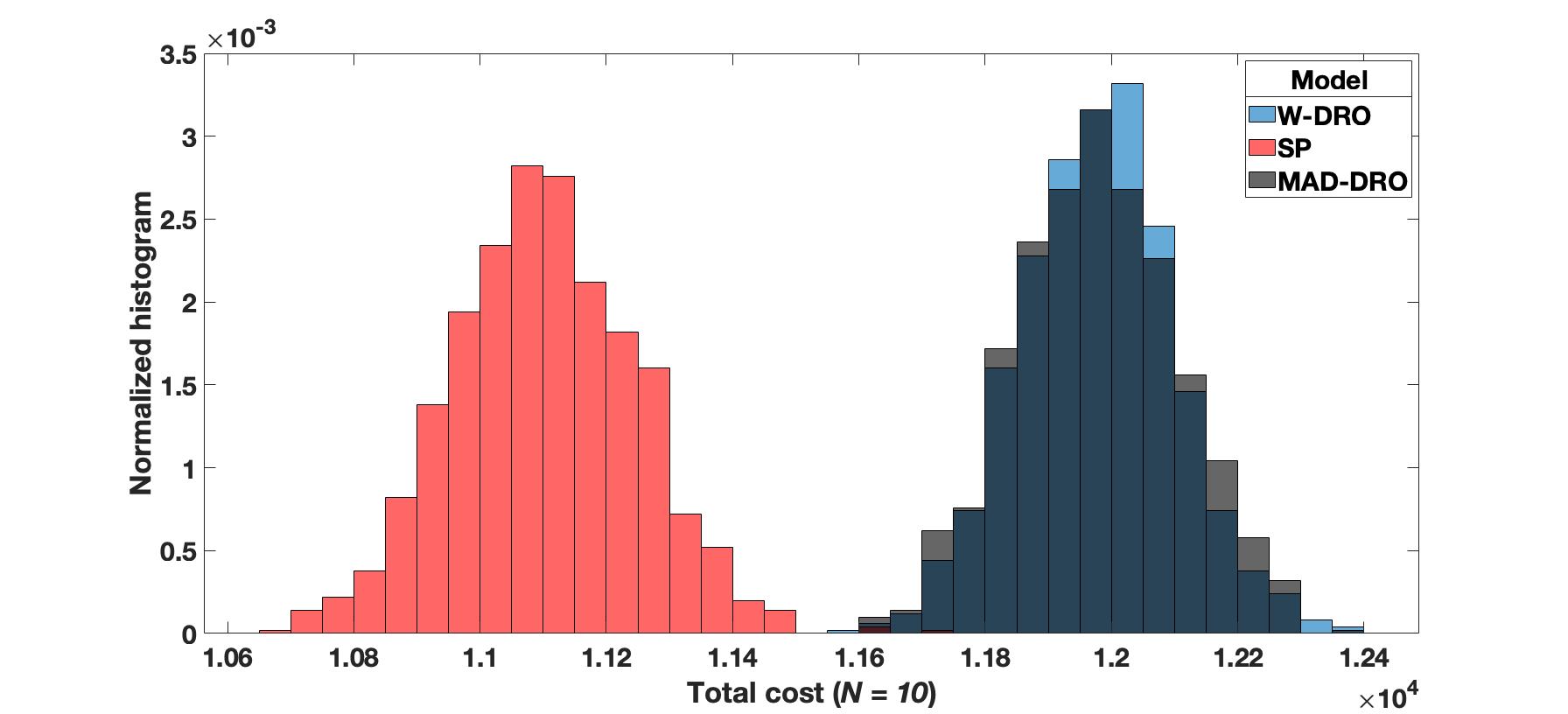}
        \caption{TC (Set 2, $\Delta=0$)}
    \end{subfigure}%
      \begin{subfigure}[b]{0.5\textwidth}
          \centering
        \includegraphics[width=\textwidth]{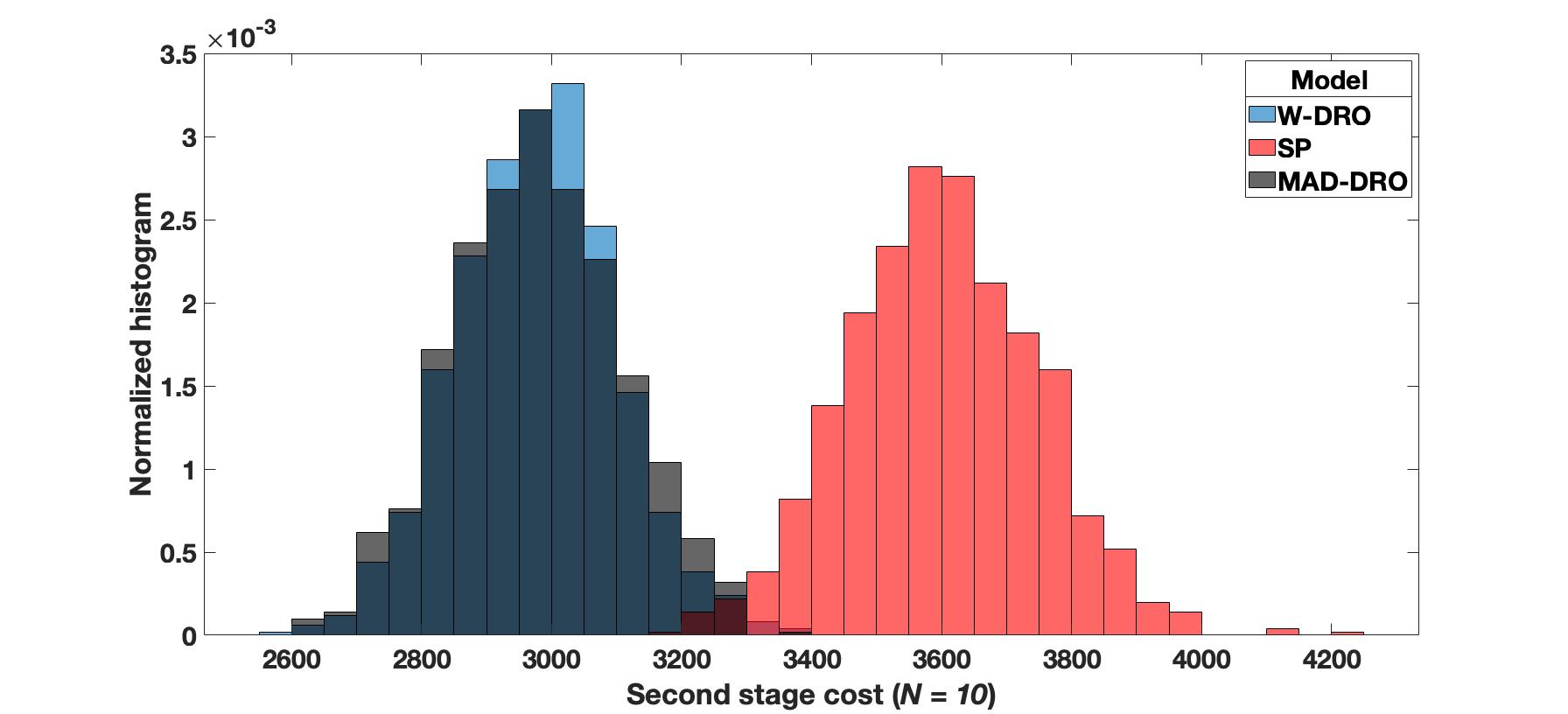}
        \caption{2nd (Set 2, $ \Delta=0$)}
    \end{subfigure}%

  \begin{subfigure}[b]{0.5\textwidth}
        \centering
        \includegraphics[width=\textwidth]{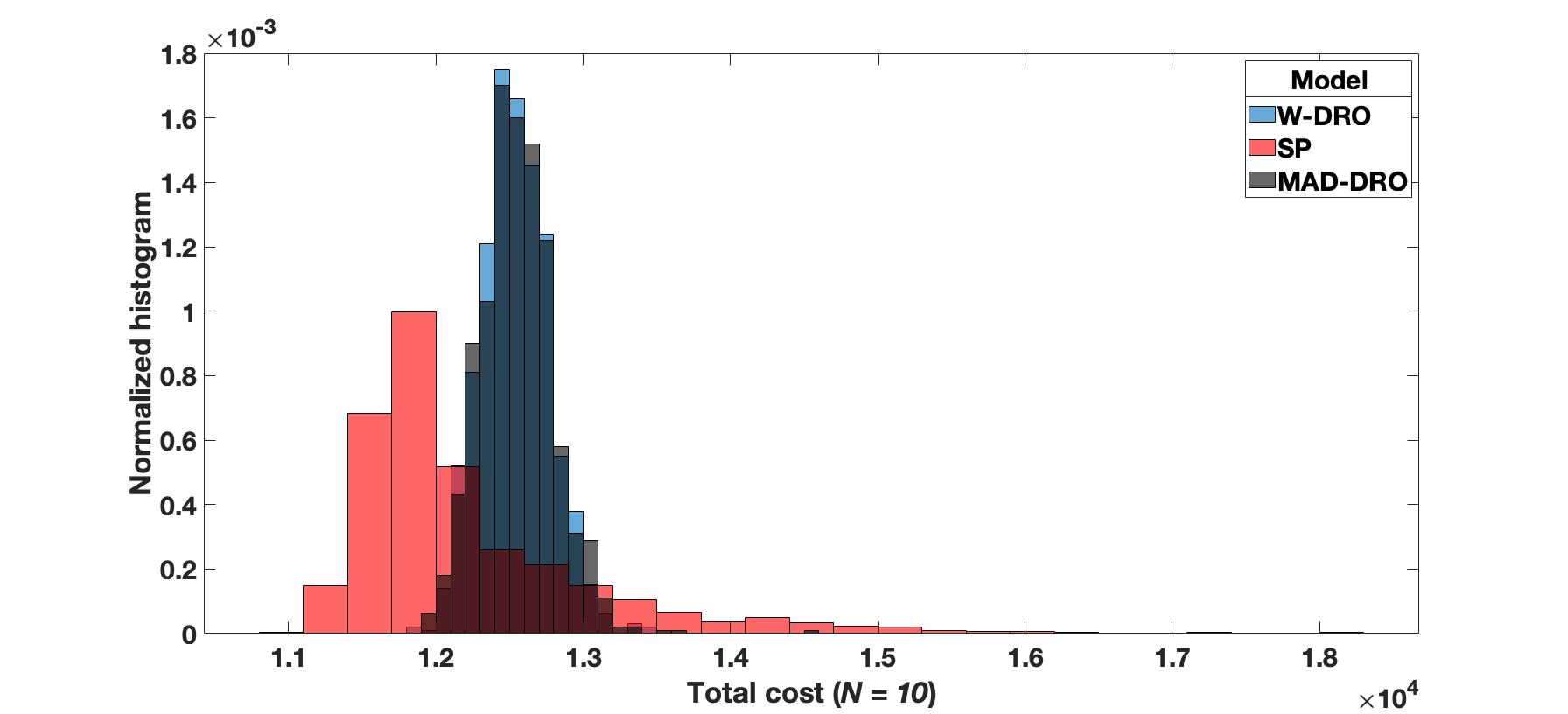}
        \caption{TC (Set 2, $ \Delta=0.25$)}
    \end{subfigure}%
      \begin{subfigure}[b]{0.5\textwidth}
        \centering
        \includegraphics[width=\textwidth]{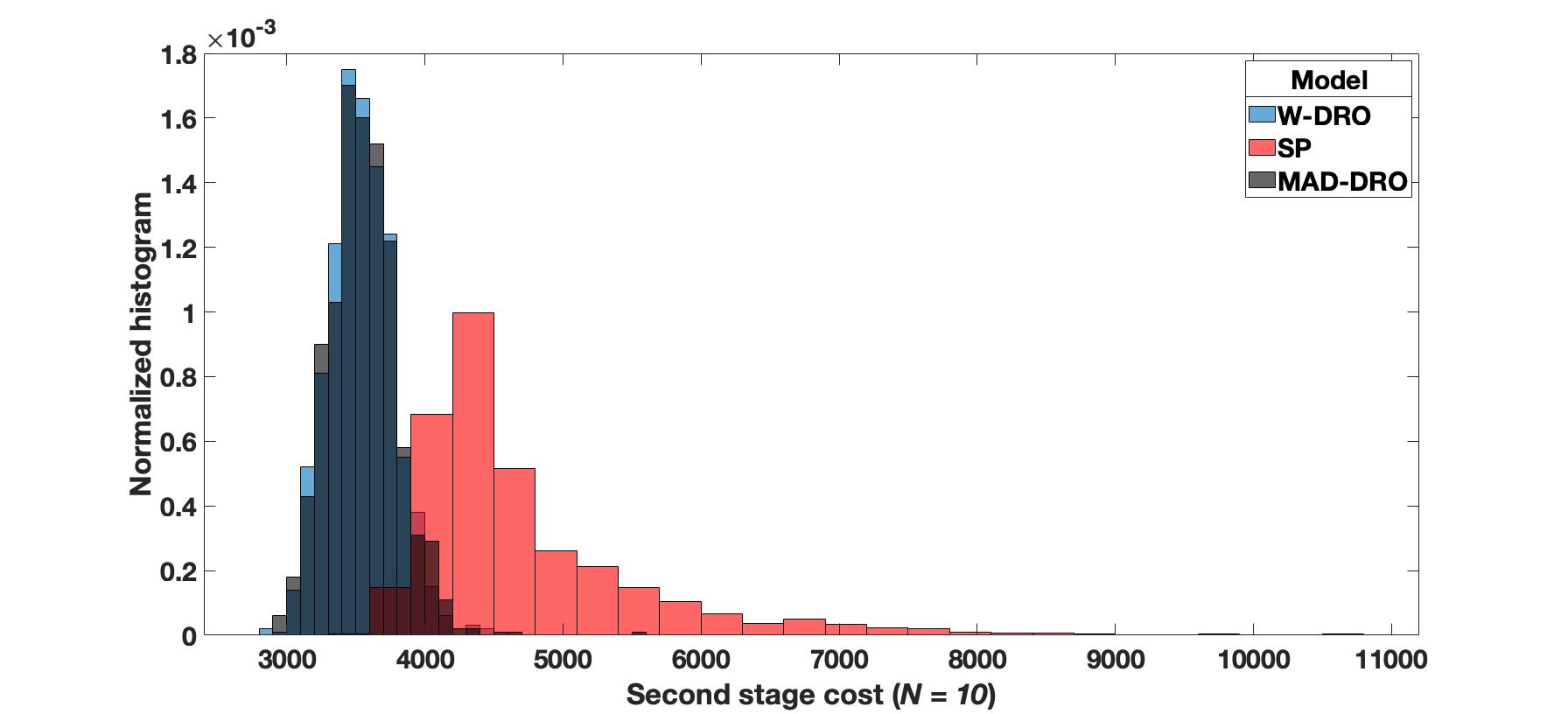}
        \caption{2nd (Set 2, $ \Delta=0.25$)}
    \end{subfigure}%

      \begin{subfigure}[b]{0.5\textwidth}
          \centering
        \includegraphics[width=\textwidth]{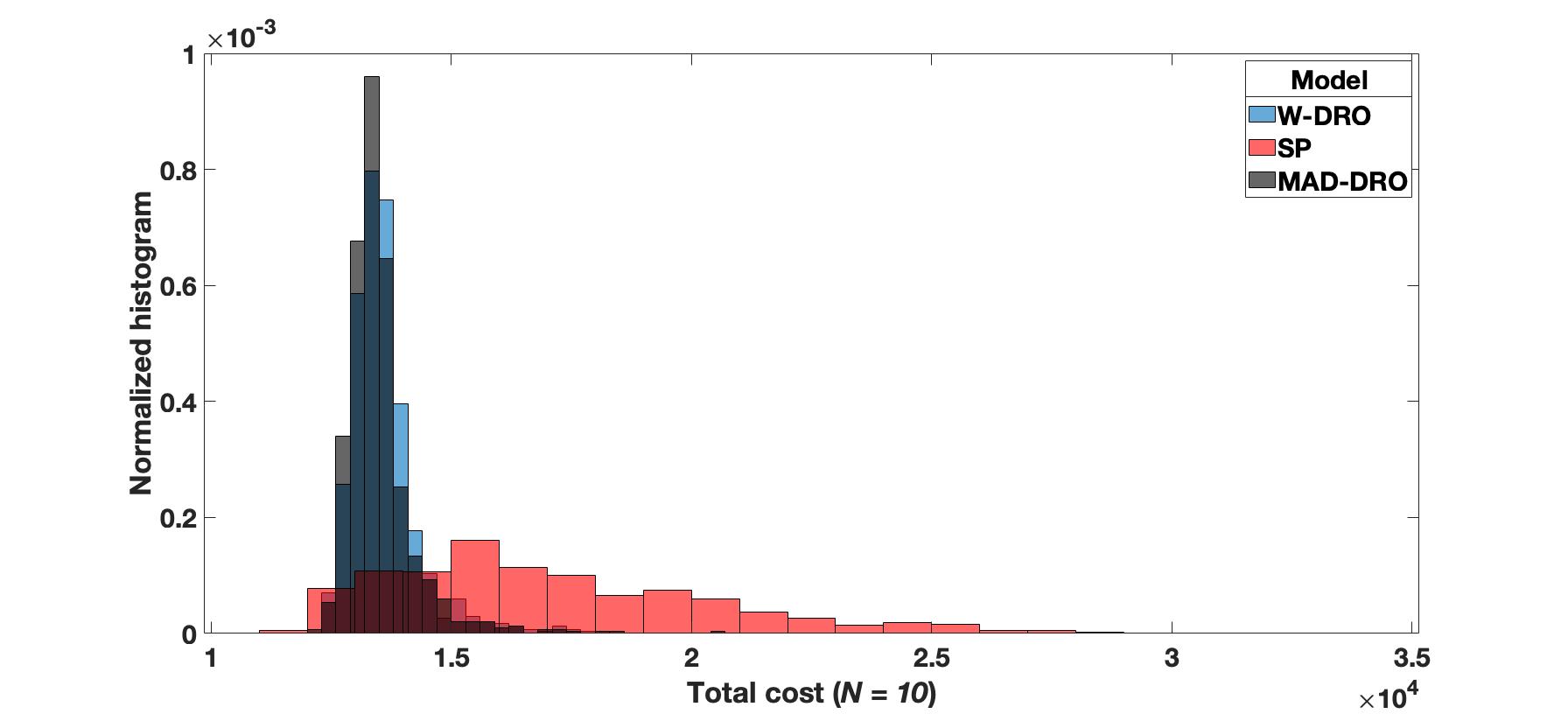}
        \caption{TC (Set 2, $ \Delta=0.5$)}
    \end{subfigure}%
    \begin{subfigure}[b]{0.5\textwidth}
          \centering
        \includegraphics[width=\textwidth]{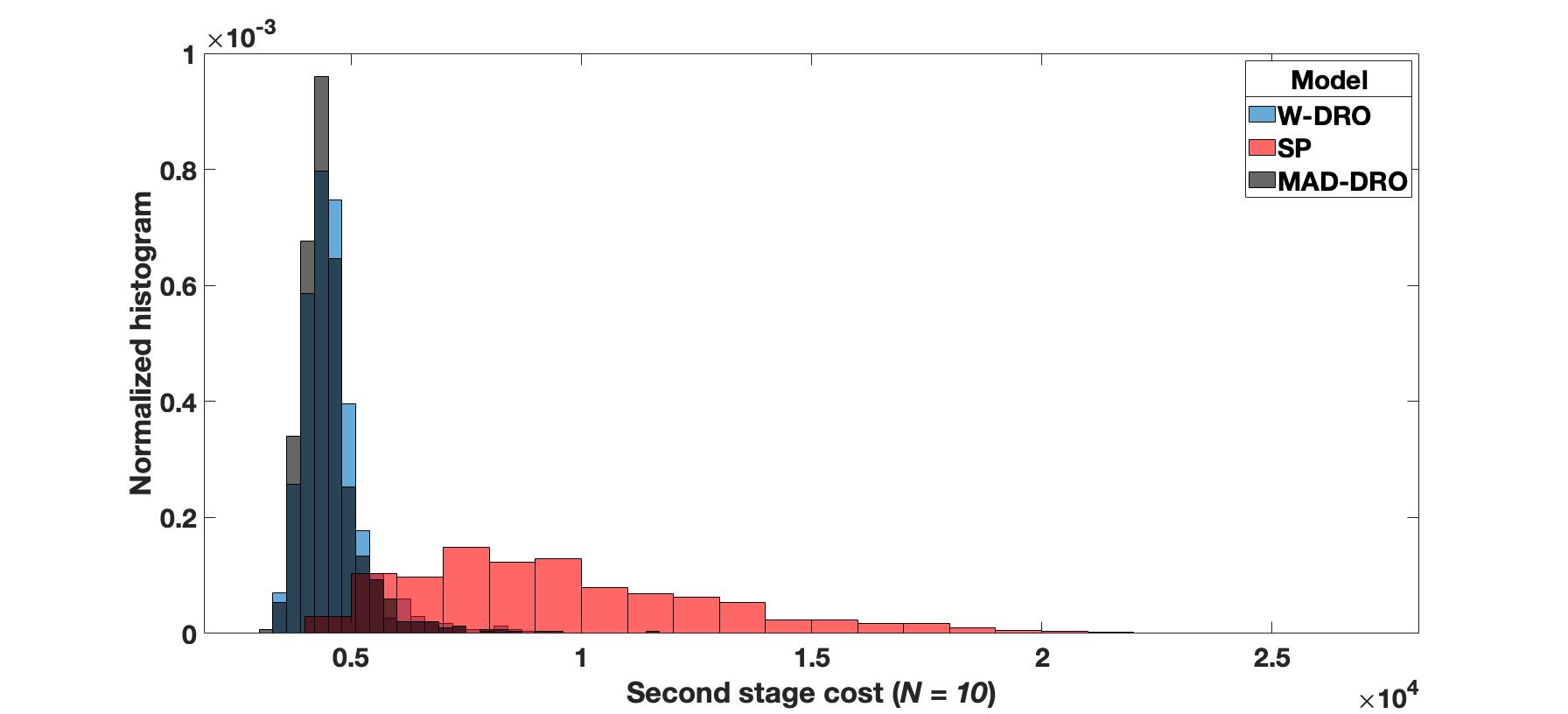}
        \caption{2nd (Set 2, $\Delta=0.5$)}
    \end{subfigure}%
\caption{Normalized histograms of out-of-sample TC and 2nd for Lehigh 2 ($N=10$) under Set 1 (LogN) and Set 2 (with $\pmb{\Delta \in \{0, 0.25, 0.5\}}$).}\label{Fig3_Uni10_Lehigh2}
\end{figure}

\begin{figure}[t!]
 \centering

  \begin{subfigure}[b]{0.5\textwidth}
          \centering
        \includegraphics[width=\textwidth]{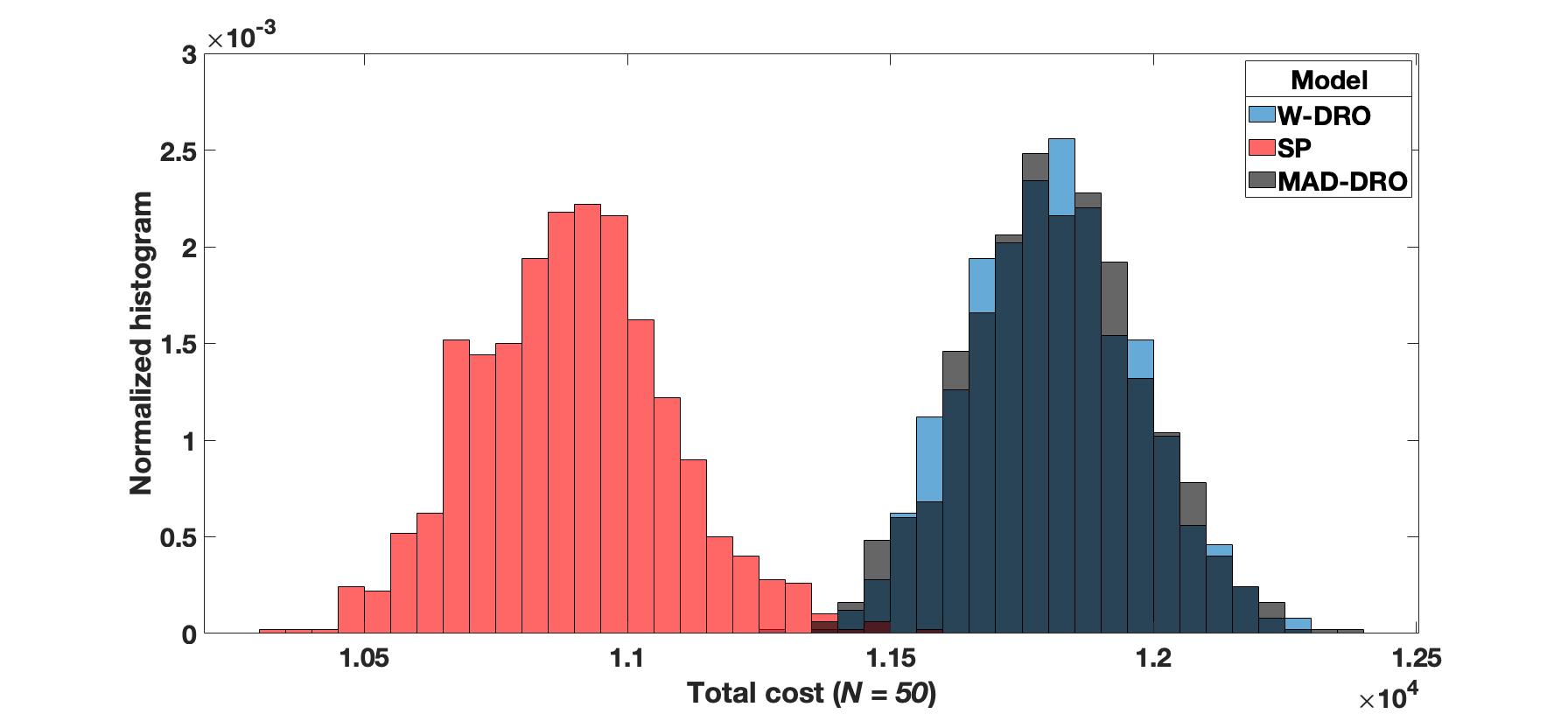}
        \caption{TC (Set 1, LogN)}
    \end{subfigure}%
      \begin{subfigure}[b]{0.5\textwidth}
          \centering
        \includegraphics[width=\textwidth]{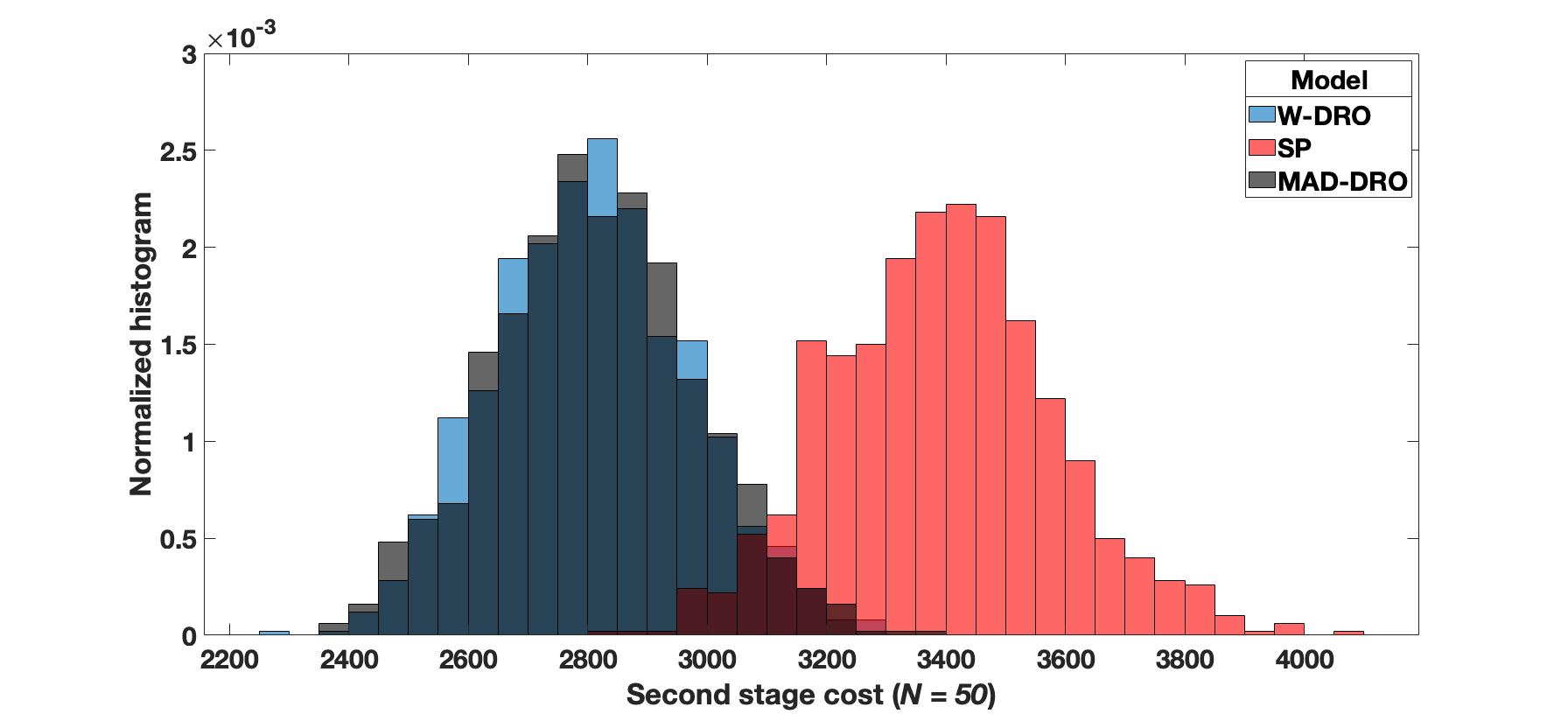}
        \caption{2nd (Set 1, LogN)}
    \end{subfigure}%
 
  \begin{subfigure}[b]{0.5\textwidth}
          \centering
        \includegraphics[width=\textwidth]{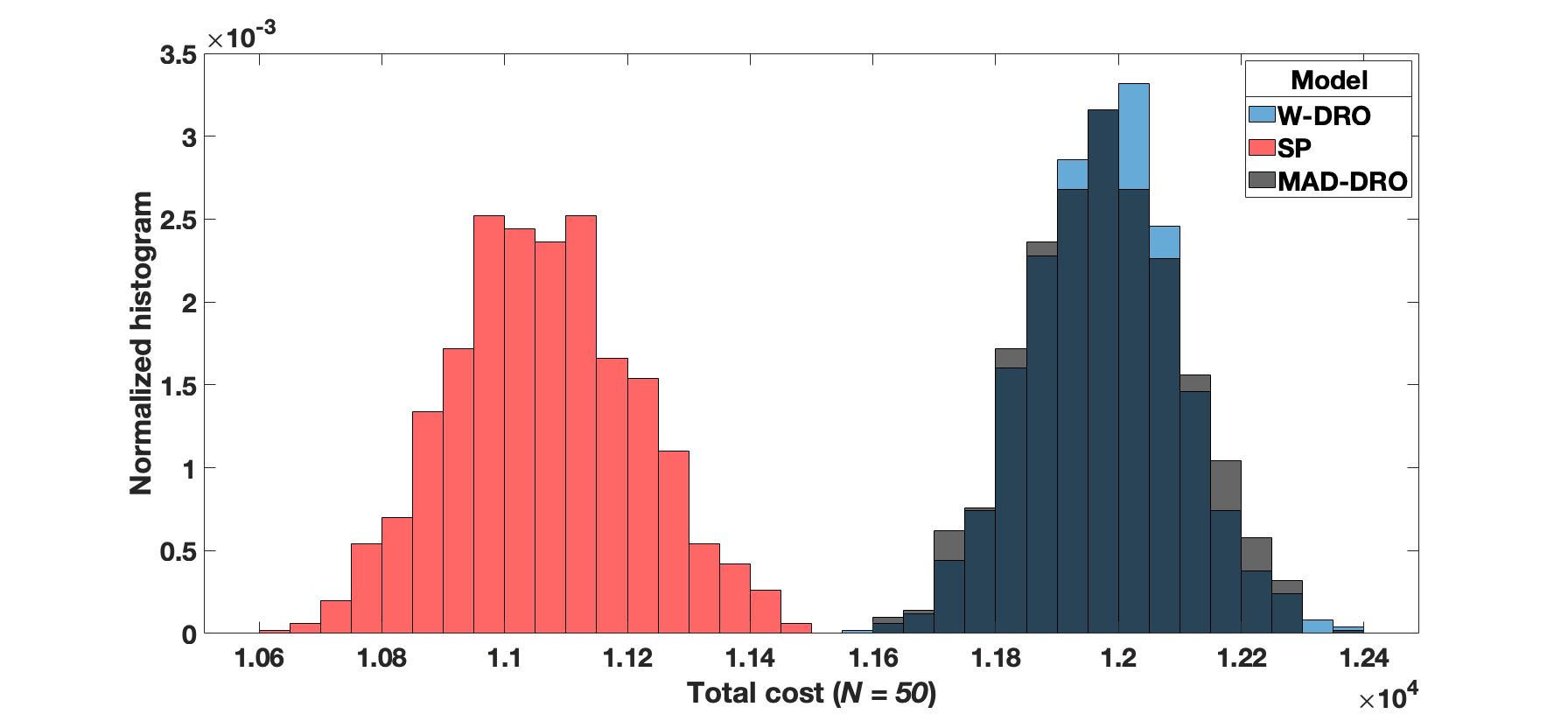}
        \caption{TC (Set 2, $ \Delta=0$)}
    \end{subfigure}%
      \begin{subfigure}[b]{0.5\textwidth}
          \centering
        \includegraphics[width=\textwidth]{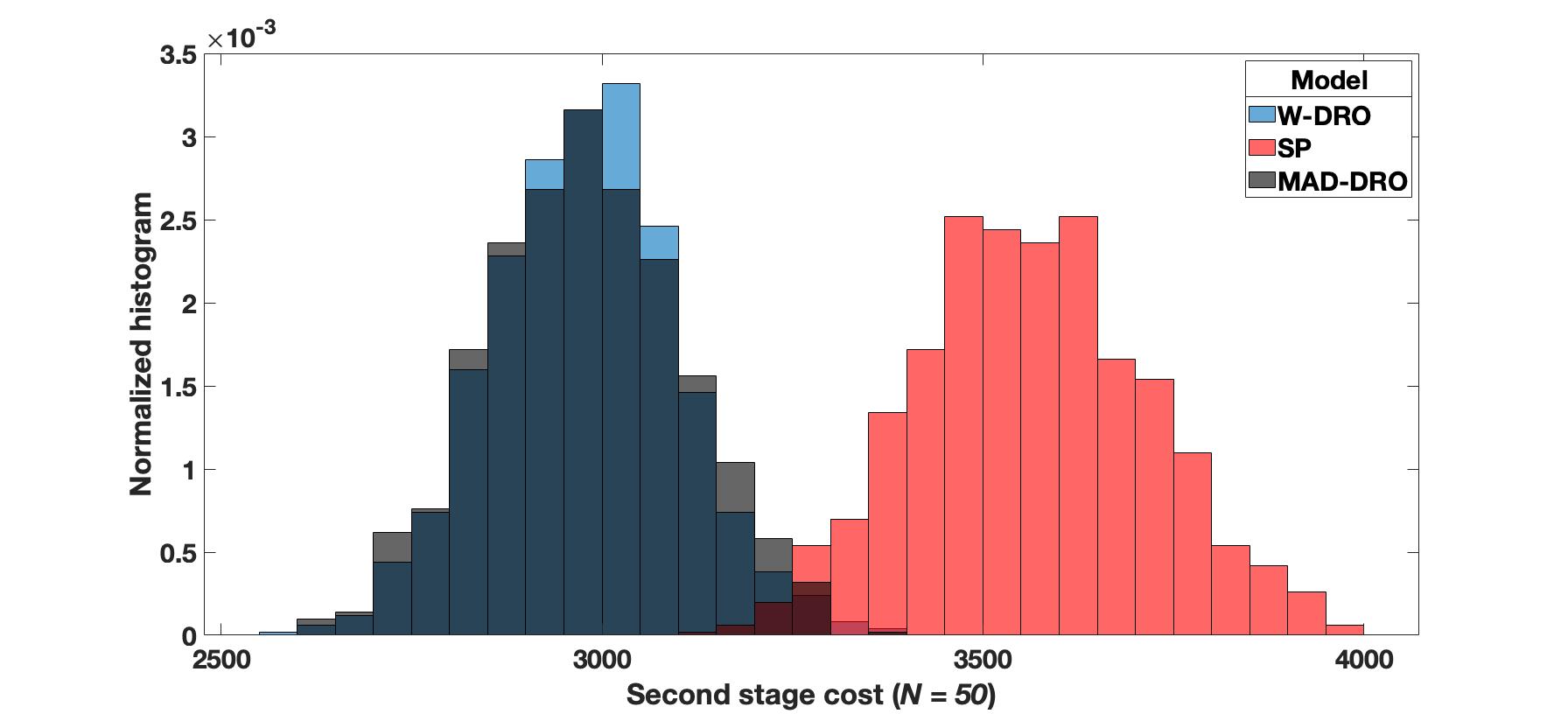}
        \caption{2nd (Set 2, $ \Delta=0$)}
    \end{subfigure}%

  \begin{subfigure}[b]{0.5\textwidth}
        \centering
        \includegraphics[width=\textwidth]{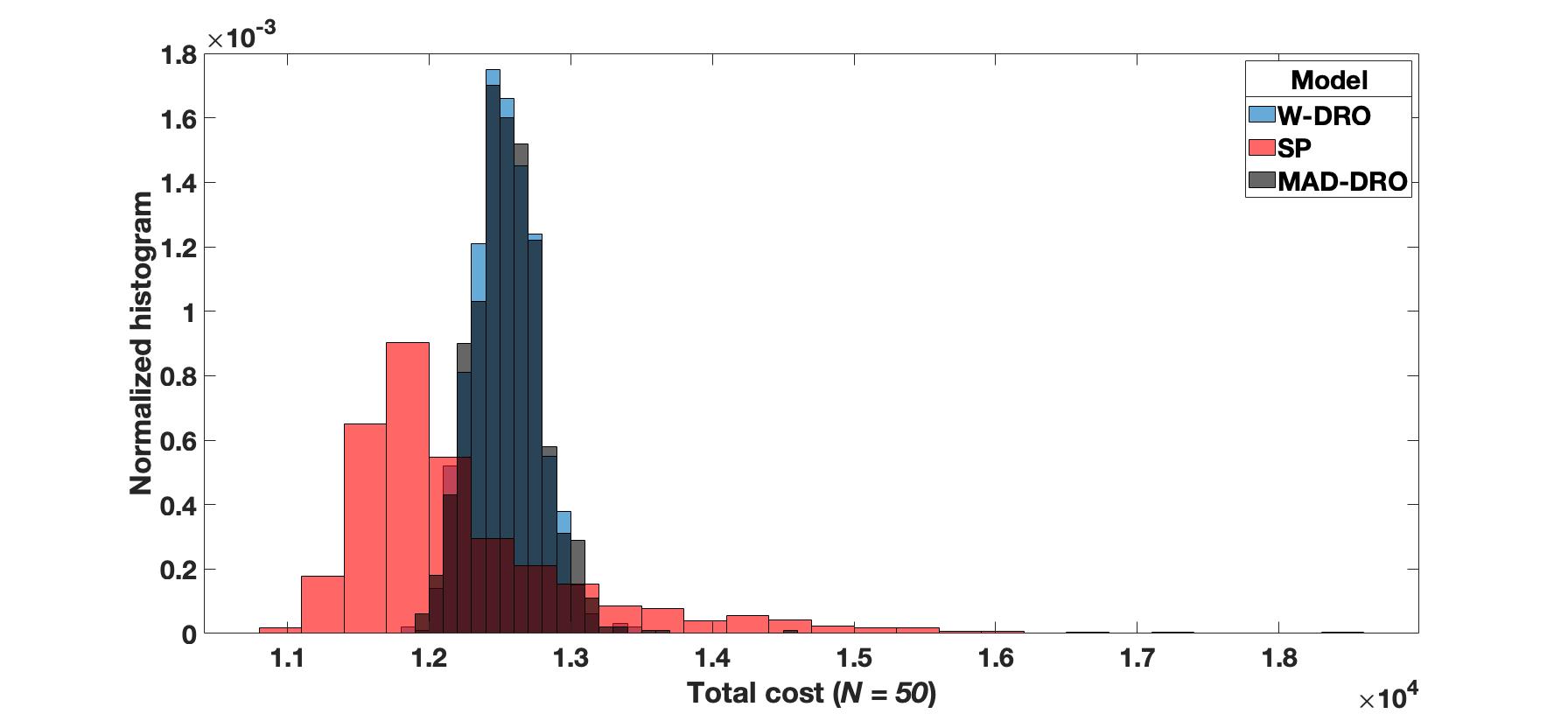}
        \caption{TC (Set 2, $\Delta=0.25$)}
    \end{subfigure}%
      \begin{subfigure}[b]{0.5\textwidth}
        \centering
        \includegraphics[width=\textwidth]{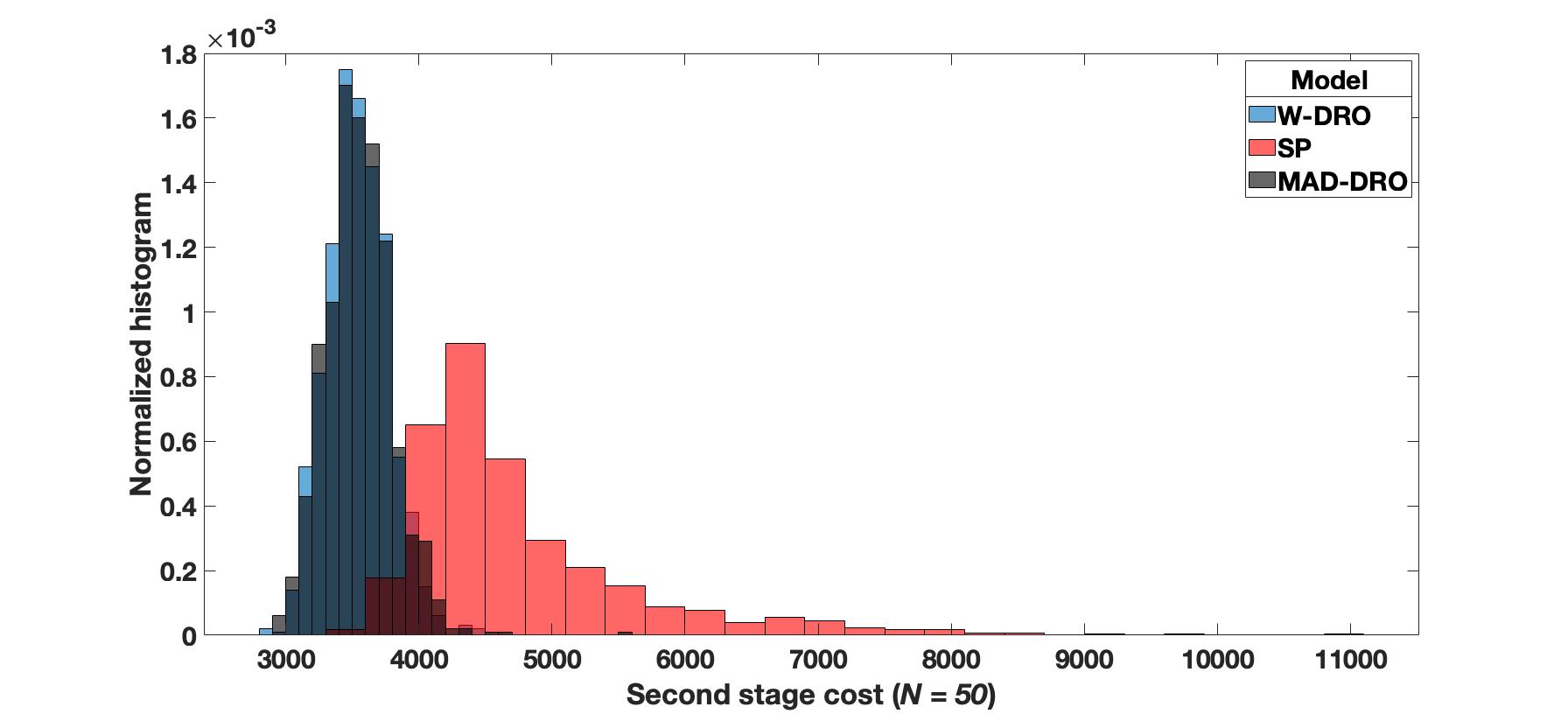}
        \caption{2nd (Set 2, $ \Delta=0.25$)}
    \end{subfigure}%

      \begin{subfigure}[b]{0.5\textwidth}
          \centering
        \includegraphics[width=\textwidth]{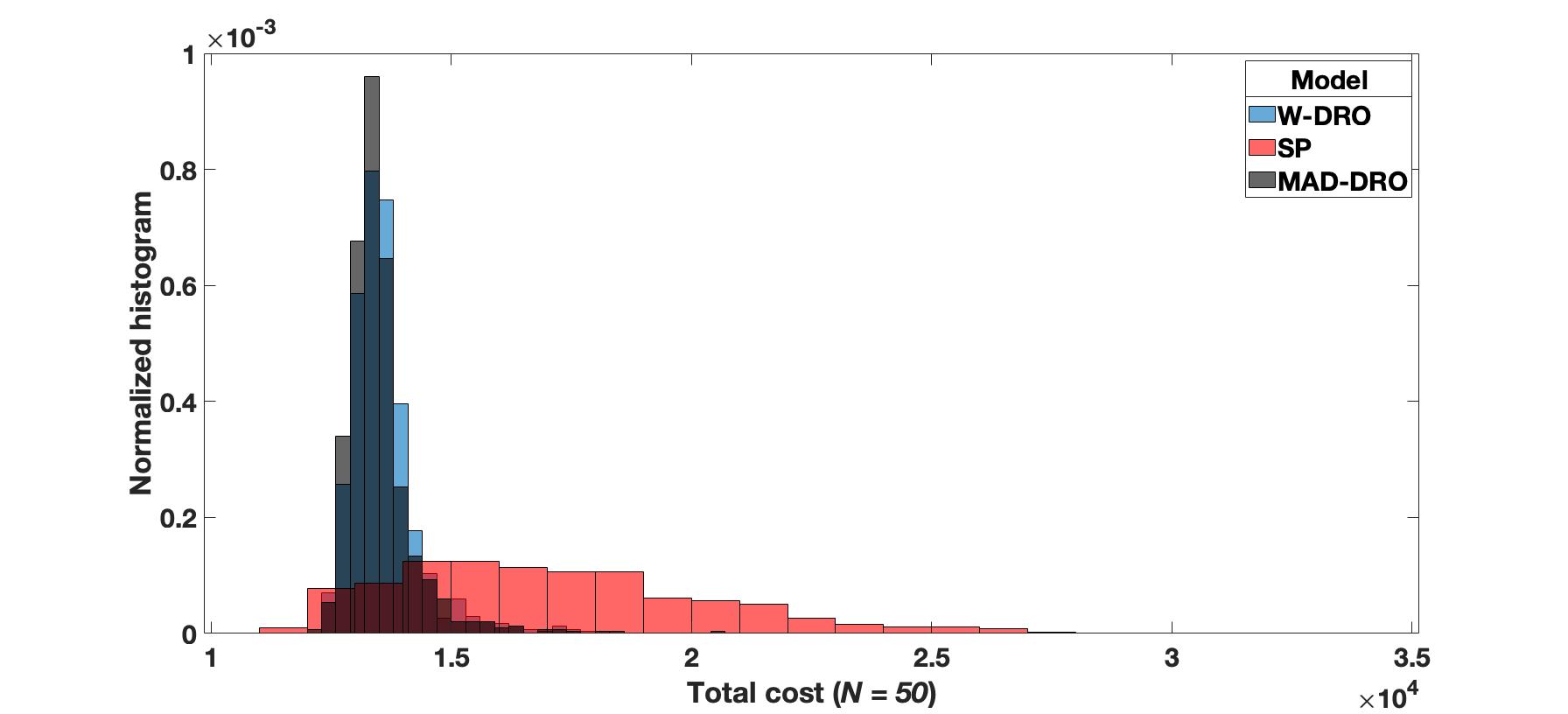}
        \caption{TC (Set 2, $\Delta=0.5$)}
    \end{subfigure}%
    \begin{subfigure}[b]{0.5\textwidth}
          \centering
        \includegraphics[width=\textwidth]{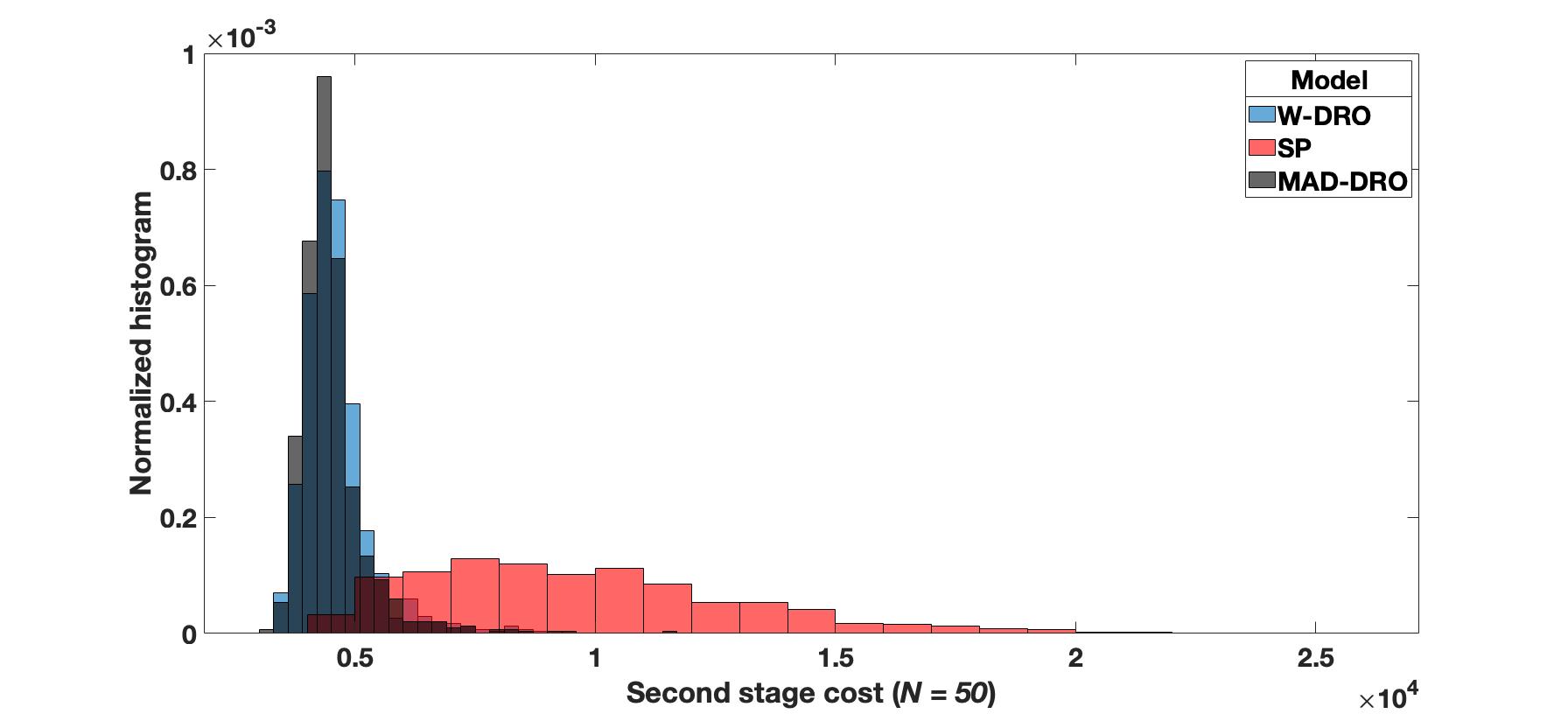}
        \caption{2nd (Set 2, $\Delta=0.5$)}
    \end{subfigure}%
\caption{Normalized histograms of out-of-sample TC and 2nd for Lehigh 2 ($N=50$) under Set 1 (LogN) and Set 2 (with $\pmb{\Delta \in \{0, 0.25, 0.5\}}$).}\label{Fig3_Uni50_Lehigh2}
\end{figure}

\begin{figure}[h!]
 \centering
  \begin{subfigure}[b]{0.5\textwidth}
          \centering
        \includegraphics[width=\textwidth]{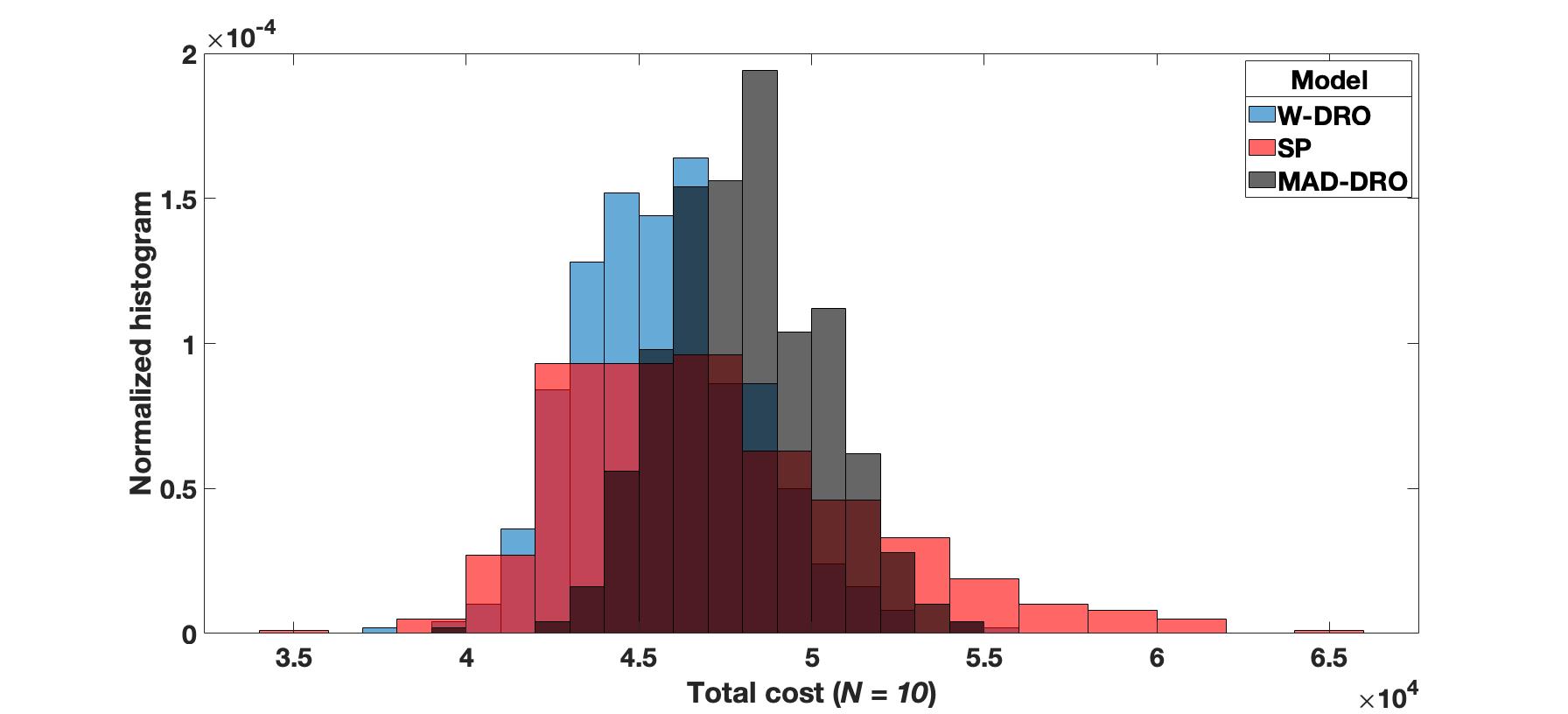}
        \caption{TC (cor=0.2)}
    \end{subfigure}%
  \begin{subfigure}[b]{0.5\textwidth}
        \centering
        \includegraphics[width=\textwidth]{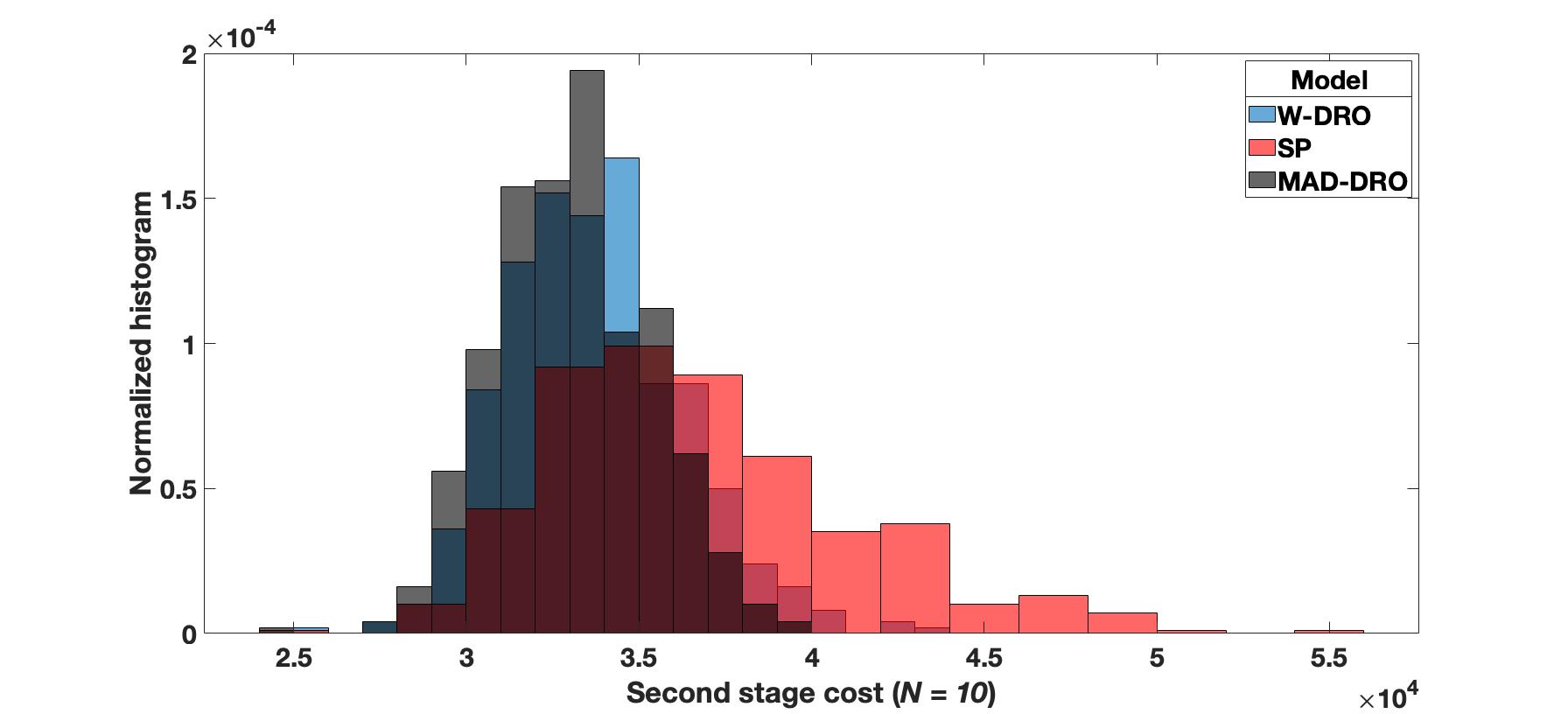}
        \caption{2nd (cor=0.2)}
    \end{subfigure}%
    
  \begin{subfigure}[b]{0.5\textwidth}
          \centering
        \includegraphics[width=\textwidth]{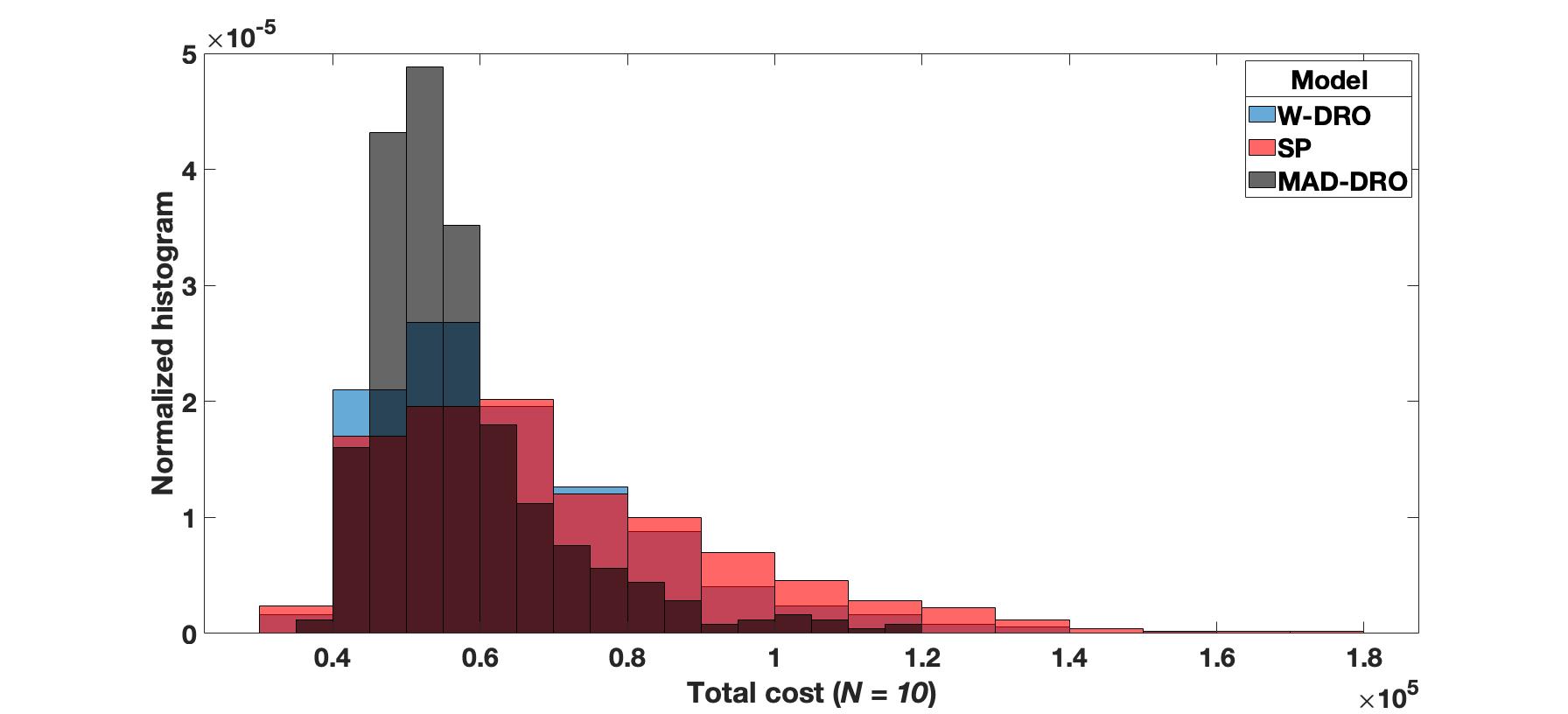}
        \caption{TC (cor=0.6) }
    \end{subfigure}%
  \begin{subfigure}[b]{0.5\textwidth}
        \centering
        \includegraphics[width=\textwidth]{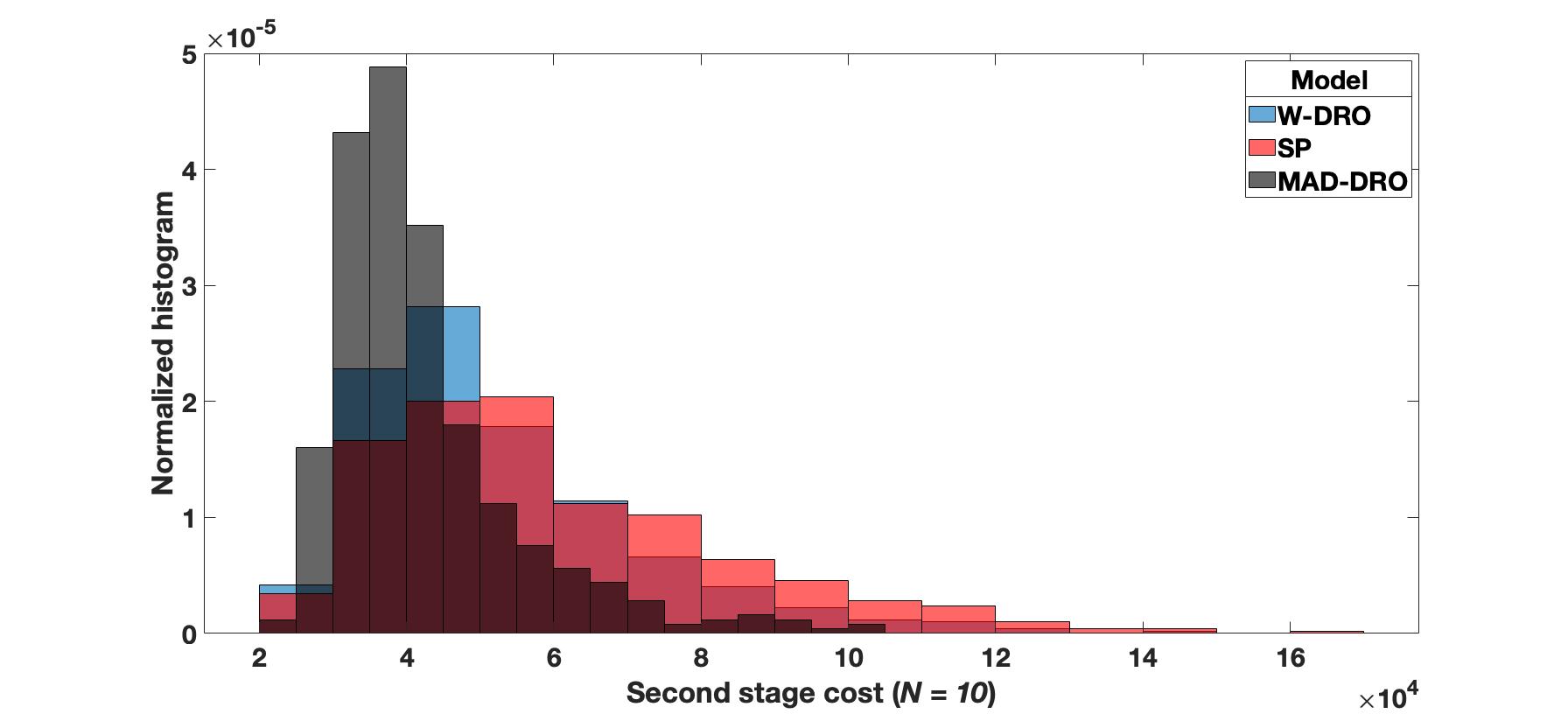}
        \caption{2nd (cor=0.6)}
    \end{subfigure}%

\caption{Out-of-sample performance under correlated data for Instance 3. Notation: cor is correlation coefficient}\label{Fig_outcorr}
\end{figure}

\clearpage
\newpage
\section{Additional Sensitivity Results}\label{Appx:sensitivity}

\begin{figure}[h!]
     \begin{subfigure}[b]{0.5\textwidth}
 \centering
        \includegraphics[width=\textwidth]{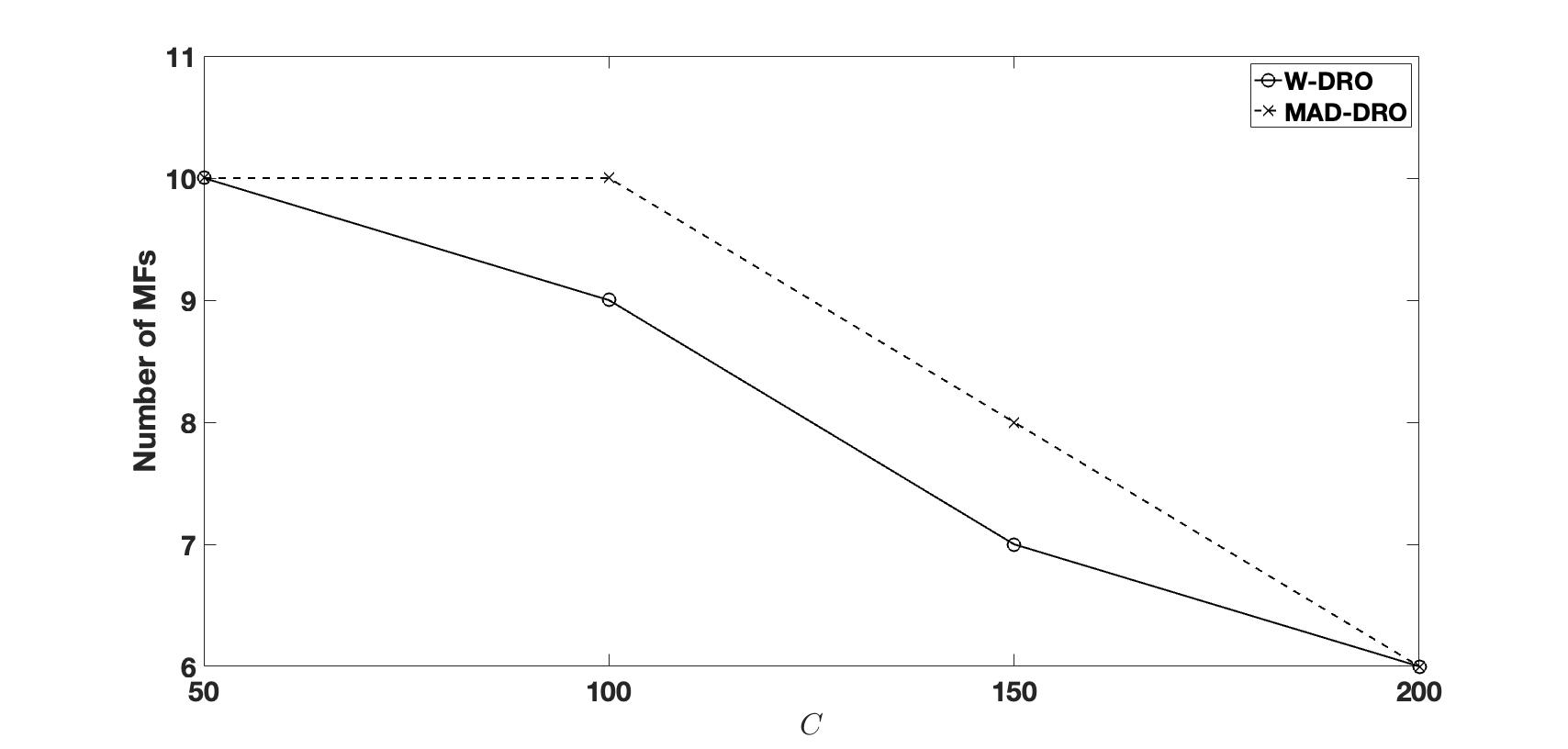}
        \caption{Number of MFs, $f=1,500$}
        \label{Fig7a}
    \end{subfigure}%
    \begin{subfigure}[b]{0.5\textwidth}
            \includegraphics[width=\textwidth]{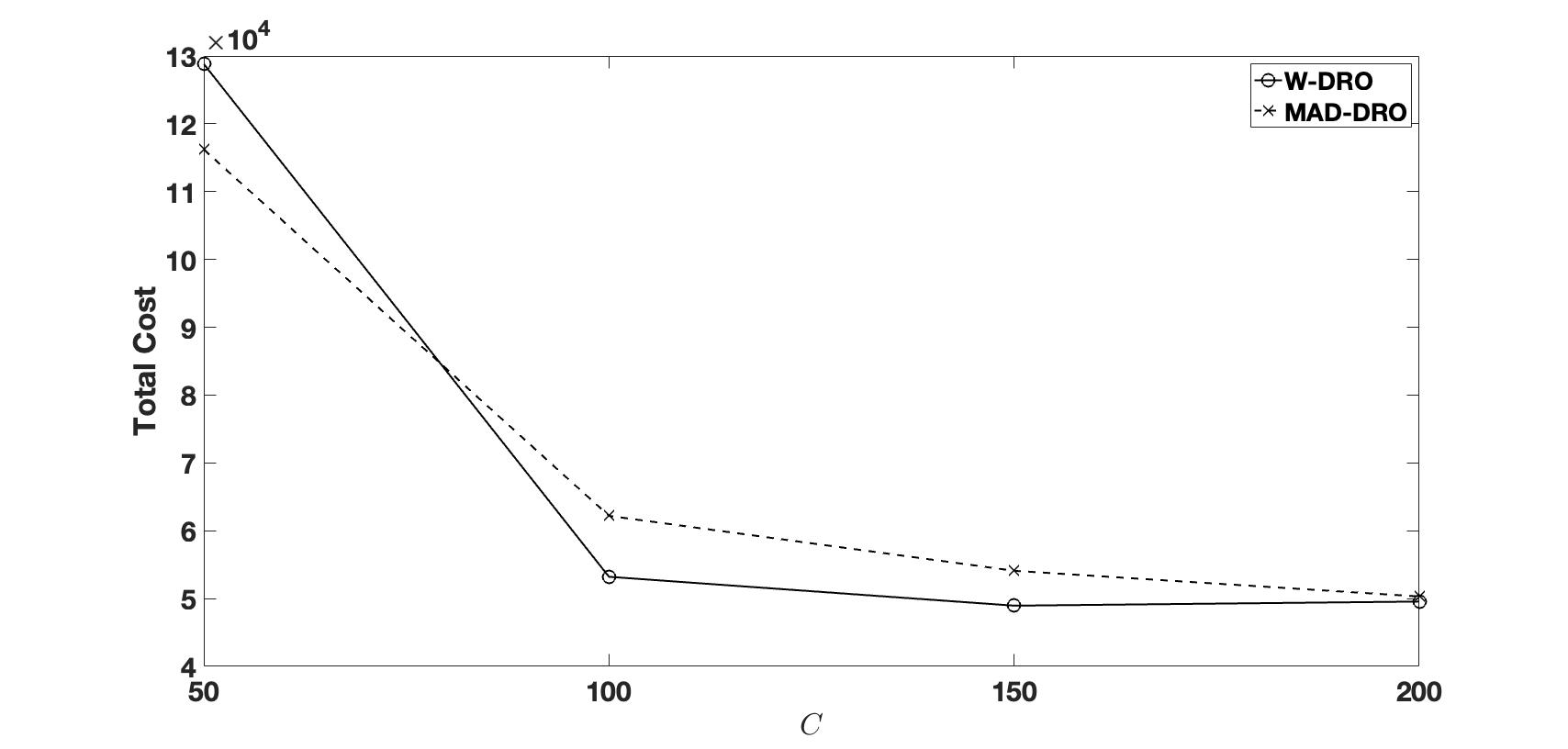}
      \caption{Total cost, $f=1,500$}
      \label{Fig7b}
    \end{subfigure}%
    
        \begin{subfigure}[b]{0.5\textwidth}
 \centering
        \includegraphics[width=\textwidth]{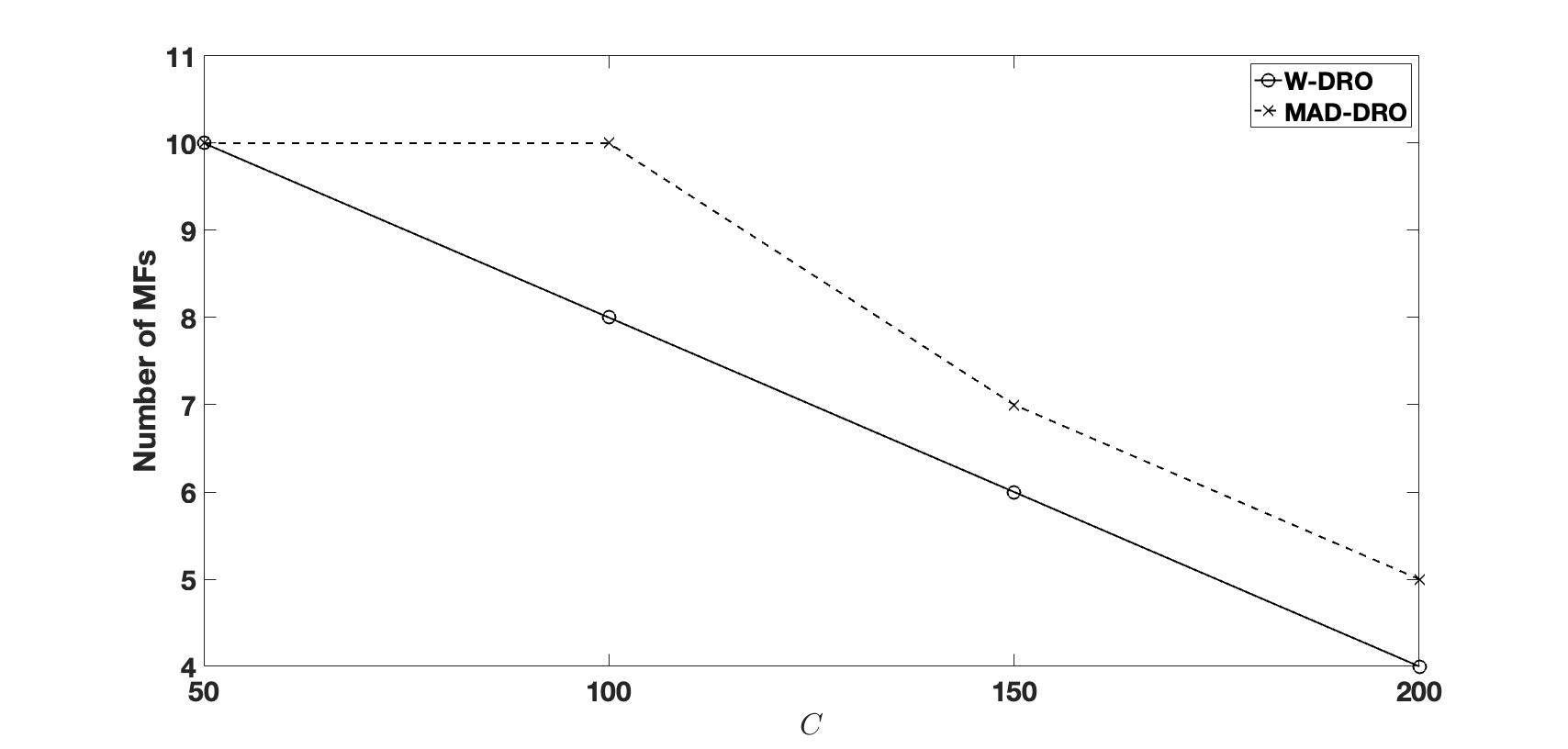}
        \caption{Number of MFs, $f=6,000$}
        \label{Fig7c}
    \end{subfigure}%
    \begin{subfigure}[b]{0.5\textwidth}
            \includegraphics[width=\textwidth]{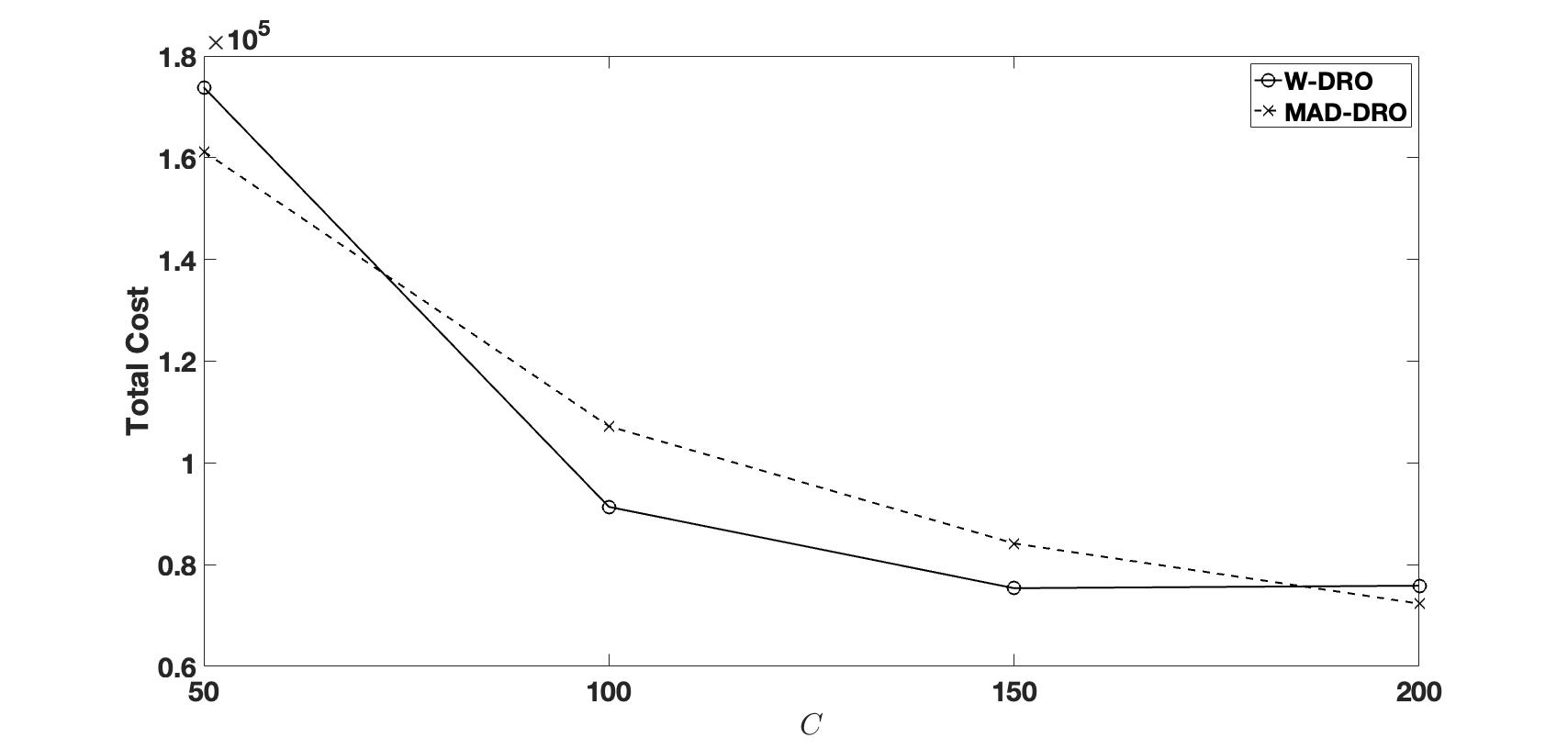}
      \caption{Total cost, $f=6,000$}
      \label{Fig7d}
    \end{subfigure}%
    
            \begin{subfigure}[b]{0.5\textwidth}
 \centering
        \includegraphics[width=\textwidth]{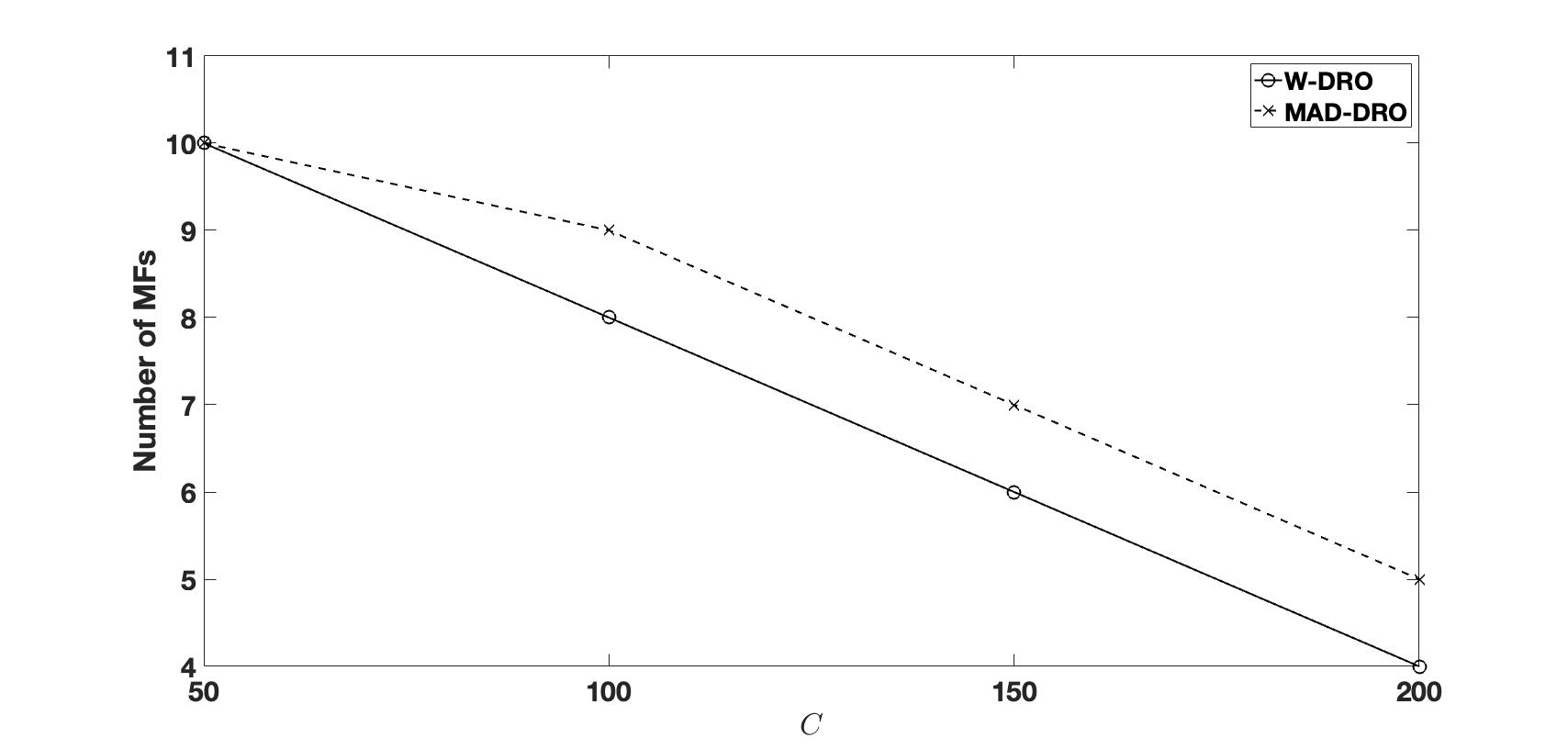}
        \caption{Number of MFs, $f=10,000$}
        \label{Fig7e}
    \end{subfigure}%
    \begin{subfigure}[b]{0.5\textwidth}
            \includegraphics[width=\textwidth]{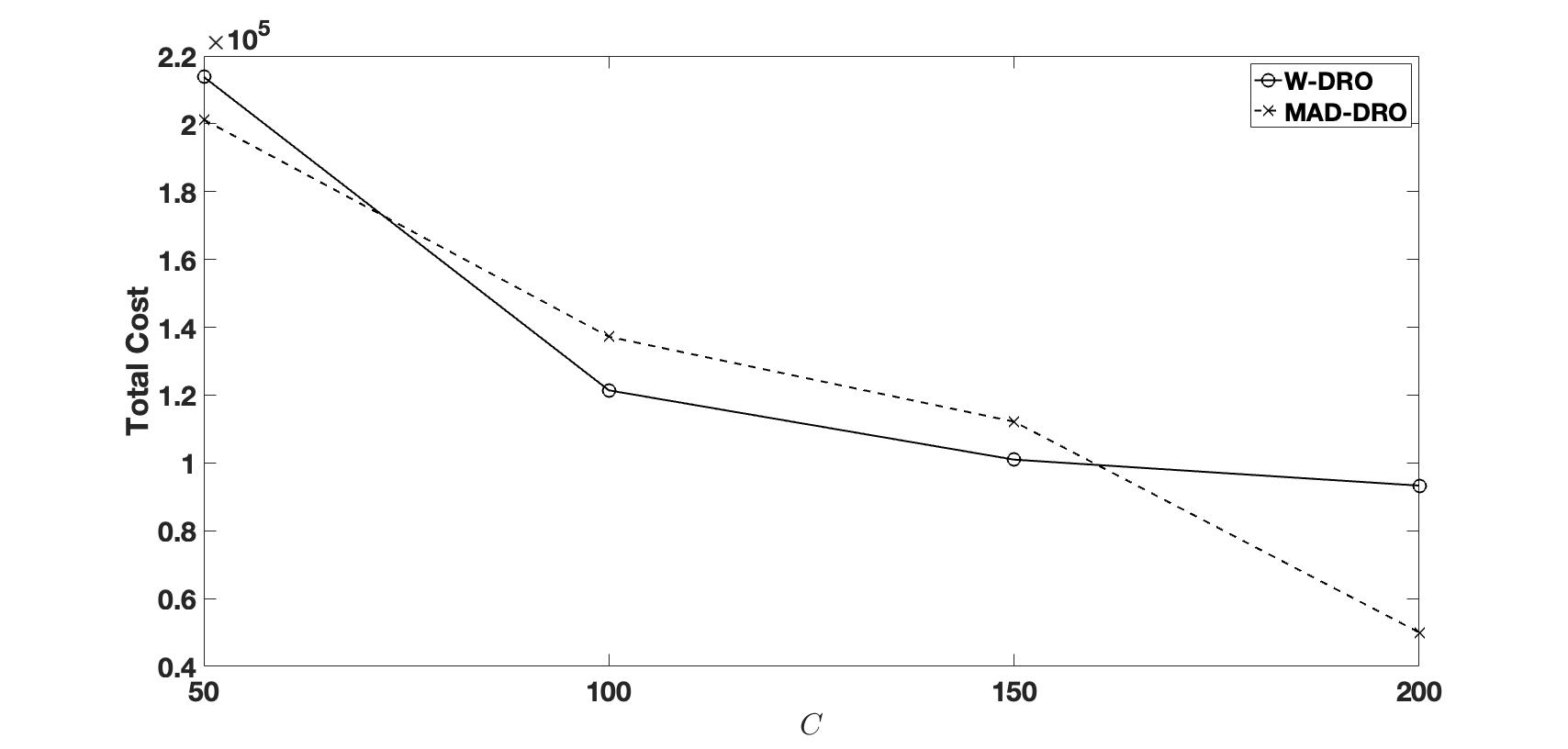}
      \caption{Total cost, $f=10,000$}
      \label{Fig7f}
    \end{subfigure}%
    \caption{Comparison of the results for different values of $C$ and $f$ under $\Wb \in [50, 100] $. Instance 1}\label{Fig7:MF_vs_C_Range2_inst1}
\end{figure}

\begin{figure}[t!]
     \begin{subfigure}[b]{0.5\textwidth}
 \centering
        \includegraphics[width=\textwidth]{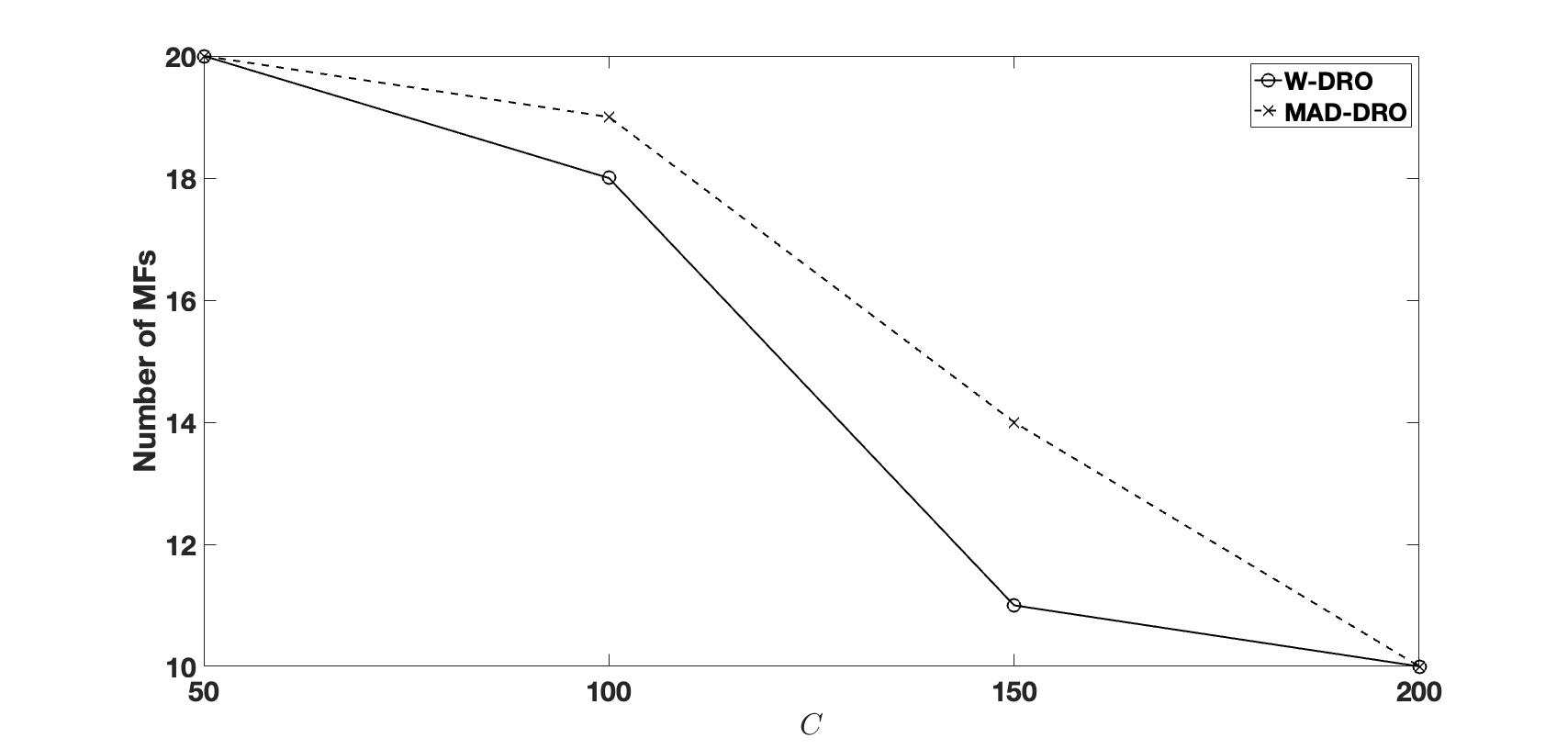}
        \caption{Number of MFs, $f=1,500$}
        \label{Fig10a}
    \end{subfigure}%
    \begin{subfigure}[b]{0.5\textwidth}
            \includegraphics[width=\textwidth]{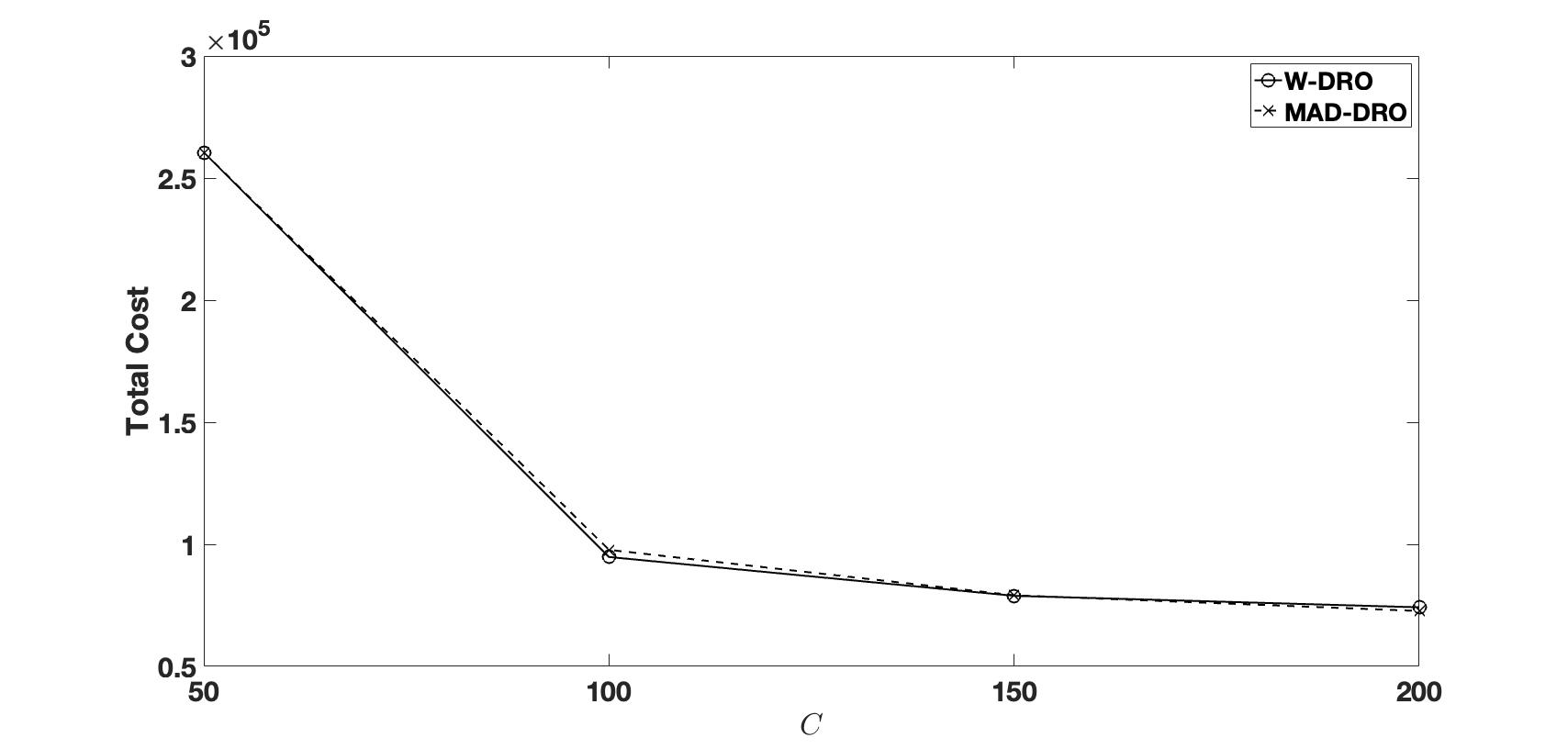}
      \caption{Total cost, $f=1,500$}
      \label{Fig10b}
    \end{subfigure}%
    
        \begin{subfigure}[b]{0.5\textwidth}
 \centering
        \includegraphics[width=\textwidth]{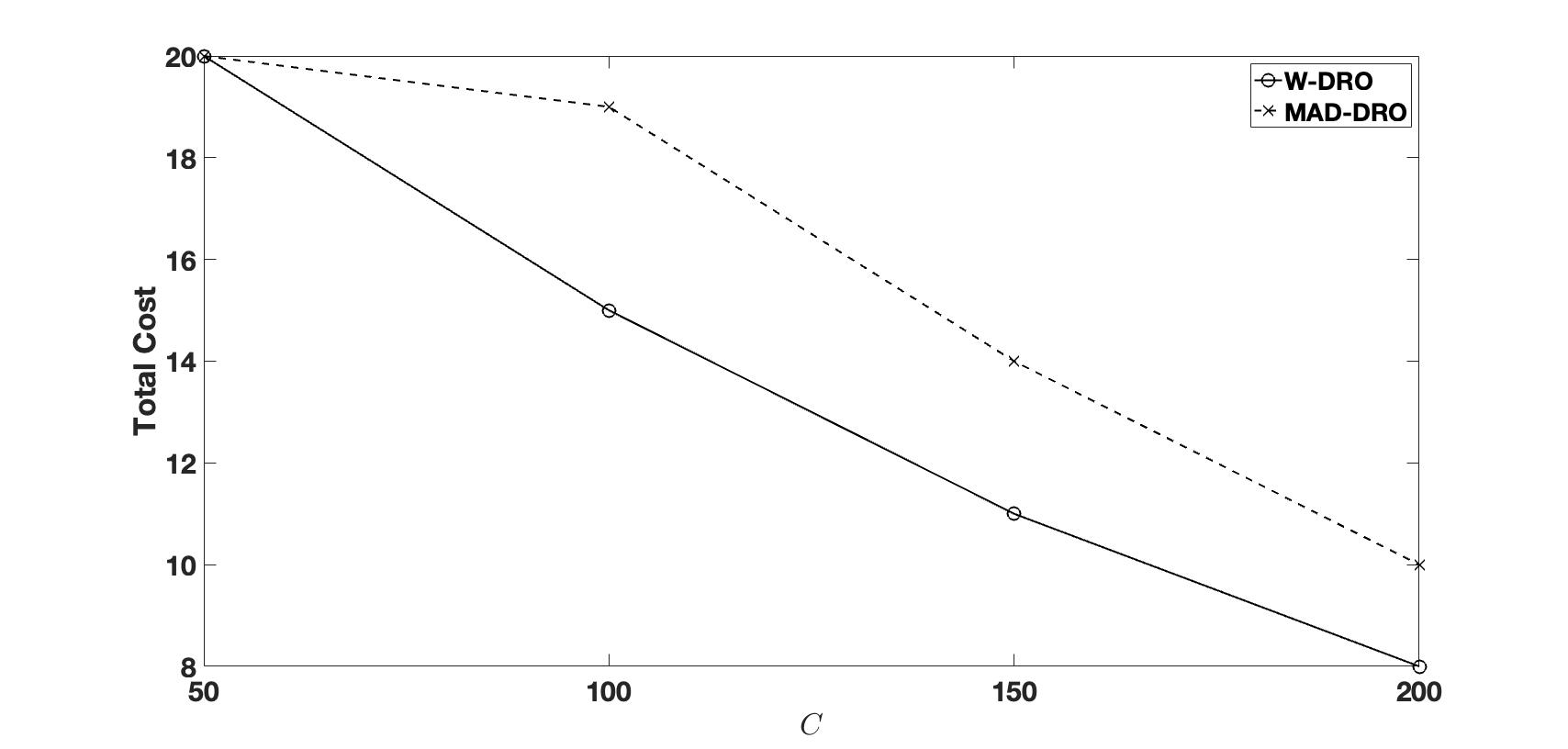}
        \caption{Number of MFs, $f=6,000$}
        \label{Fig10c}
    \end{subfigure}%
    \begin{subfigure}[b]{0.5\textwidth}
            \includegraphics[width=\textwidth]{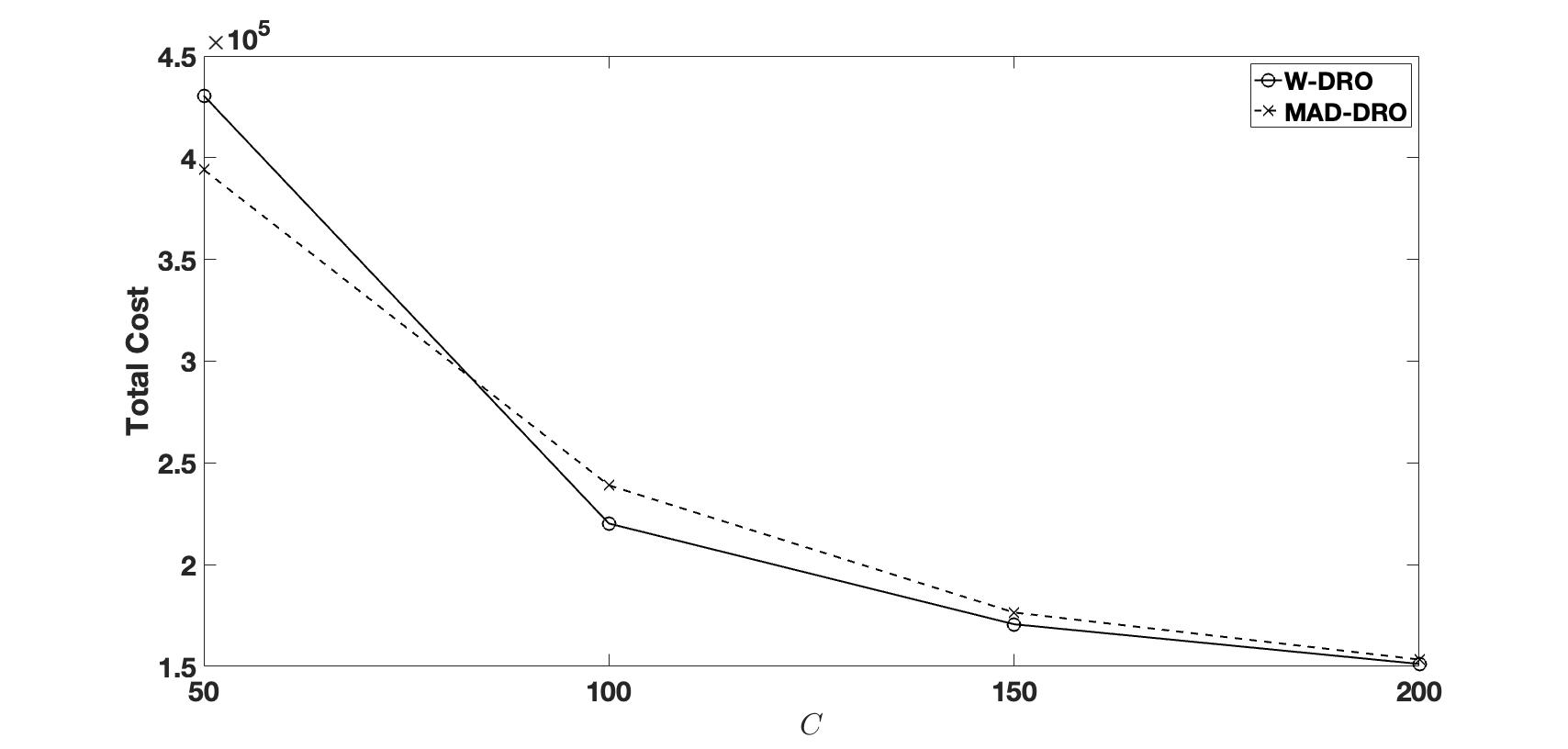}
      \caption{Total cost, $f=6,000$}
      \label{Fig10d}
    \end{subfigure}%
    
            \begin{subfigure}[b]{0.5\textwidth}
 \centering
        \includegraphics[width=\textwidth]{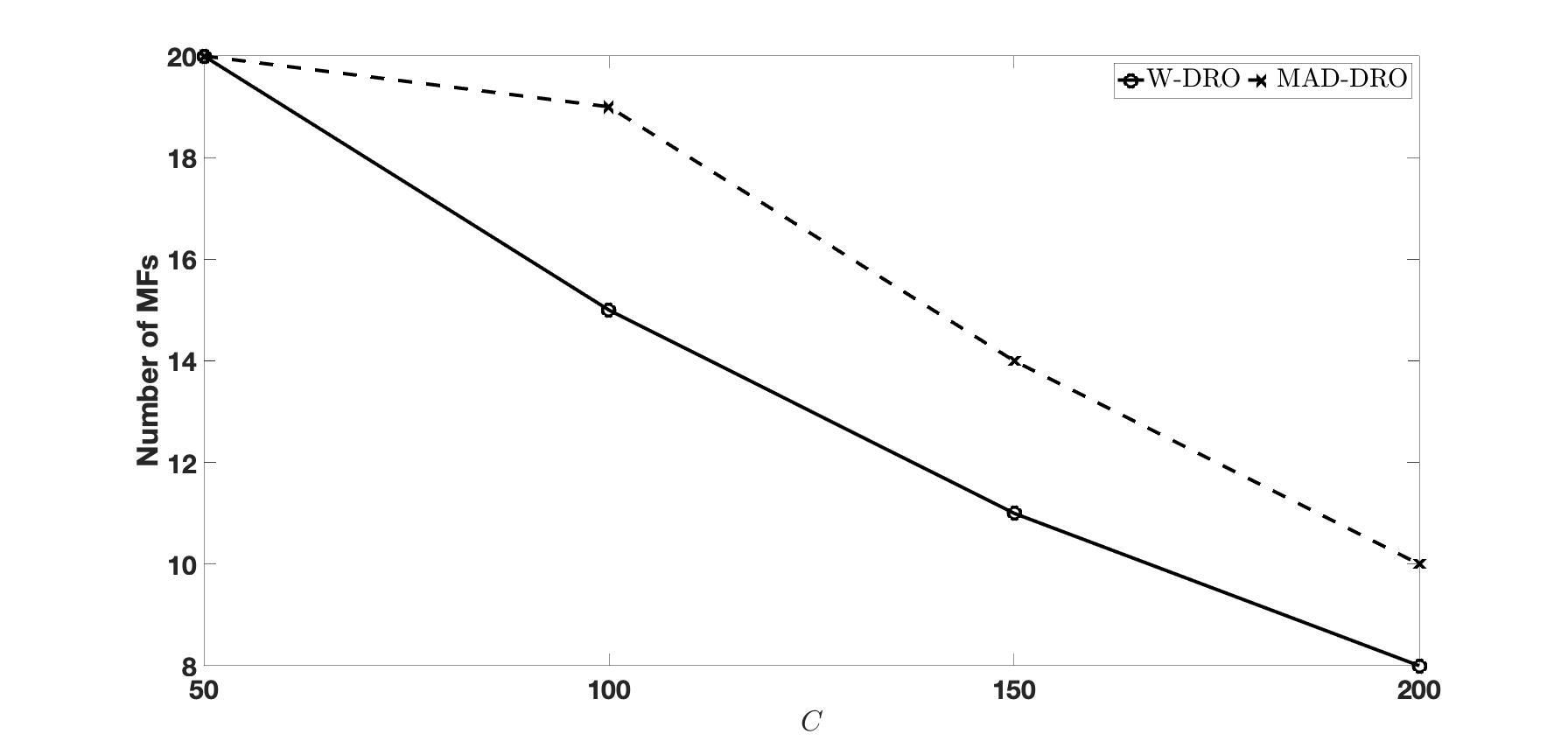}
        \caption{Number of MFs, $f=10,000$}
        \label{Fig10e}
    \end{subfigure}%
    \begin{subfigure}[b]{0.5\textwidth}
            \includegraphics[width=\textwidth]{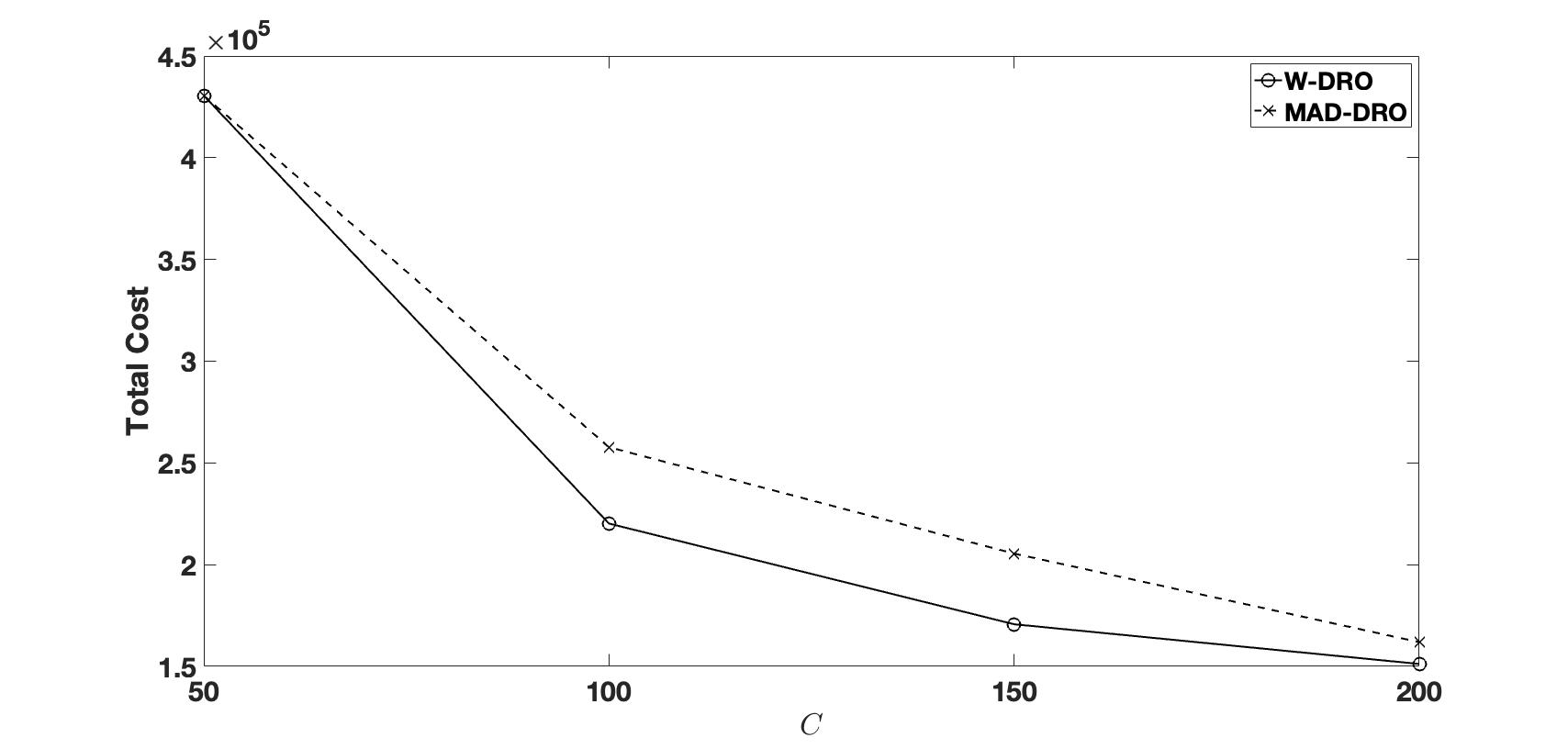}
      \caption{Total cost, $f=10,000$}
      \label{Fig10f}
    \end{subfigure}%
    \caption{Comparison of the results for different values of $C$ and $f$ under $\Wb \in [50, 100] $. Instance 5}\label{Fig10:MF_vs_C_Range2_inst5}
\end{figure}

\begin{figure}[t!]
 \centering
     \begin{subfigure}[b]{0.5\textwidth}
 \centering
        \includegraphics[width=\textwidth]{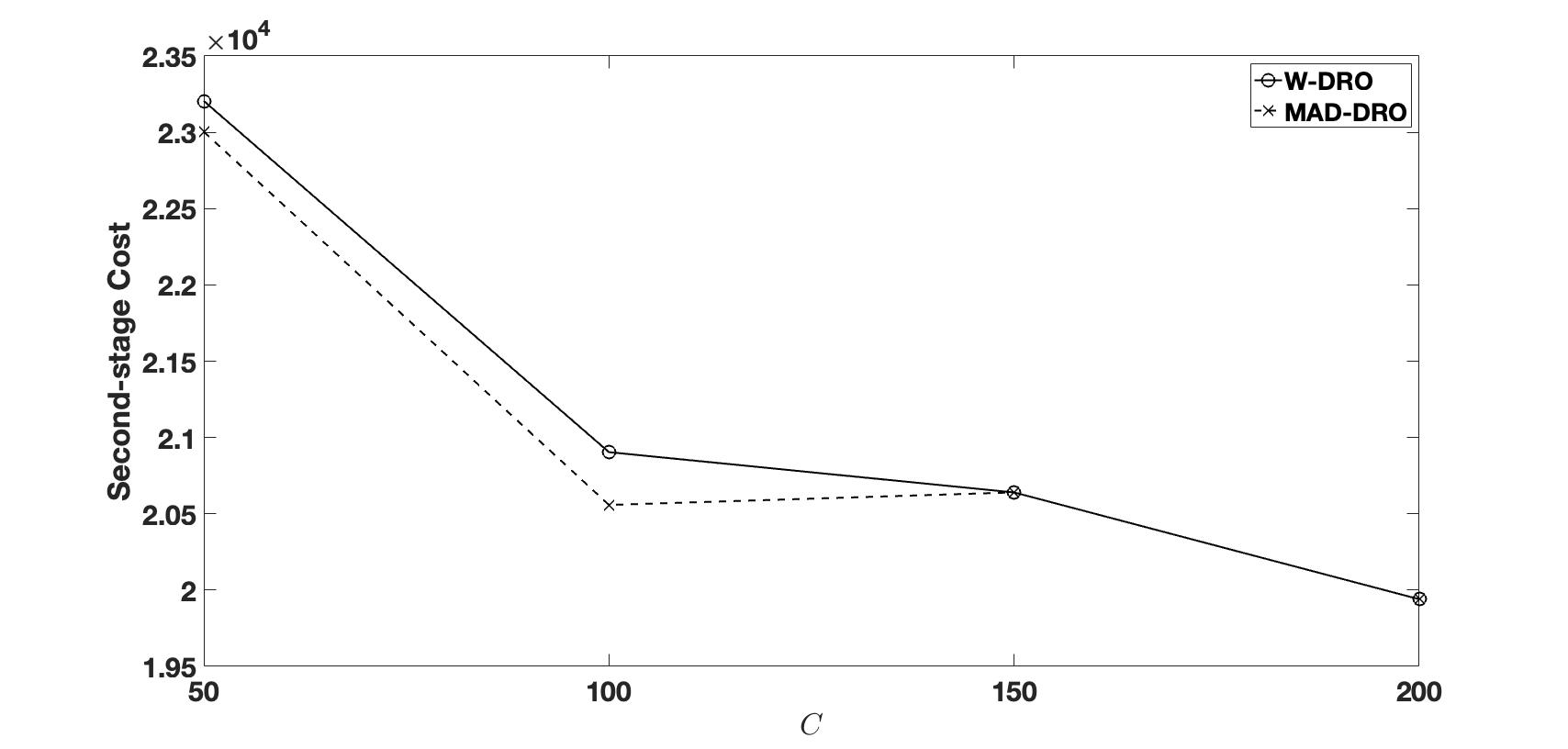}
        \caption{$f=1,500$}
        \label{Fig6a}
    \end{subfigure}%
    \begin{subfigure}[b]{0.5\textwidth}
            \includegraphics[width=\textwidth]{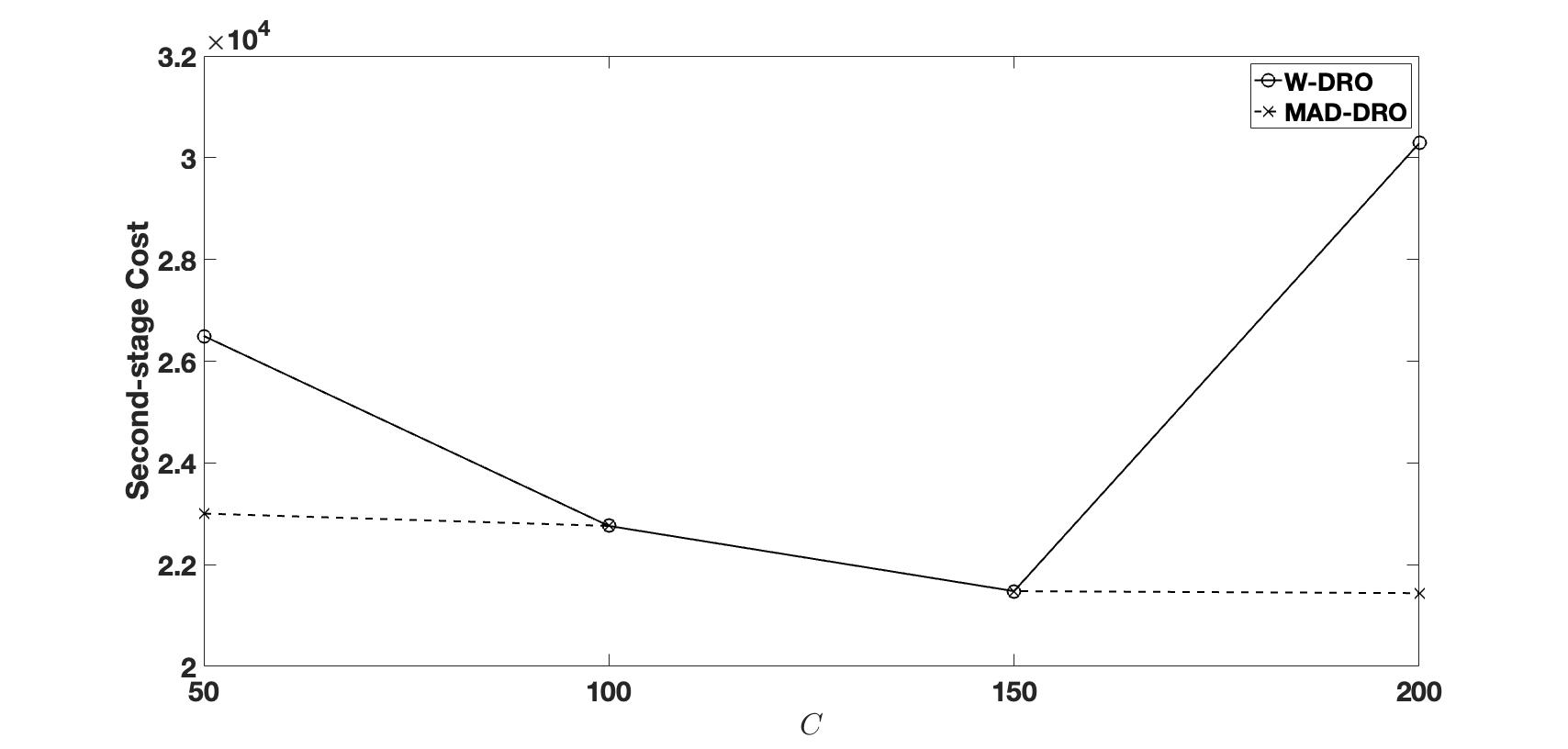}
      \caption{$f=6,000$}
      \label{Fig6b}
    \end{subfigure}%
    
        \begin{subfigure}[b]{0.5\textwidth}
 \centering
        \includegraphics[width=\textwidth]{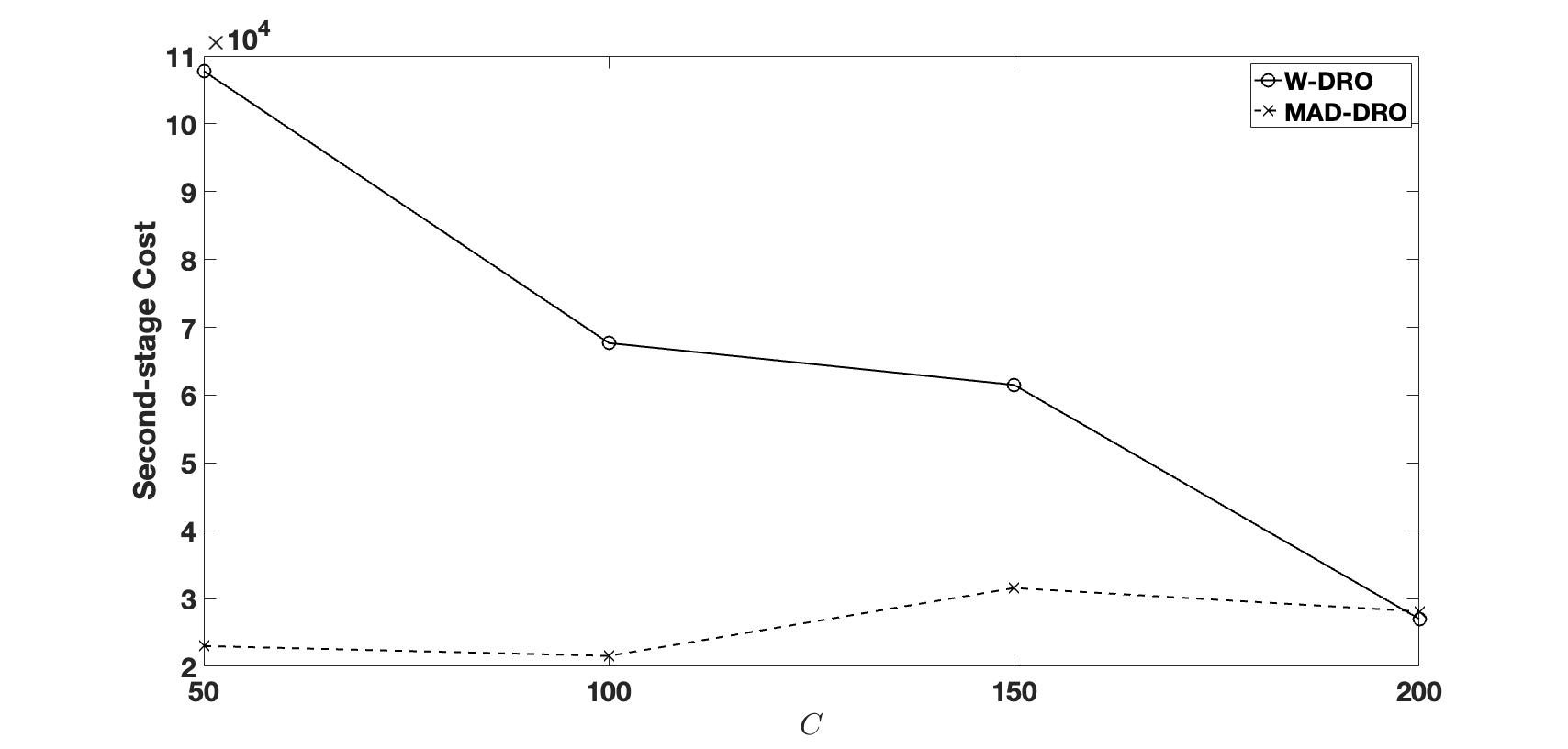}
        \caption{$f=10,000$}
        \label{Fig6c}
    \end{subfigure}%
  
    \caption{Comparison of second-stage cost for different values of $C$ and $f$ under $\Wb \in [20, 60]$. Instance 1}\label{Fig8:MF_vs_C_inst1}
\end{figure}

\begin{figure}[t!]
 \centering
     \begin{subfigure}[b]{0.5\textwidth}
 \centering
        \includegraphics[width=\textwidth]{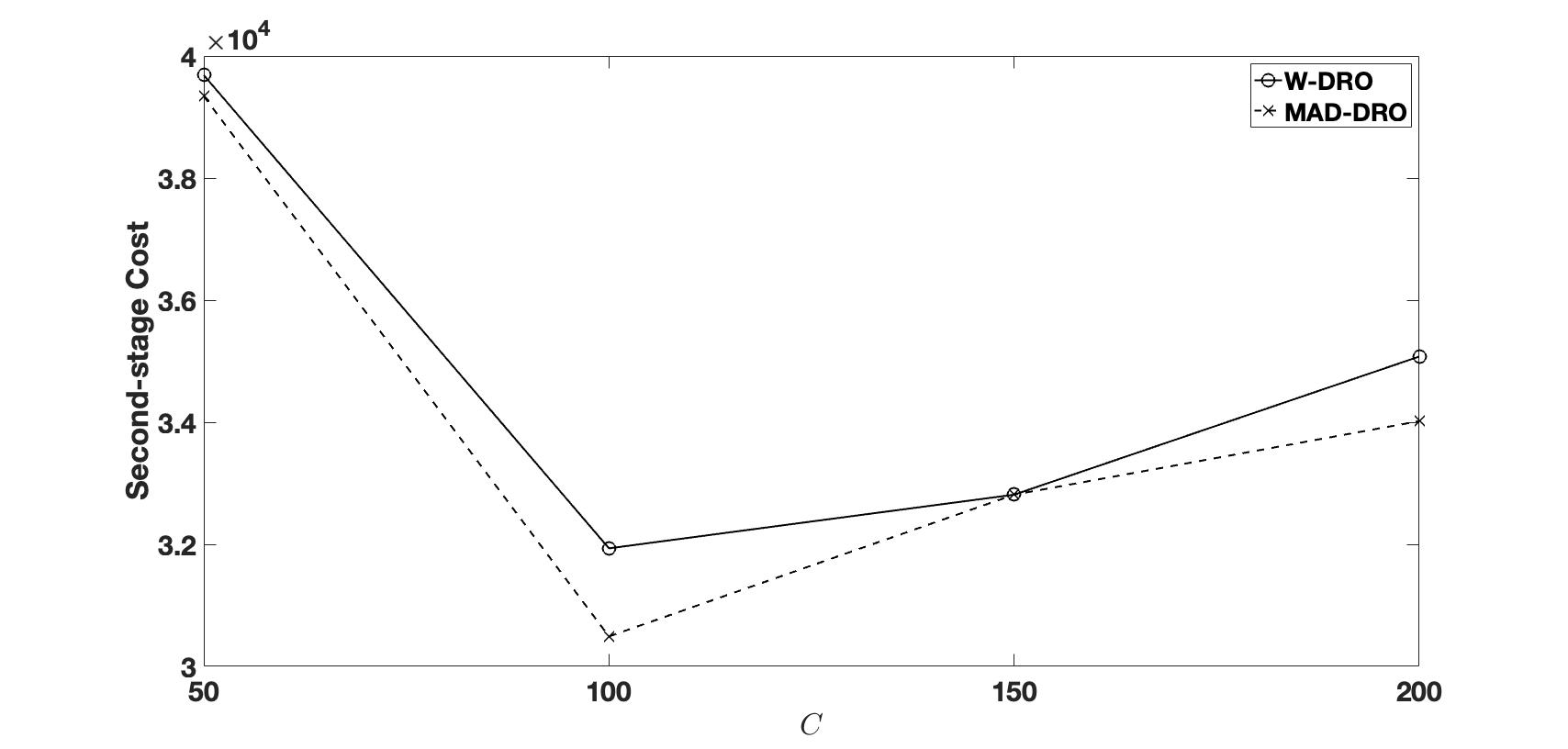}
        \caption{$f=1,500$}
        \label{Fig11a}
    \end{subfigure}%
    \begin{subfigure}[b]{0.5\textwidth}
            \includegraphics[width=\textwidth]{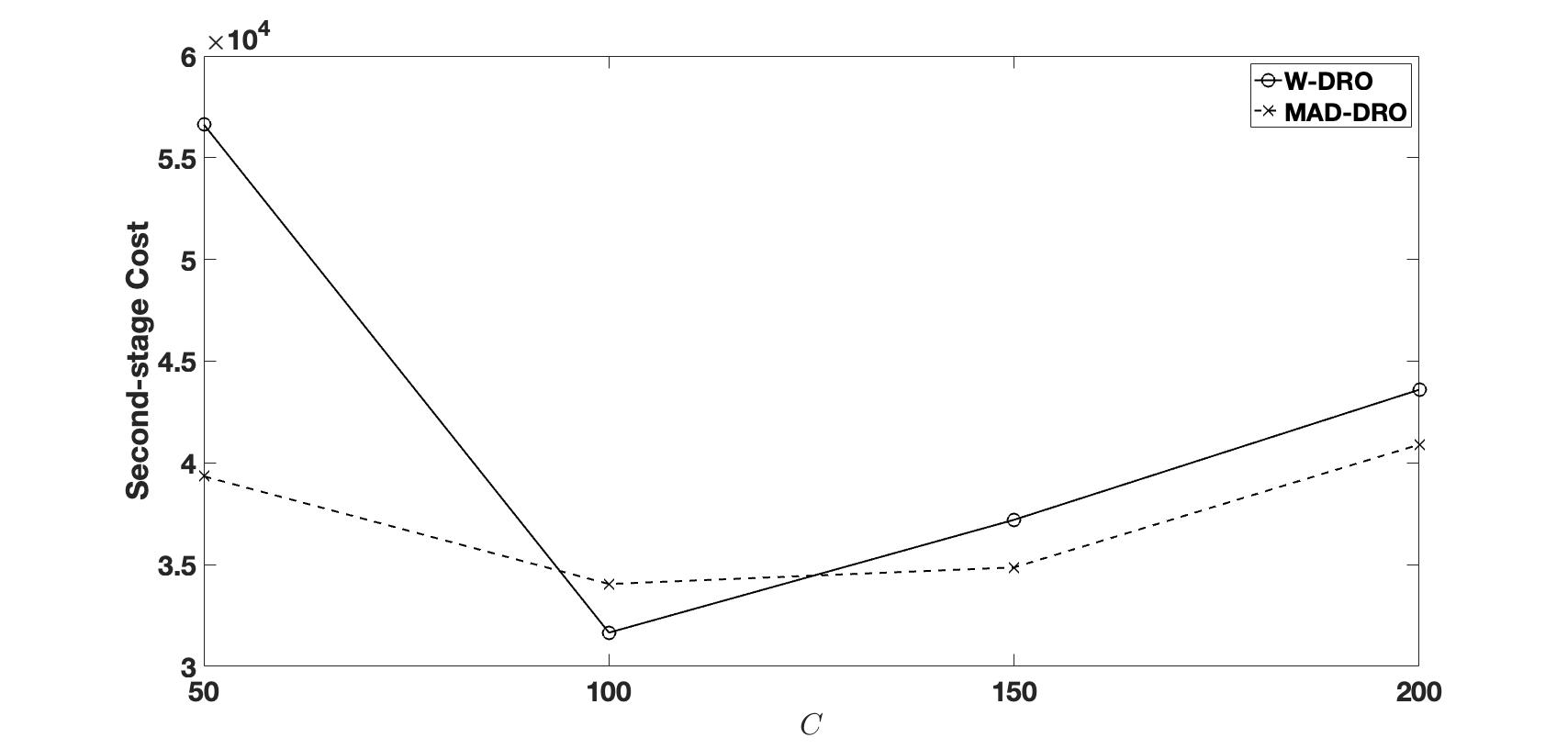}
      \caption{$f=6,000$}
      \label{Fig11b}
    \end{subfigure}%
    
        \begin{subfigure}[b]{0.5\textwidth}
 \centering
        \includegraphics[width=\textwidth]{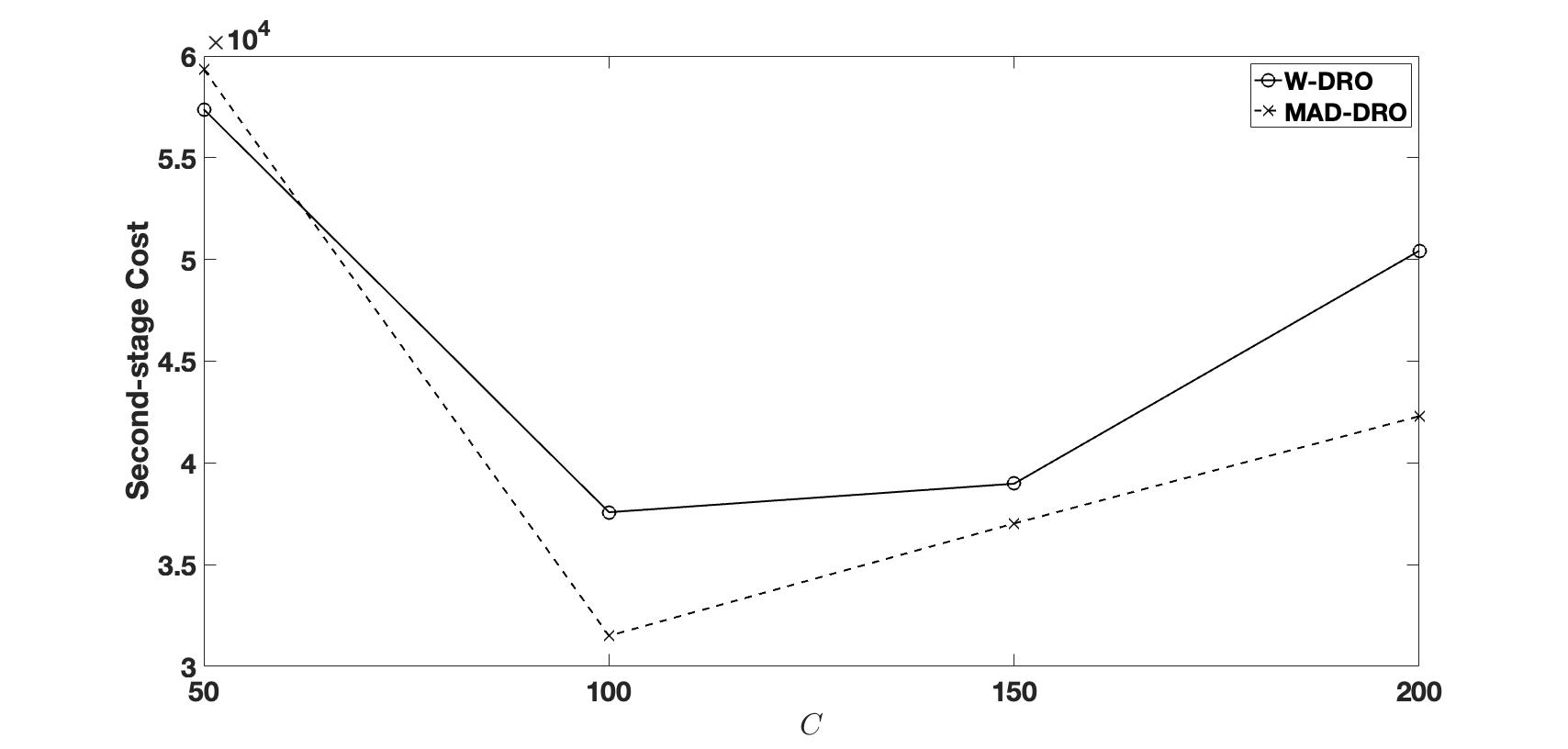}
        \caption{$f=10,000$}
        \label{Fig11c}
    \end{subfigure}%
  
    \caption{Comparison of second-stage cost for different values of $C$ and $f$ under $\Wb \in [20, 60]$. Instance 5}\label{Fig11:MF_vs_C_inst5}
\end{figure}

\begin{figure}[t!]
 \centering
     \begin{subfigure}[b]{0.5\textwidth}
 \centering
        \includegraphics[width=\textwidth]{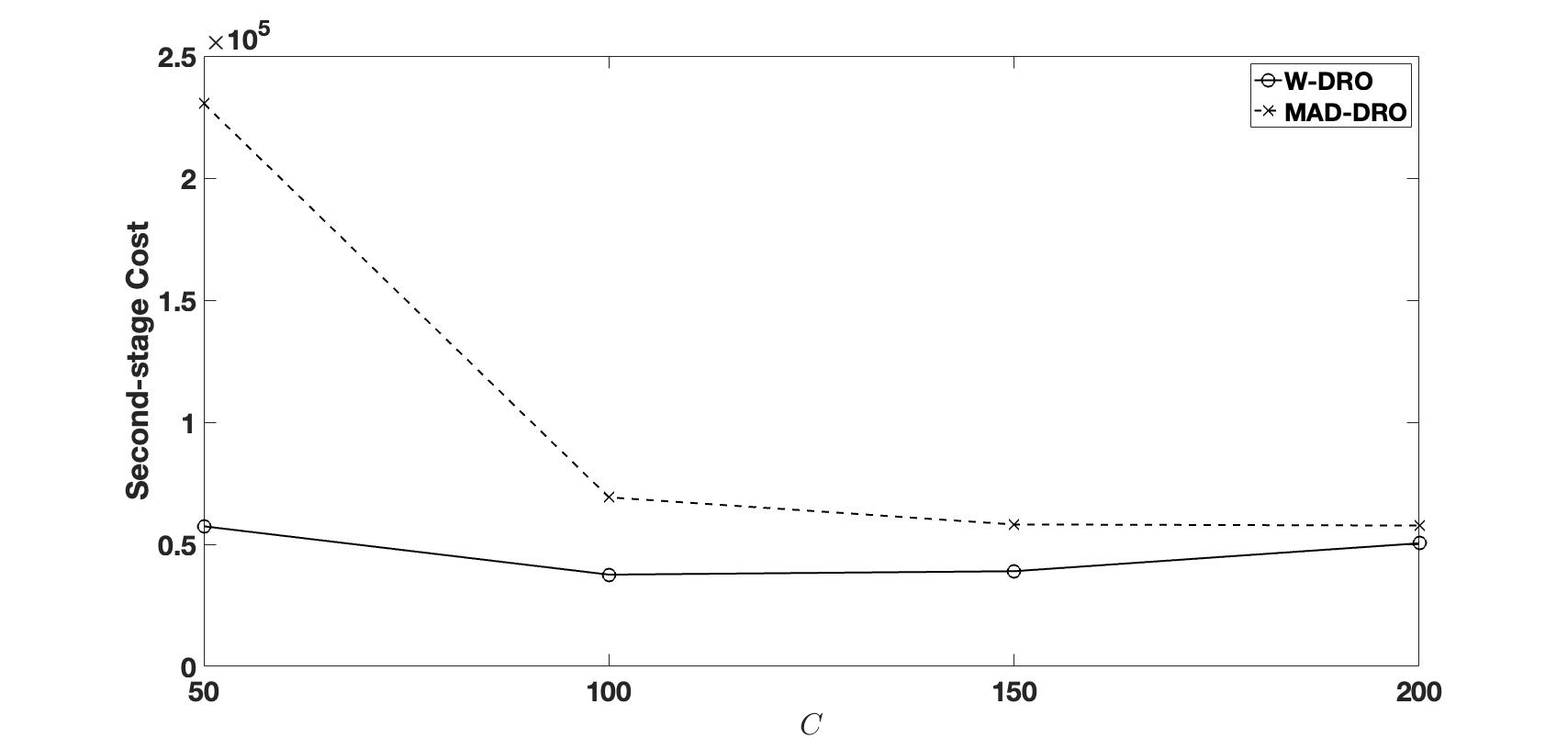}
        \caption{$f=1,500$}
        \label{Fig12a}
    \end{subfigure}%
    \begin{subfigure}[b]{0.5\textwidth}
            \includegraphics[width=\textwidth]{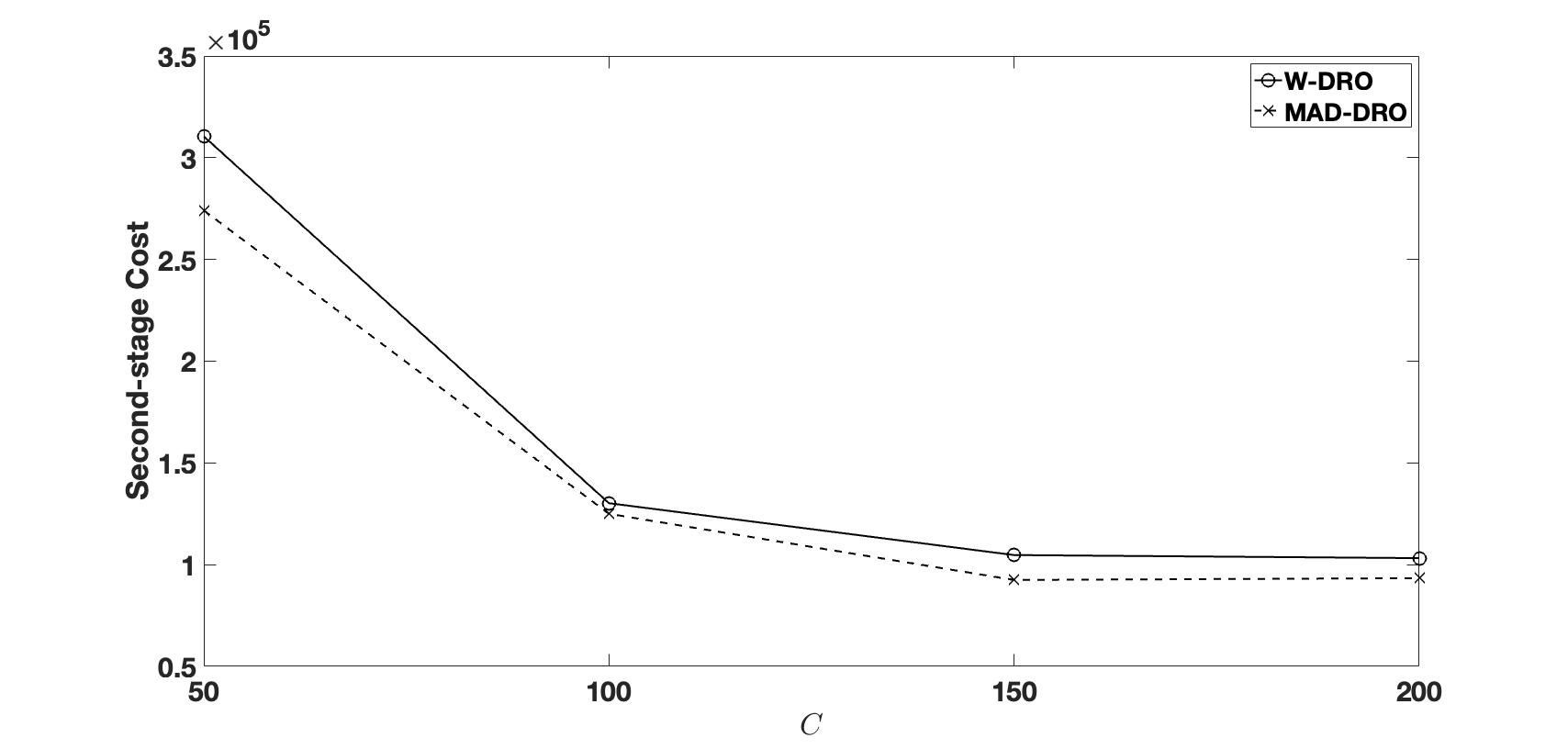}
      \caption{$f=6,000$}
      \label{Fig12b}
    \end{subfigure}%
    
        \begin{subfigure}[b]{0.5\textwidth}
 \centering
        \includegraphics[width=\textwidth]{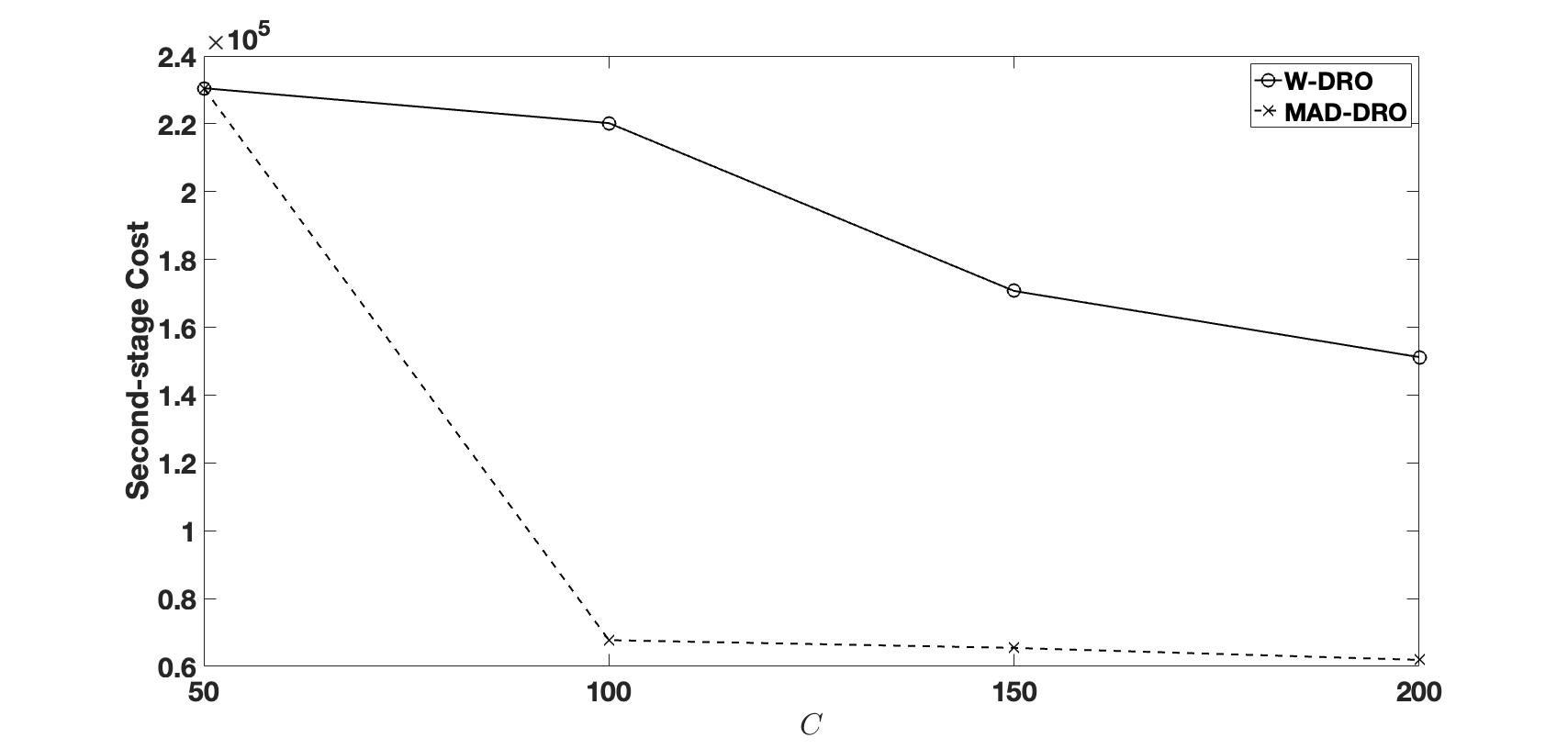}
        \caption{$f=10,000$}
        \label{Fig12c}
    \end{subfigure}%
  
    \caption{Comparison of second-stage cost for different values of $C$ and $f$ under $\Wb \in [50, 100]$. Instance 5}\label{Fig12:MF_vs_C_range2_inst5}
\end{figure}


%
%
%


\end{APPENDICES}
\clearpage
\newpage
\bibliographystyle{informs2014trsc} 
\bibliography{MFRSP_arXiv_KS} 


\end{document}